\documentclass{amsart}

\usepackage{latexsym, graphicx, psfrag, amscd, amssymb}

\newtheorem{thm}{Theorem}[section]
\newtheorem{df}[thm]{Definition}
\newtheorem{tm}[thm]{Theorem}
\newtheorem{ex}[thm]{Example}
\newtheorem{rem}[thm]{Remark}
\newtheorem{pro}[thm]{Proposition}
\newtheorem{cor}[thm]{Corollary}
\newtheorem{lm}[thm]{Lemma}
\newtheorem{con}[thm]{Conjecture}

\newcommand{\C}{\mathbb{C}}
\newcommand{\R}{\mathbb{R}}
\newcommand{\Z}{\mathbb{Z}}
\newcommand{\Q}{\mathbb{Q}}
\newcommand{\E}{\mathbb{E}}

\newcommand{\rel}{\mathrm{rel}}
\newcommand{\Aut}{\mathrm{Aut}}

\newcommand{\bB}{\mathbf{B}}
\newcommand{\bP}{\mathbf{P}}
\newcommand{\bp}{\mathbf{p}}
\newcommand{\bq}{\mathbf{q}}
\newcommand{\bs}{\mathbf{s}}
\newcommand{\bt}{\mathbf{t}}
\newcommand{\bx}{\mathbf{x}}
\newcommand{\by}{\mathbf{y}}
\newcommand{\bz}{\mathbf{z}}

\newcommand{\cA}{\mathcal{A}}
\newcommand{\cC}{\mathcal{C}}
\newcommand{\cE}{\mathcal{E}}
\newcommand{\cH}{\mathcal{H}}
\newcommand{\cK}{\mathcal{K}}
\newcommand{\cM}{\mathcal{M}}
\newcommand{\cN}{\mathcal{N}}
\newcommand{\cO}{\mathcal{O}}
\newcommand{\cP}{\mathcal{P}}
\newcommand{\cS}{\mathcal{S}}

\newcommand{\cU}{\mathcal{U}}
\newcommand{\cW}{\mathcal{W}}

\newcommand{\tU}{\tilde{U}}

\newcommand{\tg}{\tilde{g}}
\newcommand{\tn}{\tilde{n}}

\newcommand{\tsi}{\tilde{\sigma}}

\newcommand{\hC}{\hat{C}}
\newcommand{\hp}{\hat{p}}
\newcommand{\hq}{\hat{q}}
\newcommand{\hr}{\hat{r}}
\newcommand{\hs}{\hat{s}}
\newcommand{\hu}{\hat{u}}

\newcommand{\hS}{\hat{\Sigma}}
\newcommand{\hsi}{\hat{\sigma}}

\newcommand{\Si}{\Sigma}
\newcommand{\si}{\sigma}
\newcommand{\ep}{\epsilon}

\newcommand{\vm}{\vec{m}}
\newcommand{\vga}{\vec{\gamma}}

\newcommand{\domain}{ (\Si,\bB;\bp;\bq^1,\ldots,\bq^h)}
\newcommand{\another}{ (\Si',\bB';\bp';(\bq')^1,\ldots,(\bq')^h)}
\newcommand{\tdomain}{(\tilde{\Si},\tilde{\bB};\tilde{\bp};
\tilde{\bq}^1,\ldots,\tilde{\bq}^h)}
\newcommand{\bord}{(\Si;\bp;\bq)}

\newcommand{\Exone}{\mathrm{Ext}^1_{\cO_X}(\Omega_X,\cO_X)}
\newcommand{\Extwo}{\mathrm{Ext}^2_{\cO_X}(\Omega_X,\cO_X)}
\newcommand{\exone}{\mathrm{Ext}^1_{\cO_{X_0}}(\Omega_{X_0},\cO_{X_0})}
\newcommand{\extwo}{\mathrm{Ext}^2_{\cO_{X_0}}(\Omega_{X_0},\cO_{X_0})}
\newcommand{\pone}{\mathrm{Ext}^1_{\cO_X}(\Omega_X(D_\bx),\cO_X)}
\newcommand{\ptwo}{\mathrm{Ext}^2_{\cO_X}(\Omega_X(D_\bx),\cO_X)}
\newcommand{\lzero}{\cE xt^0_{\cO_X}(\Omega_X(D_\bx),\cO_X)}
\newcommand{\lone}{\cE xt^1_{\cO_X}(\Omega_X(D_\bx),\cO_X)}
\newcommand{\ltwo}{\cE xt^2_{\cO_X}(\Omega_X(D_\bx),\cO_X)}
\newcommand{\lhom}{\cH om_{\cO_X}(\Omega_X(D_\bx),\cO_X)}

\newcommand{\V}{ T_{\hat{r}_\alpha}\hat{X}\otimes
                 T_{\hat{r}_{l_0+\alpha}}\hat{X} } 
\newcommand{\Vbar}{ T_{\hat{\bar{r}}_\alpha}\hat{X}\otimes
                    T_{\hat{\bar{r}}_{l_0+\alpha}}\hat{X} }
\newcommand{\Vprime}{ T_{\hat{s}_{\alpha'}}\hat{X}\otimes
                 T_{\hat{s}_{l_0+\alpha'}}\hat{X} }
\newcommand{\W}{ H^1(\hat{C}_i,T_{\hat{C}_i}(-D_{\by^i})) }
\newcommand{\Wbar}{ H^1(\hat{\Bar{C}}_i, 
                    T_{\hat{\Bar{C}}_i}(-D_{\bar{\by}^i} )) }
\newcommand{\Wprime}{ H^1((\hat{\Si}_{i'})_\C, T_{(\hat{\Si}_{i'})_\C}
                    (-D_{\bz^{i'}}) } 

\newcommand{\PP}{\bP^1}
\newcommand{\OO}{\cO_{\PP}(-1)\oplus \cO_{\PP}(-1)}

\newcommand{\dbar}{\bar{\partial}}
\newcommand{\bS}{\partial\Sigma}
\newcommand{\inS}{\Sigma^\circ}
\newcommand{\Ce}{C_{\epsilon}}
\newcommand{\Se}{\Sigma_{\epsilon}}
\newcommand{\se}{\sigma_{\epsilon}}

\newcommand{\ra}{\rightarrow}
\newcommand{\form}{\Lambda^{0,1}}
\newcommand{\bra}{\langle}
\newcommand{\ket}{\rangle}
\newcommand{\pa}{\partial}

\newcommand{\MXL}{\bar{M}_{(g,h),(n,\vm)}(X,L\mid\beta,\vga,\mu)}
\newcommand{\WXL}{W^{k,p}_{(g,h),(n,\vm)}(X,L\mid\beta,\vga,\mu)}
\newcommand{\CXL}{C^l_{(g,h),(n,\vm)}(X,L\mid\beta,\vga,\mu)}
\newcommand{\onep}{W^{1,p}_{(g,h),(n,\vm)}(X,L \mid\beta,\vga,\mu)}

\newcommand{\half}{\frac{1}{2}}
\newcommand{\Wl}{W^{k,p}(\Sigma,\partial\Sigma,u^*TX,
                  (u|_{\partial\Sigma})^*TL)}
\newcommand{\Wless}{W^{k-1,p}(\Sigma,\Lambda^{0,1}\Sigma\otimes u^*TX)}
\newcommand{\Tl}{W^{1,p}(\Sigma,\partial\Sigma,u^*TX,
                  (u|_{\partial\Sigma})^*TL)}
\newcommand{\Tless}{L^p(\Sigma,\Lambda^{0,1}\Sigma\otimes u^*TX)}
\newcommand{\Cl}{C^l(\Sigma,\partial\Sigma,u^*TX,
                  (u|_{\partial\Sigma})^*TL)}
\newcommand{\Cless}{C^{l-1}(\Sigma,\Lambda^{0,1}\Sigma\otimes u^*TX)}
\newcommand{\Dl}{C^\infty(\Sigma,\partial\Sigma,u^*TX,
                  (u|_{\partial\Sigma})^*TL)}
\newcommand{\Dless}{C^\infty(\Sigma,\Lambda^{0,1}\Sigma\otimes u^*TX)}
\newcommand{\map}{(\Si,\bB;\bp;\bq^1,\ldots,\bq^h;u)}
\newcommand{\anothermap}{(\Si',\bB';\bp';(\bq')^1,\ldots, (\bq')^h;u')}

\newcommand{\mapaut}{\mathrm{Aut}\,\rho}
\newcommand{\domaut}{\mathrm{Aut}\,\lambda}
\newcommand{\SXL}{\overline{M}_{(g,h),(n+\hat{n},\vm+\vec{\hat{m}})}(X,L\mid 
\beta,\vga,\mu)}

\begin{document}

\title{Moduli of $J$-Holomorphic Curves with Lagrangian Boundary Conditions
and Open Gromov-Witten Invariants for an $S^1$-Equivariant Pair}
\author{Chiu-Chu Melissa Liu}
\email{ccliu@math.harvard.edu}
\address{Harvard University, Department of Mathematics, 
One Oxford Street, Cambridge, MA 02138, USA}
\date{\today} 

\begin{abstract}
Let $(X,\omega)$ be a symplectic manifold, $J$ be an $\omega$-tame
almost complex structure, and $L$ be a Lagrangian submanifold.
The stable compactification of the moduli space of parametrized $J$-holomorphic
curves in $X$ with boundary in $L$ (with prescribed topological data)
is compact and Hausdorff in Gromov's $C^\infty$-topology. 
We construct a Kuranishi structure with corners in the sense of Fukaya and
Ono. This Kuranishi structure is orientable if $L$ is spin. 
In the special case where the expected dimension of the moduli space  
is zero, and there is an $S^1$ action on the pair $(X,L)$ which
preserves $J$ and acts freely on $L$, we define the
Euler number for this $S^1$ equivariant pair and the prescribed 
topological data. We conjecture that  this rational number is 
the one computed by localization techniques using the given $S^1$ action.
\end{abstract}

\maketitle
\markboth{ }{ }

\tableofcontents

\section{Introduction}
\subsection{Background}

String theorists have been making predictions on enumerative invariants using
dualities. One of the most famous examples is the astonishing predictions of
the number of rational curves in a quintic threefold in \cite{CdGP}.
To understand these predictions, mathematicians first developed Gromov-Witten
theory to give the numerical invariants a rigorous mathematical definition so
that these predictions could be formulated as mathematical statements, and 
then tried to prove these statements. The predictions in \cite{CdGP} are
proven in \cite{G,LLY}.

Recently, string theorists have produced enumerative predictions 
about holomorphic curves with Lagrangian boundary conditions by 
studying dualities involving open strings \cite{OV,LMV,AV,AKV,MV}. 
Moreover, assuming the existence of a virtual fundamental cycle and the
validity of localization formulas, mathematicians have carried out
computations which coincide with these predictions \cite{KL,GZ}. 
In certain cases, these numbers can be reproduced by considering
relative morphisms \cite{LS}.
It is desirable to give a rigorous mathematical definition 
of these enumerative invariants, so that we may formulate physicists'
predictions as mathematical theorems, and then try to prove these theorems. 

Gromov-Witten invariants count $J$-holomorphic maps from a Riemann surface
to a fixed symplectic manifold $(X,\omega)$ together with an 
$\omega$-tame almost complex structure $J$.  These numbers can be 
viewed as intersection numbers on the moduli space of such maps. We want the
moduli space to be compact without boundary and oriented so that there 
exists a fundamental cycle which allows us to do intersection theory. The 
moduli space of $J$-holomorphic maps can be compactified by adding 
``stable maps'', whose domain is a Riemann surface which might have nodal
singularities. The stable compactification is compact and Hausdorff
in the $C^\infty$ topology defined by Gromov \cite{Gr}.
The moduli space is essentially almost complex, so it has a natural 
orientation. 
In general, the moduli space is not of the expected dimension 
and has bad singularities, but there exists a ``virtual fundamental cycle''
which plays the role of fundamental cycles \cite{LT1,BF,FO,LT2,Si}. These are now 
well-established in Gromov-Witten theory.
 
The ``open Gromov-Witten invariants'' that I want to establish shall count
$J$-holomorphic maps from a bordered Riemann surface to a symplectic manifold
$X$ as above such that the image of the boundary lies in a Lagrangian 
submanifold $L$ of $X$. To compactify the moduli space of such maps,
Sheldon Katz and I \cite{KL} introduced stable maps in
this context. The stable compactification is compact and Hausdorff in the
$C^\infty$ topology, as in the ordinary Gromov-Witten theory.
However, orientation is a nontrivial issue in the open Gromov-Witten
theory.  Moreover, the boundary is of real codimension one, so the
compactified moduli space does not ``close up'' as in the ordinary
Gromov-Witten theory, and the best we can expect is a fundamental ``chain''.

A fundamental ``chain'' is not satisfactory for intersection theory.
For example, the Euler characteristic of a compact oriented manifold without
boundary can be defined as the number of zeros of a generic vector field,
counted with signs determined by the orientation. This number is independent
of the choice of the vector field, and thus well-defined. For a compact oriented
manifold with boundary, one can still count the number of zeros of a generic
vector field with signs determined by the orientation, but the number will 
depend on the choice of the vector field. Therefore, we need to specify extra
boundary conditions to get a well-defined number.  

\subsection{Main results}

Let $(X,\omega)$ be a symplectic manifold of dimension $2N$, and 
$L$ be a Lagrangian submanifold. 
To compactify the moduli space of parametrized $J$-holomorphic curves in
$X$ with boundary in $L$, Gromov introduced cusp curves
with boundary \cite{Gr}, which I will call prestable maps.
A \emph{prestable map} to $(X,L)$ is a continuous map $f:(\Si,\bS)\ra(X,L)$ such that
$f\circ\tau:(\hat{\Si},\pa\hat{\Si})\ra (X,L)$ is $J$-holomorphic, where
$\Si$ is a prestable (i.e., smooth or nodal) bordered Riemann surface, and 
$\tau: \hat{\Si}\ra \Si$ is the normalization map \cite[Definition 3.6.2]{KL}.

A smooth bordered Riemann surface $\Si$ is of type $(g,h)$ if it is topologically
a sphere with $g$ handles and $h$ holes. The boundary of $\Si$ consists of $h$ disjoint
circles $B^1,\ldots,B^h$. We say $\Si$ has $(n,\vm)$ marked points
if there are $n$ distinct marked points in its interior and 
$m^i$ distinct marked points on $B^i$, where $\vm=(m^1,\ldots,m^h)$, $m^i\geq 0$. 
By allowing nodal singularities, we have the notion
of  a prestable Riemann surface of type $(g,h)$ with $(n,\vm)$ marked
points and an ordering $B^1,\ldots,B^h$ of the boundary components.    
An isomorphism between two such prestable bordered Riemann surfaces 
is an isomorphism of prestable bordered Riemann surfaces which
preserves the marked points and ordering of boundary components.
An isomorphism between two prestable maps $f:\Si\ra X$ and $f':\Si'\ra X$ is 
an isomorphism $\phi:\Si\ra\Si'$ in the above sense such that $f=f'\circ\phi$.
A prestable map is \emph{stable} if its automorphism group is finite. This is 
the analogue of Kontsevich's stable maps \cite{Ko} in the ordinary Gromov-Witten
theory.

For $\beta\in H^2(X,L;\Z)$, $\vga=(\gamma^1,\ldots,\gamma^h)\in H^1(L;\Z)^{\oplus h}$,
and $\mu\in\Z$, define 
\[ 
  \MXL 
\] 
to be the moduli space of isomorphism classes of stable maps 
$f:(\Si,\bS)\ra (X,L)$, where $\Si$ is a prestable
bordered Riemann surface of type $(g,h)$ with $(n,\vm)$ marked points and an ordering 
$B^1,\ldots, B^h$ of the boundary components, $f_*[\Si]=\beta$,
$f_*[B^i]=\gamma^i$, $i=1,\ldots, h$, and $\mu(f^*T_X,(f|_{\bS})^*T_L)=\mu$. Here
$\mu(f^*T_X,((f|_{\bS})^*T_L)$ is the Maslov index defined in
\cite[Definition 3.3.7]{KL}. 

\begin{tm}\label{cpthsf}
 $\MXL$ is compact and Hausdorff in the $C^\infty$ topology. 
\end{tm}

Here the $C^\infty$ topology is the one defined by Gromov's weak convergence
\cite{Gr}. The stability condition is necessary for Hausdorffness. The
compactness follows from \cite{Gr, Ye}, which will be explained in Section
~\ref{compact}. I do not claim any originality for Theorem~\ref{cpthsf}.

The boundary of the moduli space corresponds to degeneration of the domain or
blowup of the map. An interior node corresponds to a (real) codimension
$2$ stratum, while a boundary node corresponds to a codimension $1$
stratum. Blowup of the map at an interior point leads to the
well-known phenomenon of bubbling off of spheres which is codimension 2, 
while blowup at a boundary point leads to bubbling off of discs which is
codimension $1$. The intersection of two or more codimension $1$ strata
forms a corner. It is  shown in Section~\ref{constructku} that 

\begin{tm} \label{structure}
 $\MXL$ has a Kuranishi structure with corners of (real) virtual
 dimension $\mu+(N-3)(2-2g-h)+2n+m^1+\cdots+m^h$,
 where $2N$ is the (real) dimension of $X$. The Kuranishi structure
 is orientable if $L$ is spin or if $h=1$ and $L$ is relatively spin
 (i.e., $L$ is orientable and $w_2(T_L)=\alpha|_L$ for some
 $\alpha\in H^2(X,\Z_2)$). 
\end{tm}

The case for the disc with only boundary marked points ($g=n=0$, $h=1$) is
proven in \cite{FO3}.
I will describe briefly what a Kuranishi structure with corners is and 
refer to Section~\ref{defineku} for the complete definition.
A chart of a Kuranishi structure with corners is a $5$-uple
$(V,E,\Gamma,\psi,s)$, where $V$ is a smooth manifold (possibly
with corners), $\Gamma$ is a finite group acting on $V$, $E$ is 
a $\Gamma$-equivariant vector bundle over $V$, $s:V\ra E$ is 
a $\Gamma$-equivariant section, and $\psi$ maps $s^{-1}(0)/\Gamma$
homeomorphically to an open set of the moduli.
The dimension of $V$, rank of $E$, and the finite group $\Gamma$
might vary with charts, but $d=\dim V-\mathrm{rank}E$ is fixed
and is the virtual dimension of the Kuranishi structure with corners. 
$\det(TV)\otimes (\det E)^{-1}$ can be glued to an orbibundle,
the orientation bundle, and the Kuranishi structure with corners
is orientable if its orientation bundle is a trivial real line bundle.

If $s$ intersects the zero section transversally, $s^{-1}(0)$ is 
a manifold (possibly with corners) of dimension $d$.
In general, $s^{-1}(0)$ might have dimension larger than
$d$ and bad singularities due to the nontransversality of $s$.
The virtual fundamental chain can be
constructed by perturbing $s$ to a transversal section. Locally it is a
singular chain with rational coefficients in $V/\Gamma$ 
which is a rational combination of the images of
$d$-dimensional submanifolds of $V$.   

A virtual fundamental chain is not satisfactory for intersection theory.
For example, when $X$ is a Calabi-Yau threefold and $L$ is a special 
Lagrangian submanifold, 
$\bar{M}_{(g,h),(0,\vec{0})}(X,L\mid\beta,\vga,\mu)$
is empty for $\mu\neq 0$, and the expected dimension of 
$\bar{M}_{(g,h),(0,\vec{0})}(X,L\mid\beta, \vga, 0)$
is zero for any $g$, $h$, $\beta$, $\vga$.
The virtual fundamental chain is a zero chain with rational coefficients, and
we would like to define the invariant 
$\chi_{(g,h)}(X,L\mid \beta,\vga, \mu)\in\Q$
to be the degree of this zero chain. However, this number
depends on the perturbation, so we need to impose extra boundary conditions
to obtain a well-defined number. Now assume that 
\begin{itemize}
\item There is an $S^1$ action $\varrho:S^1\times X\ra X$ 
      which preserves  $J$ and $L$.
\item The restriction of $\varrho$ to $L$ is \emph{free}.
\item The virtual dimension of 
      $\bar{M}_{(g,h),(0,\vec{0})}(X,L\mid \beta,\vga,\mu)$
      is zero.
\end{itemize}
Under the above assumptions I can, using the $S^1$ action $\varrho$, impose boundary 
conditions to get a well-defined rational number
\[ 
\chi_{(g,h)}(X,L,\varrho \mid \beta,\vga,\mu) 
\]
which is an invariant of the equivariant pair $(X,L,\varrho)$, but 
\emph{not} an invariant of the pair $(X,L)$. I conjecture that
these rational numbers are the ones computed by localization techniques
using the $S^1$ action $\varrho$ \cite{KL,GZ}. The computations
in \cite{KL,LS,GZ} coincide with physicists' predictions.

\bigskip

\paragraph{\bf Acknowledgments:}
First and foremost, I would like to thank my thesis advisor Shing-Tung 
Yau for leading me to the field of mirror symmetry and providing the best
environment to learn its newest developments. 
Secondly, I thank Clifford Taubes and Gang Liu for answering my questions 
on symplectic geometry and carefully reading the draft. 
I thank Sheldon Katz for being so generous to cooperate with me. The 
cooperation \cite{KL} is a very instructive experience for me and led 
me to this project. 
I thank Cumrun Vafa for suggesting to me this fruitful problem and patiently 
explaining his works to me.
I thank Jason Starr for being an incredibly generous and patient mentor of 
algebraic geometry. 
I thank Xiaowei Wang for always being a source of mathematical knowledge 
and moral support. 
I thank Arthur Greenspoon for his numerous valuable suggestions on the first 
draft, and Chien-Hao Liu for his meticulous proofreading.
I thank Spiro Karigiannis for his great help on my English. 
In addition, it is a pleasure to thank William Abikoff, Selman Akbulut, 
Raoul Bott, Kevin Costello, Yong Fu, Kenji Fukaya, Tom Graber, Irwin Kra,
Kefeng Liu, Curtis McMullen, Maryam Mirzakhani, Yong-Geun Oh, Kaoru Ono, 
Scott Wolpert, and Eric Zaslow for helpful conversations.
Finally, special thanks go to Ezra Getzler for corrections and
refinements of the part on moduli spaces of bordered Riemann
surfaces.

\section{Surfaces with Analytic or Dianalytic Structures}

In this section, we review some definitions
and facts of surfaces with analytic or dianalytic structures,
following \cite[Chapter 1]{AG} closely. This section is 
an expansion of Section 3.1 and 3.2 of \cite{KL}.
 
The {\em marked bordered Riemann surfaces} defined in Section 
\ref{sec:marked} are directly related to open Gromov-Witten theory.  

\subsection{Analyticity and dianalyticity}

   \begin{df}
    Let $A$ be a nonempty open subset of $\C$, $f:A\ra \C$ a map.
    $f$ is \emph{analytic} on $A$ if $\frac{\pa f}{\pa \bar{z}}=0$,
    and \emph{antianalytic} on $A$ if $\frac{\pa f}{\pa z}=0$.
    $f$ is said to be \emph{dianalytic} if its restriction to each component
    of $A$ is either analytic of antianalytic.
   \end{df}

   \begin{df}
    Let $A$ and $B$ be nonempty subsets of
    $\C^+=\{z\in \C\mid\textup{Im}z\geq 0 \}$.
    A continuous function $f:A\ra B$ is \emph{analytic} (resp.
    \emph{antianalytic} on $A$ if it extends to an \emph{analytic}
    (resp. \emph{antianalytic}) function $f_\C: U\ra\C$,
    where $U$ is an open neighborhood of $A$ in $\C$. $f$ is said to be
    \emph{dianalytic} on $A$ if its restriction to each
    component of $A$ is either analytic or antianalytic.
   \end{df}

   \begin{tm}[Schwarz reflection principle]
    Let $A$ and $B$ be nonempty subsets of
    $\C^+=\{z\in \C\mid\textup{Im}z\geq0\}$.
    A continuous function $f:A\ra B$ is dianalytic (resp. analytic)
    if it is dianalytic (resp. analytic) on the interior of $A$ and satisfies
    $f(A\cap\R)\subset \R$.
   \end{tm}

    \begin{df}
    A \emph{surface} is a Hausdorff, connected, topological space $\Si$
    together with a family $\cA=\{(U_i, \phi_i)\mid i\in I\}$ such that
    $\{U_i\mid i\in I\}$ is an open covering of $\Si$ and each map
    $\phi_i:U_i\ra A_i$ is a homeomorphism onto an open subset
    $A_i$ of  $\C^+$. $\cA$ is called a \emph{topological atlas} on $\Si$,
    and each pair $(U_i, \phi_i)$ is called a \emph{chart} of $\cA$.
    The boundary  of $\Si$ is the set
    \[ \bS=\{x\in \Si\mid \exists\,i\in I \textup{ s.t. } x\in U_i,
     \phi_i(x)\in\R \} \]
    $\phi_{ij}\equiv \phi_i\circ\phi_j^{-1}:\phi_j(U_i\cap U_j)\ra
    \phi_i(U_i\cap U_j) $ are surjective homeomorphisms, called the
    \emph{transition functions} of $\cA$. $\cA$ is called a \emph{dianalytic}
    (resp. \emph{analytic}) \emph{atlas} if all its transition functions
     are dianalytic (resp. analytic).
   \end{df}

   \subsection{Various categories of surfaces}

   \subsubsection{Riemann surfaces}

    \begin{df}
     A \emph{Riemann surface} is a surface equipped with the analytic
     structure induced by an analytic atlas on it.
    \end{df}

    A Riemann surface is canonically oriented by its analytic structure.

   \subsubsection{Symmetric Riemann surfaces}

   \begin{df}
    A \emph{symmetric Riemann surface} is a Riemann surface $\Si$ together
    with an antiholomorphic involution $\si:\Si\ra \Si$.
    $\si$ is called the \emph{symmetry} of $\Si$.
   \end{df}

   \begin{df} \label{symmor}
    A \emph{morphism} between symmetric Riemann surfaces $(\Si,\si)$ and
    $(\Si', \si')$ is an analytic map $f:\Si\ra\Si'$ such that
    $f\circ\si=\si'\circ f$.
   \end{df}

   \begin{df}
   A \emph{symmetric Riemann surface with $(n,m)$ marked points} is a
   symmetric Riemann surface $(\Si,\si)$ together with $2n+m$ distinct
   points $p_1,\ldots,p_{2n+m}$ in $\Si$ such that $\si(p_i)=p_{n+i}$ for
   $i=1,\ldots,n$ and $\si(p_i)=p_i$ for $i=2n+1,\ldots,2n+m$.
   \end{df}

   \subsubsection{Klein surfaces}

    \begin{df}
     A \emph{Klein surface} is a surface equipped with the dianalytic
     structure induced by a dianalytic atlas on it.
    \end{df}

     A Riemann surface can be viewed as a Klein surface. A Klein surface
     can be equipped with an analytic structure compatible with
     the dianalytic structure if and only if it is orientable.
     In particular, an orientable Klein surface without boundary
     admits a compatible structure of a Riemann surface.

   \begin{df}
    A \emph{morphism} between Klein surfaces $\Si$ and $\Si'$ is
    a continuous map $f:(\Si, \bS)\ra (\Si', \bS')$ such that for
    any $x\in\Si$ there exist analytic charts $(U,\phi)$ and
    $(V,\psi)$ about $x$ and $f(x)$ respectively, and an analytic function
    $F:\phi(U) \ra \C$ such that the following diagram commutes:
\[
    \begin{CD}
     U &@>f>>& V\\
     @VV{\phi}V && @VV{\psi}V\\
     \phi(U) @>F>> \C @>\Phi>> \C^+
    \end{CD}
\]
    where $\Phi(x+iy)=x+i|y|$ is the \emph{folding map}.
   \end{df}

    Given a Klein surface $\Si$, there are three ways to construct
    an unramified double cover of $\Si$. We refer to
    \cite[1.6]{AG} for the precise definition
    of an unramified double cover and detailed constructions.
    The \emph{complex double}
    $\Si_\C$ is an orientable Klein surface without boundary. The
    \emph{orienting double} $\Si_\mathbf{O}$ is an orientable Klein
    surface. It is disconnected if and only if $\Si$ is
    orientable, and it has nonempty boundary if and only of $\Si$ does.
    The \emph{Schottkey double} $\Si_\mathbf{S}$ is a Klein surface
    without boundary. It is disconnected
    if and only if the boundary of $\Si$ is empty, and it is nonorientable
    if and only if $\Si$ is.

    If $\Si$ is orientable, then $\Si_\C=\Si_\mathbf{S}$, and $\Si_\mathbf{O}$
    is disconnected (the trivial double cover). If $\bS=\phi$, then
    $\Si_\C=\Si_\mathbf{O}$, and $\Si_\mathbf{S}$ is disconnected. In
    particular, if $\Si$ comes from a Riemann surface, then all three covers are
    the trivial disconnected double cover.

   \begin{ex}
    Let $\Si$ be a M\"obius strip. Then $\Si_\C$ is a torus, $\Si_\mathbf{S}$ is a
    Klein bottle, and $\Si_\mathbf{O}$ is an annulus.
   \end{ex}

 \subsubsection{Bordered Riemann surfaces}

  \begin{df}
   A \emph{bordered Riemann surface} is a compact surface with nonempty
   boundary equipped with the analytic structure induced by an
   analytic atlas on it.
  \end{df}

  \begin{rem}\label{boundary}
   A bordered Riemann surface is canonically oriented by the analytic
   (complex) structure.
   In the rest of this paper, the boundary circles $B^i$ of a bordered
   Riemann surface $\Si$ with boundary $\bS=B^1\cup\ldots\cup B^h$
   will always be given the orientation induced by the
   complex structure, which is a choice of
   tangent vector to $B^i$ such that the basis
   (the tangent vector of $B^i$, inner normal) for the real tangent space
   is consistent with the orientation of $\Si$ induced by the complex
   structure.
  \end{rem}

  \begin{df}\label{bordermor}
   A \emph{morphism} between bordered Riemann surfaces $\Si$ and $\Si'$ is
   a continuous map $f:(\Si, \bS)\ra (\Si', \bS')$
   such that for any $x\in\Si$ there exist analytic charts $(U,\phi)$ and
   $(V,\psi)$ about $x$ and $f(x)$ respectively, and an analytic function
   $F:\phi(U) \ra \C$ such that the following diagram commutes:
\[
   \begin{CD}
   U @>f>> V\\
   @VV{\phi}V  @VV{\psi}V\\
   \phi(U) @>F>> \C
   \end{CD}
\]
  \end{df}

  A bordered Riemann surface is topologically a sphere with $g\geq 0$ handles
  and with $h>0$ discs removed. Such a bordered Riemann surface is said to
  be of type $(g,h)$.

  A bordered Riemann surface can be viewed as a Klein surface.
  Its complex double and Schottkey double coincide since it is orientable.

 \subsubsection{Marked bordered Riemann surfaces}\label{sec:marked}

The following refinement of an earlier definition is suggested to the author
by Ezra Getzler.
  \begin{df}\label{marked}
   Let $h$ be a positive integer, $g,n$ be nonnegative integers, and
   $\vm=(m^1,\ldots,m^h)$ be an $h$-uple of nonnegative integers.
   A \emph{marked bordered Riemann surface of type $(g,h)$ with
   $(n,\vm)$ marked points} is an $(h+3)$-uple
$$
\domain
$$ 
whose components are described as follows.
   \begin{itemize}
    \item $\Si$ is a bordered Riemann surface of type $(g,h)$.
    \item $\bB=(B^1,\ldots,B^h)$, where $B^1,\ldots, B^h$ are connected components 
          of $\bS$, oriented as in Remark~\ref{boundary}.
    \item $\bp=(p_1,\ldots,p_n)$ is an $n$-uple of distinct points in $\inS$.
    \item $\bq^i=(q^i_1,\ldots,q^i_{m^i})$ is an $m^i$-uple of distinct points on 
          the circle $B^i$.
   \end{itemize}
   \end{df}

Let $\vec{0}=(0,\ldots,0)$.
Note that a marked bordered Riemann surface of type $(g,h)$ with
$(n,\vec{0})$ marked points is a bordered Riemann surface togther
with an ordering of the $h$ boundary components.

   \begin{df}\label{markedmor}
    A \emph{morphism} between marked bordered Riemann surfaces of type
    $(g,h)$ with $(n,\vm)$ marked points
$$
\domain \to \another
$$
is an isomorphism of bordered Riemann surface $f:\Si\ra\Si'$ such that
$f(B^i)=(B')^i$ for $i=1,\ldots,h$, $f(p_j)=p'_j$ for $j=1,\ldots,n$, and
$f(q^i_k)=(q')^i_k$ for $k=1,\ldots,m^i$.
   \end{df}

   \begin{rem}
    The category of marked bordered Riemann surfaces of type $(g,h)$
    with $(n,\vm)$ marked points is a groupoid since every morphism in
    Definition~\ref{markedmor} is an isomorphism.
   \end{rem}

   \subsection{Topological types of compact symmetric Riemann surfaces}

    A compact symmetric Riemann surface is topologically a compact orientable
    surface without boundary $\Si$ together with an orientation reversing
    involution $\si$, which is classified by the following three invariants:
    \begin{enumerate}
    \item The genus $\tg$ of $\Si$.
    \item The number $h=h(\si)$ of connected components of $\Si^\si$, the
         fixed locus of $\si$.
    \item The index of orientability, $k=k(\si)\equiv 2-$the number of
         connected components of $\Si \backslash \Si^\si$.
    \end{enumerate}

   These invariants satisfy:
   \begin{enumerate}
    \item $0\leq h\leq \tg+1$.
    \item For $k=0$, $h>0$ and $h\equiv \tg+1$ (mod $2$).
    \item For $k=1$, $0\leq h \leq \tg$.
   \end{enumerate}
  The above classification was realized already by Felix Klein (see
  e.g.  \cite{Klein}, \cite{W}, \cite{S}). This classification is
  probably better understood in terms of the quotient
  $Q(\Si)=\Si/\bra\si\ket$, where $\bra\si\ket=\{id, \si\}$ is the
  group generated by $\si$.  The quotient $Q(\Si)$ is orientable if
  $k=0$ and nonorientable if $k=1$, hence the name ``index of
  orientability''.  Furthermore, $h$ is the number of connected
  components of the boundary of $Q(\Si)$.  If $Q(\Si)$ is orientable,
  then it is topologically a sphere with $g\geq 0$ handles
  and with $h>0$ discs removed, and the invariants of
  $(\Si,\si)$ are $(\tg,h,k)=(2g+h-1,h,0)$.
  If $Q(\Si)$ is nonorientable, then it is topologically a sphere with
  $g>0$ crosscaps
  and with $h\geq 0$ discs removed, and the invariants of $\Si$ are
  $(\tg,h,k)=(g+h-1,h,1)$.

  From the above classification we see that symmetric Riemann surfaces of a
  given genus $\tg$ fall into $[\frac{3\tg +4}{2}]$ topological types.

 \section{Deformation theory of Bordered Riemann Surfaces}

In this section, we study deformation theory of bordered Riemann surfaces.
We refer to \cite[Section 3]{KL} for some preliminaries such as doubling 
constructions and the Riemann-Roch theorem for bordered Riemann surfaces.

\subsection{Deformation theory of smooth bordered Riemann surfaces}
\label{deform}

  Let $\Si$ be a bordered Riemann surface, $(\Si_\C,\si)$ be its complex double
  (see e.g. \cite[Section 3.3.1]{KL} for the definition).
  Analytically, $(\Si_\C,\si)$ is a compact symmetric Riemann surface.
  Algebraically, it is a smooth complex algebraic curve $X$ which is
  the complexification of some smooth real algebraic curve
  $X_0$, i.e., $X=X_0\times_\R\C$ (see \cite[Chapter II, Exercise 4.7]{H}).
  Alternatively, $(X,S)$ is a complex algebraic curve with a real
  structure (see \cite[I.1]{Sil}), where $S$ is a
  \emph{semi-linear} automorphism in the sense of
  \cite[Chapter II, Exercise 4.7]{H}
  which induces the antiholomorphic involution $\si$ on $\Si_\C$.

 \subsubsection{Algebraic approach}

  First order deformations of the complex algebraic curve $X$ is canonically
  identified with the complex vector space $\Exone$, where $\Omega_X^1$ is
  the sheaf of K\"{a}hler differentials on X. The obstruction lies in $\Extwo=0$.
  Similarly, the first order deformation of the real algebraic curve $X_0$
  is identified with the real vector space $\exone$, and the obstruction
  vanishes. We have  
$$
\Exone\cong\exone\otimes_{\R}\C
$$
  since $X=X_0\times_\R\C$.
  The semi-linear automorphism $S$ induces a complex conjugation
$$
S: \Exone\rightarrow\Exone.
$$
  The fixed locus $\Exone^S$ gives the
  first order deformation of $(X,S)$ as a complex algebraic curve with a
  real structure, and is naturally isomorphic to  $\exone$.

  More explicitly, $X$ can be covered by complex affine curves
  which is a complete intersection of hypersurfaces defined by
  polynomials with real coefficients. Deformation of $X$ are given by
  varying the coeffients (in $\C$). Deformation of $(X,S)$ is given by
  varying the coeffients in $\R$. The above polynomials with real coefficients
  also define the real algebraic curve $X_0$, and
  varying the coeffients in $\R$ gives the deformation of $X_0$.
  The complex conjugation of coefficients corresponds to the above
  complex conjugation S on $\Exone$.


  $X$ is a smooth algebraic variety, so $\Exone$ is isomorphic to
  the sheaf cohomology group $H^1(X,\Theta_X)$, where $\Theta_X$
  is the tangent sheaf of $X$, and 
$$
\Extwo\cong H^2(X,\Theta_X)=0.
$$
  Similarly, we have 
$$
\exone\cong H^1(X,\Theta_{X_{0}}),\ \
\extwo\cong H^2(X,\Theta_{X_{0}})=0.
$$

  We now return to the original bordered Riemann surface $\Si$.
\begin{df}
  Let $\cO_\Si$ be the sheaf of local holomorphic functions on $\Si$ 
  with real boundary values.
 
  Let $\Omega_\Si$ be a sheaf of $\cO_\Si$-modules, together
  with an $\R$ derivation $d:\cO_\Si\ra \Omega_\Si$, which
  satisfy the following universal property: for any sheaf of
  $\cO_\Si$-modules $\mathcal{F}$, and for any $\R$ derivation
  $d':\cO_\Si\ra\mathcal{F}$, there exists a unique
  $\cO_\Si$-module homomorphism $f: \Omega_\Si\ra\mathcal{F}$
  such that $d'=f\circ d$. We call $\Omega_\Si$ the {\em sheaf of
  K\"{a}hler differentials on $\Si$}. 

  Let $\Theta_\Si=\mathcal{H}om_{\cO_\Si}(\Omega_\Si,\cO_\Si)
  =\Omega_\Si^\vee$, the dual of $\Omega_\Si$ in the
  category of sheaves of $\cO_\Si$-modules. 
\end{df}

  Note that $\Omega_\Si$ and $\Theta_\Si$ are locally free sheaves of
  $\cO_\Si$-modules of rank $1$. Analytically,
  $\Omega_\Si$ is the sheaf of local holomorphic $1$-forms
  on $\Si$ whose restriction to $\bS$ are real $1$-forms,
  and $\Theta_\Si$ is the sheaf of holomorphic vector fields
  with boundary values in $T_{\bS}$. There are natural
  isomorphisms
$$
  \mathrm{Ext}^i_{\cO_\Si}(\Omega^1_\Si,\cO_\Si) \cong
   \mathrm{Ext}^i_{\cO_{X_0}}(\Omega^1_{X_0},\cO_{X_0}),\ \
  H^i(\Si,\Theta_\Si)\cong H^i(X_0, \Theta_{X_0})
$$
  for $i\geq 0$.
  Since the first order deformations of the bordered Riemann
  surface $\Si$, of the symmetric Riemann surface $(\Si_\C,\si)$,
  and of the real algebraic curve $X_0$ are identified,
  the first order deformation of $\Si$ is canonically
  identified with 
$$
\mathrm{Ext}_{\cO_\Si}^1(\Omega_\Si,\cO_\Si)\cong H^1(\Si,\Theta_\Si),
$$
and the obstruction lies in
$$
\mathrm{Ext}_{\cO_\Si}^2(\Omega_\Si,\cO_\Si)\cong H^2(\Si,\Theta_\Si)=0.
$$

 \subsubsection{Analytic approach}

  We first give the definitions of a {\em differentiable family of compact symmetric
  Riemann surfaces} and a {\em differentiable family of bordered Riemann
  surfaces}, which are modifications of \cite[Definition 4.1]{Ko}.

  \begin{df}
   Suppose given a compact symmetric Riemann surface $(M_t,\si_t)$
   for each point $t$ of a domain $B$ of $\R^m$. $\{(M_t,\si_t) | t\in B\}$ is
   called a {\em differentiable family of symmetric Riemann surfaces} if
   there are a differentiable manifold $\cM$, a surjective $C^\infty$ map
   $\pi:\cM\ra B$ and a  $C^\infty$ map
   $\si:\cM\ra \cM$
   such that
   \begin{enumerate}
    \item The rank of the Jacobian matrix of $\pi$ is equal to $m$ at every point
          of $\cM$.
    \item For each $t\in B$, $\pi^{-1}(t)$ is a compact connected subset of $\cM$.
    \item $\pi^{-1}(t)=M_t$.
    \item There are a locally finite open covering $\{\cU_i \mid i\in I\}$
          of $\cM$ and $C^\infty$ functions
          $z_i:\cU_i\ra \C$ such that
$$
          \{\left(\cU_i\cap\pi^{-1}(t),z_i|_{\cU_i\cap\pi^{-1}(t)}\right) \mid
              i\in I, \cU_i\cap\pi^{-1}(t)\neq \phi \}
$$
          is an analytic atlas for $M_t$.
    \item $\pi\circ\si=\pi$, and $\si|_{\pi^{-1}(t)}=\si_t:M_t\ra M_t$ is an
          antiholomorphic involution.

   \end{enumerate}
  \end{df}

 \begin{df} \label{borderfamily}
   Suppose given a bordered Riemann surface $M_t$
   for each point $t$ of a domain $B$ of $\R^m$. $\{M_t \mid t\in B\}$ is
   called a {\em differentiable family of bordered Riemann surfaces} if
   there are a differentiable manifold with boundary $\cM$ and
   a surjective $C^\infty$ map $\pi:\cM\ra B$ such that
   \begin{enumerate}
    \item The rank of the Jacobian matrix of $\pi$ is equal to $m$ at every point
          of $\cM$.
    \item For each $t\in B$, $\pi^{-1}(t)$ is a compact connected subset of
          $\cM$.
    \item $\pi^{-1}(t)=M_t$.
    \item There are a locally finite open covering $\{\cU_i \mid i\in I\}$
          of $\cM$ and $C^\infty$ functions
          $z_i:\cU_i\ra \C_+$ such that
$$
      \{\left(\cU_i\cap\pi^{-1}(t),z_i|_{\cU_i\cap\pi^{-1}(t)}\right) \mid
              i\in I, \cU_i\cap\pi^{-1}(t)\neq \phi \}
$$
          is an analytic atlas for $M_t$.
    \end{enumerate}

  \end{df}

  The complex double $\Si_{\C}$ of a bordered Riemann surface $\Si$
  is a complex manifold of dimension 1.  Infinitesimal deformation
  of $\Si_{\C}$ can be identified with
$$
H^1(\Si_\C,T_{\Si_\C})=H^1(X,\Theta_X),
$$ 
  and the obstruction lies in  
$$
H^2(\Si_\C,T_{\Si_\C})=H^2(X,\Theta_X)=0.
$$ 
  (See \cite{Ko}.) The differential $d\si$ of $\si$ is an antiholomorphic 
  involution on the holomorphic line bundle $T_{\Si_\C}\ra \Si_\C$
  which covers $\si:\Si_\C\ra \Si_\C$. $d\si$ induces a
  complex conjugation
$$
  \tsi:H^1(\Si_\C,T_{\Si_\C})\rightarrow H^1(\Si_\C,T_{\Si_\C})
$$
  which gets identified with the action of $S$ on
  $\mathrm{Ext}_{\cO_X}^1(\Omega_X,\cO_X)$ under the isomorphism
$$
H^1(\Si_\C,T_{\Si_\C})\cong\Exone.
$$

The pair $(T_{\Si_\C},d\si)$ is the {\em holomorphic complex double} 
of the {\em Riemann-Hilbert bundle} $(T_{\Si},T_{\bS})\ra (\Si,\bS)$,
where a Riemann-Hilbert bundle and its holomorphic complex double are
defined in \cite[Section 3.3.4]{KL}.
There is an isomorphism (see \cite[Section 3.4]{KL})
$$
  H^1(\Si_\C,T_{\Si_\C})^{\tsi}\cong H^1(\Si,\bS, T_{\Si},T_{\bS}).
$$
   From the above discussion, we know that
  $H^1(\Si_\C,T_{\Si_\C})^{\tsi}$ gives first order deformation
  of the symmetric Riemann surfaces $(\Si_\C,\si)$. Given a
  differential family of symmetric Riemann surface
  $(\tilde{\Si}_t,\si_t)$ such that
  $(\tilde{\Si}_0,\si_0)=(\Si_\C,\si)$,
  $K_t=\tilde{\Si}_t/\bra\si_t\ket$ is a family of Klein surfaces.
  Each $K_t$ is homeomorphic to $\Si$, which is orientable,
  so it admits two analytic structures compatible
  with its dianalytic structure, and one is the complex conjugate
  of the other. We get two differentiable families
  $M_t, M'_t=\Bar{M_t}$ of bordered
  Riemann surfaces, one is a deformation of $\Si$, and the other
  is a deformation of $\Bar{\Si}$. So $H^1(\Si,\bS, T_{\Si},T_{\bS})$
  should give infinitesimal deformations of $\Si$.

  Recall that an infinitesimal deformation of $\Si_\C$ determines a
  \v{C}ech 1 cocycle in
  $H^1(\tilde{\cA},\Theta_\C)\subset H^1(\Si_\C,T_{\Si_\C})$,
  where $\tilde{\cA}$ is an analytic atlas of $\Si_\C$, and
  $\Theta_\C$ is the sheaf of local holomorphic vector fields on
  $\Si_\C$ \cite{Ko}.
  (The inclusion is an isomorphism if $\tilde{\cA}$ is acyclic).
  Following argument similar to that in \cite{Ko}, we now show that
  an infinitesimal deformation of $\Si$ determines
  a \v{C}ech 1 cocycle in $H^1(\cA,\Theta_\Si)$,
  where $\cA$ is an analytic atlas of $\Si$, and
  $\Theta_{\Si}$ is the sheaf of local holomorphic vector
  fields on $T_{\Si}$ with boundary values in $T_{\bS}$.

  Let $\{M_t\mid t\in B\}$ be a differentiable family of bordered Riemann
  surfaces, $M_0=\Si$. We use the notation in Definition~\ref{borderfamily}.
  Then 
$$
\cA=\{(U_i,\phi_i)\equiv \left(\cU_i\cap\pi^{-1}(0),
              z_i|_{\cU_i\cap\pi^{-1}(t)}\right)\mid 
              i\in I, \cU_i\cap\pi^{-1}(t)\neq \phi \}
$$
  is an analytic atlas of $\Si$. Without loss of generality, we may assume
  that $\cA$ is acyclic.
  We define t-dependent transition functions $f_{ij}$ by
  $z_i=f_{ij}(z_j,t)=f_{ik}(z_k,t)$. Then
  $f_{ij}(z_j,t)\in\R$ if $z_j\in \R$ by 4. of
  Definition~\ref{borderfamily}.
  \begin{eqnarray*}
   z_i&=&f_{ik}(z_k,t)=f_{ij}(f_{jk}(z_k,t),t)\\
   \frac{\pa f_{ik}}{\pa t}&=&
   \frac{\pa f_{ij}}{\pa z_j}\frac{\pa f_{jk}}{\pa t}
   +\frac{\pa f_{ij}}{\pa t}\\
  \end{eqnarray*}

  Multiplying by $\frac{\pa}{\pa z_i}$ and noting that
  $\frac{\pa f_{ij}}{\pa z_j}=\frac{\pa z_i}{\pa z_j}$,
  we have
\[
  \frac{\pa f_{ik}}{\pa t}\frac{\pa}{\pa z_i}
 =\frac{\pa f_{jk}}{\pa t}\frac{\pa}{\pa z_j}
  +\frac{\pa f_{ij}}{\pa t} \frac{\pa}{\pa z_i}
\]

  Let $(x_i,y_i)$ be real coordinates defined by $z_i=x_i+iy_i$,
  then the boundary is defined by $\{y_i=0\}$, and the tangent line
  to the boundary is spanned by $\frac{\pa}{\pa x_i}$.
  Under the isomorphism $T_\Si \to T_{\Si}^{0,1}$ given by
  $v \mapsto (v-iJv)/2$, we have
  $\frac{\pa}{\pa x}\mapsto \frac{\pa}{\pa z}$ and
  $\frac{\pa}{\pa y}\mapsto i\frac{\pa}{\pa z}$.
  So $\theta_{ik}\equiv\left. \frac{\pa f_{ik}}{\pa t}\right|_{t=0}
  \frac{\pa}{\pa z_i}$ defines a \v{C}ech 1 cochain in
  $C^1(\mathcal{A},\Theta_\Si)$ which satisfies the
  cocycle condition $\theta_{ik}=\theta_{ij}+\theta_{jk}$.

  By exactly the same argument in \cite{Ko} we see that another system of
  coordinates will give rise to a \v{C}ech 1 cocycle
  $\theta'=\theta +\delta \alpha$, where $\alpha$ is a  \v{C}ech
  0 cochain. Therefore, the infinitesimal deformation of $\Si$ is given by
$$
   H^1(\cA,\Theta_{\Si})=H^1(\Si,\bS,T_{\Si},T_{\bS}).
$$

  \subsection{Nodal bordered Riemann surfaces}
  \label{stabbord}

   To compactify the moduli of bordered Riemann surface, we will allow
   nodal singularities. The complex double of a bordered Riemann surface
   is a complex algebraic curve with real structure. The stable
   compactification of moduli of such curves parametrizes stable complex
   algebraic curves with real structure \cite{SS, S}, or equivalently,
   stable compact symmetric Riemann surfaces. Naively, the quotient
   of a stable compact symmetric Riemann surface by its antiholomorphic
   involution will give rise to a ``stable bordered Riemann surface''.
   We will make this idea precise in this section.

   Let $\Si$ be a (smooth) bordered Riemann surface of type $(g,h)$.
   Note that if $\phi:\Si\to \Si$ is an automorphism (Definition~\ref{bordermor}),
   then its {\em complex double} (\cite[Section 3.3.2]{KL}) 
   $\phi_\C:\Si_\C\to\Si_\C$ is an automorphism of $(\Si_\C,\si)$
   (Definition~\ref{symmor}). This gives an inclusion
   $\Aut(\Si) \subset \Aut(\Si_\C,\si)$. It is easy to see that
   the following are equivalent:
   \begin{itemize}
    \item $\Si$ is stable, i.e., $\mathrm{Aut}(\Si)$ is finite.
    \item $\Si_\C$ is stable.
    \item The genus $\tg=2g+h-1$ of $\Si_\C$ is greater than one.
    \item The Euler characteristic $\chi(\Si)=2-2g-h$ of $\Si$ is negative.
   \end{itemize}

   We start with $\tg=2$. Let $\Bar{M}_2$ be the moduli of
   stable complex algebraic curves of genus $2$.
   The strata of $\Bar{M}_2$ are shown in Figure 1.

\begin{figure}\label{one}
\includegraphics[scale=0.5]{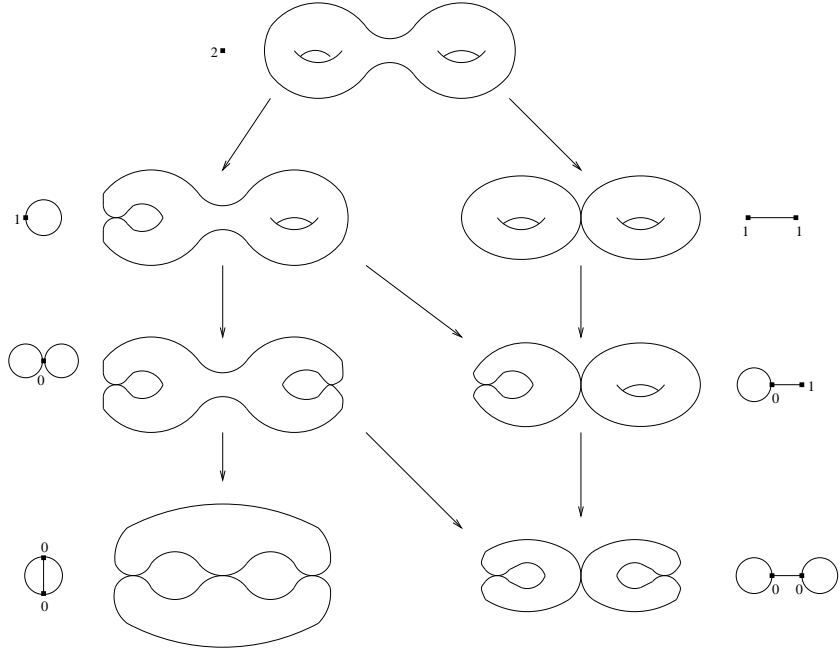}
 \caption{strata of $\Bar{M}_2$}
\end{figure}

 The dual graph and the underlying topological surface of a prestable curve
 are shown. In the dual graph of the curve $C$, each vertex corresponds to an
 irreducible component of $C$, labeled by the genus of the normalization,
 while each edge corresponds to a node of $C$, whose two end points correspond to
 the two irreducible components which intersect at this node.

If  $\tg=2$, $(g,h)$ can be $(0,3)$ (Figure 2) or $(1,1)$
(Figure 3). 

\begin{figure} \label{two}
\begin{center}
\includegraphics[scale=0.5]{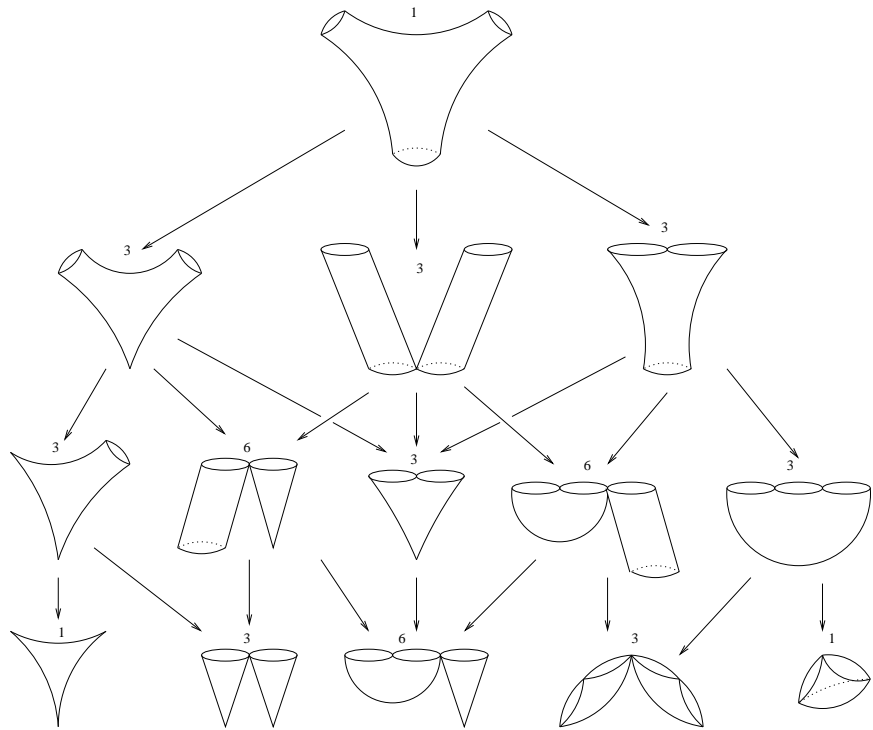}
\end{center}
\caption{strata of $\Bar{M}_{(0,3)}$. Each stable bordered
Riemann surface above represents a topological type, and
the number above it is the number of strata associated to it. 
Strata associated to the same topological type 
are related by relabelling the three boundary circles. There are 
one 3-dimensional stratum, nine 2-dimensional
strata, twenty-one 1-dimensional strata, and fourteen 0-dimensional
strata. We will see in Example \ref{ex:pants} that $\Bar{M}_{(0,3)}$ can be
identified with the associahedron $K_5$ defined by J. Stasheff \cite{St}.}
\end{figure}

\begin{figure}\label{three}
\begin{center}
\includegraphics[scale=0.42]{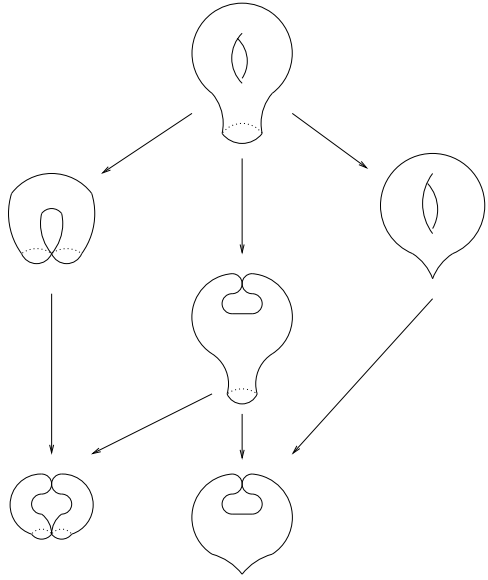}
\end{center}
\caption{strata of $\Bar{M}_{(1,1)}$}
\end{figure}

   \begin{df}
    Let $(x,y)$ be coordinates on $\C^2$, and $A(x,y)=(\bar{x},\bar{y})$ be
    the complex conjugation. A \emph{node} on a bordered Riemann surface
    is a singularity isomorphic to one of the following:
   \begin{enumerate}
    \item $(0,0)\in\{xy=0\}$ (interior node)
    \item $(0,0)\in\{x^2+y^2=0\}/A$ (boundary node of type E)
    \item $(0,0)\in\{x^2-y^2=0\}/A$ (boundary node of type H)
   \end{enumerate}

    A \emph{nodal bordered Riemann surface} is a singular bordered Riemann
    surface whose singularities are nodes.
   \end{df}

   A type E boundary node on a bordered Riemann surface corresponds
   to a boundary component shrinking to a point, while a type H
   boundary node corresponds to a boundary component intersecting
   itself or another boundary component.
   The boundary of a nodal Riemann surface is a union of points
   and circles, where one circle might intersect other circles in
   finitely many points.

   The notion of morphisms and complex doubles can be easily extended to
   nodal bordered Riemann surfaces. The complex double of a nodal bordered
   Riemann surface is a nodal compact symmetric Riemann surface.

   \begin{ex}
    Consider $\Ce=\{X^2-Y^2+\epsilon Z^2=0\}\subset\bP^2$, where $[X,Y,Z]$ are
    homogeneous coordinates on $\bP^2$, and $\epsilon\in\R$. $\Ce$ is
    invariant under the standard complex conjugation
    $A([X,Y,Z])=[\bar{X},\bar{Y},\bar{Z}]$ on $\bP^2$, so $\se=A|_{\Ce}$ is
    an antiholomorphic involution on $\Ce$. For $\epsilon\neq 0$, $(\Ce,\se)$
    is a symmetric Riemann surface of type $(0,1,0)$, and $\Ce/\bra\se\ket$ is the
    disc, which is a bordered Riemann surface. $C_0$ has two irreducible components
    $\{X+Y=0\}$ and $\{X-Y=0\}$ which are projective lines, and the
    intersection point $[0,0,1]$ of the two lines is a node on $C_0$.
    Both lines are invariant under the antiholomorphic involution $\si_0$.
    $C_0/\bra\si_0\ket$ is a nodal bordered Riemann surface: it is the union of two discs
    whose intersection is a boundary node of type H.
   \end{ex}

   \begin{ex}
    Consider $\Ce=\{X^2+Y^2+\epsilon Z^2=0\}\subset\bP^2$, where
    $\epsilon\in\R$.  $\Ce$ is invariant under the complex conjugation $A$ on
    $\bP^2$, so $\se=A|_{\Ce}$ is an antiholomorphic involution on $\Ce$.
    Set $\Se=\Ce/\bra\se\ket$. For $\epsilon>0$, $(\Ce,\se)$ is a symmetric Riemann
    surface of type $(0,0,1)$, and $\Se$ is the real projective plane; for
    $\epsilon<0$, $(\Ce,\se)$ is a symmetric Riemann surface of type
    $(0,1,0)$, and $\Se$ is the disc. $C_0$ has two irreducible components
    $\{X+\sqrt{-1}Y=0\}$ and $\{X-\sqrt{-1}Y=0\}$ which are  projective lines, 
    and their intersection point $[0,0,1]$ is a node on $C_0$. The antiholomorphic 
    involution $\si_0$ interchanges the two irreducible components of $C_0$ and leaves
    the node invariant, so $\Si_0\cong\PP$, which is a smooth Riemann surface
    without boundary. However, we would like to view it as a disc whose
    boundary shrinks to a point which is a boundary node of type E.
   \end{ex}

   \begin{df} \label{norm}
    Let $\Si$ be a nodal bordered Riemann surface. The antiholomorphic
    involution $\si$ on its complex double $\Si_\C$ can be lifted to
    $\hsi:\widehat{\Si_\C}\rightarrow\widehat{\Si_\C}$, where
    $\widehat{\Si_\C}$ is the normalization of $\Si_\C$ (viewed as a
     complex algebraic curve).
    The \emph{normalization} of $\Si=\Si_\C/\bra\si\ket$ is defined to be
    $\hS=\widehat{\Si_\C}/\bra\hsi\ket$.
   \end{df}

    From the above definition, the complex double of the normalization
    is the normalization of the complex double, i.e.,
    $\hS_\C=\widehat{\Si_\C}$.

  Let $\Si$ be a smooth bordered Riemann surfaces of type $(g,h)$.
  The following are possible degenerations of $\Si$ whose only
  singularity is a boundary node.

  \begin{enumerate}
   \item[E.] One boundary component shrinks to a point. The normalization
            is a smooth bordered Riemann surface of type $(g,h-1)$
            (Figure 4).

\begin{figure}\label{four}
\begin{center}
\psfrag{(g,h)=(1,3)}{$(g,h)=(1,3)$}
\psfrag{(g,h)=(1,2)}{$(g,h)=(1,2)$}
\includegraphics[scale=0.7]{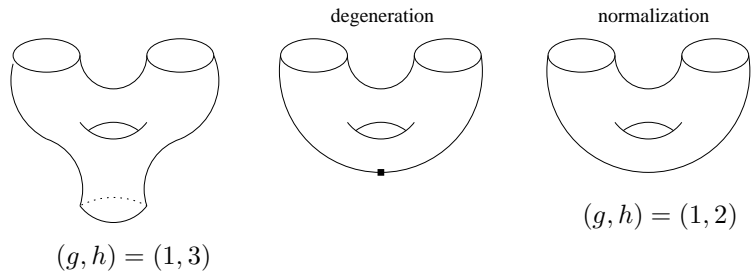}
\end{center}
\caption{boundary node of type E}
\end{figure}

   \item[H1.] Two boundary components intersect at one point. The normalization
              is a smooth bordered Riemann surface of type $(g,h-1)$, and the
              two preimages of the node are on the same boundary component
              (Figure 5).

\begin{figure}\label{five}
\begin{center}
\psfrag{(g,h)=(1,3)}{$(g,h)=(1,3)$}
\psfrag{(g,h)=(1,2)}{$(g,h)=(1,2)$}
\includegraphics[scale=0.7]{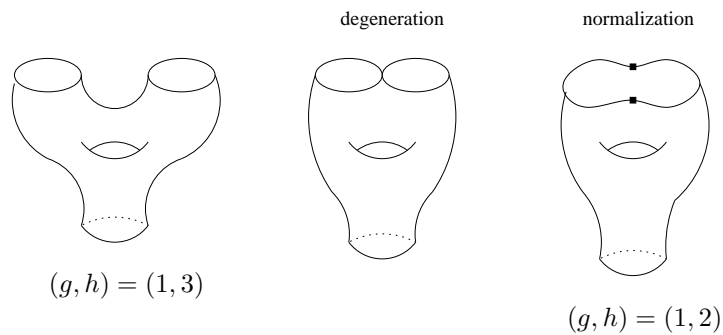}
\end{center}
\caption{boundary node of type H1}
\end{figure}

\item[H2.] One boundary component intersects itself, and the normalization
              of the surface is connected. The normalization is a smooth
              bordered Riemann surface of type $(g-1,h+1)$, and the two
              preimages of the node are on different boundary components
              (Figure 6).

\begin{figure}\label{six}
\begin{center}
\psfrag{(g,h)=(1,3)}{$(g,h)=(1,3)$}
\psfrag{(g,h)=(0,4)}{$(g,h)=(0,4)$}
\includegraphics[scale=0.6]{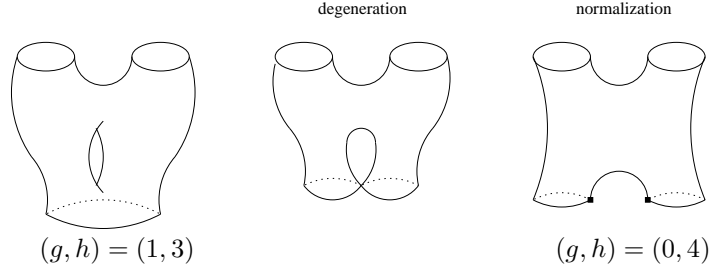}
\end{center}
\caption{boundary node of type H2}
\end{figure}

 \item[H3.] One boundary component intersects itself, and the normalization
              of the surface is disconnected. The normalization is a disjoint
              union of two smooth bordered Riemann surfaces of types
              $(g_1,h_1)$ and $(g_2,h_2)$ such that $g=g_1+g_2$ and
              $h=h_1+h_2-1$, and each connected component contains
              one of the two preimages of the node
              (Figure 7).

\begin{figure}\label{seven}
\psfrag{(g,h)=(1,3)}{$(g,h)=(1,3)$}
\psfrag{(g1,h1)=(1,2)}{$(g_1,h_1)=(1,2)$}
\psfrag{(g2,h2)=(0,2)}{$(g_2,h_2)=(0,2)$}
\includegraphics[scale=0.6]{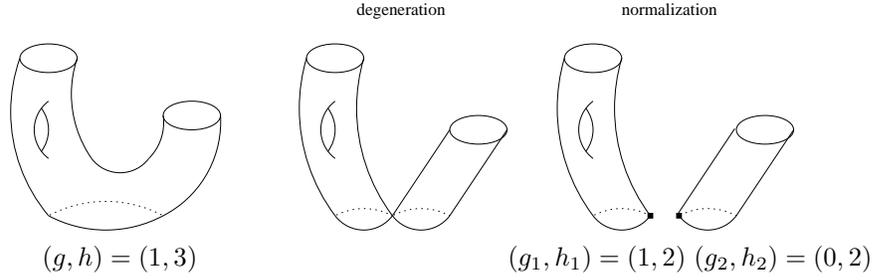}
\caption{boundary node of type H3}
\end{figure}

\end{enumerate}

 \begin{df}
  A \emph{prestable} bordered Riemann surface is either a smooth bordered
  Riemann surface or a nodal bordered Riemann surface.
 \end{df}

  Let $\Si$ be a prestable bordered Riemann surface, $\hat{\Si}$ be its
  normalization. Let $\hat{C}_1,\ldots,\hat{C}_\nu$, $\hat{\Si}_1,\ldots,
  \hat{\Si}_{\nu'}$ be the connected components of $\hat{\Si}$, where
  $\hat{C}_i$ is a smooth Riemann surface of genus $\hat{g}_i$, 
  and $\hat{\Si}_{i'}$ is a smooth bordered Riemann surface
  of type $(g_{i'},h_{i'})$. Let $\delta$ be the number of connecting
  interior nodes (Figure 8), and $\delta_E, \delta_{H1}, \delta_{H2},
   \delta_{H3}$ be the numbers of boundary nodes described in
  $E, H1, H2, H3$, respectively.

\begin{figure}\label{eight}
\psfrag{(g,h)=(2,1)}{$(g,h)=(2,1)$}
\includegraphics[scale=0.6]{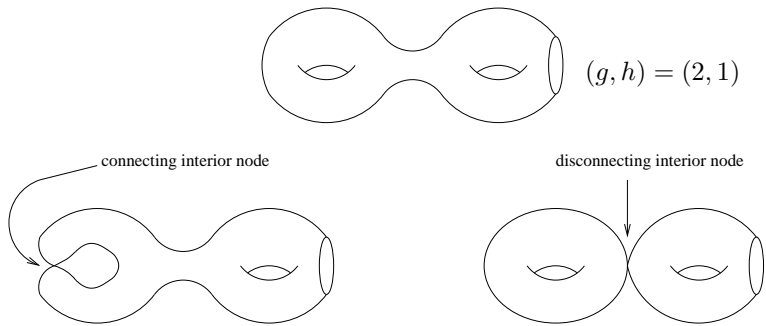}
\caption{connecting and disconnecting nodes} 
\end{figure} 

 The topological type $(g,h)$
  of $\Si$ is given by
  \begin{eqnarray*}
   g&=&\tg_1+\cdots+\tg_\nu+g_1+\cdots+g_{\nu'}+\delta+\delta_{H2}\\
   h&=&h_1+\cdots+h_{\nu'}+\delta_E+\delta_{H1}-\delta_{H2}
       -\delta_{H3}
  \end{eqnarray*}

 It is now straightforward to extend the notion of marked bordered Riemann
 surfaces to prestable bordered Riemann surfaces.

 \begin{df}\label{mark}
   Let $h$ be a positive integer, $g,n$ be nonnegative integers, and
   $\vm=(m^1,\ldots,m^h)$ be an $h$-uple of nonnegative integers.
   A \emph{prestable marked bordered Riemann surface of type $(g,h)$ with
   $(n,\vm)$ marked points} is an $(h+3)$-uple
$$\domain$$
whose components are described as follows.
   \begin{itemize}
    \item $\Si$ is a prestable bordered Riemann surface of type $(g,h)$.
    \item $\bB=(B^1,\ldots,B^h)$, where $\bS=\bigcup_{i=1}^h B^i$, and 
          each $B^i$ is an immersed circle. 
          The  circles $B^1,\ldots, B^h$ may intersect with each other
          at boundary nodes, and become $h$ disjoint embedded circles 
          under smoothing of all boundary nodes.
    \item $\bp=(p_1,\ldots,p_n)$ is an $n$-uple of distinct smooth points 
          in $\inS$.
    \item $\bq^i=(q^i_1,\ldots, q^i_{m^i})$ is an $m^i$-uple of distinct 
          smooth points on $B^i$.
   \end{itemize}
   \end{df}

   \begin{df}\label{markmor}
    A \emph{morphism} between prestable marked bordered Riemann surfaces of type
    $(g,h)$ with $(n,\vec{m})$ marked points
$$
\domain \to \another
$$
    is an isomorphism of prestable bordered Riemann surfaces $\phi:\Si\ra\Si'$ such that
    $\phi(B^i)=(B')^i$ for $i=1,\ldots,h$, $\phi(p_j)=p'_j$ for $j=1,\ldots,n$, and
    $\phi(q^i_k)=(q')^i_k$ for $k=1,\ldots,m^i$.
   A prestable marked bordered Riemann surfaces of type $(g,h)$ with
    $(n,\vm)$  marked points is \emph{stable} if its automorphism group is finite.
   \end{df}

   \begin{rem}
    The category of prestable marked bordered Riemann surfaces of type $(g,h)$
    with $(n,\vm)$ marked points is a groupoid since every morphism in
    Definition~\ref{markmor} is an isomorphism.
   \end{rem}

\begin{ex}
Consider the case $(g,h)=(0,2), n=0, \vm=(1,0)$ (annuli with one boundary
marked point). The moduli space
$\Bar{M}_{(0,2)(0,(1,0))}$ is an interval $[0,1]$.
There are three strata: $t\in (0,1), t=0, t=1$ (Figure 9)
\begin{figure}\label{nine}
\begin{center} 
\psfrag{t=0}{$t=0$}
\psfrag{t=1}{$t=1$}
\psfrag{0<t<1}{$t\in (0,1)$}
\includegraphics[scale=0.6]{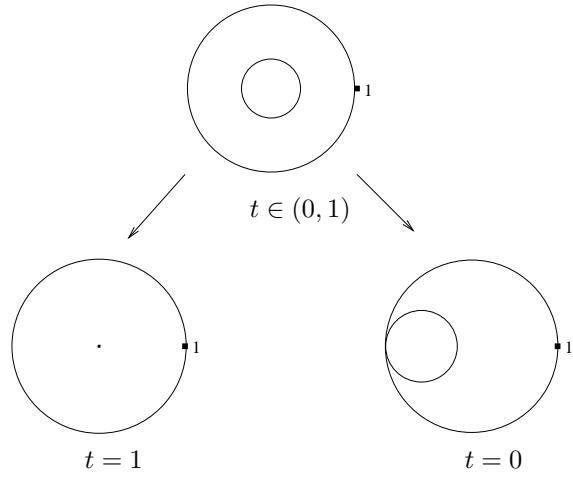}
\caption{strata of $\Bar{M}_{(0,2)(0,(1,0))}$} 
\end{center}
\end{figure} 
\end{ex}

\begin{ex}
Consider the case $(g,h)=(0,2), n=0, \vm=(2,0)$ (annuli with two boundary
marked points on the same boundary circle). The moduli space
$\Bar{M}_{(0,2)(0,(2,0))}$ is a pentagon. 
There are eleven strata (Figure 10).
\begin{figure}\label{ten}
\begin{center} 
\includegraphics[scale=0.6]{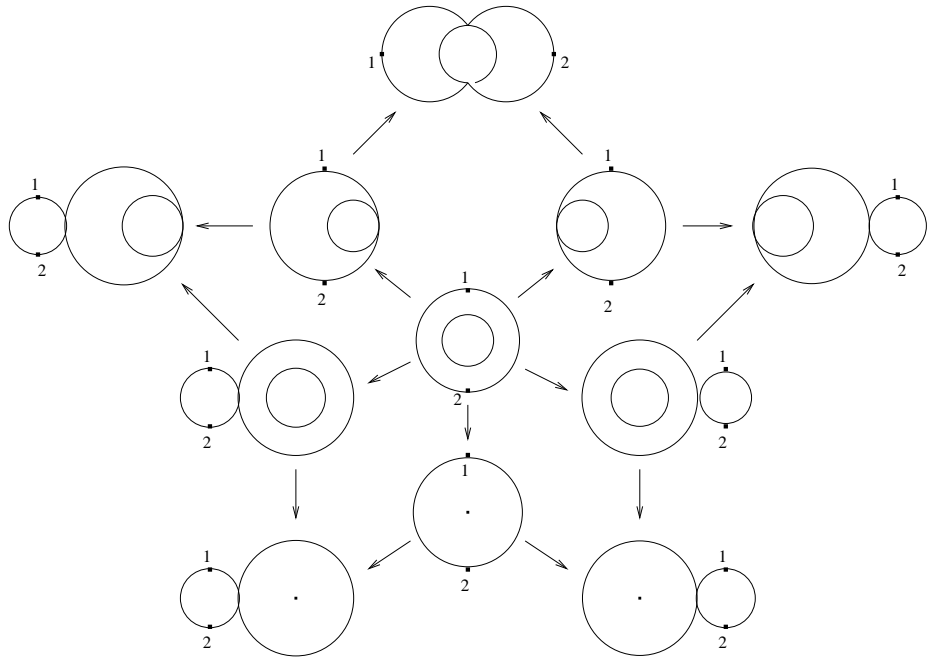}
\caption{strata of $\Bar{M}_{(0,2)(0,(2,0))}$} 
\end{center}
\end{figure} 
\end{ex}

\begin{ex}
Consider the case $(g,h)=(0,2), n=0, \vm=(1,1)$
(annuli with one marked point on each boundary circle). 
The moduli space $\Bar{M}_{(0,2)(0,(1,1))}$ is 
a disc $\{z\in\mathbb{C}\mid |z|\leq 1\}$. 
There are four strata: $0<|z|<1$, $|z|=1$ but $z\neq 1$,
$z=1$, $z=0$ (Figure 11)

\begin{figure}\label{eleven}
\begin{center}
\psfrag{0<z<1}{$0<|z|<1$}
\psfrag{z>1}{$|z|=1, z\neq 1$}
\psfrag{z=0}{$z=0$}
\psfrag{z1}{$z=1$}
\includegraphics[scale=0.6]{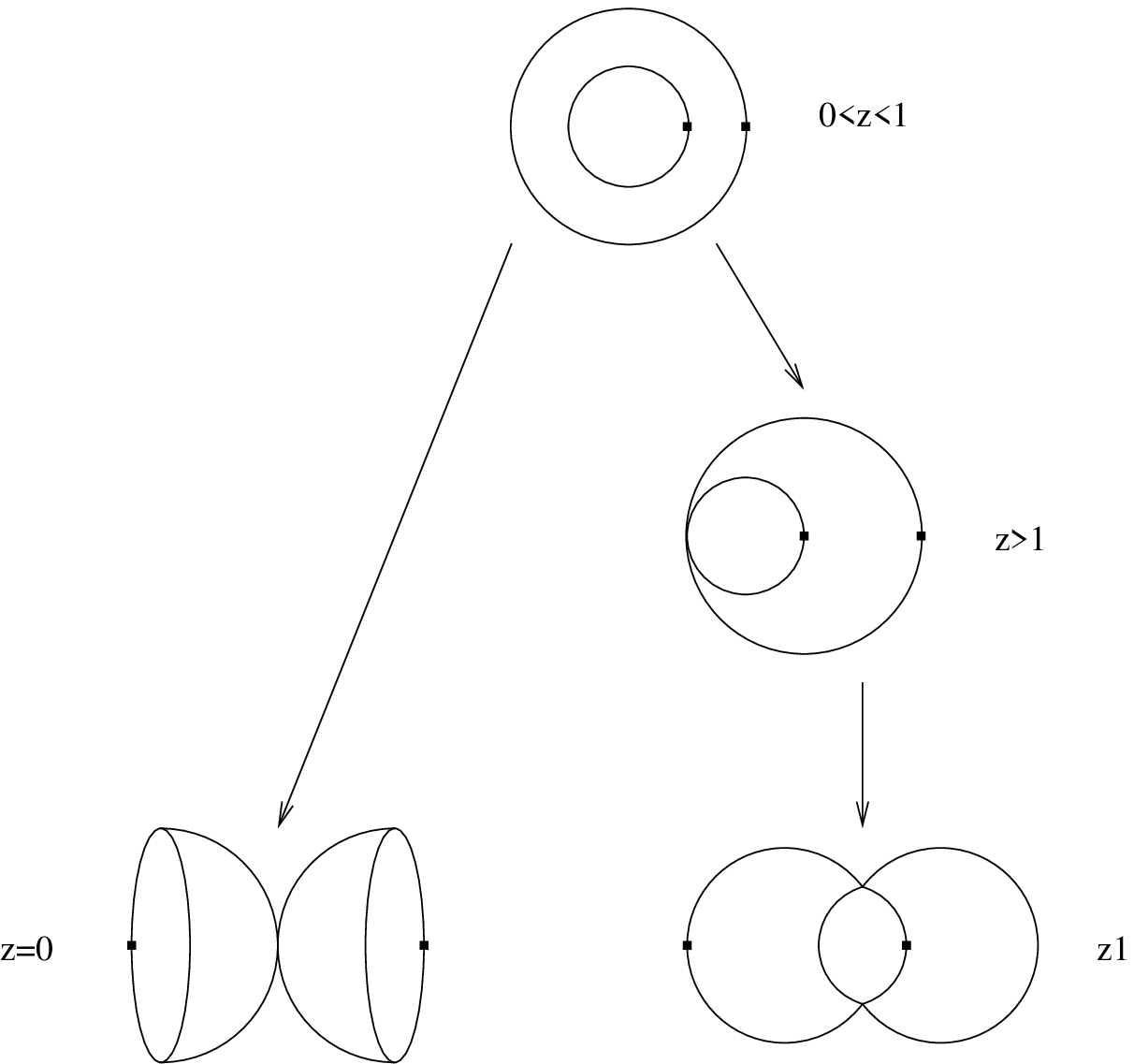}
\caption{strata of $\Bar{M}_{(0,2)(0,(1,1))}$} 
\end{center}
\end{figure} 
\end{ex}

\subsection{Deformation theory for prestable bordered Riemann surfaces}
 \label{nodaldeform}

  The algebraic approach of deformation theory for smooth bordered Riemann
  surfaces in Section~\ref{deform} can be easily extended to nodal bordered Riemann
  surfaces. We will also consider marked points.

  Let $\domain$ be a marked prestable bordered Riemann surface of type
  $(g,h)$ with $(n,\vm)$ marked points. We want to study its infinitesimal
  deformation. The ordering of the boundary circles is irrelevant to
  the infinitesimal deformation, so in this section we will ignore it
  and write $\bord$, where 
$$
\bq=(q_1,\ldots,q_{m})=
(q^1_1,\ldots,q^1_{m^1}, q^2_1,\ldots,q^2_{m^2},\ldots, q^h_1,\ldots,q^h_{m^h}),
$$
and $m=m^1+\cdots + m^h$.

  The complex double $(\Si_\C,\si)$ of $\Si$ is a nodal complex
  algebraic curve $X$ which is the complexification of some real algebraic
  curve $X_0$, i.e., $X=X_0\times_\R\C$. Alternatively, $(X,S)$ is a
  complex algebraic curve with a real structure, where $S$ is a semi-linear
  automorphism which induces the antiholomorphic involution $\si$ on $\Si_\C$.

  Algebraically, the complex double of $\bord$ is a nodal complex algebraic
  curve with $(2n+m)$ marked points $(X,\bx)$, where 
$$
\bx=(x_1,\ldots,x_{2n+m})=(p_1,\ldots,p_n, \bar{p}_1,\ldots,\bar{p}_n,q_1,\ldots,q_m).
$$
Here we identify $\Si$ with the image under the inclusion $i:\Si\to \Si_\C$, and
denote $\si(p)$ by $\bar{p}$.

Let 
$$
D_\bx=x_1+\cdots+x_{2n+m}
$$
be the divisor in $X$ associated to $\bx$.
The set of first order deformation of the pointed complex algebraic curve
$(X,\bx)$ is canonically identified with the complex vector space
\[ 
  \pone, 
\] 
  and the obstruction lies in 
\[
\ptwo. 
\]

We claim that $\ptwo=0$. Three terms in the local to global
  spectrum sequence contribute to $\ptwo$:
\begin{eqnarray*}
  &&  H^0(X,\ltwo),\\
  &&  H^1(X,\lone),\\ 
  && H^2(X,\lzero).
\end{eqnarray*}
The curve $X$ has only nodal singularities, so
\[ 
  \ltwo=0, 
\] 
thus 
\[ 
   H^0(X,\ltwo)=0.
\] 
The sheaf $\lone$ is supported on nodes, so
\[ 
H^1(X,\lone)=0. 
\] 
Finally, 
\[ H^2(X,\lzero)=0\] since $X$ is one dimensional.
So we get the desired vanishing. 

  The semi-linear automorphism $S$ induces a
  complex conjugation 
\[ 
   S: \pone\rightarrow\pone. 
\] 
   The fixed locus $\pone^S$ gives the first order deformation of 
   $(X,\bx,S)$ as a pointed complex algebraic curve with a real structure.

  We will study the group $\pone$ and the action of $S$ on it more closely.
  We first introduce some notation.

  Let $C_1, ..., C_\nu$ be the irreducible components of $\Si$ which
  are (possibly nodal) Riemann surfaces, and let $\Si_1, ..., \Si_{\nu'}$
  be the remaining  irreducible components of $\Si$, which are (possibly
  nodal) bordered Riemann surfaces.
  Then the irreducible components of $X$ are
\[
  C_1,\ldots, C_\nu,\ \Bar{C}_1,\ldots,\Bar{C}_\nu,\
  (\Si_1)_\C,\ldots,(\Si_{\nu'})_\C.
\]

  Let $\hC_i$ denote the normalization of $C_i$, $i=1,\ldots,\nu$, and
  let $\hS_{i'}$ denote the normalization of $\Si_{i'}$,
  $i'=1,\ldots,\nu'$. Then
\[
  \hC_1,\ldots, \hC_\nu,\ \hS_1,\ldots,\hS_{\nu'}
\]
  are the connected components of $\hat{\Si}$, the normalization
  of $\Si$, and
\[
  \hC_1,\ldots, \hC_\nu,\
  \hat{\Bar{C}}_1,\ldots,\hat{\Bar{C}}_\nu,\
  (\hS_1)_\C,\ldots,(\hS_{\nu'})_\C.
\]
  are the connected components of $\hat{X}$, the normalization of $X$.

Let $r_1,\ldots,r_{l_0}\in\inS$ be interior nodes of $\Si$, and
    $s_1,\ldots,s_{l_1}\in\bS$ be boundary nodes of $\Si$.
Then $X$ has $2l_0+l_1$ nodes $r_1,\ldots,r_{l_0},\
\bar{r}_1,\ldots,\bar{r}_{l_0},\ s_1,\ldots,s_{l_1}$.

Let $\hp_j\in\hS$ be the preimage of $p_j$ under the
normalization map $\pi: \hat{X}\ra X, j=1,\ldots,n$.
Define $\hat{\bar{p}}_j, \hat{q}_{j'}$ similarly.
Let $\hr_\alpha, \hr_{l_0+\alpha}$ be the preimages of
$r_\alpha, \alpha=1,\ldots, l_0$, and define
$\hs_{\alpha'}, \hs_{l_1+\alpha'}$ similarly.

Consider $\hat{X}$ with marked points
$$
  \hat{\bx}=(\hp_1,\ldots,\hp_n,\
  \hat{\bar{p}}_1,\ldots,\hat{\bar{p}}_n,\
  \hq_1,\ldots, \hq_m,\
  \hr_1,\ldots,\hr_{2l_0},\
  \hat{\bar{r}}_1,\ldots,\hat{\bar{r}}_{2l_0},\
  \hs_1,\ldots,\hs_{2l_1}),
$$
  which can be written as a disjoint union of pointed curves
$$
(\hC,\by^i),\ \  (\hat{\Bar{C}},\bar{\by}^i),\ \,
((\hS_{i'})_\C,\bz^{i'})= ((\hS_{i'})_\C,\bp^{i'},\bar{\bp}^{i'},\bq^{i'}),
$$
where
\begin{eqnarray*}
&& \by^i=(y^i_1,\ldots,y^i_{\tn_i}),\ \  
 \bar{\by}^i=(y^i_1,\ldots,\bar{y}^i_{\tn_i}),\\
&& \bp^{i'}=(p^{i'}_1,\ldots,p^{i'}_{n_{i'}}),\ \
   \bar{\bp}^{i'}= (\bar{p}^{i'}_1,\ldots,\bar{p}^{i'}_{n_{i'}}),\ \
    \bq^{i'}=(q^{i'}_1,\ldots,q^{i'}_{m_{i'}}),
\end{eqnarray*}
$i=1,\ldots,\nu$, $i'=1,\ldots,\nu'$,
$p^{i'}_1,\ldots,p^{i'}_{n_{i'}}\in \hS_{i'}^\circ$,
$q^{i'}_1,\ldots,q^{i'}_{m_{i'}}\in \pa \hS_{i'}$,  and
$$
   \sum_{i=1}^\nu \tn_i+\sum_{i'=1}^{\nu'}n_i=n+2l_0,\ \
   \sum_{i=1}^\nu m_i=m+2l_1.
$$

  From the local to global spectrum sequence, we have an exact sequence
\begin{eqnarray*}
   0 &\to & H^1(X,\lzero)\to \pone \\
     &\to & H^0(X,\lone)\to 0,
\end{eqnarray*}
where
$$
\lzero=\lhom=\Omega_X^1(D_\bx)^{\vee}.
$$

 We have the following elementary fact:

\begin{lm}
\[
  \Omega_X^1(D_\bx)^{\vee}=\pi_*T_{\hat{X}}(-D_{\hat{\bx}}),
\]
  where $\pi:\hat{X}\ra X$ is the normalization map, and $D_{\hat{\bx}}$ is the divisor
  corresponding to the marked points $\hat{\bx}$ in $\hat{X}$.
\end{lm}
 \paragraph{Proof.} The equality obviously holds for smooth points. It
 suffices to show that $\Omega^1_Y=\pi_* T_{\hat{Y}}$ for
 $Y=\mathrm{Spec}\C[x,y]/(xy)$, which follows from a local
 calculation. $\Box$

\bigskip

The map $\pi:\hat{X}\to X$ is an affine morphism, so by 
\cite[Chapter III, Exercise 4.1]{H} we have
$$
 H^1(X,\lzero)= H^1(X,\pi_*T_{\hat{X}}(-D_{\hat{\bx}})) 
   =\bigoplus_{i=1}^\nu(W_i\oplus\Bar{W}_i)\oplus\bigoplus_{i'=1}^{\nu'}W'_{i'}
$$
where the vector spaces
\begin{eqnarray*}
&& W_i= \W, \ \ \Bar{W}_i =\Wbar,\\
&& W'_{i'}=\Wprime
\end{eqnarray*}
correspond to deformations of pointed curves  
$(\hC_i,\by^i)$, $(\hat{\Bar{C}}_i,\bar{\by}^i)$, $((\hS_{i'})_\C,\bz^{i'})$,
respectively.

Another local calculation shows that
\[
  H^0(X,\lone)\cong\bigoplus_{\alpha=1}^{l_0}(V_\alpha\oplus \Bar{V}_\alpha)
              \oplus\bigoplus_{\alpha'=1}^{l_1}V'_{\alpha'}\ ,
\]
where 
\[
V_\alpha=\V, \ \ \Bar{V}_\alpha=\Vbar,\ \ V'_{\alpha'}=\Vprime
\]
correspond to smoothing of the nodes 
$r_\alpha, \bar{r}_\alpha, s_{\alpha'}$ in X, respectively.

Now the action of $S$ on $\pone$ is clear: 
it acts on $W_i\oplus\Bar{W}_i$ and
$V_\alpha\oplus\Bar{V}_\alpha$ by $(a,b)\mapsto(\bar{b},\bar{a})$,
and it acts on $W'_{i'}$ and  $V'_{\alpha'}$ by
$a \mapsto \bar{a}$. The real vector space
${W'_{i'}}^S$ corresponds to deformation of the pointed
symmetric Riemann surface $((\hS_{i'})_\C, \bz^{i'},\si)$.
We also have ${W'_{i'}}^S\cong \hat{W}_{i'}$, where
\[
  \hat{W}_{i'}=H^1(\hS_{i'}, \pa \hS_{i'},
  T_{\hS_{i'}}(-p^{i'}_1-\cdots-p^{i'}_{n_{i'}}),
  T_{\pa \hS_{i'}}(-q^{i'}_1-\ldots-q^{i'}_{m_{i'}}))
\]
  corresponds to deformation of the pointed bordered Riemann surface
$(\hS_{i'},\bp^{i'},\bq^{i'})$.

  The action of $S$ on $V'_{\alpha'}$ can be understood by studying local
  models
 \begin{eqnarray*}
  \{x^2+y^2=0\}/A && \textup{(type E) },\\
  \{x^2-y^2=0\}/A && \textup{(type H) },
 \end{eqnarray*}
  where the antiholomorphic involution is $A(x,y)=(\bar{x},\bar{y})$.
  The deformation of $\{x^2\pm y^2=0\}$ is given by
  $\{ x^2\pm y^2=\ep\}$, where $\ep\in\C$ is small.
  $A$ acts on deformation by $\ep\mapsto \bar{\ep}$,
  so the deformation of  $(\{x^2\pm y^2=0\},A)$ is given by
  $(\{x^2\pm y^2=\ep\},A)$, where $\ep\in\R$ is small.
  $\{x^2\pm y^2=\ep\}/A$, $\ep\in\R$ small gives a deformation
  of the boundary nodes, and there is a topological transition from
  $\ep>0$ to $\ep<0$ (Figure 12, 13, 14, 15).
  $\ep\in \C$ corresponds to $V'_{\alpha'}$,
  $\ep\in \R$ corresponds to
  ${V'_{\alpha'}}^S=\hat{V}_{\alpha'}$,
  and $\ep\geq 0$ corresponds to $\hat{V}_{\alpha'}^+$.

\begin{figure}\label{twelve}
\begin{center}
\psfrag{(tg,h,k)=(4,3,0)}{$(\tg,h,k)=(4,3,0)$}
\psfrag{(tg,h,k)=(4,2,1)}{$(\tg,h,k)=(4,2,1)$}
\includegraphics[scale=0.55]{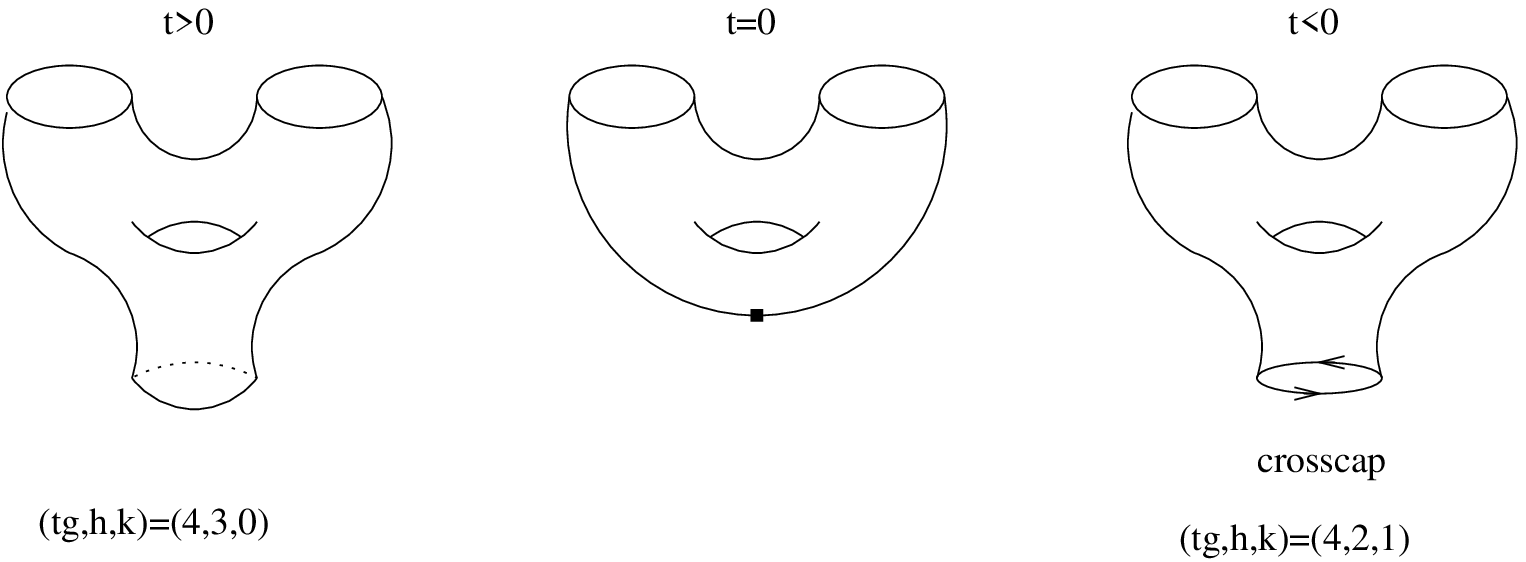}
\end{center}
\caption{type E}
\end{figure}

\begin{figure}\label{thirteen}
\psfrag{(tg,h,k)=(4,3,0)}{$(\tg,h,k)=(4,3,0)$}
\psfrag{(tg,h,k)=(4,2,1)}{$(\tg,h,k)=(4,2,1)$}
\includegraphics[scale=0.48]{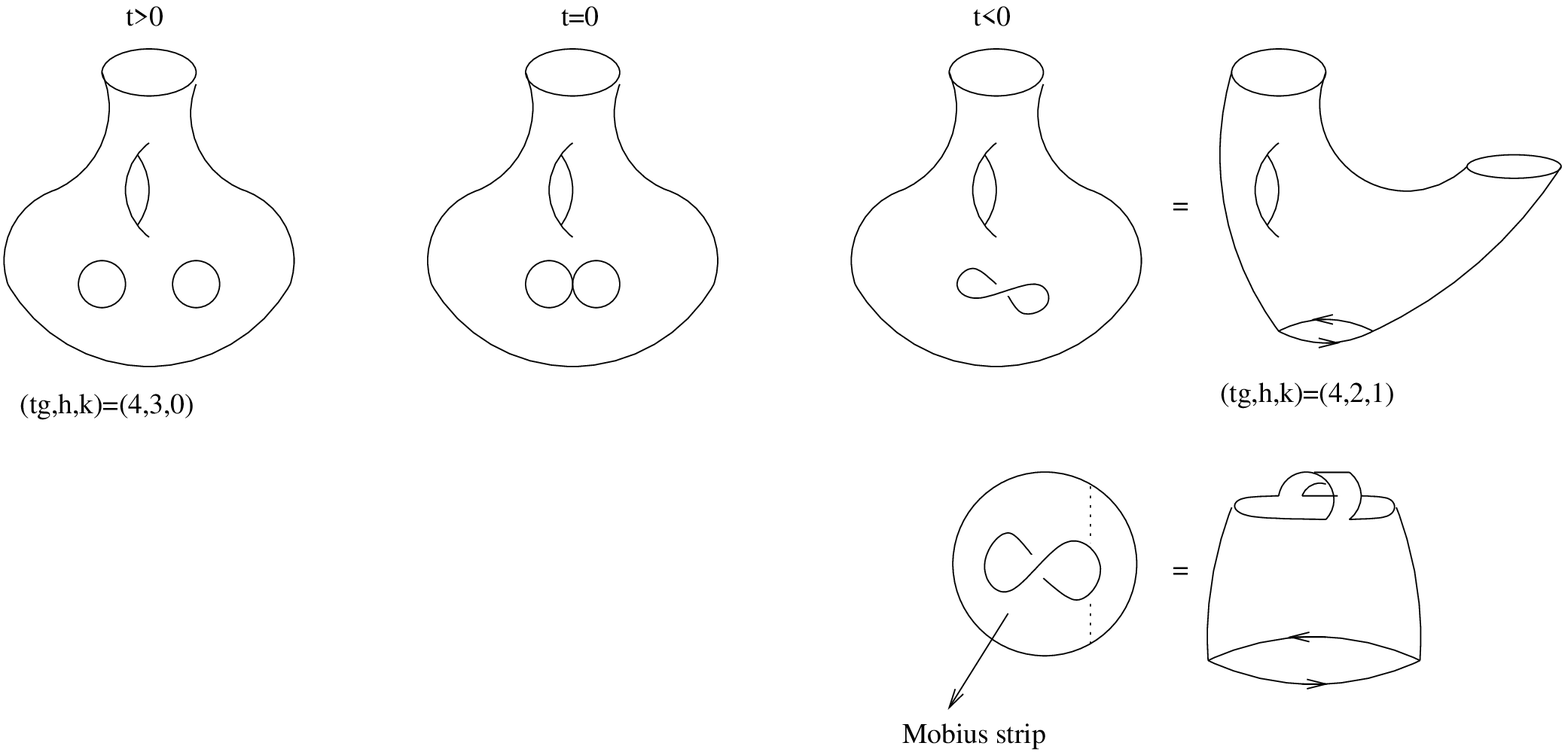}
\caption{type H1}
\end{figure}

\begin{figure}\label{fourteen}
\psfrag{(tg,h,k)=(4,3,0)}{$(\tg,h,k)=(4,3,0)$}
\psfrag{(tg,h,k)=(4,3,1)}{$(\tg,h,k)=(4,3,1)$}
\includegraphics[scale=0.48]{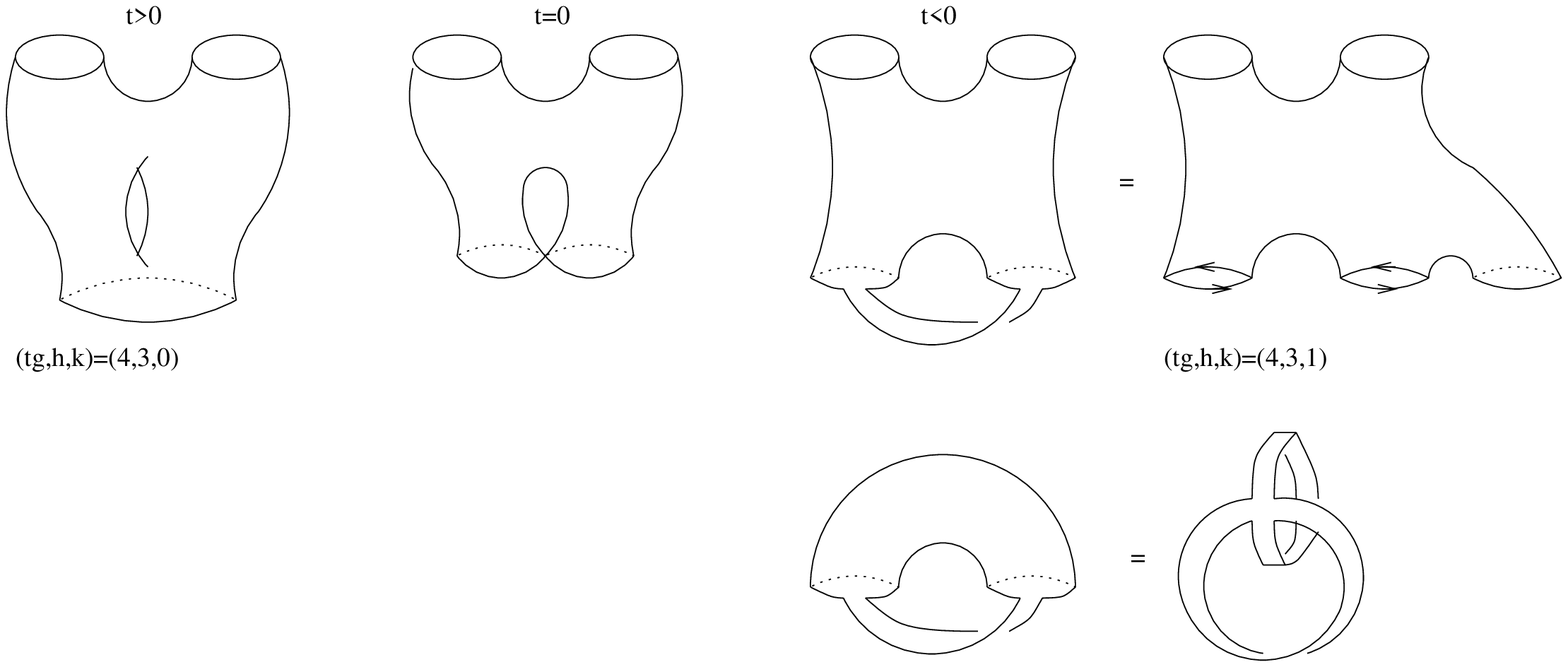}
\caption{type H2}
\end{figure}

\begin{figure}\label{fifteen}
\psfrag{(tg,h,k)=(4,3,0)}{$(\tg,h,k)=(4,3,0)$}
\includegraphics[scale=0.49]{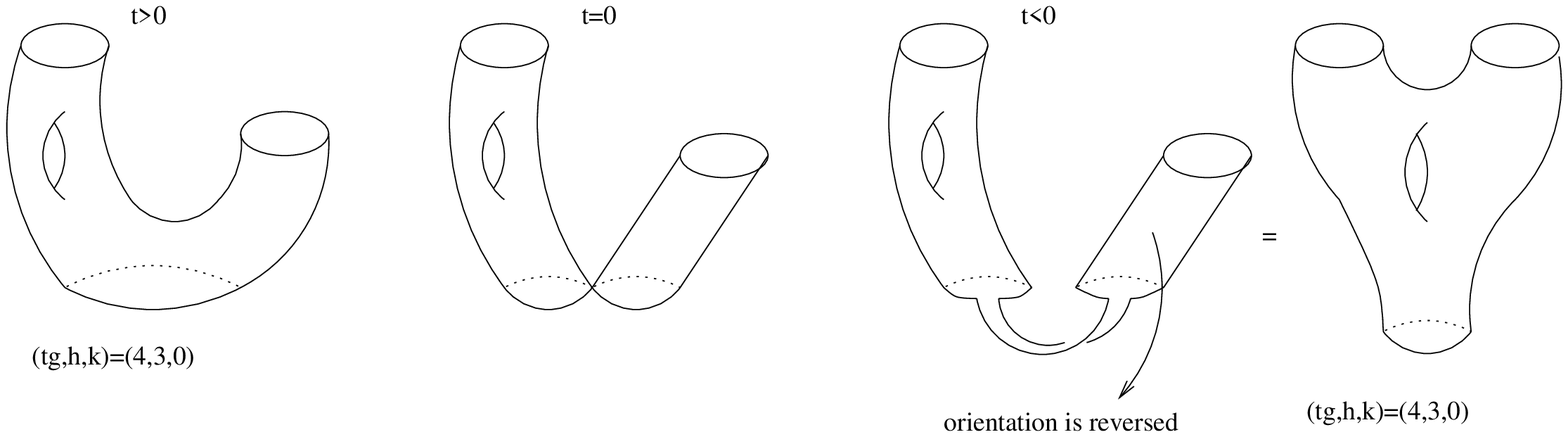}
\caption{type H3}
\end{figure}

 Above discussion can be summarized as follows.

\begin{enumerate}
 \item The infinitesimal deformation of the pointed complex algebraic curve
       $(X,\bx)$ is given by
    \[ \pone\cong
       \bigoplus_{i=1}^\nu(W_i\oplus\Bar{W}_i)\oplus\bigoplus_{i'=1}^{\nu'}W'_{i'}
       \oplus\bigoplus_{\alpha=1}^{l_0}(V_\alpha\oplus \Bar{V}_\alpha)
       \oplus\bigoplus_{\alpha'=1}^{l_1}V'_{\alpha'}.
    \]
 \item The infinitesimal deformation of the pointed complex algebraic curve with a real
       structure $ (X,\bx,S) $ is given by
     \[
       \pone^S\cong
       \bigoplus_{i=1}^\nu (W_i\oplus\Bar{W}_i)^S\oplus\bigoplus
        _{i'=1}^{\nu'} {W'_{i'}}^S
       \oplus\bigoplus_{\alpha=1}^{l_0}(V_\alpha\oplus \Bar{V}_\alpha)^S
       \oplus\bigoplus_{\alpha'=1}^{l_1}{V'_{\alpha'}}^S.
     \]

 \item The infinitesimal deformation of the pointed prestable bordered Riemann surface
       $(\Si;\bp;\bq)$  is given by
     \[
       \bigoplus_{i=1}^\nu W_i\oplus\bigoplus_{i'=1}^{\nu'}\hat{W}_{i'}
       \oplus\bigoplus_{\alpha=1}^{l_0} V_\alpha
       \oplus\bigoplus_{\alpha'=1}^{l_1}\hat{V}_\alpha^+.
     \]
\end{enumerate}

\section{Moduli of Bordered Riemann Surfaces}

For stable complex algebraic curves, the moduli of
complex structures can be identified with the moduli of
hyperbolic structures.  Under this identification,
analytic methods are applicable to the study of moduli
of stable curves.
In \cite{A}, Abikoff constructed a topology on $\Bar{M}_{g,n}$, the moduli of
stable complex algebraic curves of genus $g$ with $n$ marked points, and
showed that $\Bar{M}_{g,n}$ is compact and Hausdorff in this topology.
Similarly, Sepp\"{a}l\"{a} \cite{S} constructed a topology for the moduli space
of stable real algebraic curves of a given genus $g>1$, and showed that
the moduli space is compact and Hausdorff in the topology.
In this section, we describe how the above works can be modified to study
stable bordered Riemann surfaces.

\subsection{Various moduli spaces and their relationships} \label{model}

Let
$$
\domain
$$
be a marked stable bordered Riemann surface of type
$(g,h)$ with $(n,\vm)$ marked points. The moduli of the ordering
of the circles $B^1,\ldots,B^h$ is given by the permutation group of $h$ elements.
Let $(\Si_\C, \si)$ be the complex double of $\Si$, and let
$$
\bx=(x_1,\ldots,x_{2n+m})=(p_1,\ldots,p_n, \bar{p}_1,\ldots,\bar{p}_n, q_1,\ldots, q_m)
$$
$(\Si_\C,\si, p_1,\bx)$ is a stable symmetric Riemann surface
of genus $\tg=2g+h-1$ with $(n,m)$ marked points.
Removing $x_1,\ldots,x_{\tn}$ from $\Si_\C$, where $\tn=2n+m$, we obtain $(S,\si)$,
a stable symmetric Riemann surface of genus $\tg$ with $\tn$ punctures.
Let $S'$ be the complement of nodes
in $S$. There is a one to one correspondence between connected components
of $S'$ and irreducible components of $S$. Each connected components of $S'$
is a smooth punctured Riemann surface. The stability condition is equivalent
to the statement that each connected component of $S'$ has negative
Euler characteristic. Therefore, there is a unique complete hyperbolic
metric in the conformal class of Riemannian metrics on $S'$ determined by
the complex structure.

Let $\Bar{M}_{\tg,\tn}$ be the moduli of stable compact Riemann surfaces of genus $\tg$ with
$\tn$ marked points, or equivalently, the moduli of stable complex algebraic
curves of genus $\tg$ with $\tn$ marked points.
Let $\Bar{P}_{\tg,\tn}$ be the moduli of stable oriented hyperbolic surfaces of
genus $\tg$ with $\tn$ punctures. From the above discussion we know that
there is a surjective map $\tilde{\pi}:\Bar{M}_{\tg,\tn}\ra
\Bar{P}_{\tg,\tn}$. $\tilde{\pi}$ is
generically $\tn!$ to one since the marked points are ordered, while the
punctures are not. The fiber over a point in $\Bar{P}_{\tg,\tn}$ represented
by the surface $S$ consists of less than $\tn!$ points if and only if there
is an automorphism of $S$ permuting its punctures.

Let $\Bar{M}_{\tg,(n,m)}^\R$ be the moduli of stable symmetric compact
Riemann surface of genus $\tg$ with $(n,m)$ points, and let
$M_{(\tg,h,k),(n,m)}^\R$ be the moduli of smooth symmetric compact
Riemann surface of type $(\tg,h,k)$ with $(n,m)$ marked points.
Note that $M_{(\tg,0,k),(n,m)}^\R$ is empty if $m>0$.
$M_{(\tg,h,k),(n,m)}^\R$ are disjoint subsets of
$\Bar{M}_{\tg,(n,m)}^\R$, and their closures $\Bar{M}_{(\tg,h,k),(n,m)}^\R$
(in the topology defined later) cover $\Bar{M}_{\tg,(n,m)}^\R$.

There is an involution $A:\Bar{M}_{\tg,\tn}\ra\Bar{M}_{\tg,\tn}$, given by
\[
[ (\Si,x_1,\ldots,x_{\tn}) ] \mapsto
[ (\Bar{\Si},\si(x_{n+1}),\ldots,\si(x_{2n}),
      \si(x_1),\ldots,\si(x_n), \si(x_{2n+1}),\ldots,\si(x_{2n+m}))]
\]
where $\si:\Si\ra\Bar{\Si}$ is the canonical anti-holomorphic map from
$\Si$ to its complex conjugate $\Bar{\Si}$. Let
$\Bar{M}_{\tg,\tn}^A$ denote the fixed locus of $A$. Then there is a
surjective map $\Bar{M}_{g,(n,m)}^\R\ra\Bar{M}_{\tg,\tn}^A$, given by
forgetting the symmetry $\si$. This map is generically injective.
It fails to be injective exactly when the automorphism
group of $(\Si,\bx)$ is larger than that of
$(\Si,\si, \bx)$.

Let $\Bar{P}_{\tg,\tn}^\R$ be the moduli of stable symmetric oriented hyperbolic
surfaces of genus $\tg$ with $\tn$ punctures, and let
$P_{(\tg,h,k),\tn}^\R$  the moduli of smooth symmetric oriented hyperbolic
surfaces of type $(\tg,h,k)$ with $\tn$ punctures.
$P_{(\tg,h,k),\tn}^\R$ are disjoint subsets of $\Bar{P}_{\tg,\tn}^\R$, and
their closures $\Bar{P}_{(\tg,h,k),\tn}^\R$ (in the topology to be defined
later) cover $\Bar{P}_{\tg,\tn}^\R$.

There is an involution $A':\Bar{P}_{\tg,\tn}\ra\Bar{P}_{\tg,\tn}$, given by
$[S]\mapsto [\bar{S}]$.
There is a surjective map $\Bar{P}_{\tg,\tn}^\R\ra\Bar{P}_{\tg,\tn}^{A'}$, given by
forgetting the symmetry $\si$. This map is generically injective.
It fails to be injective exactly when the automorphism
group of $S$ is larger than that of $(S,\si)$.

We have the following commutative diagrams:

\[
\begin{array}{lll}
{ \begin{CD}
\Bar{M}_{\tg,\tn} @>A>> \Bar{M}_{\tg,\tn} \\
@V\tilde{\pi}VV @V\tilde{\pi}VV  \\
\Bar{P}_{\tg,\tn} @>A'>> \Bar{P}_{\tg,\tn}
\end{CD} }
& &
{ \begin{CD}
\Bar{M}_{\tg,(n,m)}^\R @>>> \Bar{M}_{\tg,\tn}^A\\
@V{\pi^\R}VV @V{\pi^A}VV\\
\Bar{P}_{\tg,(n,m)}^\R @>>> \Bar{P}_{\tg,\tn}^{A'}
\end{CD}}
\end{array}
\]
where the generic fiber of $\pi^\R$ consists of $2^n n!m!$ points.
The factor $2^n$ corresponds to the permutation of the two points
in each of the $n$ conjugate pairs. The factor $n!$ corresponds to
the permutation of the $n$ conjugate pairs. The factor
$m$ corresponds to the permutation of the $m$ marked points fixed
by the symmetry. Similarly, the generic fiber of
$\pi^A$ consists of $2^n n! m!$ points.
For a generic point in $\Bar{P}_{\tg,\tn}^{A'}$, its preimage
under $\tilde{\pi}$ consists of $\tn !$ points, but only
$2^n n!m!$ lie in the fixed locus of $A$.

Neither $\pi^\R$ nor $\pi^A$ is surjective because the number of punctures fixed by
the symmetry can be any integer between $0$ and $\tn$, not only $m$.

We are interested in the moduli space $\Bar{M}_{(g,h),(n,\vm)}$ of stable bordered 
Riemann surfaces of type $(g,h)$ with $(n,\vm)$ marked points. There is a finite to
one map $\Bar{M}_{(g,h),(n,\vm)}\ra\Bar{M}_{\tg,(n,m)}^\R$ given by complex double,
so there is a finite to one map $\Bar{M}_{(g,h),(n,\vm)}\ra
\Bar{P}_{(\tg,h,0),\tn}^\R\subset \Bar{P}_{\tg,\tn}^\R$.
We will first study $\Bar{P}_{\tg,\tn}$ and $\Bar{P}_{\tg,\tn}^\R$,
following \cite{A,S}.

\subsection{Decomposition into pairs of pants}

A \emph{pair of pants} $P$ is a sphere from which three
disjoint closed discs (or points) have been removed.
It is the interior a stable bordered Riemann surface of type $(0,3)$.
There is a unique hyperbolic structure compatible
with the complex structure of $P$ such that the
boundary curves are geodesics. Conversely, given
a hyperbolic structure on $P$ such that the boundary curves
are geodesics, the conformal structure is determined up to
conformal or anti-conformal equivalence by the lengths
$l_1$, $l_2$, and $l_3$ of the three boundary curves
(\cite[Chapter II (3.1), Theorem]{A}).

\subsubsection{Riemann surfaces with punctures}
A Riemann surface $S$ of genus $\tg$ with $\tn$ punctures can be decomposed
into pairs of pants. More precisely, there are $3\tg-3+\tn$ disjoint curves
$\alpha_1,\ldots, \alpha_{3\tg-3+\tn}$ on $\Si$, each of which is either a
closed geodesic (in the hyperbolic metric) or a node, such that the
complement of $\cup_{i=1}^{3\tg-3+\tn} \alpha_i$ is a disjoint union of
$2\tg-2+\tn$ pairs of pants $P_1,\ldots, P_{2\tg-2+\tn}$. A boundary
component of the closure of a pair of pants in this decomposition is either
a decomposing curve or a puncture. We call
\[\cP=\{ P_1,P_2,\ldots,P_{2\tg-2+\tn} \} \]
a {\em geodesic decomposition of $S$ into pairs of pants}.

Suppose that there exists an anti-holomorphic involution $\si: S\ra S$.
Then $(S,\si)$ is a stable symmetric Riemann surface of genus $\tg$ with
$\tn$ punctures, and $\si$ is an isometry of the hyperbolic metric. Let
\[
  \cP=\{ P_1, P_2, \ldots, P_{2\tg-2+\tn} \}
\]
be a geodesic decomposition of $S$ into pairs of pants. Then
\[
  \si(\cP)=\{ \si(P_1), \si(P_2), \ldots, \si(P_{2\tg-2+\tn}) \}
\]
is another geodesic decomposition of $S$ into pairs of pants.
$\cP$ is said to be \emph{$\si$-invariant} if
$\si(\cP)=\cP$.
The argument in \cite[Section 4]{S}, combined with
\cite[Chapter II (3.3), Lemma 3]{A}, shows that
\begin{tm}\label{sympant}
Let $(S,\si)$ be a stable symmetric Riemann surface of genus $\tg$ with $\tn$ punctures.
There exists a $\si$-invariant geodesic decomposition of $\Si$ into pairs of
pants such that the decomposing curves are simple closed geodesics of length less than
$C(\tg,\tn)$, where $C(\tg,\tn)$ is a constant depending only on $\tg$, $\tn$.
\end{tm}

\subsubsection{Riemann surfaces with boundary and punctures}

Let $(\Si,\bB;\bp)$ be a stable bordered Riemann surface
of type $(g,h)$ with $(n,\vec{0})$ marked points, and suppose that $\Si$ has no
boundary nodes. Let $S$ be the complement of marked points in $\Si$, then
$S$ is a surface of type $(g,h,n)$ in the sense of \cite{A}, namely,
$S$ is obtained by removing $h$ open discs and $n$ points from a compact
(possibly nodal) Riemann surface of genus $g$, where the discs and points
are all disjoint. $S$ is stable in the sense that its automorphism group
is finite. Let $S'$ be the complement of nodes in $S$.
Each connected component of $S'$ is a smooth Riemann surface with boundary or
punctures. The stability condition on $S$ is equivalent to the statement that
each connected component of $S'$ has a negative Euler characteristic, so
there exists a unique hyperbolic metric on $S'$ in the conformal class
determined by the complex structure such that the boundary circles are geodesics.

Let $S$ be a stable surface of type $(g,h,n)$. The $S$ can be decomposed
into pairs of pants. More precisely, there are $3g+h-3+n$ disjoint curves
$\alpha_1,\ldots, \alpha_{3g+h-3+n}$ on $S$, each of which is either a
closed geodesic (in the hyperbolic metric) or a node, such that the
complement of $\cup_{i=1}^{3g+h-3+n} \alpha_i$ is a disjoint union of
$2g+h-2+n$ pairs of pants $P_1,\ldots, P_{2g+h-2+n}$. A boundary
component of the closure of a pairs of pants in this decomposition is
a decomposing curve, a boundary component or a puncture. We call
\[\cP=\{ P_1,P_2,\ldots,P_{2g+h-2+n} \} \]
a {\em geodesic decomposition of $S$ into pairs of pants}.
We have the following result (\cite[Chapter II (3.3), Lemma 3]{A}):

\begin{tm}\label{bordpant}
Let $S$ be a stable surface of type $(g,h,n)$.
There is a geodesic decomposition of $S$ into pairs of pants such that
the decomposition curves are simple closed geodesic with length
less than $C(g,h,n,L_1,\ldots,L_h)$, where $C(g,h,n,L_1,\ldots,L_h)$ is
a constant depending only on $g$, $h$, $n$, and the lengths of the
$h$ border curves $L_1,\ldots,L_h$.
\end{tm}

We will see later that the moduli of stable surfaces of type $(g,h,n)$
is of (real) dimension $6g+3h-6+2n$, and $L_1,\ldots,L_h$ are among
the $6g+3h-6+2n$ real parameters. They are not good coordinates for
compactness since $L_1,\ldots,L_h$ can be arbitrarily large. Actually,
the length of some border curve tends to infinity as $S$ acquires a 
type H boundary node. To deal with boundary nodes and boundary marked points, 
we go to the complex double, where the local coordinates of the moduli can be 
chosen to be bounded.

\subsection{Fenchel-Nielsen coordinates}
\begin{df}\label{order}
Let $S$ be a stable symmetric Riemann surface of genus $\tg$ with $\tn$ punctures.
A geodesic pants decomposition $\cP$ is \emph{oriented} if
\begin{enumerate}
 \item The pairs of pants in $\cP$ is ordered.
 \item The boundary components of each pair of pants in $\cP$ is ordered.
 \item Any decomposing curve which is not a node is oriented.
\end{enumerate}
\end{df}

\begin{rem}
The orientability of a geodesic decomposition of a surface of type
$(g,h,n)$ can be defined similarly, with the additional assumption
that the boundary components are ordered.
\end{rem}

Let $P$ be a pair of pants with hyperbolic structure
with ordered boundary components $\alpha_1$, $\alpha_2$, $\alpha_3$.
The base points $\xi_i$ on $\alpha_i$, $i=1,2,3$, are shown in
Figure 16.
\begin{figure}\label{sixteen}
\begin{center}
\psfrag{a1}{$\alpha_1$} \psfrag{a2}{$\alpha_2$}
\psfrag{a3}{$\alpha_3$} \psfrag{x1}{$\xi_1$} \psfrag{x2}{$\xi_2$}
\psfrag{x3}{$\xi_3$} \psfrag{r1}{$\gamma_{2,3}$}
\psfrag{r2}{$\gamma_{3,1}$} \psfrag{r3}{$\gamma_{1,2}$}
\includegraphics{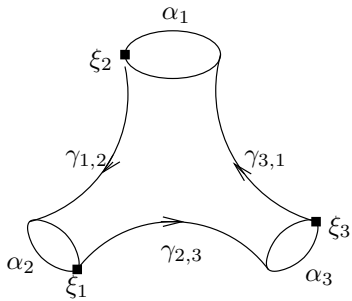}
\end{center}
\caption{a pair of pants}
\end{figure}
In Figure 16, $\gamma_{1,2}$ is the geodesic which realizes the distance between
$\alpha_1$, $\alpha_2$, etc.

If $\cP$ is oriented, then each boundary curve of a pair of pants
in $\cP$ has a base point. Each decomposing curve has
two distinguished points since it appears twice as a boundary of some
pair of pants in $\cP$. The orderings in 1. and 2.
of Definition~\ref{order} determine
an ordering of the decomposing curves
$\alpha_j$, $j=1,\ldots, 3\tg-3+\tn$ and an ordering on
the two distinguished points  $\xi_j^1, \xi_j^2$ on each decomposing curve.
We define $l_j(\cP)$ to be the length of $\alpha_j$, and
$\tau_j(\cP)$ be the distance one travels from $\xi_j^1$ to $\xi_j^2$
along $\alpha_j$, in the direction determined by 3. in
Definition~\ref{order}. Set $\theta_j(\cP)=2\pi\frac{\tau_j}{l_j}$ if $l_j\neq 0$.
We have $l_j\geq 0$, $0\leq\tau_j<l_j$, and $0\leq \theta_j <2\pi$.

\begin{df}
 Let $\tilde{S}$, $S$ be two stable Riemann surfaces of genus $\tg$ with
 $\tn$ punctures. Let $\tilde{S}'$, $S'$ be $1$-dimensional complex
 manifolds obtained from $\tilde{S}$, $S$ by removing the nodes.
 A {\em strong deformation} $\kappa:\tilde{S}\ra S$
 is a continuous map such that
\begin{enumerate}
 \item If $\tilde{r}$ is a node on $\tilde{S}$, then
       $\kappa(\tilde{r})$ is a node on $S$.
 \item If $r$ is a node on $S$, then $\kappa^{-1}(r)$ is a node or
       an embedded circle on a connected component of $S'$.
 \item $f|_{\kappa^{-1}(S')}:\kappa^{-1}(S')\ra S'$ is a diffeomorphism.
\end{enumerate}
\end{df}

Note that 1. implies that $\kappa^{-1}(S')\subset \tilde{S}'$, so
3. makes sense. There is a strong deformation $\kappa:\tilde{S}\ra S$
if and only if $\tilde{S}$ can be obtained by deforming $S$
as a quasiprojective variety over $\C$.

\begin{rem}
A strong deformation $\kappa:(\tilde{S},\tilde{\si})\ra (S,\si)$
between stable symmetric Riemann surfaces can be defined similarly,
with the additional assumption that $\si\circ\kappa=\kappa\circ \tilde{\si}$.
A strong deformation between two surfaces of type $(g,h,n)$ can also be
defined similarly.
\end{rem}

We now describe {\em Fenchel-Nielsen coordinates}
for various category of surfaces.
\begin{enumerate}
\item
Let $S$ be a stable Riemann surface of genus $\tg$ with $\tn$
punctures. Let $\cP$ be an oriented geodesic decomposition of $S$ into
pairs of pants, and let
$\alpha_1,\alpha_2,\ldots,\alpha_{3\tg-3+\tn}$ be the decomposing
curves of $\cP$. Suppose that there is a strong
deformation $\kappa:\tilde{S}\ra S$. Let $\tilde{\alpha}_j$ be the
closed geodesic homotopic to $\kappa^{-1}(\alpha_j)$. There exists
another strong deformation $\kappa'$ such that
$\kappa'({\tilde{\alpha}_j})=\alpha_j$, so $\cP$ is pulled
back under $\kappa'$ to an oriented geodesic decomposition
$\cP_S$.

\[
  [S]\mapsto (l_1(\cP_S),\theta_1(\cP_S),\ldots,
  l_{3\tg-3+\tn}(\cP_S),\theta_{3\tg-3+\tn}(\cP_S))
\]
defines local coordinates
$(l_1,\theta_1,\ldots,l_{3\tg-3+\tn},\theta_{3\tg-3+\tn})$
on $\Bar{P}_{\tg,\tn}$. Therefore, both $\Bar{M}_{\tg,\tn}$ and
$\Bar{P}_{\tg,\tn}$ are $(6\tg-6+2\tn)$ dimensional.

\item Let $(S,\si)$ be a stable symmetric surface of genus $\tg$ with
$\tn$ punctures, and let $\cP$ be an oriented geodesic
decomposition of $S$ into pairs of pants which is invariant under
$\si$. By considering $\si$-invariant decompositions into pants in
a neighborhood of $(S,\si)$ in $\mathcal{P}_{\tg,\tn}^\R$ we
obtain local coordinates
$(l_1,\theta_1,\ldots,l_{3\tg-3+\tn},\theta_{3\tg-3+\tn})$.
However, these parameters are not independent.
If $\si(\alpha_i)=\alpha_j$, $i\neq j$, then
$l_i=l_j$, and $\theta_i= c\pm \theta_j$  for some constant $c$.
If $\si(\alpha_i)=\alpha_i$, then $\theta_i=0$. So there are
$3\tg-3+\tn$ independent parameters, and the dimension of
$\mathcal{P}_{\tg,\tn}^\R$ is $3\tg-3+\tn$.

\item Let $S$ be a stable surface of type $(g,h,n)$, and let
$R_1,\ldots,R_h$ be its border curves.
Let $\cP$ be an oriented geodesic decomposition of $S$ into
pairs of pants, and let $\alpha_1,\ldots,\alpha_{3g+h-3+n}$ be
the decomposing curves. We have local coordinates
$(l_1,\theta_1,\ldots,l_{3g+h-3+n},\theta_{3g+h-3+n},L_1,\ldots,L_h)$,
where $l_j$ is the length of $\alpha_j$, $\theta_j$ is the angle of
gluing along $\alpha_j$, and $L_i$ is the length of the border
curve $R_i$. Therefore, the moduli of stable surfaces of type
$(g,h,n)$ is $6g+3h-6+2n = 3\tg -3 +\tn$, where $\tg=2g+h-1$, and
$\tn=2n$. This is consistent with the previous paragraph since
the complex double of $S$ is a stable symmetric Riemann surface
with genus $\tg=2g+h-1$ and $\tn=2n$ punctures.

\item Let $\domain$ be a stable marked bordered Riemann
surface of type $(g,h)$ with $(n,\vm)$ marked points.
$(\Si_\C,\si,\bx)$ be its complex double, and
$(S,\si)$ be the symmetric Riemann surface obtained from $\Si$
by removing the marked points, as in Section~\ref{model}.
Choose a $\si$-invariant geodesic decomposition of $S$ into pants,
we have local coordinates
$(l_1,\theta_1,\ldots,l_{3\tg-3+\tn},\theta_{3\tg-3+\tn})$,
as described in 1., where 
\begin{eqnarray*}
&& \tg=2g+h-1,\  \ \tn=2n+m^1+\cdots+m^h,\\  
&& 3\tg-3+\tn=6g+3h-3+2n+m^1+\cdots +m^h.
\end{eqnarray*}
We have seen that half of the
$2(6g+3h-6+2n+m^1+\cdots+m^h)$ parameters are independent, so the dimension
of $\Bar{M}_{(g,h),(n,\vm)}$ is 
$$
6g+3h-6+2n+m^1+\cdots+m^h.
$$
\end{enumerate}

In the following example, we describe the Frenchel-Nelson coordintates
of the moduli space $\Bar{M}_{0,3}$ of a pair of pants explicitly.
\begin{ex}\label{ex:pants}
The hexigon in Figure 17 is obtained by cutting the pair of pants
in Figure 16 along the geodesics $\gamma_{1,2}, \gamma_{2,3},
\gamma_{3,1}$. $\beta_1$ is the geodesic which realizes the distance
between $\gamma_1$ and $\gamma_{2,3}$, etc.

Let $l_1,l_2,l_3, l_4, l_5, l_6, l_7, l_8, l_9$
be twice the lengths of
$\gamma_1,\gamma_2,\gamma_3,\gamma_{2,3}, \gamma_{3,1}, \gamma_{1,2},
\beta_1,\beta_2,\beta_3$, respectively. The degeneration $l_i=0$ corresponds to 
a real codimension one stratum $V_i$ of $\Bar{M}_{0,3}$. Let $V_{ij}=V_i\cap V_j$ and
$V_{ijk}=V_i\cap V_j \cap V_k$. $\Bar{M}_{0,3}$ can be identified
with the associahedron $K_5$ defined by J. Stasheff \cite{St}.
The configuration of the strata in $\Bar{M}_{0,3}\cong K_5$ is
shown in Figure 17. There is one 3-dimensional stratum.
There are nine 2-dimensional strata:
$$
V_1,V_2,V_3,V_4,V_5,V_6,V_7,V_8,V_9.
$$ 
There are twenty-one 1-dimensional strata:
$$
V_{12}, V_{13}, V_{23},V_{45},V_{56},V_{46}, V_{14},V_{25},V_{36},
V_{27},V_{37},V_{57},V_{67},
V_{18},V_{38},V_{48},V_{68},
V_{19},V_{29},V_{49},V_{59}.
$$
There are fourteen 0-dimensional strata:
$$
V_{123},V_{456}, 
V_{237},V_{257},V_{367},V_{567}
V_{138},V_{148},V_{368},V_{468}
V_{129},V_{149},V_{259},V_{459}.
$$

There is a one-to-one correspondence between the 0-dimensional
strata and Frenchel-Nelson coordinate charts of $\Bar{M}_{0,3}$:
the Frenchel-Nelson coordinates near $V_{ijk}$ are $l_i,l_j,l_k$.  

\end{ex}
\begin{figure}\label{seventeen}
\begin{center}
\psfrag{r1}{$\gamma_1$}
\psfrag{r2}{$\gamma_2$}
\psfrag{r3}{$\gamma_3$}

\psfrag{r12}{$\gamma_{1,2}$}
\psfrag{r23}{$\gamma_{2,3}$}
\psfrag{r31}{$\gamma_{3,1}$}
\psfrag{b1}{$\beta_1$}
\psfrag{b2}{$\beta_2$}
\psfrag{b3}{$\beta_3$}
\psfrag{V12}{\tiny $V_{12}$}
\psfrag{V13}{\tiny $V_{13}$}
\psfrag{V23}{\tiny $V_{23}$}
\psfrag{V45}{\tiny $V_{45}$}
\psfrag{V56}{\tiny $V_{56}$}
\psfrag{V46}{\tiny $V_{46}$}
\psfrag{V14}{\tiny $V_{14}$}
\psfrag{V25}{\tiny $V_{25}$}
\psfrag{V36}{\tiny $V_{36}$}
\psfrag{V27}{\tiny $V_{27}$}
\psfrag{V37}{\tiny $V_{37}$}
\psfrag{V57}{\tiny $V_{57}$}
\psfrag{V67}{\tiny $V_{67}$}
\psfrag{V18}{\tiny $V_{18}$}
\psfrag{V38}{\tiny $V_{38}$}
\psfrag{V48}{\tiny $V_{48}$}
\psfrag{V68}{\tiny $V_{68}$}
\psfrag{V19}{\tiny $V_{19}$}
\psfrag{V29}{\tiny $V_{29}$}
\psfrag{V49}{\tiny $V_{49}$}
\psfrag{V59}{\tiny $V_{59}$}
\psfrag{K5}{$K_5$}
\includegraphics{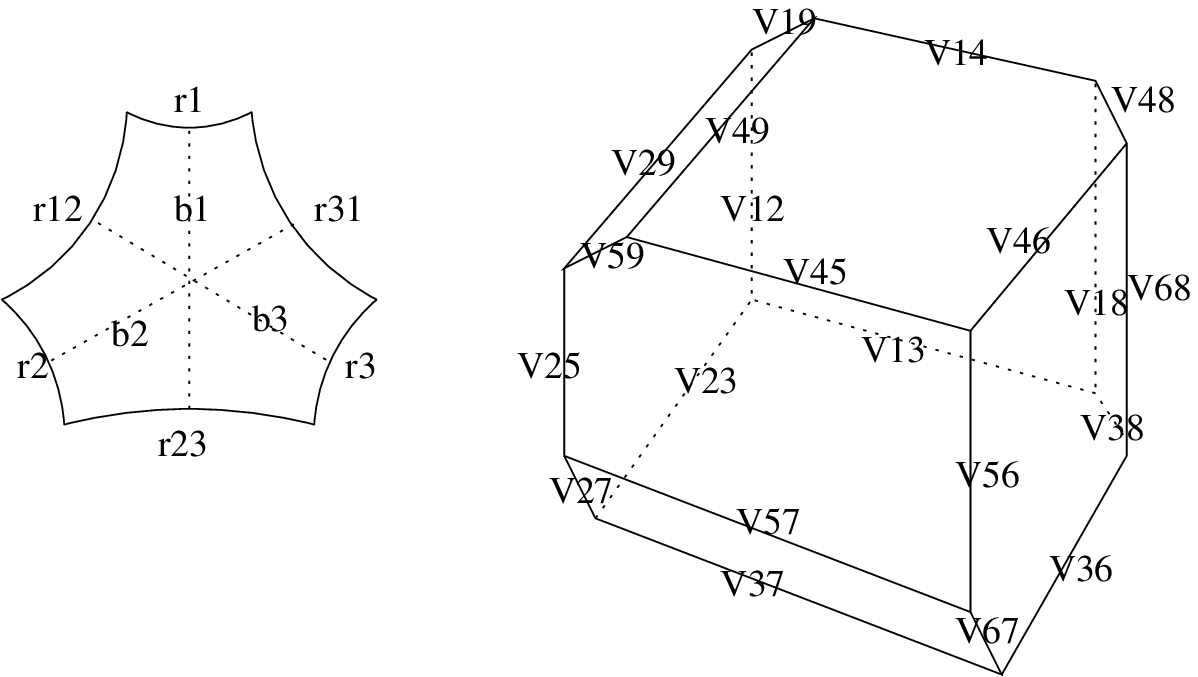}
\end{center}
\caption{}
\end{figure}

\subsection{Compactness and Hausdorffness}

We first define a topology on $\Bar{M}_{(g,h),(n,\vm)}$, following
\cite[Chapter II (3.4)]{A}, and \cite[Section 5]{S}.
We will call it \emph{Fenchel-Nielsen topology}.

\begin{df}\label{contract}
 A {\em strong deformation} between two stable marked bordered Riemann
 surfaces 
 $\tdomain$ and $\domain$ of type $(g,h)$ with $(n,\vm)$ marked points
 is a continuous map $\kappa:\tilde{\Si}\ra \Si$ such that
\begin{enumerate}
 \item $\kappa(\tilde{B}^i)=B^i$, $\kappa(\tilde{q}^i_k)=q^i_k$, $\kappa(\tilde{p}_j)=p_j$.
 \item If $\tilde{r}$ is an interior node on $\tilde{\Si}$, then
       $\kappa(\tilde{r})$ is an interior node on $\Si$.
 \item If $\tilde{s}$ is a boundary node of type E (H) on $\tilde{\Si}$, then
       $\kappa(\tilde{s})$ is a boundary node of type E (H) on $\Si$.
 \item If $r$ is an interior node on $\Si$, then $\kappa^{-1}(r)$ is an
       interior node or a circle.
 \item If $s$ is a boundary node of type E, then $\kappa^{-1}(r)$ is a
       boundary node of type E or a border circle.
 \item If $s$ is a boundary node of type H, then
       $\kappa^{-1}(s)$ is a boundary node or type H
       or  an arc with ends in $\bS$.
 \item $f|_{\kappa^{-1}(S')}:\kappa^{-1}(S')\ra S'$ is a diffeomorphism.
\end{enumerate}
\end{df}

Let
\[
\kappa:\tdomain\ra\domain
\]
be a strong deformation.
Let $(\tilde{S},\si)$, $(S,\si)$ be the symmetric Riemann surfaces obtained
by removing marked points from $\tilde{\Si}_\C$, $\Si_\C$, respectively.
Define $\kappa_\C:\tilde{\Si}_\C \ra \Si_\C$ by
\[
  \kappa_\C(z)=\left\{
  \begin{array}{ll}
   \kappa(z)& \textup{ if }z\in\Si\\
   \tilde{\si}\circ\kappa\circ\si(z)& \textup{ if }z\in\Bar{\Si}
   \end{array} \right.
\]
Let $\hat{\kappa}$ denote the restriction of $\kappa_\C$ to
$\tilde{S}$.  Then $\hat{\kappa}:(\tilde{S},\tilde{\si})\ra(S,\si)$
is a strong deformation.

Given $\ep, \delta >0$ and 
$$
\rho=[\domain]\in \Bar{M}_{(g,h),(n,\vm)},
$$
we will define a neighborhood $U(\ep,\delta,\rho)$ of $\rho$ in
$\Bar{M}_{(g,h),(n,\vm)}$. Let $(S,\si)$ be the symmetric Riemann
surface obtained by removing marked points from $\Si_\C$. Let
$$
\tilde{\rho}=[\tdomain]\in \Bar{M}_{(g,h),(n,\vm)},
$$
and
$(\tilde{S},\tilde{\si})$ be the associated symmetric Riemann surface.
Then $\tilde{\rho}\in\Bar{M}_{(g,h),(n,\vm)}$ if
\begin{enumerate}
\item There exists a $\si$-invariant oriented geodesic decomposition
   of $S$ into pairs of pants.
\item There exists a strong deformation
$$ 
\kappa:\tdomain\ra\domain
$$ i
n the sense of Definition~\ref{contract}. So we have a strong deformation
  $\hat{\kappa}:\tilde{S}\ra S$ as above.
\item Let $l_j,\theta_j$ and $\tilde{l}_j,\tilde{\theta}_j$ be the
     Fenchel-Nielsen coordinates for $\mathcal{P}$ and
     $\hat{\kappa}^*(\mathcal{P})$, respectively.
     Set $d=6g+3h-6+2n+m$. We have
     $|l_j-\tilde{l}_j|<\ep$ for $j=1,\ldots,d$, and
     $|\theta_j-\tilde{\theta}_j|<\delta$ if $l_j>0$.
\end{enumerate}

$$ 
\{ U(\ep,\delta,\rho)\mid \ep,\delta>0, \rho\in\Bar{M}_{(g,h),(n,\vm)}\}
$$
form a basis of the \emph{Fenchel-Nielsen topology}.

\bigskip

$U(\ep,\delta,\rho)$ can be described more precisely.
Set $z_j=l_je^{i\theta_j}$, then up to permutation
and complex conjugation of some $z_k$ we have
\[
\si(z_1,z_2,\ldots,z_{2d_1-1},z_{2d_1},z_{2d_1+1},\ldots,z_{d})
=(\bar{z_2},\bar{z_1},\ldots,\bar{z}_{2d},\bar{z}_{2d-1},
\bar{z}_{2d_1+1},\ldots,\bar{z}_d),
\]
so the fixed locus of $\si$ consists of points of the form
\[ (\bar{z}_2,z_2,\ldots,\bar{z}_{2d_1-1},\bar{z}_{2d_1},x_1,\ldots,x_{d_2}), \]
where $2d_1+d_2=d$, $z_2, z_4,\ldots,z_{2d_1}\in \C$, and
$x_1,\ldots,x_{d_2}\in\R$.
The coordinates take values in the fixed locus of $\si$, and $x_i$ are
nonnegative on $\Bar{M}_{(g,h),(n,m)}$
because negative values correspond to nonorientable surfaces, as we have
seen in Section~\ref{nodaldeform}.
We conclude that $U(\ep,\delta,\rho)$ is homeomorphic to $\tU/\Gamma$, where
$\tU$ is an open subset of $\C^{d_1}\times [0,\infty)^{d_2}$,
and $\Gamma$ is the automorphism group of $\tau$. The transition
functions between charts are real analytic (\cite[Appendix]{Wo}), so
Fenchel-Nielsen coordinates give $\Bar{M}_{(g,h),(n,\vm)}$ the structure
of an orbifold with corners. The topology determined by
the structure of an orbifold with corners coincides with
Fenchel-Nielsen topology. Therefore, we may equip
$\Bar{M}_{(g,h),(n,\vm)}$ with a metric which induces
Fenchel-Nielson topology. In particular, the topology
is Hausdorff, and compactness is equivalent to
sequential compactness. A straightforward generalization of
the argument in \cite[Section 6]{S} shows that
$\Bar{M}_{(g,h),(n,\vm)}$  is sequentially compact
in Fenchel-Nielsen topology. Therefore,

\begin{tm} \label{domaincpt}
$\Bar{M}_{(g,h),(n,\vm)}$ is Hausdorff and compact in
Fenchel-Nielson topology.
\end{tm}



\subsection{Orientation}

$\bar{M}_{(g,h),(n,\vm)}$ is an orbifold with corners, so we may ask
if it is orientable as an orbifold. By Stasheff's results in \cite{St},
we have
\begin{tm}\label{thm:contractible}
$\bar{M}_{(0,1),(0,(m))}$ has $(m-1)!$ isomorphic connected components,
which correspond to the cyclic ordering of the $m$ boundary
marked points. Each connected component of  $\bar{M}_{(0,1),(0,(m))}$ 
is homeomorphic to $\R^{m-3}$. 
\end{tm}

\begin{lm}\label{thm:induction}
Suppose that $(g,h,n)\neq (0,1,0)$, and $m^i>0$. If
$\Bar{M}_{(g,h),(n,\vm)}$ is orientable,
then $\Bar{M}_{(g,h),(n,(m^1,\ldots, m^i+1,\ldots, m^h))}$ is orientable. 
\end{lm}
\paragraph{Proof}
Assume that  $\Bar{M}_{(g,h),(n,(m^1,\ldots, m^i,\ldots,m^h)))}$
is orientable. Consider the map
\begin{equation}\label{eqn:forget}
F:  \Bar{M}_{(g,h),(n,(m^1,\ldots, m^i+1,\ldots, m^h))}
\ra \Bar{M}_{(g,h),(n,\vm)},
\end{equation}
given by forgetting the last boundary marked point on the
$i$-th boundary circle. Under our assumption, the fiber
of $F$ over $[\domain]$ is a union of $m^i$ 
intervals and inherits the orientation of $B^i$.
So  $\Bar{M}_{(g,h),(n,(m^1,\ldots, m^i+1,\ldots,m^h)))}$
is orientable. $\Box$

\begin{lm}\label{thm:reduction}
Supposet that $(g,h,n)\neq (0,1,0)$, and $m^i=0$. If\\
$\Bar{M}_{(g,h),(n,(m^1,\ldots, m^i+1,\ldots,m^h)))}$ is orientable,
then $\Bar{M}_{(g,h),(n,\vm)}$ is orientable.
\end{lm}
\paragraph{Proof}
Assume that $\Bar{M}_{(g,h),(n,(m^1,\ldots, m^i+1,\ldots,m^h)))}$ 
is orientable. Let $T$ be the tangent bundle of $\Bar{M}_{(g,h),(n,\vm)}$,
which is an orbibundle over $\Bar{M}_{(g,h),(n,\vm)}$.
To show that $\Bar{M}_{(g,h),(n,\vm)}$ is orientable,
it suffices to show that the restriction of $T$ to every
loop in $\Bar{M}_{(g,h),(n,\vm)}$ is orientable.
Let $N_{(g,h),(n,\vm)}$ be the interior of $\bar{M}_{(g,h),(n,\vm)}$. 
More precisely, $N_{(g,h),(n,\vm)}$  corresponds to surfaces with no
boundary nodes. Since every loop in $\Bar{M}_{(g,h),(n,\vm)}$ is 
homotopic to a loop in $N_{(g,h),(n,\vm)}$, it suffices to show that
$N_{(g,h),(n,\vm)}$ is orientable. 

Suppose that $\rho=[\domain]\in N_{(g,h),(n,\vm)}$. Then $B^i$
is an embedded circle in $\Si$, oriented as in Remark \ref{boundary}, 
and the fiber of the map $F$ in (\ref{eqn:forget})
over $\rho$ can be identified with $B^i$.
So $N_{(g,h),(n,\vm)}$ is orientable. $\Box$ 

\bigskip

It is shown in \cite{IS2} that 
\begin{tm}\label{thm:complex}
Suppose that $(g,h,n)\neq (0,1,0)$. Then
$\Bar{M}_{(g,h),(n,(1,\ldots,1))}$ is a complex orbifold.
\end{tm}

Theorem \ref{thm:contractible}, Lemma \ref{thm:induction},
Lemma \ref{thm:reduction}, and  Theorem \ref{thm:complex} 
imply that
\begin{tm}\label{thm:Morientable}
$\bar{M}_{(g,h),(n,\vm)}$ is orientable.
\end{tm}

Let $\Bar{Q}_{(g,h),n}$ be the moduli space of stable bordered
Riemann surfaces of type $(g,h)$ with $n$ interior points. There is
an $h!$ to one map $\Bar{M}_{(g,h),(n,\vec{0})}\ra \Bar{Q}_{(g,h),n}$,
given by forgetting the ordering of boundary components.
Then $\Bar{Q}_{(g,h),n}$ is nonorientable.

For example, consider $\rho\in\Bar{M}_{(1,2),(0,0)}$ represented by
the surface as show in Figure 18.
\begin{figure}\label{eighteen}
\begin{center}
\psfrag{a1}{$\alpha_1$} \psfrag{a2}{$\alpha_2$}
\psfrag{a3}{$\alpha_3$} \psfrag{a4}{$\alpha_4$}
\includegraphics{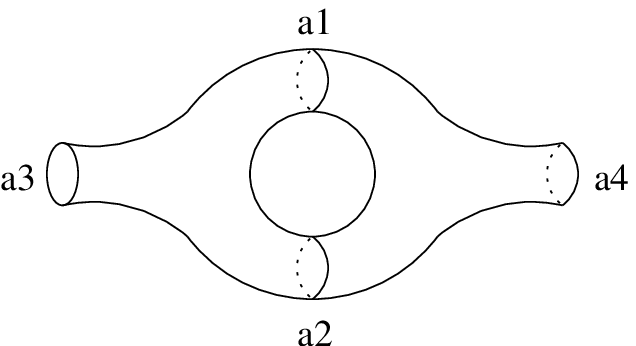}
\end{center}
\caption{$\Bar{Q}_{(g,h),n}$ is nonorientable}
\end{figure}
The local coordinates are  $(l_1,\theta_1,l_2,\theta_2,l_3,l_4)$, where
$l_j$ is the length of $\alpha_j$ for $j=1,\ldots,4$, and
$\theta_1,\theta_2$ are gluing angles for $\alpha_1,\alpha_2$,
respectively. There is an automorphism $\phi$ of order $2$
of $\rho$ which rotates the above Figure~14 by $180^\circ$.
$\phi(\alpha_1)=\alpha_2$, $\phi(\alpha_3)=\alpha_4$, and
$\alpha(l_1,\theta_1, l_2,\theta_2,l_3,l_4)=
 (l_2,-\theta_2,l_1,-\theta_1,l_4,l_3)$, which is orientation
reversing.

In general, an automorphism induces permutation of $d_1$ decomposing curves
and permutation of $d_2$ bordered curves. The former corresponds to permutation
of pairs $(l_j,\theta_j)$, $j=1,\ldots,d_1$, which is orientation
preserving. The later corresponds to permutation of $(l_{d_1+1},\ldots, l_{d_1+d_2})$
which is orientation preserving if and only if it is an even permutation.
When we consider $\Bar{M}_{(g,h),(n,\vec{0})}$, automorphisms permuting border
curves are not allowed.

\section{Moduli Space of Stable Maps} 

\subsection{Prestable and stable maps}

 Let $(X,\omega)$ be a compact symplectic manifold, and let 
 $L$ be a Lagrangian submanifold. Let $J$ be an $\omega$-tame
 almost complex structure.   

 \begin{df}
  A \emph{prestable map} is a continuous map
  $u:(\Si,\bS)\ra (X,L)$ such that $J\circ d\hu=d\hu\circ j$,
  where $\Si$ is a prestable bordered Riemann surface, $\hu=u\circ \tau$,
  $\tau:\hS\ra \Si$ is the normalization map 
  (Definition~\ref{norm}). 
 \end{df} 
 
 \begin{df}
  A \emph{prestable map} of type $(g,h)$ with $(n,\vm)$ marked points  
  consists of a prestable marked bordered Riemann surface of type $(g,h)$
  with $(n,\vm)$ marked points $\domain$
  and a prestable map $u:(\Si,\bS)\ra (X,L)$. 
 \end{df}

 \begin{df}\label{mapmor}
  A \emph{morphism} between prestable maps of type $(g,h)$ with $(n,\vm)$
  marked points 
\[
\map \to \anothermap
\]
  is an isomorphism  
\[  
  \phi: \domain \to \another
\]
 between prestable marked bordered Riemann surfaces
  of type $(g,h)$ with $(n,\vm)$ points such that $u=u'\circ\phi$.
 \end{df}

\begin{df}
  A prestable map of type $(g,h)$ with $(n,\vm)$ marked points
  is \emph{stable} if its automorphism group is finite.
\end{df}

\subsection{$C^\infty$ topology}

Let $\beta\in H_2(X,L;\Z)$, $\vga=(\gamma^1,\ldots,\gamma^h)\in H_1(L;\Z)^{\oplus h}$
be such that $\gamma_1+\cdots+\gamma_h =\pa\beta$, where
$\pa: H_2(X,L;\Z)\ra H_1(L;\Z)$ is the connecting map in the long
exact sequence for relative homology groups. Let $h$ be a positive
integer, $g,n$ be nonnegative integers, $\vm=(m^1,\ldots,m^h)$
be an $h$-uple of nonnegative integers, and $\mu$ be an integer.
Given the above data, define 
\[ 
\MXL 
\] 
to be the moduli space of isomorphism classes of 
stable maps of type $(g,h)$ with $(n,\vm)$ marked points 
\[ 
\map 
\]
such that $u_*[\Si]=\beta$, $u_*[B^i]=\gamma^i$ for $i=1,\ldots,h$,
and $\mu(u^*TX,u^*TL)=\mu$. Here $\mu(u^*TX,u^*TL)$ is the
\emph{Maslov index} defined in \cite[Definition 3.3.7, Definition 3.7.2]{KL}.
From now on, we assume that $L$ is oriented, so 
$\mu(u^*TX,u^*TL)$ is even, and we may restrict ourselves to even $\mu$.
We will also assume that {\em none of the $\gamma^i$ is zero, so the domain
cannot have boundary nodes of type E}.

Let $\Bar{M}_{(g,h),(n,\vm)}(X,L)$ be the moduli space of isomorphism classes
of stable maps of type $(g,h)$ with $(n,\vm)$ marked points. Then
$\MXL$ are disjoint subsets of $\Bar{M}_{(g,h),(n,\vm)}(X,L)$
for different $(\beta,\vga,\mu)$.

Let $\Si$ be a prestable bordered Riemann surface, and let
$\tau:\hat{\Si}\ra\Si$ be the normalization. Let 
$u:(\Si,\bS)\ra (X,L)$ be a continuous map such that
$\hat{u}=u\circ\tau: (\hat{\Si},\pa\hat{\Si})\ra (X,L)$ is 
$C^\infty$ w.r.t. $g_0$ on $X$ and some Hermitian metric 
$h$ on $\hat{\Si}$, $l\geq 1$. Define
\[
  a(u)=\int_{\hat{\Si}}\sqrt{
  \left|\frac{\pa \hat{u}}{\pa x}\right|^2 \left|\frac{\pa \hat{u}}{\pa y}\right|^2
  -\left\bra \frac{\pa \hat{u}}{\pa x}, \frac{\pa \hat{u}}{\pa y}\right\ket ^2}
   g(x,y) dx \wedge dy,
\]
  where $(x,y)$ are local isothermal coordinates on $\hat{\Si}$, and $g(x,y) dx \wedge dy$
  is the volume form for the metric $h$. If $u$ is an embedding,
  $a(u)$ is the area of $u(\Si)$ w.r.t. $g_0$. If $u$ is a prestable map,
  then $a(u)=\half\parallel du \parallel_{L^2}=(u_*[\Si])\cap [\omega]$,
  where
\[ \parallel du \parallel_{L^2}^2 
  =\int_{\hat{\Si}}\left( \left|\frac{\pa \hat{u}}{\pa x}\right|^2
      +\left|\frac{\pa \hat{u}}{\pa y} \right|^2 \right)g(x,y)dx\wedge dy.
\]
  $\left( u_*[\Si] \right) \cap [\omega]$ only depends on the relative
  homology class $u_*[\Si]\in H_2(X,L;\Z)$, so we have a function
  $a:\Bar{M}_{(g,h),(n,m)}\ra [0,\infty)$, which takes the constant value 
  $\beta\cap[\omega]$ on $\MXL$.  

With the above definition, $\MXL$ is a set. We will equip it with the structure
of a topological space, and show that it is sequentially compact 
and Hausdorff in this topology. This topology was introduced by Gromov \cite{Gr}.

We want to say two stable maps are close if the complex structures on the
domain are close, and the maps are close. To measure the closeness, we use
metrics on the domain and on the target.
For the target $X$, $J$ is an $\omega$-tame complex structure, so 
$g_0(X,Y)=\half(\omega(X,JY)+\omega(Y,JX))$
is a Riemannian metric on $X$ such
that $J$ is an isometry. For the domain, by a {\em Hermitian metric} $h$ on a prestable 
bordered Riemann surface $\Si$ we mean a Hermitian metric on $\hat{h}$
on $\hS$, the normalization of $\Si$. 

We now introduce some notation. Let $\tau:\hS\ra\Si$ be the normalization map.
Given a node $r\in\Si$ and a small positive number $\ep$,
let $B_\ep(r)=\tau(B_\ep(r_1)\cup B_\ep(r_2))$, where
$\tau^{-1}(r)=\{ r_1,r_2\}$, and $B_\ep(r_\alpha)$ is the geodesic ball
of radius $\ep$ for $\alpha=1,2$. Let $\ep$ be sufficiently small
so that $B_\ep(r)$ are disjoint for $r\in \Si_\mathrm{sing}$, where 
$\Si_\mathrm{sing}$ denotes the set of nodes on $\Si$. Set
$N_\ep(\Si)=\bigcup_{r\in\Si_\mathrm{sing}}B_\ep(r)$,
$K_\ep(\Si)=\Si-N_\ep(\Si)$. 

\begin{df}[$C^\infty$ topology] \label{Cl}
Let $\rho=\map$ be a prestable map of type $(g,h)$ with $(n,\vm)$ marked points.
For a Hermitian metric $h$ on $\Si$ and $\ep_1, \ldots, \ep_4>0$,
a neighborhood $U(\rho,h,\ep_1,\ldots,\ep_4)$ of $u$ in 
$\Bar{M}_{(g,h),(n,\vm)}(X,L)$ in the 
\emph{$C^\infty$ topology} is defined as follows. A prestable map 
\[ \rho=\anothermap \] belongs to $U(\rho,h,\ep_1,\ldots,\ep_4)$ if
\begin{enumerate}
 \item There is a strong deformation 
   \[ \kappa:\another\ra\domain \]
       such that $\kappa^{-1}$ is defined on $K_{\ep_1}(\Si)$. 
 \item $\parallel j-(\kappa^{-1})^* j'\parallel_{C^\infty(K_{\ep_1}(\Si))}<\ep_2$,
       where $j$, $j'$ are complex structures on $\Si$, $\Si'$, respectively.
 \item $\parallel u-u'\circ \kappa^{-1}\parallel_{C^\infty(K_{\ep_1}(\Si))} <\ep_3$.
 \item $|a(u)-a(u')|<\ep_4$. 
\end{enumerate}
\end{df}
 
(1) says that $\Si'$ can be obtained by deforming $\Si$, or equivalently, 
$\Si'$ is in the same or a higher stratum in $\widetilde{M}_{(g,h),(n,\vm)}$,
the moduli space  of $\emph{prestable}$ marked bordered
Riemann surfaces of type $(g,h)$ with $(n,\vm)$ marked points.
(2) says that \[ \domain,\  \another\] are close
in the $C^\infty$ topology (Definition~\ref{Cldomain}).
(3) says that the maps $u$, $u'$ are $C^\infty$ close away from the nodes.
(4) implies that $a:\Bar{M}_{(g,h),(n,\vm)}(X,L)\ra [0,\infty)$ is a continuous function.

\[ \map,\ \ \anothermap \] 
represent the same point in 
$\Bar{M}_{(g,h),(n,\vm)}(X,L)$ if and only if there is a Hermitian metric 
$h$ on $\Si$ such that

\begin{enumerate}
 \item There is a homeomorphism 
       \[ \kappa:\another\ra\domain \]
       which induces a diffeomorphism $\hS'\ra\hS$. 
 \item $\parallel j-(\kappa^{-1})^* j'\parallel_{C^\infty(\Si)}=0$,
       where $j$, $j'$ are complex structures
       on $\Si$, $\Si'$, respectively.
 \item $\parallel u-u'\circ \kappa^{-1}\parallel_{C^\infty(\Si)}=0$.
 \item $|a(u)-a(u')|=0$.
\end{enumerate}

So the $C^\infty$ topology is actually a topology on the moduli space 
$\Bar{M}_{(g,h),(n,\vm)}(X,L)$. 
$\MXL$ is a closed subspace of $\Bar{M}_{(g,h),(n,\vm)}(X,L)$ and is
equipped with the subspace topology. 

\bigskip
 
When $X$ is a point, we get the $C^\infty$ topology of 
$\Bar{M}_{(g,h),(n,\vm)}$.

\begin{df} \label{Cldomain}
Let $\lambda=\domain$ be a stable (prestable) bordered Riemann surface of 
type $(g,h)$ with $(n,\vm)$ marked points.
For a Hermitian metric $h$ on $\Si$ and $\ep_1,\ep_2>0$,
a neighborhood $U(\rho,h,\ep_1,\ep_2)$ of $u$ in 
$\Bar{M}_{(g,h),(n,\vm)}$ in the  
\emph{$C^\infty$ topology} is defined as follows. A stable (prestable) bordered Riemann
surface $\lambda'=\another$ belongs to $U(\rho,h,\ep_1,\ep_2)$ if
\begin{enumerate}
 \item There is a strong deformation $\kappa:\lambda'\ra\lambda$
       such that $\kappa^{-1}$ is defined on $K_{\ep_1}(\Si)$. 
 \item $\parallel j-(\kappa^{-1})^*
       j'\parallel_{C^\infty(K_{\ep_1}(\Si))}<\ep_2$,
       where $j$, $j'$ are complex structures on $\Si$, $\Si'$, respectively.
\end{enumerate}
\end{df}

\subsection{Compactness and Hausdorffness} \label{compact}

The following is the main theorem of this section.

\begin{tm}
$\MXL$ is Hausdorff and sequentially compact in the $C^\infty$ topology.
\end{tm}

The Hausdorffness can be proven as in the case of curves without boundary,
see e.g. \cite[Proposition 3.8]{Si}.     
The compactness is a consequence of the following theorem.

\begin{tm}[Gromov's Compactness Theorem] \label{gromov}
Let $\{\rho_l \}$ be a sequence in \\
$\Bar{M}_{(g,h),(n\vm)}(X,L)$ such that
$a(\rho_l)<C$ for all $l\in\cN$. Then there is a subsequence
of $\{\rho_l\}$ convergent in the $C^\infty$ topology. 
\end{tm}

Gromov's compactness theorem \cite[1.5]{Gr} for $J$-holomorphic curves
without boundary was carried out in details in \cite{PW, Ye}.  
The case with boundary was proved in \cite{Ye} (see also \cite{IS1, IS2}).
In \cite{Ye}, the moduli space is compactified by the moduli space of
\emph{cusp curves}, or prestable maps in this paper. 
We will describe how the proof in \cite{Ye} gives
Theorem~\ref{gromov}.
 
The $C^\infty$ topology can be equivalently defined as follows.
\begin{df}[$C^\infty$ Topology]\label{converge}
A sequence 
\[
\rho_l=(\Si_l,\bB_l;\bp_l;\bq^1_l,\ldots,\bq^h_l;u_l)
\]
converges to $\rho=\map$ in the $C^\infty$ topology 
if for each $\ep_1,\ldots, \ep_4 >0$, there is an integer $N$ such that for $l\geq N$,
\begin{enumerate}
 \item There is a strong deformation $\kappa_l:\Si_l \ra\Si$
       such that $\kappa_l^{-1}$ is defined on $K_{\ep_1}(\Si)$. 
 \item $\parallel j-(\kappa_l^{-1})^* j_l\parallel_{C^\infty(K_{\ep_1}(\Si))}<\ep_2$.
 \item $\parallel u-u_l\circ \kappa_l^{-1}\parallel_{C^\infty(K_{\ep_1}(\Si))} <\ep_3$.
 \item $|a(u)-a(u_l)|<\ep_4$.
\end{enumerate}
\end{df}

Recall that $\widetilde{M}_{(g,h),(n,\vm)}$ denotes the moduli space of prestable
bordered Riemann surfaces of type $(g,h)$ with $(n,\vm)$ marked points.
There is a map $F: \MXL\ra\widetilde{M}_{(g,h),(n,\vm)}$, given by forgetting
the map. $\widetilde{M}_{(g,h),(n,m)}$ has infinitely many strata since one 
can keep on going to lower and lower strata by adding non-stable components
-- spheres and discs.  

We claim that the image of $F$ is covered by only
finitely many strata, or equivalently:
\begin{lm}\label{finite}
The domains in $\MXL$ have only finitely many topological types. 
\end{lm}

\paragraph{Proof.}
There is a map
$\widetilde{M}_{(g,h),(n,\vm)}\ra \Bar{M}_{(g,h), (n,\vm)}$, given by 
contracting non-stable components. Since a stable bordered Riemann surfaces
of type $(g,h)$ with $(n,\vm)$ points can have only finitely many possible
topological types, it suffices to get an upper bound for
the number of non-stable irreducible components. The restriction of a stable 
map to a non-stable irreducible component is nonconstant, so there is a lower
bound $\ep>0$ for the area of the restriction of the map to each non-stable
component by \cite[Lemma 4.3, Lemma 4.5]{Ye}. Therefore, the number of non-stable
irreducible components cannot exceed $(\beta\cap[\omega])/\ep$. $\Box$  

\bigskip

Let $\{\rho_l\}$ be a sequence
in $\Bar{M}_{(g,h),(n,\vm)}(X,L)$ such that $a(\rho)<C$ for all $l\in \mathbb{N}$.
By Lemma~\ref{finite}, there is a subsequence of $\{\rho_l\}$ such that
the domains are of the same topological type. 
By normalization we obtain several sequences of stable
maps with smooth domains of the same topological type and with uniform area
bound. Note that each node gives rise to two marked points
on the normalization. It suffices to
show that each sequence has a subsequence convergent in the $C^\infty$  topology.
So we may assume that the domain is a smooth marked bordered Riemann surface
or a smooth curve with marked points. In this case, it is proven in \cite{Ye}
that there is a subsequence which converges to a prestable map in the $C^\infty$
topology. However, it is straightforward to check that the limit produced in
\cite{Ye} is actually a stable map.

\section{Construction of Kuranishi Structure} \label{constructku}

\subsection{Kuranishi structure with corners} \label{defineku}

We first quote the following definition from \cite[A2.1.1-A2.1.4]{FO3}, which is 
a slight modification of \cite[Definition 5.1]{FO}.
\begin{df}[Kuranishi neighborhood] \label{nbd}
Let $M$ be a Hausdorff topological space. A {\em Kuranishi neighborhood (with corners)}
of $p\in M$ is a $5$-uple $(V_p,E_p,\Gamma_p,\psi_p,s_p)$ such that
\begin{enumerate}
\item $V_p$ is a smooth manifold (with corners), and $E_p$ is a smooth vector
      bundle on it.
\item $\Gamma_p$ is a finite group which acts smoothly on $E_p\ra V_p$.
\item $s_p$ is a $\Gamma_p$-equivariant continuous section of $E_p$.
\item $\psi_p: s_p^{-1}(0)\ra M$ is a continuous map which induces a
      homeomorphism from $s_p^{-1}(0)/\Gamma_p$ to a
      neighborhood of $p$ in $M$.
\end{enumerate}
We call $E_p$  the {\em obstruction bundle} and $s_p$ the 
{\em Kuranishi map}.
\end{df}

The following equivalence relation is weaker than the one in 
\cite[Definition 5.2]{FO}, so the resulting equivalence class is larger. 
\begin{df}\label{same}
 Let $M$ be a Hausdorff topological space.
 Two Kuranishi neighborhoods (with corners)
 $(V_{1,p},E_{1,p},\Gamma_{1,p},\psi_{1,p},s_{1,p})$ and 
 $(V_{2,p},E_{2,p},\Gamma_{2,p},\psi_{2,p},s_{2,p})$ of $p\in M$ 
 are \emph{equivalent} if
 \begin{enumerate}
  \item $\dim V_{1,p}-\mathrm{rank}E_{1,p}=
         \dim V_{2,p}-\mathrm{rank}E_{2,p}\equiv d$.
  \item There is another Kuranishi neighborhood (with corners)
        $(V_p,E_p,\Gamma_p,\psi_p,s_p)$ of $p$ 
        such that $\dim V_p-\mathrm{rank}E_p=d$. 
  \item There are homomorphisms $h_i:\Gamma_{i,p}\ra \Gamma_p$ for $i=1,2$.
  \item For $i=1,2$, there is a $\Gamma_{i,p}$-invariant open neighborhood
        $V_i$ of $\psi^{-1}_{i,p}(p)$, an $h_i$-equivariant  
        embedding $\phi_i:V_i\ra V_p$, and an $h_i$-equivariant
        embedding of vector bundles 
        $\hat{\phi}_i: E_{i,p}|_{V_i} \ra E_p$ which covers $\phi_i$.
  \item $\hat{\phi}_i\circ s_{i,p}=s_p\circ\phi_i$ for $i=1,2$.
  \item $\psi_{i,p}=\psi_p\circ\phi_i$ for $i=1,2$.
 \end{enumerate} 
 In this case, we write  
 $(V_{1,p},E_{1,p},\Gamma_{1,p},\psi_{1,p},s_{1,p})\sim 
  (V_{2,p},E_{2,p},\Gamma_{2,p},\psi_{2,p},s_{2,p})$ 
 \end{df}

The following definition is a combination of 
\cite[A2.1.5-A2.1.11]{FO3} and \cite[Definition 5.3]{FO}.
\begin{df}[Kuranishi structure]\label{Kuranishi}
 Let $M$ be a Hausdorff topological space. A {\em Kuranishi structure (with
 corners)} on $M$ assigns a Kuranishi neighborhood (or a Kuranishi neighborhood with
 corners) $(V_p,E_p,\Gamma_p,\psi_p,s_p)$ to each $p\in M$ and a $4$-uple
 $(V_{pq}, \hat{\phi}_{pq},\phi_{pq},h_{pq})$ to each pair $(p,q)$ where 
 $p\in M, q\in\psi_p(s_p^{-1}(0))$ such that
 \begin{enumerate}
  \item $V_{pq}$ is an open subset of $V_q$ containing $\psi_q^{-1}(q)$.
  \item $h_{pq}$ is a homomorphism $\Gamma_q\ra \Gamma_p$.
  \item $\phi_{pq}:V_{pq}\ra V_p$ is an $h_{pq}$-equivariant embedding.
  \item $\hat{\phi}_{pq}: E_q|_{V_{pq}}\ra E_p$ is 
        an $h_{pq}$-equivariant embedding of vector bundles
        which covers $\phi_{pq}$.
  \item $\hat{\phi}_{pq}\circ s_q=s_p\circ\phi_{pq}$.
  \item $\psi_q=\psi_p\circ\phi_{pq}$.
  \item If $r\in \psi_q(s_q^{-1}(0)\cap V_{pq})$, then 
        $\hat{\phi}_{pq}\circ \hat{\phi}_{qr}=\hat{\phi}_{pr}$ in a 
        neighborhood of $\psi_r^{-1}(r)$.
  \item $\dim V_p -\mathrm{rank}\,E_p$ is independent of $p$ and is called the 
        {\em virtual dimension} of the Kuranishi structure (with corners).
  \end{enumerate}
$(V_{pq},\hat{\phi}_{pq},\phi_{pq},h_{pq})$ 
is called a {\em transition function} from $(V_q,E_q,\Gamma_q,\psi_q,s_q)$
to $(V_p,E_p,\Gamma_p,\psi_p,s_p)$. 
\end{df}

\begin{rem}\label{boundaryku}
Let $M$ be a Hausdorff space with a Kuranishi structure with corners
\[ 
  \mathcal{K}=\left\{(V_p,E_p,\Gamma_p,\psi_p,s_p):p\in M,
  (V_{pq}, \hat{\phi}_{pq},\phi_{pq},h_{pq}):q\in \psi_p(s_p^{-1}(0)) \right\}
\]
of virtual dimension $d$. Let $\pa M=\cup_{p\in M}\psi_p(s_p^{-1}(0)\cap\pa
V_p)$, where $\pa V_p$ is the union of corners in $V_p$. Then
\[ 
  \pa \mathcal{K}=\left\{(\pa V_p,E_p|_{\pa V},\Gamma_p,\psi_p,s_p):p\in \pa M,
  (\pa V_{pq}, \hat{\phi}_{pq},\phi_{pq},h_{pq}):q\in
  \psi_p(s_p^{-1}(0)\cap\pa V_p) \right\}
\]
is a Kuranishi structure with corners of virtual dimension $d-1$ on $\pa M$.
\end{rem}

\begin{df} \label{sameku}
Let $M$ be a Hausdorff topological space.
Two Kuranishi structures 
\[
  \mathcal{K}_1=\left\{ 
  (V_{1,p},E_{1,p},\Gamma_{1,p},\psi_{1,p},s_{1,p}):p\in M,
  (V_{1,pq}, \hat{\phi}_{1,pq},\phi_{1,pq},h_{1,pq}):
   q\in \psi_{1,p}(s_{1,p}^{-1}(0)) \right\}
\]
and 
\[
  \mathcal{K}_2=\left\{ 
  (V_{2,p},E_{2,p},\Gamma_{2,p},\psi_{2,p},s_{2,p}):p\in M,
  (V_{2,pq}, \hat{\phi}_{2,pq},\phi_{2,pq},h_{2,pq}):
  q\in \psi_{2,p}(s_{2,p}^{-1}(0)) \right\}
\]
on $M$ are \emph{equivalent} if there is another Kuranishi structure
\[
  \mathcal{K}=\left\{
  (V_p,E_p,\Gamma_p,\psi_p,s_p): p\in M,
  (V_{pq},\hat{\phi}_{pq},\phi_{pq},h_{pq}): 
  q\in \psi_p(s_p^{-1}(0)) \right \}
\]
on $M$ such that for all $p\in M$, 
$(V_{1,p},E_{1,p},\Gamma_{1,p},\psi_{1,p},s_{1,p})$,
$(V_{2,p},E_{2,p},\Gamma_{2,p},\psi_{2,p},s_{2,p})$, and
$(V_p,E_p,\Gamma_p,\psi_p,s_p)$ satisfy the relation described in
Definition~\ref{same}. In this case, we write
$\mathcal{K}_1\sim\mathcal{K}_2$.
\end{df} 

\bigskip

Let $(V_p,E_p,\Gamma_p,\psi_p,s_p)$ be a Kuranishi neighborhood (with corners)
of $p$. If $s_p$ intersects the zero section of $E_p$ transversally, then 
$\tilde{M}_p=s_p^{-1}(0)$ is a smooth submanifold (with corners) of $V_p$ of dimension 
$\dim V_p -\mathrm{rank} E_p$,
and there is an exact sequence of smooth vector bundles
\[
0\ra T\tilde{M}_p\ra TV_p|_{\tilde{M}_p} \stackrel{ds_p}{\ra} E_p|_{\tilde{M}_p}
\]
over $\tilde{M}_p$. So $T\tilde{M}_p$ is equivalent to the two term complex
$[TV_p|_{\tilde{M}_p} \stackrel{ds_p}{\ra} E_p|_{\tilde{M}_p}]$
as an element of the Grothendieck group $KO(\tilde{M}_p)$.

Both $TV_p$ and $E_p$ are $\Gamma_p$-equivariant vector bundles
over $V_p$, so $TV_p/\Gamma_p$, $E_p$ are orbibundles over
the orbifold (with corners) $U_p=V_p/\Gamma_p$. We call
$TU_p\equiv TV_p/\Gamma_p$ the tangent bundle of the orbifold (with corners)
$U_p$, and $TM_p\equiv T\tilde{M}_p/\Gamma_p$ is the tangent bundle of the orbifold
(with corners) $M_p=\tilde{M}_p/\Gamma_p$.

In general, $\tilde{M}_p$ might be singular, so $T\tilde{M}_p$ does not
exist. Nevertheless, $\tilde{M}_p$ is a topological space, so  
$KO(\tilde{M}_p)$ makes sense. We define
\[
T^{\mathrm{vir}}\tilde{M}_p = 
[TV_p|_{\tilde{M}_p} \stackrel{ds_p}{\ra} E_p|_{\tilde{M}_p}]\in KO(\tilde{M}_p)
\]
to be the {\em virtual tangent bundle} of $\tilde{M}_p$, and 
\[
T^{\mathrm{vir}}M_p = 
[TU_p|_{M_p} \stackrel{ds_p}{\ra} (E_p/\Gamma_p)|_{M_p}]
\]
to be the {\em virtual tangent bundle} of $M_p$.
We have $T^\mathrm{vir}\tilde{X}_p=T\tilde{M}_p$ and 
$T^\mathrm{vir}M_p=TM_p$ when $s_p$ intersects the zero section transversally.
The transition functions $(V_{pq},\hat{\phi}_{pq},\phi_{pq},h_{pq})$
in Definition~\ref{Kuranishi} enable
us to glue $T^{\mathrm{vir}}M_p$ to obtain the {\em virtual tangent bundle}
$T^{\mathrm{vir}} M$ of the Kuranishi structure on $M$.
 
$\det TV_p\otimes (\det E_p)^{-1}$ glue to a real line orbibundle
$\det (T^{\mathrm{vir}} M)$, the {\em orientation bundle}
of the Kuranishi structure (with corners). 
It is a real line bundle if the action of each $\Gamma_p$ on
$\det TV_p\otimes (\det E_p)^{-1}$ is orientation preserving. 
We say a Kuranishi structure is \emph{orientable} if its orientation
bundle is a trivial real line bundle. 
If $\mathcal{K}$ and $\mathcal{K}'$ are equivalent Kuranishi structures
(with corners), then $\mathcal{K}$ is orientable if and only if 
$\mathcal{K}'$ is.

In the ordinary Gromov-Witten theory, there is an algebraic approach to
define Gromov-Witten invariants when the target is a smooth projective
variety \cite{BF}. In the algebraic approach, the moduli of stable maps is a
{\em Deligne-Mumford stack}, which is locally \'{e}tale covered by affine schemes.
In Definition~\ref{Kuranishi}, $s_p^{-1}(0)$ is the analogue of an affine scheme
-- an affine scheme is the zero locus of polynomials, while $s_p^{-1}(0)$ is
the zero locus of smooth functions.

The moduli space of stable maps admits a {\em perfect obstruction theory}, which
is an element in the derived category locally isomorphic to a two term complex 
of vector bundles $[E_{-1}\ra E_0]$. Given a perfect obstruction theory, the
{\em virtual dimension}
is defined to be $\mathrm{rank}\,E_0-\mathrm{rank}\,E_{-1}$, and a virtual 
fundamental class of the virtual dimension can be constructed.
The two term complex  $[TV_p\stackrel{ds_p}{\ra} E_p]$ is the analogue of 
$[E_0^\vee\ra E_{-1}^\vee ]$. 

A Kuranishi structure can be viewed as the
analytic counterpart of a Deligne-Mumford stack together with a perfect
obstruction theory. We will show that 
 
\begin{tm} \label{kstructure}
$\MXL$ has a Kuranishi structure of virtual dimension
$$
\mu+(N-3)(2-2g-h)+2n+m^1+\cdots +m^h,
$$ 
where $2N$ is the dimension of $X$. The Kuranishi structure is orientable if 
$L$ is spin or if $h=1$ and $L$ is relatively spin (i.e., $L$ is orientable
and $w_2(TL)=\alpha|_L$ for some $\alpha\in H^2(X,\Z_2)$).
\end{tm}

\subsection{Stable $W^{k,p}$ maps} \label{virtual}

Let $(X,\omega)$ be a compact symplectic manifold together with an
$\omega$-tame almost complex structure $J$, and let $L$ be a compact
Lagrangian submanifold of $X$, as before. To construct a Kuranishi
structure on $\MXL$, we need to enlarge the category of stable maps.
We first specify metrics on the target (Section~\ref{metrictarget}) and on the domain
(Section~\ref{metricdomain}) so that we can define norms on relevant Banach spaces.
The definition of stable $W^{k,p}$ maps is given in Section~\ref{wkp}.
The virtual dimension in Theorem~\ref{kstructure}
is computed in Section~\ref{virtualdim}.

\subsubsection{Metric on the target} \label{metrictarget}
 Let $g_0$ be the Riemannian 
metric on $X$ defined by $g_0(v,w)=\half(\omega(v,Jw)+\omega(w,Jv))$. 
We will modify $g_0$ to obtain a Riemannian metric $g_1$ such that $L$ 
is totally geodesic w.r.t. $g_1$.

\begin{lm} \label{metric}
Given a Riemannian vector bundle $(V,h)$ over a compact Riemannian manifold
$(M,g)$, there is a Riemannian metric $\tg$ on the total space of $V$
such that
\begin{enumerate}
\item For any $x\in M$, the restriction of $\tg$ to the fiber $V_x$ over $x$
      is $h(x)$.
\item The zero section $i_0: (M,g)\ra(V,\tg)$ is an isometric embedding.
\item $i_0(M)$ is totally geodesic in $(V,\tg)$.
\end{enumerate}
\end{lm}

\paragraph{Proof.}Let $\pi:V\ra M$ be the canonical projection.
Choose a connection on $V$ which is compatible with $h$. This gives
a decomposition $TV= \pi^*V\oplus H$, where $H\cong \pi^* TM$.
Let $x\in M$, $w\in V_x$, so that $(x,w)\in V$.
Given $\xi\in T_{(x,w)}V$, there is a unique decomposition
$\xi=\xi_v+\xi_h$, where $\xi_v\in \pi^*V$, and $\xi_h\in H$. Define
a quadratic form $Q$ on $T_{(x,w)}V$ by $Q(\xi,\xi)= h(x)(\xi_v,\xi_v)
+ g(x)(\pi_*(\xi_h),\pi_*(\xi_h))$. Then $Q$ determines an inner product
$\tg(x,w)$ on $T_{(x,w)}V$. $\tg$ is a Riemannian metric on $V$ which
clearly satisfies 1. and 2.

For $(3)$, let $x_0,x_1\in i_0(M)$ be close enough such that there is a unique
length minimizing geodesic $\gamma:[0,1]\ra (V,\tg)$
such that $\gamma(0)=x_0$, $\gamma(1)=x_1$.  
It suffices to show that this geodesic lies in $i_0(M)$. 
We have
\begin{eqnarray*}
l(\gamma)&=&\int_0^1\tg(\gamma',\gamma')dt\\
&=&\int_0^1( h(\pi\circ\gamma)(\gamma'_v,\gamma'_v) 
       + g(\pi\circ\gamma)((\pi\circ\gamma)',(\pi\circ\gamma)'))dt\\
&\geq& \int_0^1 g(\pi\circ\gamma)((\pi\circ\gamma)',(\pi\circ\gamma)'))dt\\
&=&\int_0^1 \tg((i_0\circ\pi\circ\gamma)',(i_0\circ\pi\circ\gamma)') dt\\
&=& l(i_0\circ\pi\circ\gamma),
\end{eqnarray*}
The equality holds since $\gamma$ is length minimizing.
So $\gamma'_v\equiv 0$, and $\gamma=i_0\circ\pi\circ\gamma:[0,1]\ra L$. 
$\Box$

\bigskip

The Riemannian metric $g_0$ on $X$ gives an orthogonal decomposition
\[ 
TX|_L= TL\oplus N_{L/X}, 
\]
where $N_{L/X}$ is the normal bundle of $L$ in X. Let 
$\exp^0$ denote the exponential map $TX\ra X$ determined by $g_0$. For
$R>0$, let $B_R(TX)$ denote the ball bundle of radius $R$ in $TX$. There
exists $R>0$ such that $\exp^0$ maps $B_R(N_{L/X})$ diffeomorphically to 
its image in $X$. For $r\leq R$, let $N_r(L)$ denote the image of $B_r(N_{L/X})$
under $\exp^0$. We have a diffeomorphism $G:N_R(L)\ra B_R(N_{L/X})$ which
is the inverse of $\exp^0|_{B_R(N_{L/X})}$.
Construct a Riemannian metric $\tg$ on $N_{L/X}$ as in Lemma~\ref{metric}.
Let $\chi$ be a smooth cut-off function defined on $X$ such that
$\chi=1$ on $N_\frac{R}{3}(L)$ and $\chi=0$ on $X-N_\frac{2R}{3}(L)$.
Define $g_1=\chi G^* \tg  + (1-\chi)g_0$ on $N_R(L)$. Then
$g_1=g_0$ on $N_R(L)-N_\frac{2R}{3}(L)$, so $g_1$ extends to a Riemannian metric
on $X$ such that $g_1=g_0$ on $X-N_R(L)$ and on $TX|_L$.

$X$ is compact, so there exists some constant $C_0>0$
such that $C_0^{-1}g_0 \leq g_1 \leq C_0 g_0$.
The $C^\infty$ topology (Definition~\ref{Cl}) defined by $g_1$ is equivalent
to that defined by $g_0$. From now on, all the parallel transports,
exponential maps, and norms are defined by $g_1$ instead of $g_0$.

\subsubsection{Metric on the domain} \label{metricdomain}

Let $\lambda=[\domain]\in \widetilde{M}_{(g,h),(n,\vm)}$, and let
$\MXL_\lambda$ denote the fiber of 
\[ 
F:\MXL\ra \widetilde{M}_{(g,h),(n,\vm)} 
\]
over $\lambda$. Choose a Hermitian metric $\tilde{h}$ on $\Si_\C$ which is compact, flat
near nodes, and invariant under the antiholomorphic involution 
$\si:\Si_\C\ra\Si_\C$. Let $h$ be the restriction of $\tilde{h}$ to
$\Si$. Then the border curves of $\Si$ are geodesics in the Riemannian
metric determined by $h$. We further require that
\begin{enumerate}
\item If $s$ is an interior node, then there is an isometric holomorphic
      embedding $B_\ep(r)\ra\C^2$, where $\C^2$ is equipped with the 
      standard metric, such that the image is 
      $ \{(x,y)\in\C^2\mid xy=0,|x|<\ep,|y|<\ep \}$.
\item If $s$ is a boundary node of type H, then there is an isometric holomorphic
      embedding $B_\ep(r)\ra\C^2/A$, where $A(x,y)=(\bar{x},\bar{y})$,
      such that the image is 
      $ \{(x,y)\in\C^2\mid xy=0,|x|<\ep,|y|<\ep\}/A$.
\item $h$ is invariant under $\mathrm{Aut}\,\rho$.
\end{enumerate}
We call such a Hermitian metric an \emph{admissible metric}. 

\subsubsection{$W^{k,p}$ maps and $C^l$ maps} \label{wkp}

\begin{df} Let $\Si$ be a prestable bordered Riemann surface.
A continuous map $u:(\Si,\bS)\ra (X,L)$ is a {\em $W^{k,p}$ map on $\lambda$}
if $\hat{u}=u\circ\tau:(\hS,\pa\hS)\ra (X,L)$ is of class $W^{k,p}$ in the
sense of \cite[Appendix B]{MS}, where $\tau:\hS\ra \Si$ is the 
normalization map. 
\end{df}

In the above definition, we assume that $kp>2$, so the embedding
$W^{k,p}\subset C^0$ is compact.

 \begin{df}
 A {\em prestable $W^{k,p}$ map} of type $(g,h)$ with $(n,\vm)$ marked points  
  consists of a prestable marked bordered Riemann surface of type $(g,h)$
  with $(n,\vm)$ marked points $\domain$
  and a prestable $W^{k,p}$ map $u:(\Si,\bS)\ra (X,L)$. 
 \end{df}

 \begin{df}
  A \emph{morphism} between prestable $W^{k,p}$ maps of type $(g,h)$ with $(n,\vm)$
  marked points 
\[ 
\map \to \anothermap
\]
 is an isomorphism  
\[  
  \phi: \domain\to \another
\]
 between prestable bordered Riemann surfaces
 of type $(g,h)$ with $(n,\vm)$ points such that $u=u'\circ\phi$.
 \end{df}

\begin{df}
  A prestable $W^{k,p}$ map of type $(g,h)$ with $(n,\vm)$ marked points
  is \emph{stable} if its automorphism group is finite.
\end{df}

$C^l$ maps and stable $C^l$ maps are defined similarly. 

\bigskip

Let $\WXL, \CXL$ be the moduli space of isomorphism classes of 
stable $W^{k,p}$, $C^l$ maps of type $(g,h)$ with $(n,\vm)$ marked points
satisfying the topological conditions as in the definition of $\MXL$,
respectively. There are maps
\[
  F^{k,p}:\WXL\ra\widetilde{M}_{(g,h),(n,\vm)}
\]
and
\[
  C^l:\CXL\ra\widetilde{M}_{(g,h),(n,\vm)}
\] 
given by forgetting the map. Recall that forgetting the map also gives
\[
  F:\MXL\ra\widetilde{M}_{(g,h),(m,n)}. 
\]
Given $\lambda\in\widetilde{M}_{(g,h),(m,n)}$, let
$\MXL_\lambda$, $\WXL_\lambda$, $\CXL_\lambda$ denote the fiber 
of $F$, $F^{k,p}$, $F^l$ over $\lambda$, respectively.
From now on, we will write
\begin{eqnarray*}
M_\lambda&=&\MXL_\lambda\\
W^{k,p}_\lambda&=&\WXL_\lambda\\
C^l_\lambda&=&\CXL_\lambda
\end{eqnarray*}
for convenience.

Let $\exp$ denote the exponential map of the Riemannian metric $g_1$ on X,
and let $h$ be an admissible metric on $\Si$.
For a stable $W^{k,p}$ map $u$ on $\lambda$ and $\ep>0$, define
\[
  U^{k,p}(u,\ep)=\{\exp_u(w)\mid w\in\Wl, \parallel w\parallel_{W^{k,p}}<\ep\},
\]
where $\Wl$ will be defined in Section~\ref{virtualdim}, and 
the norms are defined by $g_1$, $h$.
Similarly, for a stable $C^l$ map $u$ on $\lambda$ and $\ep>0$, define
\[
  U^l(u,\ep)=\{\exp_u(w)\mid w\in\Cl, \parallel w\parallel_{C^l}<\ep\},
\]  
where $\Cl$ will be defined in Section~\ref{virtualdim}.
Note that $\exp_u(w)|_{\bS}\subset L$ because $L$ is totally geodesic
w.r.t. $g_1$.
Then
\[ 
\{ U^{k,p}(u,\ep)\mid u\textup{ is a stable }W^{k,p} 
   \textup{ map on }\lambda,\ep>0 \}
\]
generate the $W^{k,p}$ topology on $W^{k,p}_\lambda$, while
\[ 
\{ U^l(u,\ep)\mid u\textup{ is a stable }C^l
   \textup{ map on }\lambda,\ep>0 \}
\]
generate the $C^l$ topology on $C^l_\lambda$.

Define
\[
W^{k,p}(u,\ep)=\{ w\in\Wl\mid\; \parallel w\parallel_{W^{k,p}}<\ep\}
\]
and
\[
C^l(u,\ep)=\{ w\in\Cl\mid\; \parallel w\parallel_{C^l}<\ep\}.
\]
For sufficiently small $\ep>0$, $w\mapsto \exp_u(w)$ gives 
\[
U^{k,p}(u,\ep)\cong W^{k,p}(u,\ep)/\mathrm{Aut}(\lambda,u)
\]
and
\[
U^l(u,\ep)\cong C^l(u,\ep)/\mathrm{Aut}(\lambda,u).
\]
So $W^{k,p}_\lambda$ and $C^l_\lambda$ are Banach orbifolds.
They are Banach manifolds if $\mathrm{Aut}\lambda$ is trivial.
$M_\lambda$ is contained in both $W^{k,p}_\lambda$ and $C^l_\lambda$
since stable maps are $C^\infty$ (\cite[Theorem 2.1]{Ye}).

\subsubsection{Virtual dimension} \label{virtualdim}
Let $\rho=[\map]\in M_\lambda$.
Let $C_1,\ldots,C_\nu$ be irreducible components of $\Si$ which are 
(possibly nodal) Riemann surfaces, and let $\Si_1,\ldots,\Si_{\nu'}$
be the remaining irreducible components of $\Si$, which are 
(possibly nodal) bordered Riemann surfaces. 
Let $\hC_i$ denote the normalization of $C_i$ for $i=1,\ldots, \nu$,
let $\hS_{i'}$ denote the normalization of $\Si_{i'}$ for  
$i'=1,\ldots, \nu'$, and let $\tau:\hS\ra\Si$ be the normalization map.

Let $r_1,\ldots,r_{l_0}\in\inS$ be interior nodes of $\Si$, and 
$s_1,\ldots,s_{l_1}\in\bS$ be boundary nodes of $\Si$. 
Let $\hp_j\in\hS$ be the preimage of $p_j$ under $\mu$ for
$j=1,\ldots,n$, $\hq_{j'}$ be the preimage of $q_{j'}$ under
$\tau$ for $j'=1,\ldots,m$, $\hr_\alpha,\hr_{l_0+\alpha}$
be the preimages of $r_\alpha$ under $\mu$ for $\alpha=1,\ldots, l_0$, 
and $\hs_{\alpha'},\hat{s}_{l_1+\alpha'}$ be the preimages
of $s_{\alpha'}$ under $\tau$ for $\alpha'=1,\ldots,l_1$.
  
Set $\hu=u\circ\tau$. Then $\hu:(\hS,\pa \hS)\ra (X,L)$ is $J$-holomorphic. 
For $l>1$, let 
\[ 
   C^l(\hS,\pa\hS,\hu^*TX, (\hu|_{\pa\hS})^*TL) 
\]
   denote the vector space of $C^l$ sections 
   of $\hu^*TX$ with boundary values in $(\hu|_{\pa\hS})^*TL$. Then
  \begin{eqnarray*}
  && C^l(\hS,\pa\hS,\hu^*TX, (\hu|_{\pa\hS})^*TL)\\ 
  &=&\bigoplus_{i=1}^\nu C^l(\hC_i,(\hu|_{\pa\hC_i})^*TX)\oplus
     \bigoplus_{i'=1}^{\nu'} C^l(\hS_{i'},\pa\hS_{i'},
      (\hu|_{\hS_{i'}})^*TX, (\hu|_{\pa\hS_{i'}})^*TL)
   \end{eqnarray*}
  
   Let $\Cl$ be the kernel of
\[
   C^l(\hS,\pa\hS,\hu^*TX, (\hu|_{\pa\hS})^*TL)\longrightarrow
   \bigoplus_{\alpha=1}^{l_0} T_{r_\alpha}X \oplus 
   \bigoplus_{\alpha'=1}^{l_1}T_{s_{\alpha'}}L 
\]
\[
  s \mapsto \left(\{ s(\hr_\alpha)-s(\hr_{l_0+\alpha})\}_{\alpha=1}^{l_0},
               \{ s(\hs_{\alpha'})-s(\hs_{l_1+\alpha'}) \}_{\alpha'=1}^{l_1} \right),
\]
 and set 
\begin{eqnarray*}
 && \Cless\\
 &=& \bigoplus_{i=1}^\nu C^{l-1} (\hC_i,\form\hC_i\otimes (\hu|_{\hC_i})^*TX) \oplus
    \bigoplus_{i'=1}^{\nu'} C^{l-1} (\hS_{i'},\form\hS_{i'}\otimes (\hu|_{\hS_{i'}})^*TX)
\end{eqnarray*} 

 The linearization of $\dbar_{J,\Si}$ at the map $u$ gives rise to operators
\[  
    C^l(\hC_i,(\hu|_{\hC_i})^*TX)\ra
    C^{l-1} (\hC_i,\form\hC_i\otimes (\hu|_{\hC_i})^*TX),
    \;\;i=1,\ldots,\nu,
\]

\[
    C^l(\hS_{i'},\pa\hS_{i'}, (\hu|_{\hS_{i'}})^*TX, (\hu|_{\pa\hS_{i'}})^*TL)\ra
    C^{l-1} (\hS_{i'},\form\hS_{i'}\otimes (\hu|_{\hS_{i'}})^*TX)\;\;i'=1,\ldots,\nu',
\]
  and thus the operator
\[
  D_u:\Cl \ra\Cless
\]
  Similarly, we have
\[
  D_u:\Wl\ra\Wless.
\]
  
The following proposition follows from straightforward computations.

\begin{pro}\label{linear}
 Let $\nabla$ be a connection on $TX$. Then   
\[ 
  D_u(\xi)=\half\left(\nabla\xi\circ du+J\circ\nabla \xi\circ du\circ j
  +\nabla_\xi J\circ du \circ j + T(\xi, du)+ J T(\xi,du\circ j)\right)
\] 
where $T(v,w)=\nabla_v w-\nabla_w v- [v,w]$ is the torsion of $\nabla$. 
\end{pro}

\begin{lm}\label{find}
 Let $\Si$ be a smooth bordered Riemann surface of type
 $(g,h)$, and $u:(\Si,\bS)\ra (X,L)$ be a $J$-holomorphic map. Then
\[
  D_u:\Dl\ra\Dless 
\]
 is a Fredholm operator of index $\mu+N(2-2g-h)$, where
 $\mu=\mu(u^*TX,(u|_{\bS})^*TL)$, and $2N=\dim_\R X$.
\end{lm}

\paragraph{Proof.} Let $\nabla$ be the Levi-Civita connection.
 $\nabla$ is torsion-free, so by Proposition~\ref{linear}, 
\[
  D_u(\xi)=\half\left(\nabla\xi\circ du+J\circ\nabla \xi\circ du \circ j
  +\nabla_\xi J\circ du \circ j \right).
\]
We have $D_u=D''+R$, where
\begin{eqnarray*}
D''(\xi)&=&\half(\nabla\xi\circ du+J\circ\nabla\xi\circ du\circ j)\\
R(\xi)&=&=\half \nabla_\xi J\circ du\circ j
\end{eqnarray*}
 
$D''$ defines a holomorphic structure on $u^* TX$ such that $D''=\dbar$. By
\cite[Theorem 3.4.2]{KL}, $D''$ is a Fredholm operator of index $\mu +N(2-2g-h)$. 
$R$ is a compact operator, so $D_u=D''+R$ is Fredholm of index 
$\mu+N(2-2g-h).\;\;\Box$

\begin{pro}\label{singfind}
 Let $u:(\Si,\bS)\ra(X,L)$ be a prestable map. Then
\[
  D_u:\Dl\ra\Dless 
\]
 is a Fredholm operator of index $\mu + N(1-\tg)$,
 where $\mu=\mu (u^*TX,(u|_{\bS})^*TL)$, and $\tg$ is the 
 arithmetic genus of $\Si_\C$.
\end{pro}
\paragraph{Proof.} 
 We use the above notation. Set
\begin{eqnarray*}
 C^0&=&\Dl\\
 C^1&=&\Dless\\
 \tilde{C}^0_i&=&C^\infty(\hC_i,(\hu|_{\hC_i})^*TX)\\
 \tilde{C}^1_i&=&C^\infty(\hC_i,\form\hC_i\otimes (\hu|_{\hC_i})^*TX)\\
 C^0_{i'}&=&C^\infty(\hS_{i'},\pa\hS_{i'}, (\hu|_{\hS_{i'}})^*TX,(\hu|_{\pa\hS_{i'}})^*TL)\\
 C^1_{i'}&=&C^\infty(\hS_{i'},\form\hS_{i'}\otimes (\hu|_{\hS_{i'}})^*TX)
\end{eqnarray*}
where $i=1,\ldots,\nu$, $i'=1,\ldots,\nu'$. The linearization of
$\dbar_{J,\Si}$ gives rise to Fredholm operators
$\tilde{D}_i:\tilde{C}^0_i\ra \tilde{C}^1_i$ for $i=1,\ldots,\nu$ and
$D_{i'}:C^0_{i'}\ra C^1_{i'}$ for $i'=1\ldots, \nu'$.
We have the following commutative diagram:
\[
\begin{CD}
0 @>>> C^0 @>>> 
\bigoplus_{i=1}^\nu \tilde{C}^0_i \oplus \bigoplus_{i'=1}^{\nu'} {C}^0_{i'}        
@>>> \bigoplus_{\alpha=1}^{l_0} T_{r_\alpha}X\oplus 
   \bigoplus_{\alpha'=1}^{l_1}T_{s_{\alpha'}}L @ >>> 0\\
 & & @V{D_u}VV @VDVV @VVV\\
0 @>>> C^1 @>>> 
\bigoplus_{i=1}^\nu \tilde{C}^1_i \oplus \bigoplus_{i'=1}^{\nu'} {C}^1_{i'}
@>>> 0 @>>> 0
\end{CD}
\]
where $D=\bigoplus_{i=1}^\nu \tilde{D}_i\oplus \bigoplus_{i'=1}^{\nu'}
D_{i'}$, and the rows are exact. So $D_u$ is Fredholm.

Given a Fredholm operator $D$, let 
$\mathrm{Ind}(D)$ denote the virtual real vector space
\[ 
\mathrm{Ker}(D) -\mathrm{Coker}(D), 
\]
whose dimension
\[ 
\dim \mathrm{Ind}(D)=\dim \mathrm{Ker}(D) -\dim \mathrm{Coker}(D) 
\]
is the Fredholm index of $D$. With the above notation, we have
\[
\dim\mathrm{Ind}(D_u)=\sum_{i=1}^\nu \dim\mathrm{Ind}(\tilde{D}_i) 
                    + \sum_{i'=1}^{\nu'}\dim\mathrm{Ind}(D_{i'})-2N l_0-N l_1.
\]

Suppose that $\hC_i$ is of genus $\hat{g}_i$, and $\hS_{i'}$ is of type
$(g_{i'},h_{i'})$. We have
\begin{eqnarray*}
\dim\mathrm{Ind}(\tilde{D}_i)
 &=& 2\deg((\hu|_{\hC_i})^*TX) + 2N(1-\hat{g}_i),\\ 
\dim\mathrm{Ind}(D_{i'})
 &=&\mu((\hu|_{\hS_{i'}})^*TX,(\hu|_{\pa\hS_i})^*TL)+N(2-2g_{i'}-h_{i'}),
\end{eqnarray*}
where the second equality follows from Proposition~\ref{find}.
$\Si_\C$ has $2l_0+l_1$ nodes and $2\nu+\nu'$ irreducible components
\[  
  C_1,\ldots, C_\nu,\ \Bar{C}_1,\ldots,\Bar{C}_\nu,\ 
  (\Si_1)_\C,\ldots,(\Si_{\nu'})_\C,
\]
where the genus of $(\hS_{\nu'})_\C$ is $\tg_{i'}=2g_{i'}+h_{i'}-1$, so
the arithmetic genus of $\Si_\C$ is
\begin{eqnarray*}
\tg&=& 2\sum_{i=1}^\nu \hat{g}_i
   +\sum_{i'=1}^{\nu'}\tg_{i'}+2l_0+l_1-2\nu-\nu' +1 \\
   &=& 2\sum_{i=1}^\nu (\hat{g}_i -1) +\sum_{i'=1}^{\nu'}(2g_{i'}+h_{i'}-2)
       +2l_0+l_1+1
\end{eqnarray*}
by \cite[(3.1)]{HM}. Finally,
\[
\mu=\mu(u^*TX,(u|_{\bS})^* TL)= 2\sum_{i=1}^\nu \deg((\hu|_{\hC_i})^*TX
     +\sum_{i'=1}^{\nu'}\mu((\hu|_{\hS_{i'}})^*TX,(\hu|_{\pa\hS_{i'}})^*TL),
\]
so we conclude that
\[ \dim \mathrm{Ind}(D_u)=\mu(u^*TX,(u|_{\bS})^* TL)+N(1-\tg). \;\;\;\Box
\]

\begin{rem}
Corollary~\ref{singfind} remains true for  
\[ D_u:\Wl\ra\Wless \]
\[ D_u:\Cl\ra\Cless \]
\end{rem}

\bigskip
 
Set $\cE^\infty_u=\Dless$ for $u\in C^\infty_\lambda$. The $\cE^\infty_u$
fit together to form a Banach orbibundle 
$\cE^\infty_\lambda \ra C^\infty_\lambda$. There is a section
$s_J:C^\infty_\lambda \ra \cE^\infty_\lambda$, defined by $u\mapsto \dbar_{J,\Si}u$ ,
and $M_\lambda$ is the zero locus of $s_J$.
If $\lambda$ has no nontrivial automorphism, then
$C^\infty_\lambda$ is a Banach manifold, and $\cE^\infty_\lambda\ra C^\infty_\lambda$ is a 
Banach bundle. In this case, if $M_\lambda$ is nonempty and
$D_u$ is surjective for all $u\in M_\lambda$, then
$M_\lambda$ is a smooth manifold of dimension
$\mu+N(2-2g-h)$ by the implicit function theorem.
We call $\mu+N(2-2g-h)$ the \emph{virtual} dimension
of $M_\lambda$. In general, $M_\lambda$ is singular, and
the actual dimension of $M_\lambda$ can be larger
than the virtual dimension.

The dimension of $\widetilde{M}_{(g,h),(n,\vm)}$ is 
$6g+3h-6+2n+m^1+\cdots+ m^h$, so the expected (or virtual) dimension of $\MXL$ is  
\[ \mu + (N-3)(2-2g-h)+2n+m^1+\cdots +m^h, \]
which will also be the virtual dimension of the Kuranishi structure
on \\$\MXL$. 

Similarly, we have 
\begin{eqnarray*}
 \cE^l_\lambda\ra C^l_\lambda, &&  s_J: C^l_\lambda \ra \cE^l_\lambda,\\
 \cE^{k,p}_\lambda \ra W^{k,p}_\lambda,&&\ s_J: W^{k,p}_\lambda \ra \cE^{k,p}_\lambda,
\end{eqnarray*}
and $M_\lambda$ is the zero locus of $s_J$ in the above spaces.

\subsection{Deformation of the domain}\label{deformdomain}

Let 
\[ 
\rho=[\map] 
\] 
be a point in $\MXL$. We have seen in Section~\ref{nodaldeform}
that the infinitesimal deformation of the domain
\[ 
\lambda=[\domain]
\] 
is given by
\[ H_{\rho,\mathrm{domain}}=
 \bigoplus_{i=1}^\nu W_i\oplus\bigoplus_{i'=1}^{\nu'}\hat{W}_{i'}
       \oplus\bigoplus_{\alpha=1}^{l_0} V_\alpha 
       \oplus\bigoplus_{\alpha'=1}^{l_1}\hat{V}_{\alpha'}^+,
\]
where
\begin{eqnarray*}
 W_i&=&\W,\\
 \hat{W}_{i'}&=&H^1(\hS_{i'}, \pa \hS_{i'},
   T_{\hS_{i'}}(-p^{i'}_1-\cdots-p^{i'}_{n_{i'}}),
   T_{\pa \hS_{i'}}(-q^{i'}_1-\ldots-q^{i'}_{m_{i'}})),\\
 V_\alpha&=&T_{\hat{r}_\alpha} \hS \otimes T_{\hat{r}_{l_0+\alpha}} \hS \cong \C,\\
 \hat{V}_{\alpha'}&=&T_{\hat{s}_{\alpha'}} \pa\hS \otimes
 T_{\hat{s}_{l_1+\alpha'}}\pa\hS\cong \R,\\
 \hat{V}_{\alpha'}^+&\cong&[0,\infty)\subset \hat{V}_{\alpha'}
\end{eqnarray*}
are defined as in Section~\ref{nodaldeform}.  $\hat{V}_{\alpha'}^+$ is only a
semigroup, while the others are vector spaces. Set
\begin{eqnarray*}
 H_{\rho,\mathrm{deform}}&=&
  \bigoplus_{i=1}^\nu W_i\oplus\bigoplus_{i'=1}^{\nu'}\hat{W}_{i'}\\
 H_{\rho,\mathrm{interior}}&=&\bigoplus_{\alpha=1}^{l_0} V_\alpha\\
 H_{\rho,\mathrm{boundary}}&=&
 \hat{V}_1^+ \times\cdots\times \hat{V}_{l_1}^+\\
 H_{\rho,\mathrm{smooth}}&=&
 H_{\rho,\mathrm{interior}}\times H_{\rho,\mathrm{boundary}}
\end{eqnarray*}
Then $H_{\rho, \mathrm{deform}}$ corresponds to tangent directions of the
stratum that $\lambda$ belongs to, while
$H_{\rho,\mathrm{smooth}}$ corresponds to normal directions to this
stratum. $H_{\rho,\mathrm{interior}}$ corresponds to smoothing
of interior nodes, and   $H_{\rho,\mathrm{boundary}}$ corresponds
to smoothing of boundary nodes. 

Let $(\mathrm{Aut}\,\lambda)_0$ denote the identity component of
$\mathrm{Aut}\,\lambda$.
$(\mathrm{Aut}\,\lambda)_0$ is a normal subgroup of $\mathrm{Aut}\,\lambda$,
and the quotient $\mathrm{Aut}'\lambda=\mathrm{Aut}\,\lambda/(\mathrm{Aut}\,\lambda)_0$
is a  finite group. $\mathrm{Aut}\,\lambda$ acts on $H_{\rho,\mathrm{deform}}$, 
and $(\mathrm{Aut}\,\lambda)_0$ acts trivially, so 
$\mathrm{Aut}'\lambda$ acts on $H_{\rho,\mathrm{deform}}$.

We choose an admissible metric $h$ on $\Si$ in the sense of  
Section~\ref{metricdomain}.
Let $\ep_1$ be a small positive number, and
define $N_{\ep_1}(\Si)$, $K_{\ep_1}(\Si)$ as in the paragraph right before
Definition~\ref{Cl}. Then $N_{\ep_1}(\Si)$, $K_{\ep_1}(\Si)$ are invariant
under $\mathrm{Aut}\,\rho$.
We may choose a subspace $\tilde{H}$ of the space of smooth Beltrami differentials
such that the elements in $\tilde{H}$ vanish on $N_{\ep_1}(\Si)$, 
and the natural map $\tilde{H}\ra H_{\rho,\mathrm{deform}}$
is an isomorphism. We may further assume that
$\tilde{H}$ is invariant under the action of $\mathrm{Aut}'\rho$,
so that $\mathrm{Aut}'\rho$ acts on
$\tilde{H}$ and the isomorphism $\tilde{H}\ra H_{\rho,\mathrm{deform}}$
is  $\mathrm{Aut}'\rho$-equivariant. From now on, we will identify 
$\tilde{H}$ with $H_{\rho,\mathrm{deform}}$.

\subsubsection{Deformation within the stratum}\label{within}

Let $j(\xi)$ be the complex structure on $\Si$ determined by $\xi\in \tilde{H}$, 
and let $\Si_{(\xi,0,0)}$ be the prestable bordered Riemann surface
corresponding to $(\Si,j(\xi))$. In particular, $j(0)$ is the original
complex structure $j$ on $\Si$, and $\Si_{(0,0,0)}=\Si$. Let 
\[
\lambda_{(\xi,0,0)}=(\Si_{(\xi,0,0)},\bB;\bp;\bq^1,\ldots,\bq^h).
\] 
Let $\kappa_{(\xi,0,0)}:\Si_{(\xi,0,0)}\ra\Si$ be the identity map. Then 
\begin{enumerate}
\item $\kappa_{(\xi,0,0)}:\lambda_{(\xi,0,0)}\ra \lambda$ is a strong deformation 
      in the sense of Definition~\ref{contract}.  
      $\kappa_{(\xi,0,0)}: \Si_{(\xi,0,0)}\ra\Si$ is a homeomorphism.
\item $j=\left(\kappa_{(\xi,0,0)}^{-1}\right)^*j(\xi)$ on $N_{\ep_1}(\Si)$. 
\item $ \parallel j-\left(\kappa_{(\xi,0,0)}^{-1}\right)^*j(\xi)
        \parallel_{C^\infty(K_{\ep_1}(\Si))} < C|\xi|$, where 
      $|\xi|$ is the Weil-Petersson
      norm of the Beltrami differential $\xi$.
\end{enumerate}

Note that any two norms on $H_{\rho,\mathrm{deform}}$ are equivalent
since $H_{\rho,\mathrm{deform}}$ is finite dimensional.
Let $B_{\delta_2}\subset H_{\rho,\mathrm{deform}}$ be the ball of
radius $\delta_2>0$ centered at the origin. From the above discussion, we
see that there is a family
of prestable bordered Riemann surfaces of type $(g,h)$ with
$(n,\vm)$ marked points 
$\{ \lambda_{(\xi,0,0)} \mid \xi\in B_{\delta_2} \}$.
More precisely, we have
\[ 
(\pi:\cC\ra B_{\delta_2};\bs;\bt^1,\ldots,\bt^h)), 
\]
where $\bs=(s_1,\ldots,s_n)$, $\bt^i=(t^i_1,\ldots,t^i_{m^i})$,
and a contraction $\kappa:\cC\ra \Si$. 
Diffeomorphically, $\cC=B_{\delta_2}\times \Si$, $\pi$ is the projection
to the first factor, $\kappa$ is the projection of the second factor,
$s_j,t^i_k:B_{\delta_2}\to\cC$ are constant sections corresponding to 
marked points $p_j$, $q^i_k$, respectively.
Holomorphically, $\pi^{-1}(\xi)=\Si_{(\xi,0,0)}$. 

Let $\omega_0$ be the volume form on $\Si$ determined by $h$. Then
$\omega_0$ is a K\"{a}hler form. Let $h_{(\xi,0,0)}$ be the Hermitian 
metric on $\Si_{(\xi,0,0)}$ determined by $\kappa_{(\xi,0,0)}^*\omega_0$
and $j(\xi)$. Then $\kappa_{(\xi,0,0)}$ is an isometry
near the boundary and nodes, so $h_{(\xi,0,0)}$ is an admissible
metric on $\Si_{(\xi,0,0)}$. 
 
There is a map $i:B_{\delta_2}\ra \widetilde{M}_{(g,h),(n,\vm)}$,
given by $\xi\mapsto\lambda_{(\xi,0,0)}$. Given $\ep_2>0$, there
exists $\delta_{2}>0$ such that 
$i(B_{\delta_2})\subset U(\lambda,\ep_1,\ep_2)$, where
$U(\lambda,\ep_1,\ep_2)$ is the neighborhood of 
$\lambda$ in $\widetilde{M}_{(g,h),(n,\vm)}$ in the $C^\infty$ topology
defined in Definition~\ref{Cldomain}.
$i(B_{\delta_2})$ is a neighborhood of $\lambda$ in the stratum
of $\widetilde{M}_{(g,h),(n,\vm)}$ which $\lambda$ belongs to. 

\subsubsection{Smoothing of interior nodes}\label{interior}

Let $s$ be an interior node of $\Si$, and let $i:B_{\ep_1}(s)\ra \C^2$
be a holomorphic isometry such that 
$i(B_{\ep_1}(s))=\{(x,y)\in\C^2\mid xy=0, |x|<\ep_1, |y|<\ep_1 \}$.
Let $s_1,s_2\in\hS$ be the preimages of $s$ under $\tau:\hS\ra\Si$. 
Up to permutation of $s_1$, $s_2$, there exist unique
$e_1\in T_{s_1}\hS$, $e_2\in T_{s_2}\hS$  
such that $(i\circ\tau)_*(e_1)=(1,0)\in \C^2$, 
and $(i\circ\tau)_*(e_2)=(0,1)\in \C^2$.

Given $v\in T_{s_1}\hS\otimes T_{s_2}\hS$, we have
$v=t e_1\otimes e_2$ for some $t\in \C$. Suppose
that $t=r^2 e^{i\phi}$, where $0<r<\frac{\ep_1}{3}$.
Let $\Si_t$ be the bordered Riemann surface obtained from $\Si$ by
replacing 
\[
B_{\ep_1,0}=\{(x,y)\in\C^2\mid xy=0, |x|<\ep_1, |y|<\ep_1 \}
\]
with
\[
B_{\ep_1,t}=\{(x,y)\in\C^2\mid xy=t, |x|<\ep_1, |y|<\ep_1 \}
\]
More precisely, $B_{\ep_1,t}$ is obtained by identifying
\[
 x=re^{i\theta}\in \{x\in\C\mid |x|=r \}\subset \{x\in\C\mid r\leq|x|<\ep_1\} 
\]
with 
\[
 \frac{t}{x}= r e^{i(\phi-\theta)}\in \{y\in\C\mid |x|=r \}\subset \{y\in\C\mid r\leq|y|<\ep_1\}. 
\]
There is a strong deformation 
\[
\kappa_t:
B_{\ep_1,t}=\left\{\left(x,\frac{t}{x}\right)\left|\;
\frac{r^2}{\ep_1}<|x|<\ep_1\right.\right\}\ra B_{\ep_1,0}
\]
given by
\begin{eqnarray*}
\kappa_t\left(x,\frac{t}{x}\right)=&
\left(\chi(\left|\frac{x}{r}\right|^2)x,0\right) & \mathrm{if}\; r\leq|x|<\ep_1\\
 &\left(0,\chi(\left|\frac{t}{x r}\right|^2)\frac{t}{x}\right)& \mathrm{if}\;
\frac{r^2}{\ep_1}<|x|\leq r
\end{eqnarray*}
where $\chi:\R\ra [0,1]$ is a smooth function such that
\begin{enumerate}
 \item $0\leq \chi'(s) \leq 1$.
 \item $\chi(s)=0$ for $s\leq 1$.
 \item $\chi(s)>0$ for $s>1$.
 \item $\chi(s)=1$ for $s\geq 4$.
\end{enumerate}

\begin{lm}\label{glue}
Let $f, g$ be smooth functions on $\C$ such that
$f(0)=g(0)$. Let $F$ be the continuous function 
on $\{z\in\C \mid \frac{r^2}{\ep_1}<|z|<\ep_1 \}$ defined by
\[
F(z)=\left\{
\begin{array}{ll}
f(\chi(\left|\frac{z}{r}\right|^2)z) & \textup{ if } r\leq|z|<\ep_1\\
g\left(\chi(\left|\frac{t}{z r}\right|^2)\frac{t}{z}\right) & 
\textup{ if }\frac{r^2}{\ep_1}<|z|\leq r
\end{array}\right.
\]
Then $F$ is smooth.
\end{lm}

\paragraph{Proof.}
 The lemma follows from
\begin{enumerate}
\item Both $h_1(z)=f(\chi(\left|\frac{z}{r}\right|^2)z)$, 
 $h_2(z)=g\left(\chi(\left|\frac{t}{z r}\right|^2)\frac{t}{z}\right)$
 are smooth for $\frac{r^2}{\ep_1}<|z|<\ep_1$.
 \item The derivatives of $h_1$ of any order vanish when $|z|=r$, and
       the same is true for $h_2$. $\Box$  
\end{enumerate}

Let $A(r,R)$ denote the annulus $\{(u,v)\in\R^2\mid
r^2\leq u^2+v^2 \leq R^2\}$, and let $\AA (r,R)$ denote its interior.
Let $D_R$ denote the closed disc  $\{(u,v)\in\R^2\mid u^2+v^2\leq R^2\}$.
\begin{lm} \label{derivative}
 Let $f$, $g$, $F$ be defined as in Lemma~\ref{glue}.
 Define $h:\AA(\frac{r^2}{\ep_1},\ep_1)\ra \C$ by
 $h(u,v)=F(u+iv)$. Then
\begin{eqnarray*}
\max_{A(\frac{r}{2},2r)}|h|&=&\max\{ \max_{D_{2r}}|f|,\max_{D_{2r}}|g|, \}\\
\max_{A(\frac{r}{2},2r)}|\nabla h| 
 &\leq& 9\sqrt{2}\max\{\max_{D_{2r}}|f'|,4\max_{D_{2r}}|g'|\},
\end{eqnarray*}
where $|\nabla h|^2=|h_u|^2+|h_v|^2$.
\end{lm}

\paragraph{Proof.}
We have
\[
h(u,v)=\left\{\begin{array}{ll}
f(\chi(\frac{u^2+v^2}{r^2})(u+iv))&\textup{ if } r^2\leq
u^2+v^2 <\ep_1^2,\\
g\left(\chi(\frac{r^2}{u^2+v^2})\frac{t}{u+iv}\right)& \textup{ if }\
\frac{r^4}{\ep_1^2}< u^2+v^2 \leq r^2,
\end{array}\right.
\]
so
\[
\sup_{A(\frac{r}{2},2r)}|h| =\max\{ \max_{D_{2r}}|f|,\max_{D_{2r}}|g| \}.
\]

We also have
\[
h_u(u,v)=\left\{\begin{array}{ll}
f'(\chi(\frac{u^2+v^2}{r^2})(u+iv))(\chi'(\frac{u^2+v^2}{r^2})
 \frac{2u(u+iv)}{r^2} +\chi(\frac{u^2+v^2}{r^2}))
&\textup{ if } r^2<u^2+v^2 <\ep_1^2,\\
g'(\chi(\frac{r^2}{u^2+v^2})\frac{t}{u+iv})
(\chi'(\frac{r^2}{u^2+v^2})\frac{2r^2 u}{(u^2+v^2)^2}\frac{t}{u+iv}
 -\chi(\frac{r^2}{u^2+v^2})\frac{t}{(u+iv)^2})
&\textup{ if }\frac{r^4}{\ep_1^2}< u^2+v^2 < r^2,
\end{array}\right.
\]
and $h_u(u,v)=0$ for $u^2+v^2=r^2$, so 
\[
\sup_{A(\frac{r}{2},2r)} |h_u|\leq
\max\{9\max_{D_{2r}}|f'|,36\max_{D_{2r}}|g'|.\}
\]
Similarly,
\[
\sup_{A(\frac{r}{2},2r)} |h_v|\leq
\max\{9\max_{D_{2r}}|f'|,36\max_{D_{2r}}|g'|.\} \;\;\;\Box
\]

\bigskip

In particular, let $(f(x),g(y))=(x,0),(0,y)$. We see that
$\kappa_t$ is smooth as a map to $\C^2$.
$\kappa_t$ is a diffeomorphism when $|x|\neq r$, and
$\kappa_t^{-1}(0,0)=\{\left(x,\frac{t}{x})\mid |x|=r\right\}$.
Choose a Hermitian metric $h_t$ on $B_{\ep_1,t}$ such that it
is induced by inclusion in $\C^2$ on $B_{2r,t}$ and $\kappa_t$ is an 
isometry outside $B_{3r,t}$. 

We now have a family of prestable bordered Riemann surfaces of type
$(g,h)$ with $(n,\vm)$ marked points 
\[ 
\lambda_t=(\Si_t, \bB;\bp;\bq^1,\ldots,\bq^h)
\]
together with a family of admissible metrics $h_t$ on $\Si_t$ such that

\begin{enumerate}

\item There are strong deformations $\kappa_t:\lambda_t\ra\lambda$ such
      that on $K_{3r}(\Si)$, where $r=\sqrt{|t|}$,
      $\kappa_t^{-1}$ is defined and is an isometry.
\item $j= \left(\kappa_t^{-1}\right)^* j_t$ on $K_{3r}(\Si)$, where
      $j_t$ is the complex structure on $\Si_t$. 
\end{enumerate}

Let $D_{\ep_1^2/9}=\{te_1\otimes e_2\mid |t|<\ep_1^2/9\}\subset
T_{s_1}\hS\otimes T_{s_2}\hS$. The map
$D_{\ep_1^2/9}\ra\widetilde{M}_{(g,h),(n,\vm)}$ given by
$te_1\otimes e_2 \mapsto \lambda_t$ 
defines a parametrized curve in $\widetilde{M}_{(g,h),(n,\vm)}$
whose tangent line at $\lambda_0=\lambda$ is 
$T_{s_1}\hS\otimes T_{s_2}\hS\subset H_{\rho,\mathrm{interior}}$.

Let $\eta=(v_1,\ldots,v_{l_0})\in H_{\rho,\mathrm{interior}}$, where 
$v_\alpha\in V_\alpha= T_{r_\alpha}\hS\otimes T_{r_{l_0+\alpha}}\hS$.
Applying the above construction to each interior node on $\Si_{(\xi,0,0)}$, we
obtain 
\[
\lambda_{(\xi,\eta,0)}=(\Si_{(\xi,\eta,0)},\bB;\bp;\bq^1,\ldots,\bq^h).
\]
Given 
$0<d_1,\ldots,d_{l_0}<\ep_1^2/9$, 
let $D(d_1,\ldots,d_{l_0})$ denote the polydisc 
$D_{d_1}\times\ldots\times D_{d_{l_0}}$ in
$H_{\rho,\mathrm{interior}}$,
where $D_{d_\alpha}$ is the disc of radius $d_\alpha$ centered at
the origin in $V_\alpha$. We have a family
\[
(\pi:\cC\to B_{\delta_2}\times D(d_1,\ldots,d_{l_0}) ;
  \bs;\bt^1,\ldots,\bt^h)
\]
of prestable bordered Riemann surfaces of type $(g,h)$ with
$(n,\vm)$ marked points together with a family of admissible 
metrics $h_{(\xi,\eta,0)}$ on $\Si_{(\xi,\eta,0)}$. There is
a contraction $\kappa:\cC\ra \Si$ whose restriction to
$\pi^{-1}(\xi,\eta)=\lambda_{(\xi,\eta,0)}$ is a strong deformation 
$\kappa_{(\xi,\eta,0)}:\lambda_{(\xi,\eta,0)}\ra\lambda$ such that
\begin{enumerate}
 \item $\kappa_{(\xi,\eta,0)}^{-1}$ is defined on $K_{3\sqrt{|\eta|}}(\Si)$, and
       $\kappa_{(\xi,\eta,0)}^{-1}\circ\kappa_{(\xi,0,0)}$ is an isometry on
       $K_{3\sqrt{|\eta|}}(\Si_{(\xi,0,0)})$.  
 \item $\parallel j- \left(\kappa_{(\xi,\eta,0)}^{-1}\right)^* j(\xi,\eta)
        \parallel_{C^\infty(K_{\ep_1}(\Si))}
       = \parallel j- \left(\kappa_{(\xi,0,0)}^{-1}\right)^* j(\xi)
        \parallel_{C^\infty(K_{\ep_1}(\Si))} <C|\xi|$, where
         $j(\xi,\eta)$ is the complex structure on $\Si_{(\xi,\eta,0)}$.
\end{enumerate} 

\subsubsection{Smoothing of boundary nodes}\label{bnodes}

Let $s$ be a boundary node of $\Si$, and let $i:B_{\ep_1}(s)\ra \C^2$
be a holomorphic isometry such that 
\[
i(B_{\ep_1}(s))=\{(x,y)\in\C^2\mid xy=0, |x|<\ep_1,|y|<\ep_1,
\mathrm{Im}\,x\geq 0, \mathrm{Im}\,y\leq 0 \},
\]

Let $s_1,s_2\in\hS$ be the preimages of $s$ under $\tau:\hS\ra\Si$. 
Up to permutation of $s_1$, $s_2$, there exist unique
$e_1\in T_{s_1}\hS$, $e_2\in T_{s_2}\hS$  
such that $(i\circ\tau)_*(e_1)=(1,0)\in \C^2$, 
and $(i\circ\tau)_*(e_2)=(0,1)\in \C^2$.

Given $v\in T_{s_1}\pa\hS\otimes T_{s_2}\pa\hS$, we have
$v=t e_1\otimes e_2$ for some $t\in \R$.
We construct $\lambda_t$ as in Section~\ref{interior}. We have seen
in Section~\ref{nodaldeform} that there is
topological transition when $t$ changes sign.
We may assume that $\lambda_t\in \widetilde{M}_{(g,h),(n,\vm)}$
for $t\geq 0$.

Let $\eta'=(v'_1,\ldots,v'_{l_0})\in H_{\rho,\mathrm{boundary}}$, where
$v'_{\alpha'}\in \hat{V}_{\alpha'}^+$. Apply the above construction to 
each boundary node on $\Si_{(\xi,\eta,0)}$, we obtain
$\Si_{(\xi,\eta,\eta')}$ and $\lambda_{(\xi,\eta,\eta')}$.
Given $0<d'_1,\ldots,d'_{l_1} < \ep_1^2/9$, let
$D'(d'_1,\ldots,d'_{l_0})=[0,d'_1)\times\ldots\times [0,d'_{l_1})
\subset \hat{V}_1\times\ldots\times \hat{V}_{l_1}$. 
We have a universal family
\[ (\pi:\cC\ra B_{\delta_2}\times D(d_1,\ldots,d_{l_0})\times 
   D'(d'_1,\ldots,d'_{l_1});\bs;\bt^1,\ldots,\bt^h )
\]
of prestable bordered Riemann surfaces of type $(g,h)$ with $(n,\vm)$ marked
points together with a family of admissible metrics $h_{(\xi,\eta,\eta')}$
on $\Si_{(\xi,\eta,\eta')}$. There is a contraction
$\kappa:\cC\ra\Si$ whose restriction to $\pi^{-1}(\xi,\eta,\eta')$ is a
strong deformation $\kappa_{(\xi,\eta,\eta')}:
\lambda_{(\xi,\eta,\eta')}\ra\lambda$ such that
\begin{enumerate}
 \item $\kappa_{(\xi,\eta,\eta')}^{-1}$ is defined on
       $K_{3\sqrt{|\eta|+|\eta'|}}(\Si)$, and
       $\kappa_{(\xi,\eta,\eta')}^{-1}\circ\kappa_{(\xi,0,0)}$ is an isometry on
       $K_{3\sqrt{|\eta|+|\eta'|}}(\Si_{(\xi,0,0)})$.  
 \item $\parallel j- \left(\kappa_{(\xi,\eta,\eta')}^{-1}\right)^* j'
         \parallel_{C^\infty(K_{\ep_1}(\Si))} <C |\xi|$.
\end{enumerate} 

If we embed $\Si_\C$ in a complex projective space $\bP^N$, then
$\kappa_{(\xi,\eta,\eta')}\circ \tau$ is smooth as a map to $\bP^N$,
where $\tau:\hS_{(\xi,\eta,\eta')}\ra \Si_{(\xi,\eta,\eta')}$ is the 
normalization map.

Given $\ep_1,\ep_2>0$, choose $\delta_2$ as before. Suppose that
both $\max\{\sqrt{|d_\alpha|}\mid\alpha=1,\ldots,l_0\}$
and  $\max\{\sqrt{|d'_{\alpha'}|}\mid\alpha'=1,\ldots,l_1\}$ are less
than $\frac{\ep_1}{3}$. Then the image of the map
\[
 i: B_{\delta_2}\times D(d_1,\ldots,d_{l_0})\times D'(d'_1,\ldots,d'_{l_1})
 \ra \widetilde{M}_{(g,h),(n,\vm)}
\]
given by $(\xi,\eta,\eta')\mapsto \lambda_{(\xi,\eta,\eta')}$ lies in
the neighborhood $U(\lambda,\ep_1,\ep_2)$ of $\lambda$ in the $C^\infty$ topology. 

\subsubsection{Action of the automorphism group} \label{actauto}

We write $D_d$ for $D(d_1,\ldots,d_{l_0})$ and $D'_{d'}$ for
$D'(d'_1,\ldots,d'_{l_1})$. In this section, we study the action of
$\domaut$ on $B_{\delta_2}\times D_d\times D'_{d'}$ and the universal 
family over it.

We first consider deformation within the stratum. 
Let $\pi_{\delta_2}:\cC_{\delta_2}\ra B_{\delta_2}$ be the universal family, so
that $\pi_{\delta_2}^{-1}(\xi)=\lambda_{(\xi,0,0)}$.  
$\domaut$ acts on $B_{\delta_2}$ by
$j(\phi\cdot\xi)=(\phi^{-1})^* j(\xi)$. Therefore, it acts on the 
the universal family. Given $\phi\in\domaut$,
we have the following commutative diagram: 
\[
\begin{CD}
\cC_{\delta_2} @>\phi>> \cC_{\delta_2}\\
@V{\pi_{\delta_2}}VV @VV{\pi_{\delta_2}}V\\
B_{\delta_2}  @>\phi>> B_{\delta_2} 
\end{CD}
\]
$\phi:\lambda_{(\xi,0,0)}\ra \lambda_{(\phi\cdot\xi,0,0)}$ is an isomorphism.
In particular, if $\phi\in (\domaut)_0$, then
$\phi:B_{\delta_2}\ra B_{\delta_2}$ is the identity map, and
we have the following commutative diagram:
\[
\begin{CD}
\lambda_{(\xi,0,0)} @>{\phi}>> \lambda_{(\xi,0,0)}\\
@V\kappa_{(\xi,0,0)}VV @VV\kappa_{(\xi,0,0)}V\\
\lambda @>{\phi}>> \lambda
\end{CD}
\]

We now consider smoothing of nodes. Let 
$\pi_{\delta_2,d,d'}:\cC_{\delta_2,d,d'}\ra
B_{\delta_2}\times D_d\times D'_{d'}$ be the universal family,
so that $\pi_{\delta_2,d,d'}^{-1}(\xi,\eta,\eta')=\lambda_{(\xi,\eta,\eta')}$.
$\domaut$ acts on $D_d\times D'_{d'}$ by
$\phi\cdot(te_1\otimes e_2)=t\phi_* e_1\otimes \phi_* e_2$.
Given $\phi\in \domaut$, we have the following commutative diagram
\[
\begin{CD}
\cC_{\delta_2,d,d'} @>\phi>> \cC_{\delta_2,d,d'}\\
@V{\pi_{\delta_2,d,d'}}VV @VV{\pi_{\delta_2,d,d'}}V\\
B_{\delta_2}\times D_d\times D'_{d'}  @>\phi>> 
B_{\delta_2}\times D_d\times D'_{d'}  
\end{CD}
\]
$\phi:\lambda_{(\xi,\eta,\eta')}\ra 
\lambda_{(\phi\cdot\xi,\phi\cdot\eta,\phi\cdot\eta')}$ is an isomorphism.
For example, the prestable bordered Riemann surface $\Si$ in Figure~15 has
two interior nodes. The smoothing of the two interior nodes is parametrized
by $\eta=(\eta_1,\eta_2)$. Let $\Si_{(\eta_1,\eta_2)}$ be the corresponding
bordered Riemann surfaces obtained by smoothing the two interior nodes
on $\Si$. 

$\Si$ has an automorphism $\phi$ of order two which 
rotates Figure~19 by $180^\circ$. It acts on $\eta$ by 
$\phi\cdot(\eta_1,\eta_2)=(\eta_2,\eta_1)$ and gives an isomorphism
$\Si_{(\eta_1,\eta_2)}\ra \Si_{(\eta_2,\eta_1)}$ by rotating $180^\circ$. 
The case $\eta_2=0$ is shown in Figure~20.

\begin{figure}\label{nineteen}
\begin{center}
\psfrag{e1}{$\eta_1$}
\psfrag{e2}{$\eta_2$}
\includegraphics[scale=0.8]{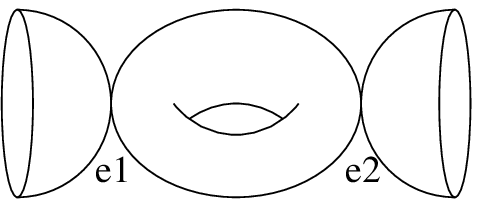}
\end{center}
\caption{$\Si=\Si_{(0,0)}$}
\end{figure}

\begin{figure}\label{twenty}
\begin{center}
\psfrag{S1}{$\Si_{(t,0)}$}
\psfrag{S2}{$\Si_{(0,t)}$}
\includegraphics[scale=0.8]{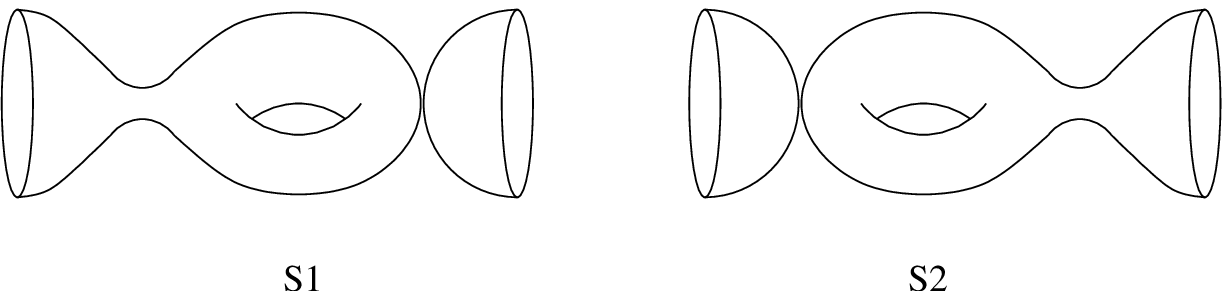}
\end{center}
\caption{$\Si_{(t,0)}$ and $\Si_{(0,t)}$}
\end{figure}
 
\subsection{Local Charts} \label{localchart}

Let \[ \rho=[\map] \]
be a point in $\MXL$, as at the beginning of  Section \ref{deformdomain}.
Consider
\[
  D_u:\Tl\ra\Tless
\]
for a large integer $p$.
$\Tless$ is a complex vector space, but $\Tl$ is only a real
vector space. $\mapaut$ acts on both $\Tl$ and $\Tless$, and 
$D_u$ is $\mapaut$-equivariant.

\begin{lm} 
 Let $\mathrm{Im}\,D_u$ denote the image of $D_u$.
 We can choose a subspace $E_\rho$ of $\Tless$ such that
 \begin{enumerate}
  \item $\mathrm{Im}\,D_u + E_\rho = \Tless$.
  \item $E_\rho$ is finite dimensional.
  \item Elements in $E_\rho$ are smooth sections supported in $K_{\ep_1}(\Si)$.
  \item $E_\rho$ is a complex subspace of $\Tless$.
  \item $E_\rho$ is $\mapaut$-invariant.
 \end{enumerate}
\end{lm}

\paragraph{Proof.}
We first claim that given any $\alpha\in \Tless$, there exists
$g\in \Tl$ such that $\alpha'=\alpha-D_u g$ has support in
$K_{\ep_1}(\Si)$. Actually, 
there exists $g'$ defined on $N_{2\ep_1}(\Si)$ such that
$D_u g'= \alpha$ on $N_{2\ep_1}(\Si)$. Let $\chi$ be a smooth function on $\Si$
which is $1$ on $N_{\ep_1}(\Si)$ and $0$ on $K_{2\ep_1}(\Si)$.
Let $g=\chi g'$ on $N_{2\ep_1}(\Si)$, and $0$ on
$K_{2\ep_1}(\Si)$. Then $\alpha'=\alpha-D_u g$ is 
supported in $K_{\ep_1}(\Si)$.

By Corollary~\ref{singfind}, we can find a finite dimensional subspace
$E'$ of \\$\Tless$ such that $\mathrm{Im}\,D_u \oplus E' =\Tless$.
We may assume that $E'$ consists of smooth sections since
any section in \\$\Tless$ can be approximated by sections in $\Dless$.
The above claim shows that we may assume that all sections in $E'$
have compact support in $K_{\ep_1}(\Si)$. 
Let $E_\rho$ be the smallest $\mapaut$-invariant
complex subspace which contains $E'$. Then
$E_\rho$ satisfies (1)--(5). $\Box$ 

\bigskip

Let $F_\rho=\Tless\cap E_\rho^\bot$, where $E_\rho^\bot$ is the orthogonal 
complement of $E_\rho$ in $L^2(\Sigma,\Lambda^{0,1}\Sigma\otimes u^*TX)$.
Then $F_\rho$ is a closed subspace of\\ $\Tless$, thus a Banach subspace of
$\Tless$. $F_\rho\cong \Tless/E_\rho$.
Let
\[ 
  \pi:\Tless\ra F_\rho
\]
be the $L^2$-orthogonal projection. (1) implies that $\pi\circ D_u$ is surjective.
Multiplication by $i$ preserves the $L^2$ inner product, so 
(4) implies that $F_\rho$ is a complex vector space. 
The action of $\mapaut$ also preserves the $L^2$ inner product,
so (5) implies that $\mapaut$ acts on 
$F_\rho$, and $\pi\circ D_u:\Tl\ra F_\rho$ is 
$\mapaut$-equivariant.

Set $H_{\rho,\mathrm{map}}=\mathrm{Ker}(\pi\circ D_u)$. We have 
\[
\dim H_{\rho,\mathrm{map}}=\mu+N(1-\tg)+\dim E_\rho,
\]
where $\mu=\mu(u^*TX,(u|_{\bS})^*TL)$, $2N=\dim_\R X$, and
$\tg$ is the arithmetic genus of $\Si_\C$, as before.

The infinitesimal deformation of the domain is given by 
\[ H_{\rho,\mathrm{aut}}=
 \bigoplus_{i=1}^\nu U_i \oplus \bigoplus_{i'=1}^{\nu'}\hat{U}_{i'},
\]
where 
\begin{eqnarray*} 
U_i&=&H^0(\hat{C}_i,T_{\hat{C}_i}(-x^i_1-\cdots-x^i_{\tn_i}))\\
\hat{U}_{i'}&=&H^0(\hS_{i'}, \pa \hS_{i'},
   T_{\hS_{i'}}(-p^{i'}_1-\cdots-p^{i'}_{n_{i'}}),
   T_{\pa \hS_{i'}}(-q^{i'}_1-\ldots-q^{i'}_{m_{i'}}))
\end{eqnarray*}

$U_i=0$ if and only $(\hat{C}_i,x^i_1,\cdots,x^i_{\tn_i})$ is stable,
and $\hat{U}_{i'}=0$ if and only if
$(\hS_{i'},(p^{i'}_1,\cdots,p^{i'}_{n_{i'}}),
 (q^{i'}_1,\ldots,q^{i'}_{m_{i'}}))$ is stable.

\bigskip

Let $\lambda=\domain$ as before. Then $H_{\rho,\mathrm{aut}}$
is the tangent space to $\domaut$ at the identity map.
$u$ is nonconstant on nonstable components, so 
$\phi\in \mathrm{Aut}\,\lambda\mapsto u\circ\phi^{-1}$ induces an inclusion
of vector spaces $H_{\rho,\mathrm{aut}}\subset H_{\rho,\mathrm{map}}$.
Let $H'_{\rho,\mathrm{map}}$ be the $L^2$-orthogonal
complement of $H_{\rho,\mathrm{aut}}$ in $H_{\rho,\mathrm{map}}$.
Set
\begin{eqnarray*}
H_\rho&=&H_{\rho,\mathrm{domain}}\times H_{\rho,\mathrm{map}},\\ 
H'_\rho&=&H_{\rho,\mathrm{domain}}\times H'_{\rho,\mathrm{map}}.
\end{eqnarray*}

With the above definitions, we are ready to state the main theorem of this section:

\begin{tm} \label{chart}
Let $\MXL$ be equipped with the $C^\infty$ topology.
There are a neighborhood $V'_\rho$ of $0$ in $H'_\rho$ such that 
$\mapaut$ acts on $V'_\rho$, an $\mapaut$-equivariant map
$s_\rho: V'_\rho \ra E_\rho$ 
such that $s_\rho(0)=0$, and a continuous map 
$\psi_\rho:s_\rho^{-1}(0)\ra \MXL$ such that 
$ s_\rho^{-1}(0)/\mapaut\ra \MXL $ gives a homeomorphism onto
a neighborhood of $\rho$ in $\MXL$.                                    
\end{tm}

$(V'_\rho,E_\rho,\mathrm{Aut}\,\rho,\psi_\rho, s_\rho)$ is a Kuranishi
neighborhood of $\rho$.

\subsubsection{Pregluing: construction of approximate solutions}
 \label{pregluing}

In this section, we will modify $u$ near nodes to obtain approximate
$J$-holomorphic maps 
\[
   u_{\eta,\eta'}:(\Si_{(0,\eta,\eta')},\pa\Si_{(0,\eta,\eta')})\ra (X,L),
\]
where the notation $\Si_{(\xi,\eta,\eta')}$ was introduced in Section 
\ref{within}, \ref{interior}, \ref{bnodes} and will be used repeatedly
in the rest of Section \ref{localchart}. 

We first consider a neighborhood of an interior node.
We will follow the construction in \cite[Appendix A]{MS} closely.
Recall that 
\[ B_{\ep_1,t}=\{(x,y)\in \C^2\mid xy=t,|x|<\ep_1,|y|<\ep_1\}. \]
When $t\neq 0$, we have
\[ B_{\ep_1,t}=\left\{\left(x,\frac{t}{x}\right)\left|x\in\C,
  \frac{r^2}{\ep_1}<|x|<\ep_1\right.\right\},
\]
where $|t|=r^2$. We assume that $r<\frac{\ep_1^2}{4}$.

Let $B_{\ep_1}=\{x\in\C \mid |x|<\ep_1\}$.
Two nonconstant $J$-holomorphic maps $f,g: B_{\ep_1}\ra X$ 
such that $f(0)=g(0)=p$ determine a stable map 
$u:B_{\ep_1,0}\ra X$ defined by $u(x,0)=f(x)$, $u(0,y)=g(y)$.
Suppose that the images of $f$, $g$ are contained
in the geodesic ball $B_{r_1}(x)$ with respect to $g_1$, where
$r_1$ is the injectivity radius of $(M,g_1)$. 
We define $u_t:B_{\ep_1,t}\ra X$ by
\[
u_t\left(z,\frac{t}{z}\right)
=\exp_p\left(\chi_1\left(\frac{z}{\sqrt{r}}\right)
           (\exp_p)^{-1}(f(z))
       +\chi_1\left(\frac{r\sqrt{r}}{z}\right)
           (\exp_p)^{-1}(g\left(\frac{t}{z}\right)) \right),  
\]
where $\chi_1:\C\ra [0,1]$ is a smooth cutoff function such that
\begin{eqnarray*}
\chi_1(z)&=&\left\{ \begin{array}{ll}1,& \textup{ if } |z|\geq 2\\
                0,& \textup{ if } |z|\leq 1\end{array}\right.\\
|\nabla\chi_1|&\leq& 2
\end{eqnarray*}
Then
\[
u_t\left(z,\frac{t}{z}\right)=\left\{\begin{array}{ll}
 u\left(0,\frac{t}{z}\right)=g(\frac{t}{z})& \textup{ if }
\frac{r^2}{\ep_1}<|z|<\frac{r\sqrt{r}}{2},\\
 p &\textup{ if } r\sqrt{r}\leq |z| \leq \sqrt{r},\\
 u(z,0)=f(z)& \textup{ if } 2\sqrt{r}<|z|<\ep_1.
\end{array}\right.
\]

Define $f_t,g_t:B_{\ep_1}\ra X$ by
\begin{eqnarray*} 
f_t(x)=\exp_p\left(\chi_1\left(\frac{x}{\sqrt{r}}\right)
           (\exp_p)^{-1}(f(x))\right)\\
g_t(y)=\exp_p\left(\chi_1\left(\frac{y}{\sqrt{r}}\right)
           (\exp_p)^{-1}(g(y))\right)\\
\end{eqnarray*}
Then 
\[
f_t(x)=\left\{\begin{array}{ll}
 f(x) & \textup{ if } 2\sqrt{r}<|x|<\ep_1\\
 p & \textup{ if }|x|<\sqrt{r}
\end{array}\right.
\]
\[
g_t(y)=\left\{\begin{array}{ll}
 g(y) & \textup{ if } 2\sqrt{r}<|x|<\ep_1\\
 p & \textup{ if }|x|<\sqrt{r}
\end{array}\right.
\]
$f_t,g_t$ determine a map $v_t:B_{\delta_1,0}\ra X$
defined by $v_t(x,0)=f_t(x)$, $v_t(0,y)=g_t(y)$.
Define $F_t,F:B_{\ep_1}\ra X$ by
$f_t(x)=\exp_p(F_t(x))$, $f(x)=\exp_p(F(x))$.
Then $F_t(0)=F(0)=0$, and 
\[
F_t(x)=\chi_1\left(\frac{x}{\sqrt{r}}\right)F(x)
\]
\begin{lm}
\[
\parallel F_t-F \parallel_{W^{1,p}(B_{\ep_1})} \leq C 
\max_{\bar{B}_{2\sqrt{r}}}|\nabla F|r^{\frac{1}{p}},
\]
where $C$ is a universal constant, and 
$\bar{B}_{2\sqrt{r}}=\{x\in\C\mid |x|\leq 2\sqrt{r}\}$.
\end{lm}

\paragraph{Proof.}
$F_t(x)-F(x)=0$ for $|x|\geq 2\sqrt{r}$, and for $|x|< 2\sqrt{r}$,
\begin{eqnarray*}
|F_t(x)-F(x)|&=& |(\chi_1\left(\frac{x}{\sqrt{r}}\right)-1)F(x)|
\leq |F(x)|\leq (\max_{\bar{B}_{2\sqrt{r}}} |\nabla F |)|x|\\ 
|\nabla(F_t(x)-F(x))|&=&
\left|\frac{1}{\sqrt{r}}\nabla\chi_1\left(\frac{x}{\sqrt{r}}\right)F(x)
+ \chi_1\left(\frac{x}{\sqrt{r}}\right)\nabla F(x)\right| \\
&\leq&\frac{2}{\sqrt{r}}|F(x)|+|\nabla F(x)|\leq 5\max_{\bar{B}_{2\sqrt{r}}}
|\nabla F|
\end{eqnarray*}
The desired estimate is obtained by integrating over $B_{2\sqrt{r}}$, and 
noting that $r<1$. $\Box$

\bigskip

We embed $(X,g_1)$ isometrically in $\R^l$ for some large $l$. Maps to 
$X$ can be viewed as maps to $\R^l$, so we may subtract one map from
another, and define their $L^{1,p}$ norms.

\begin{cor}\label{gluemap}
\[
\parallel v_t-u\parallel_{W^{1,p}(B_{\ep_1})} \leq C 
(\max_{\bar{B}_{2\sqrt{r}}}|\nabla u|)r^{\frac{1}{p}},
\]
where $C$ depends on the $C^\infty$ norms of
$\exp|_{B_{r_1}(N_{L/X})}$ and its inverse. 
\end{cor}

\begin{lm}\label{bar}
\[
\parallel \dbar_J u_t \parallel_{L^p} \leq C r^\frac{1}{p}
\]
where $C$ depends on the $C^1$ norm of $u$, the $C^1$ norm of $J$,
the $C^\infty$ norms of $\exp|_{B_{r_1}(N_{L/X})}$ and its inverse.
\end{lm}

\paragraph{Proof.}

For $r\leq|z|<\ep_1$, we have
\[
u_t\left(z,\frac{t}{z}\right)= f_t(z),
\]
so $\dbar_J u_t(z)=0$ for $|z|>\sqrt{r}$. 
For $r\leq|z|\leq\sqrt{r}$, we have
\[
\dbar_J u_t\left(z,\frac{t}{z}\right)
=\dbar_J f_t(z)= \dbar_J (f_t-f)(z).
\]
Let $z=x+iy$ ,then
\begin{eqnarray*}
&&2\dbar_J (f_t-f)\left(\frac{\pa}{\pa x}\right)\\
&=&\frac{\pa f_t}{\pa x} + J(f_t)\frac{\pa f_t}{\pa y}
   -\frac{\pa f}{\pa x} - J(f)\frac{\pa f}{\pa y}\\ 
&=&\frac{\pa}{\pa x}(f_t-f) + (J(f_t)-J(f))\frac{\pa f}{\pa y}
   +J(f)\frac{\pa}{\pa y}(f_t-f)+ (J(f_t)-J(f))\frac{\pa}{\pa y}(f_t-f)\\
&&\left|\dbar_J (f_t-f)\left(\frac{\pa}{\pa x}\right)\right| \\
&\leq& C_1(1+\sup|J|+\sup|\nabla J||f_t-f|)|\nabla(f_t-f)|+C_2
       \sup|\nabla J||f_t-f||\nabla f|\\
&\leq& C_3(1+\sup|J|+\sup|\nabla J|\sup|\nabla f||z|)|\nabla(f_t-f)|
      +C_4\sup|\nabla J|\sup|\nabla f||z|\sup|\nabla f|\\
&\leq& C_5|\nabla(f_t-f)|+C_6 |z|
\end{eqnarray*} 
Similarly,
\[
\left|\dbar_J (f_t-f)\left(\frac{\pa}{\pa y}\right)\right|\leq
C_5|\nabla(f_t-f)|+C_6|z| 
\]
The case $\frac{r^2}{\ep_1}<|z|<r$ can be estimated similarly.
Therefore,
\[
\parallel \dbar_J u_t \parallel_{L^p}\leq
C_7(\parallel \nabla(v_t-u) \parallel_{L^p} + \sqrt{r}r^\frac{1}{p})
\leq C r^\frac{1}{p} 
\]
where $C$ depends on the $C^1$ norms of $u$ and of $J$, 
 $C^\infty$ norms of $\exp^1|_{B_{r_1}(N_{L/X})}$ and its inverse.
$\Box$

\bigskip
We next consider boundary nodes. For $t\in [0,\infty)$, define
\[
  B^+_{\ep_1,t}=\{(x,y)\in \C^2\mid xy=t,|x|<\ep_1, |y|<\ep_1, 
  \mathrm{Im}x\geq 0,\mathrm{Im}y\leq 0 \}.
\]
Then for $t=r^2>0$,
\[ B^+_{\ep_1,t}=\left\{\left(x,\frac{r^2}{x}\right)\left|\; x\in\C,
  \frac{r^2}{\ep_1}<|x|<\ep_1, \mathrm{Im}x\geq 0 \right.\right\}
\]
Let $B^+_{\ep_1}=\{x\in\C \mid |x|<\ep_1,\mathrm{Im}\,x\geq0\}$.
Let $I_{\ep_1}=B^+_{\ep_1}\cap \R$, and
$I_{\ep_1,t}=B^+_{\ep_1,t}\cap \R\times \R$.
Two nonconstant $J$-holomorphic maps 
$f,g: (B^+_{\ep_1},I_{\ep_1})\ra (X,L)$ 
such that $f(0)=g(0)=p$ determine a stable map 
$u:(B^+_{\ep_1,0},I_{\ep_1,0})\ra (X,L)$ defined by 
$u(x,0)=f(x)$, $u(0,y)=g(y)$. One can construct
$u_t:(B^+_{\ep_1,t},I_{\ep_1,t})\ra(X,L)$ and
$v_t:(B^+_{\ep_1,0},I_{\ep_1,0})\ra(X,L)$ as before.

\bigskip
Applying the above construction to each node, we obtain
\begin{eqnarray*}
&&u_{\eta,\eta'}:(\Si_{(0,\eta,\eta')},\pa\Si_{(0,\eta,\eta')})\ra (X,L)\\
&&v_{\eta,\eta'}: (\Si,\bS)\ra (X,L)
\end{eqnarray*}

Lemma~\ref{bar} implies

\begin{lm}\label{almosthol}
\[
 \parallel \dbar_J u_{\eta,\eta'}\parallel_{L^{p}} \leq 
 C(|\eta|+|\eta'|)^\frac{1}{2p}
\]
where $C$ depends on the $C^1$ norms of $u$ and $J$, $C^\infty$ norms 
of $\exp|_{B_{r_1}}(N_{L/X})$ and its inverse.
\end{lm}

The linearization of $\dbar_{J,\Si}$ at the stable $W^{1,p}$ map
$v_{\eta,\eta'}$ is
\[
 D_{v_{\eta,\eta'}}:
W^{1,p}(\Si,\bS,v_{\eta,\eta'}^*TX,(v_{\eta,\eta'}|_{\bS})^*TL)
\ra L^p(\Si,\Lambda^{0,1}\Si\otimes v_{\eta,\eta'}^*TX)
\]
\begin{lm}\label{tend}
\[
\lim_{(\eta,\eta')\ra 0}\parallel D_{v_{\eta,\eta'}} \parallel
=\parallel D_u \parallel 
\] 
\end{lm}
\paragraph{Proof.}

We have a bundle isomorphism
\[
P_0:(u^*TX, (u|_{\bS})^* TL) \ra (v_{\eta,\eta'}^*TX, 
(v_{\eta,\eta'}|_{\bS})^* TL)
\]
given by parallel transport along the unique length minimizing
geodesic from $u(z)$ to $v_{\eta,\eta'}(z)$. This also gives 
\[
P_1:\Lambda^1\Si\otimes u^*TX \ra \Lambda^1\Si\otimes v_{\eta,\eta'}^*TX
\]
and
\[
P'_1=\pi\circ P_1\circ i:
\Lambda^{0,1}\Si\otimes u^*TX\ra\Lambda^{0,1}\Si\otimes v_{\eta,\eta'}^*TX,
\]
where $i:\Lambda^{0,1}\Si\otimes u^*TX\ra\Lambda^1\Si\otimes u^*TX$ is the 
inclusion, and $\pi:\Lambda^1\Si\otimes v_{\eta,\eta'}^*TX\ra
 \Lambda^{0,1}\Si\otimes v_{\eta,\eta'}^*TX$ is the projection. $P_0,P'_1$
induces
\begin{eqnarray*}
&&\tilde{P}_{\eta,\eta',0}^{-1}:
W^{1,p}(\Si,\bS,v_{\eta,\eta'}^*TX,(v_{\eta,\eta'}|_{\bS})^*TL)
 \ra \Tl\\
&&\tilde{P}_{\eta,\eta',1}:\Tless
\ra L^p(\Si,\Lambda^{0,1}\Si\otimes v_{\eta,\eta'}^*TX)
\end{eqnarray*}
Define 
\begin{eqnarray*}
D'_{\eta,\eta'}=\tilde{P}_{\eta,\eta',1}\circ D_u\circ
\tilde{P}_{\eta,\eta',0}^{-1}:&&
W^{1,p}(\Si,\bS,v_{\eta,\eta'}^*TX,(v_{\eta,\eta'}|_{\bS})^*TL)\\
&\ra& L^p(\Si,\Lambda^{0,1}\Si\otimes v_{\eta,\eta'}^*TX),
\end{eqnarray*}
then
\[
  \lim_{(\eta,\eta')\ra(0,0)}
 \parallel D'_{\eta,\eta'} \parallel = \parallel D_u \parallel.
\]
From Lemma~\ref{linear}, we see that
\begin{eqnarray*}
\parallel(D_{v_{\eta,\eta'}}-D'_{\eta,\eta'})w\parallel 
&\leq& C_2\parallel u - v_{\eta,\eta'}\parallel_{C^0} 
       \parallel \nabla w \parallel_{L^p}+
       \parallel w\parallel_{C^0}
       \parallel du-dv_{\eta,\eta'}\parallel_{L^p} \\
&\leq& C_3 \parallel u-v_{\eta,\eta'}\parallel _{W^{1,p}}
 \parallel w \parallel_{W^{1,p}}\\
\parallel D_{v_{\eta,\eta'}}-D'_{\eta,\eta'}\parallel 
&\leq& C_3\parallel u-v_{\eta,\eta'}\parallel_{W^{1,p}}
\end{eqnarray*}
which tends to $0$ as $(\eta,\eta')\ra (0,0)$. $\Box$

\subsubsection{Gluing: construction of exact solutions}\label{gluing}

 The goal of this section is to construct a local parametrization of
 solutions to  $\pi\circ \dbar_J v=0$ near the approximate $J$-holomorphic
 map $u_{\eta,\eta'}$ constructed in Section \ref{pregluing}.
 The main result in this Section is Proposition \ref{surjective}.

 Let $B_{\delta_2}\times D(d_1,\dots, d_{l_0})\times D'(d'_1,\ldots,d'_{l_1})$
 be the neighborhood of the origin in $H_{\rho,\mathrm{domain}}$ as 
 in Section~\ref{bnodes}. We write $D_d$ for $D(d_1,\dots, d_{l_0})$, and
 $D'_{d'}$ for $D'(d'_1,\ldots,d'_{l_1})$, as in Section~\ref{actauto}.
 We have seen that there is a family of prestable bordered Riemann surfaces 
\[
 (\pi:\cC\ra B_{\delta_2}\times D_d\times D'_{d'};\bs;\bt^1,\ldots,\bt^h)
\]   
of type $(g,h)$ with $(n,m)$ marked points, together with
a family of admissible metrics, such that $\pi^{-1}(0)=\Si$.
There is a map $\cC\ra \Si$ whose restriction to each fiber
of $\pi$ is a smooth strong deformation
$\kappa_{(\xi,\eta,\eta')}:\lambda_{(\xi,\eta,\eta')}\ra \lambda$.

Let $B$ be the image of the map $i:B_{\delta_2}\times D_d\times D'_{d'}
\ra \widetilde{M}_{(g,h),(n,\vm)}$ given by $(\xi,\eta,\eta')\mapsto
\lambda_{(\xi,\eta,\eta')}$. Then $B$ is a neighborhood of $\lambda$ in
the $C^\infty$ topology. Using the family of admissible metrics, we
define $W_B=\cup_{\lambda'\in B}W_{\lambda'}^{1,p}$. 
Let $M_B=\MXL\cap W_B$. Then $\rho\in M_B$.

We have a Cartesian diagram
\[
\begin{CD}
&&W_{\tilde{B}} @>>> W_B\\
&&@V\tilde{\pi}VV  @VV\pi V\\
B_{\delta_2}\times D_d\times D'_{d'}&=&\tilde{B} @>i>> B 
\end{CD}
\]
Let $S:\tilde{B}\ra W_{\tilde{B}}$ be given by
$(\xi,\eta,\eta')\mapsto u_{\eta,\eta'}$.

We first extend $E_\rho\subset 
 L^p(\Si,\Lambda^{0,1}\Si\otimes u^* TX)$ to a trivial bundle
over $S(\{0\}\times D_d\times D'_{d'})$.
Recall that elements in $E_\rho$ are supported on $K_{\ep_1}(\Si)$.
Since $v_{\eta,\eta'}=u$ on $K_{\ep_1}(\Si)$, $E_\rho$ can
be viewed as a subspace  $E_{\eta,\eta'}$ of 
$L^p(\Si,\Lambda^{0,1}\Si\otimes v_{\eta,\eta'}^*TX)$.
Let $F_{\eta,\eta'}$ be the $L^2$-orthogonal complement
of $E_{\eta,\eta'}$ in 
$L^p(\Si,\Lambda^{0,1}\Si\otimes v_{\eta,\eta'}^*TX)$, and
let 
\[
  \pi_{\eta,\eta'}:
L^p(\Si,\Lambda^{0,1}\Si\otimes v_{\eta,\eta'}^*TX)
\ra F_{\eta,\eta'}
\]
be the $L^2$-orthogonal projection. Then for sufficiently small
$\eta,\eta'$, 
\[
  \pi_{\eta,\eta'}\circ D_{v_{\eta,\eta'}}:
  W^{1,p}(\Si,\bS,v_{\eta,\eta'}^*TX,(v_{\eta,\eta'}|_{\bS})^*TL)
  \ra F_{\eta,\eta'}
\] 
is a surjection. Let $H_{(\eta,\eta')}=
\mathrm{Ker}(\pi_{\eta,\eta'}\circ D_{v_{\eta,\eta'}})$, and
let $H_{(\eta,\eta')}^\bot$ be its $L^2$-orthogonal complement
in $W^{1,p}(\Si,\bS,v_{\eta,\eta'}^*TX,(v_{\eta,\eta'}|_{\bS})^*TL)$.
Then we have an isomorphism
\[
  \pi_{\eta,\eta'}\circ D_{v_{\eta,\eta'}}:
 H_{\eta,\eta'}^\bot\ra F_{\eta,\eta'}
\]
  whose inverse is 
\[
  \hat{Q}_{\eta,\eta'}: F_{\eta,\eta'} \ra H_{\eta,\eta'}^\bot 
\]
We have
\[
  Q_{\eta,\eta'}=i\circ \hat{Q}_{\eta,\eta'}: F_{\eta,\eta'}
  \ra W^{1,p}(\Si,\bS,v_{\eta,\eta'}^*TX,(v_{\eta,\eta'}|_{\bS})^*TL)
\]
where 
\[
i: H_{(\eta,\eta')}^\bot\ra 
W^{1,p}(\Si,\bS,v_{\eta,\eta'}^*TX,(v_{\eta,\eta'}|_{\bS})^*TL)
\]
is the inclusion.
$Q_{\eta,\eta'}$ is a right inverse of 
$\pi_{\eta,\eta'}\circ D_{v_{\eta,\eta'}}$.

By Lemma~\ref{tend}, we may choose $d,d'$ sufficiently small such that
$\hat{Q}_{\eta,\eta'}$ exists 
and $\parallel Q_{\eta,\eta'}\parallel \leq M$ for all
$(\eta,\eta')\in D_d\times D'_{d'}$, where $M$ is a constant. 

We now extend $E\ra S(\{0\}\times D_d\times D'_{d'})$ to a neighborhood 
$U_{\tilde{B}}$ of $S(\tilde{B})$ in $W_{\tilde{B}}$. 
Let $\rho'=(\lambda_{(\xi,\eta,\eta')}, f)$, where
$f:(\Si_{(\xi,\eta,\eta')},\bS_{(\xi,\eta,\eta')})\ra (X,L)$  
is a stable $W^{1,p}$ map such that
\[
\sup_{\Si_{(\xi,\eta,\eta')}}d_1(u_{\eta,\eta'}(z),f(z))
\]
is less than the injectivity radius of $g_1$, where $d_1$ is the geodesic
distance of $g_1$.

We have a bundle isomorphism
\[ 
  P_0:(u_{\eta,\eta'}^*TX,(u_{\eta,\eta}|_{\bS_{(\xi,\eta,\eta')}})^*TL) 
 \ra (f^* TX,(f|_{\bS_{(\xi,\eta,\eta')}})^* TL) 
\]
given by the parallel transport along the unique length minimizing
geodesic from $u_{\eta,\eta'}(z)$ to $f(z)$, which
gives  
\begin{eqnarray*}
 P_1:&&\Lambda^1\Si_{(\xi,\eta,\eta')}\otimes u_{\eta,\eta'}^*TX 
  \ra\Lambda^1\Si_{(\xi,\eta,\eta')}\otimes f^*TX\\
 P'_1=\pi\circ P_1\circ i:&&
\Lambda^{0,1}\Si_{(\xi,\eta,\eta')}\otimes u_{\eta,\eta'}^*TX 
  \ra\Lambda^{0,1}\Si_{(\xi,\eta,\eta')}\otimes f^*TX
\end{eqnarray*}
where $i:\Lambda^{0,1}\Si_{(\xi,\eta,\eta')}\otimes u_{\eta,\eta'}^*TX
\ra\Lambda^1\Si_{(\xi,\eta,\eta')}\otimes u_{\eta,\eta'}^*TX$
is the inclusion, and $\Lambda^1\Si_{(\xi,\eta,\eta')}\otimes f^*TX\ra
\Lambda^{0,1}\Si_{(\xi,\eta,\eta')}\otimes f^*TX$ is the projection.
We have  
\[
\tilde{P}:L^p(\Si_{(\xi,\eta,\eta')},
\Lambda^{0,1}\Si_{(\xi,\eta,\eta')}\otimes u_{\eta,\eta'}^*TX)
\ra L^p(\Si_{(\xi,\eta,\eta')},\Lambda^{0,1}\Si_{(\xi,\eta,\eta')}\otimes f^*TX).
\]
Let $E_{\rho'}=\tilde{P}E_\rho\cong E_\rho$. Then we have a trivial
bundle $E\ra U_{\tilde{B}}$ together with a trivialization
$\Phi: E\cong U_{\tilde{B}}\times E_\rho$.

Let $F_{\rho'}$ be the $L^2$-orthogonal complement of $E_{\rho'}$ in 
$L^p(\Si_{(\xi,\eta,\eta')},\Lambda^{0,1}\Si_{(\xi,\eta,\eta')}\otimes
f^*TX)$, and let 
\[
\pi_{\rho'}:
L^p(\Si_{(\xi,\eta,\eta')},\Lambda^{0,1}\Si_{(\xi,\eta,\eta')}\otimes f^*TX)
\ra F_{\rho'}
\]
be the $L^2$-orthogonal projection.

We will use $Q_{\eta,\eta'}$, the right inverse of 
$\pi_{\eta,\eta'}\circ D_{v_{\eta,\eta'}}$, to construct an approximate right 
inverse $Q'_{(\xi,\eta,\eta')}$  of 
$\pi_{(\xi,\eta,\eta')}\circ D_{(\xi,\eta,\eta'),u_{\eta,\eta'}}$, where
$\pi_{(\xi,\eta,\eta')}=
\pi_{(\lambda_{(\xi,\eta,\eta')}, u_{\eta,\eta'})}$, and 
$D_{(\xi,\eta,\eta'),u_{\eta,\eta'}}$ is the linearization of 
$\dbar_{J,\Si_{(\xi,\eta,\eta')}}$ at $u_{\eta,\eta'}$.
We will use the cutoff function constructed in
\cite[A.1]{MS}. The construction in this section 
works for $W^{1,p}$ but not for general $W^{k,p}$.

\begin{lm} \label{cutoff}
For any $r\in (0,1)$, there is a smooth cutoff function
$\chi_r:\C\ra [0,1]$ such that 
\[
\chi_r(z)=\left\{ \begin{array}{ll}
   1&\textup{if } |z|\leq r\sqrt{r}\\
   0&\textup{if } |z|\geq r
  \end{array}\right.
\]
\[
\int_{|z|\leq r}\chi_r \leq \frac{4\pi}{|\log r|}.
\]
\end{lm} 

\paragraph{Proof.} We will follow the proof of \cite[Lemma A.1.1]{MS}.
We first define a cutoff function of class $W^{1,2}$ by
\[ \beta_r(z)=\left\{\begin{array}{ll}
         1& \textup{ for }|z|\leq r\sqrt{r}\\
   2\left(\frac{\log|z|}{\log r}-1\right)& 
   \textup{ for } r\sqrt{r}\leq|z|\leq r\\
   0 & \textup{ for } |z|\geq r
   \end{array} \right. 
\]
Then we have
\[
  |\nabla\beta_r(z)|=\frac{2}{|z||\log r|}
\]
for $r\sqrt{r}\leq|z|\leq r$, so 
\[
  \int_{r\sqrt{r}\leq|z|\leq r}|\nabla \beta_r(z)|^2=\frac{4\pi}{|\log r|}. 
\]
To obtain a smooth function $\chi_r$, take the convolution with
$\phi_N(z)=N^2\phi(Nz)$ where $N$ is large and $\phi:\C\ra \R$ is a 
smooth function with support in the unit ball and mean value $1$.
$\Box$

\bigskip

Let $p>2$ be fixed as before.
\begin{lm}\label{cutsection}
For every $r\in(0,1)$, there exists a smooth cutoff function
$\chi_r:\C\ra[0,1]$ as in Lemma~\ref{cutoff} such that
\[
  \parallel (\nabla\chi_r) w \parallel_{L^p}\leq 
 C\parallel w \parallel_{W^{1,p}}|\log r|^{\frac{1}{p}-1} 
\]
for any $w\in W^{1,p}(\C)$ with $w(0)=0$.
\end{lm}

\paragraph{Proof.} It follows from the proof of \cite[Lemma A.1.2]{MS}. $\Box$

\bigskip

We now look at the local model of an interior node.
Let $u:B_{\ep_1,0}\ra X$ be a stable map, and construct smooth maps
$u_t:B_{\ep_1,t}\ra X$, $v_t:B_{\ep_1,0}\ra X$ as before.
We now define linear maps
\begin{eqnarray*}
&& e_t:L^p(B_{\ep_1,t},\Lambda^{0,1}B_{\ep_1,t}\otimes u_t ^* TX)
\ra L^p(B_{\ep_1,0},\Lambda^{0,1}B_{\ep_1,0}\otimes v_t ^* TX)\\
&& g_t: W^{1,p}(B_{\ep_1,0},v_t^* TX)\ra W^{1,p}(B_{\ep_1,t},u_t^* TX)
\end{eqnarray*}

Given $s\in L^p(B_{\ep_1,t}, \Lambda^{0,1} B_{\ep_1,t}\otimes u_t^*TX)$, 
we define
\[
  e_t(s)\in L^p(B_{\ep_1,0}, \Lambda^{0,1} B_{\ep_1,0}\otimes v_t^*TX)
\]
by
\[
e_t(s)(x,0)=\left\{\begin{array}{ll} 
s(x,\frac{t}{x}) & \textup{ if }r\leq|x|\leq \ep_1\\
0 & \textup{ if } |x|<r \end{array}\right.
\]
\[
e_t(s)(0,y)=\left\{\begin{array}{ll} 
s(\frac{t}{y},y) & \textup{ if }r\leq|y|\leq \ep_1\\
0 & \textup{ if } |y|<r \end{array}\right.
\]
The above definition is valid since
\[
 u_t\left(z,\frac{t}{z}\right)=\left\{\begin{array}{ll}
 v_t(0,\frac{t}{z})& \textup{ if } \frac{r^2}{\ep_1}\leq |z|\leq r\\
 v_t(z,0) & \textup{ if }r\leq|x|\leq \ep_1
 \end{array}\right.
\]

Given $w\in W^{1,p}(B_{\ep_1,0},v_t^* TX)$, we define
\[
g_t(w)\in W^{1,p}(B_{\ep_1,t},u_t^* TX)
\]
by
\[
g_t(w)\left(z,\frac{t}{z}\right)=\left\{\begin{array}{ll}
w(z,0) & \textup{ if } \sqrt{r}\leq |z|\leq \ep_1 \\
w(z,0)+ (1-\chi_r(\frac{t}{z}))(w(0,\frac{t}{z})-w(0,0)) &
\textup{ if } r\leq|z|\leq\sqrt{r}\\
w(0,\frac{t}{z})+(1-\chi_r(z))(w(z,0)-w(0,0)) &
\textup{ if } r\sqrt{r}\leq |z|\leq r\\
w(0,\frac{t}{z}) &\textup{ if }
\frac{r^2}{\ep_1}\leq |z|\leq r\sqrt{r}
\end{array}\right.
\]

\begin{lm}\label{cutin}
If $s\in L^p(B_{\ep_1,t},\Lambda^{0,1} B_{\ep_1,t}\otimes u_t^*TX),\ 
w\in W^{1,p}(B_{\ep_1,0},v_t^* TX)$ 
satisfy $D_{v_t} w= e_t(s)$, then
\[
\parallel D_{u_t}g_t(w)-s\parallel_{L^p}\leq C
\parallel w \parallel_{W^{1,p}}(|\log r|^{\frac{1}{p}-1}+
\parallel u \parallel_{W^{1,p}} r^\frac{1}{p})
\]
\end{lm}

\paragraph{Proof.}
We have
\[
D_{u_t}\circ g_t(w)\left(z,\frac{t}{z}\right)
=\left\{ \begin{array}{ll}
D_{v_t}w(z,0)=e_t(s)(z,0)=s(z,\frac{t}{z})
&\textup{ for } \sqrt{r}\leq |z| \leq \ep_1,\\
D_{v_t}w(0,\frac{t}{z})=e_t(s)(0,\frac{t}{z})
=s(z,\frac{t}{z})
&\textup{ for } \frac{r^2}{\ep_1}\leq|z|\leq r\sqrt{r}.
\end{array}\right. 
\]
We now consider $r\leq|z|\leq\sqrt{r}$. We have
\begin{eqnarray*}
&& D_{u_t}\circ g_t(w)\left(z,\frac{t}{z}\right)\\
&=&D_{v_t}w(z,0)+ D(1-\chi_r\left(\frac{t}{z}\right))
   (w(0,\frac{t}{z})-w(0,0))\\
&&+(1-\chi_r\left(\frac{t}{z}\right))(D_{v_t}w (0,\frac{t}{z})
  -D_{v_t}w(0,0))
\end{eqnarray*}
where $D$ is a derivation, so 
\[
|D(1-\chi_r\left(\frac{t}{z}\right))|\leq C_1
\left|\frac{t}{z^2}\nabla\chi_r\left(\frac{t}{z}\right)\right|.
\]
We also have
\begin{eqnarray*}
D_{v_t}w(z,0)&=&e_t(s)(z,0)=s(z,\frac{t}{z}),\\
D_{v_t}w(0,\frac{t}{z})&=&e_t(s)(0,\frac{t}{z})=0,
\end{eqnarray*}
and
\begin{eqnarray*}
|D_{v_t}w(0,0)|&\leq& C_2\parallel v_t \parallel_{W^{1,p}}|w(0,0)|\\
               &\leq& C_2\parallel v_t \parallel_{W^{1,p}}
                       \parallel w \parallel_{C^0}\\
               &\leq& C_3 \parallel v_t \parallel_{W^{1,p}}
                          \parallel w \parallel_{W^{1,p}} 
\end{eqnarray*}
Let $h_0$ be the metric on $A(r,\sqrt{r})$ given by 
$A(r,\sqrt{r})\subset\C$, and $h_1$ be the metric
on $A(r,\sqrt{r})$ given by the embedding 
$z\in A(r,\sqrt{r})\mapsto (z,\frac{t}{z})\in \C^2$. Then
\[
h_1=\left(1+ \left(\frac{r}{|z|}\right)^4\right)h_0,
\]
so $h_0\leq h_1\leq 2h_0$. To estimate the $L^{1,p}$ norm
defined by any metric which is an interpolation of $h_0$ and
$h_1$, it suffices to calculate in $h_0$.
\begin{eqnarray*}
&&\int_{r\leq|z|\leq\sqrt{r}}
|D_{u_t}\circ g_t(w)(z,\frac{t}{z})-s(z,\frac{t}{z})|^p
\frac{i}{2}dz\wedge d\bar{z}\\
&\leq& C_1^p \int_{\sqrt{r}\leq|z|\leq r} 
     \left|\nabla\chi_r\left(\frac{t}{z}\right)\right|^p
     \left|w(0,\frac{t}{z})-w(0,0)\right|^p
     \left(\frac{r}{|z|}\right)^{2p}\frac{i}{2}dz\wedge d\bar{z}\\
   &&+C_4 \parallel v_t \parallel_{W^{1,p}}^p
     \parallel w \parallel_{W^{1,p}}^p r
\end{eqnarray*}
where
\begin{eqnarray*}
&&\int_{r\leq|z|\leq \sqrt{r}}
     \left|\nabla\chi_r\left(\frac{t}{z}\right)\right|^p
     |w(0,\frac{t}{z})-w(0,0)|^p\left(\frac{r}{|z|}\right)^{2p}
          \frac{i}{2}dz\wedge d\bar{z}\\
&=&\int_{r\sqrt{r}\leq |y|\leq r}    
   |\nabla\chi_r(y)|^p|w(0,y)-w(0,0)|^p
\left(\frac{|y|}{r}\right)^{2p-4}\frac{i}{2}dy\wedge d\bar{y}\\
&=&\int_{r\sqrt{r}\leq |y|\leq r}    
   |\nabla\chi_r(y)|^p|w(0,y)-w(0,0)|^p
  \frac{i}{2}dy\wedge d\bar{y}\\
&\leq& C_5 \parallel w\parallel_{W^{1,p}}^p |\log r|^{1-p}
\end{eqnarray*}

We finally consider the case $r\sqrt{r}\leq |z|\leq r$. Let
$y=\frac{t}{z}$, then $r\leq |y| \leq\sqrt{r}$, and
\begin{eqnarray*}
g_t(w)\left(z,\frac{t}{z}\right)&=& g_t(w)\left(\frac{t}{y},y\right)\\
&=&w(0,y)+(1-\chi_r\left(\frac{t}{y}\right))(w(\frac{t}{y},0)-w(0,0)),
\end{eqnarray*}
which is the same as the case $r\leq |z|\leq \sqrt{r}$.
So we conclude that
\[
\parallel D_{u_t}g_t(w)-s\parallel_{L^p}\leq C
\parallel w \parallel_{W^{1,p}}(|\log r|^{\frac{1}{p}-1}+
\parallel u \parallel_{W^{1,p}} r^\frac{1}{p})
\]
since $v_t$ converges to $u$ in $W^{1,p}$ norm. $\Box$

\bigskip

We next look at the local model of a boundary node.
Let $u:(B^+_{\ep_1,0}, I_{\ep_1,0})\ra (X,L)$ be a stable map, and construct smooth maps
$u_t:(B^+_{\ep_1,t},I_{\ep_1,t})\ra (X,L)$, 
$v_t:(B^+_{\ep_1,0},I_{\ep_1,0})\ra (X,L)$ as before.
We define linear maps
\begin{eqnarray*}
&& e_t:L^p(B^+_{\ep_1,t},\Lambda^{0,1}B_{\ep_1,t}\otimes u_t ^* TX)
\ra L^p(B^+_{\ep_1,0},\Lambda^{0,1}B_{\ep_1,0}\otimes v_t ^* TX)\\
&& g_t: W^{1,p}(B^+_{\ep_1,0},I_{\ep_1,0},v_t^* TX,(v_t|_{I_{\ep_1,0}})^*TL)
   \ra  W^{1,p}(B^+_{\ep_1,t},I_{\ep_1,t},u_t^* TX,(u_t|_{I_{\ep_1,t}})^*TL)
\end{eqnarray*}
in exactly the same way as for $B_{\ep_1,t}$, $B_{\ep_1,0}$. Then we have
\begin{lm}\label{cutbd}
If $s\in L^p(B^+_{\ep_1,t},\Lambda^{0,1} B_{\ep_1,t}\otimes u_t^*TX),\ 
w\in W^{1,p}(B^+_{\ep_1,0},I_{\ep_1,0},v_t^* TX, (v_t|_{I_{\ep_1,0}})^* TL)$ 
satisfy $D_{v_t} w= e_t(s)$, then
\[
\parallel D_{u_t}g_t(w)-s\parallel_{L^p}\leq C
\parallel w \parallel_{W^{1,p}}(|\log r|^{\frac{1}{p}-1}+
\parallel u \parallel_{W^{1,p}} r^\frac{1}{p}).
\]
\end{lm}

\bigskip
We now apply  above construction to each node to obtain linear maps
\begin{eqnarray*}
&e_{\eta,\eta'}:&
 L^p(\Si_{(0,\eta,\eta')},\Lambda^{0,1}\Si_{(0,\eta,\eta')}\otimes 
    u_{\eta,\eta'}^* TX)
\ra L^p(\Si,\Lambda^{0,1}\Si\otimes v_{\eta,\eta'}^* TX)\\
&g_{\eta,\eta'}:& 
W^{1,p}(\Si,\bS,v_{\eta,\eta'}^* TX,(v_{\eta,\eta'}|_{\bS})^*TL)\\
 &&  \ra  W^{1,p}(\Si_{(0,\eta,\eta')},\bS_{(0,\eta,\eta')},
   u_{\eta,\eta'}^* TX,(u_{\eta,\eta'}|_{\bS_{(0,\eta,\eta')}})^*TL)
\end{eqnarray*}

Let 
\begin{eqnarray*}
&& Q'_{(0,\eta,\eta')}=g_{\eta,\eta'}\circ Q_{\eta,\eta'}
\circ \pi_{\eta,\eta'}\circ\left(e_{\eta,\eta'}|_{F_{\eta,\eta'}}\right):
F_{\eta,\eta'}\\
&\ra& W^{1,p}(\Si_{(0,\eta,\eta')},\bS_{(0,\eta,\eta')},
u_{\eta,\eta'}^* TX,(u_{\eta,\eta'}|_{\bS_{(0,\eta,\eta')}})^*TL).
\end{eqnarray*}
The operator norm of $Q'_{(0,\eta,\eta')}$ has a uniform bound independent
of $\eta,\eta'$ since that of $Q_{\eta,\eta'}$ has a uniform bound
independent of $\eta,\eta'$. We now show that
$Q'_{(0,\eta,\eta')}$ is an approximate right inverse of
$\pi_{(0,\eta,\eta')}\circ D_{(0,\eta,\eta'),u_{\eta,\eta'}}$.  

\begin{pro} \label{closer}
\[
 \parallel\left(\pi_{(0,\eta,\eta')}\circ D_{(0,\eta,\eta'),u_{\eta,\eta'}}
  \circ Q'_{(0,\eta,\eta')}\right) s-s
  \parallel_{L^p}
\leq C(|\log(|\eta|+|\eta'|)|^{\frac{1}{p}-1})
 \parallel s \parallel_{L^p}
\]
where $C$ depends on $\parallel u \parallel_{W^{1,p}}$.
\end{pro}

\paragraph{Proof.} 
Let $\rho(\eta,\eta')=(\lambda_{(0,\eta,\eta')},u_{\eta,\eta'})$.
Given $s\in F_{\rho(\eta,\eta')}$, let 
\begin{eqnarray*}
s_1&=&\pi_{\eta,\eta'}\circ e_{\eta,\eta'}(s)\in F_{\eta,\eta'},\\
t_1&=&e_{\eta,\eta'}(s)-s_1\in E_{\eta,\eta'},\\ 
w&=&Q_{\eta,\eta'}(s_1)\in W^{1,p}(\Si,\bS,v_{\eta,\eta'}^*TX,
(v_{\eta,\eta'}|_{\bS})^* TL),
\end{eqnarray*}
then
\[
\pi_{\eta,\eta'}\circ D_{v_{\eta,\eta'}} w=s_1=
e_{\eta,\eta'}(s)-t_1.
\]
so $D_{v_{\eta,\eta'}}w= e_{\eta,\eta'}(s)+t_2$ for
some $t_2\in E_{\eta,\eta'}$. There is a unique 
$t\in E_{\rho(\eta,\eta')}$ such that 
$e_{\eta,\eta'}(t)=t_2$. We have
\[
D_{v_{\eta,\eta'}}w=e_{\eta,\eta'}(s+t),
\]
so by Lemma~\ref{cutin}, \ref{cutbd},
\begin{eqnarray*}
&&\parallel D_{(0,\eta,\eta'),u_{\eta,\eta'}}\circ g_{\eta,\eta'}(w)
 -(s+t)\parallel_{L^p\left(\kappa_{\eta,\eta'}^{-1}(N_{\ep_1}(\Si))\right)}\\
&\leq& C_1 \parallel w\parallel_{W^{1,p}}
\left(|\log(|\eta|+|\eta'|)|^{\frac{1}{p}-1}+ (|\eta|+|\eta'|)^\frac{1}{2p}\right)
\end{eqnarray*}
where $C_1$ depends on $\parallel u \parallel_{W^{1,p}}$.

We have
\[
\pi_{(0,\eta,\eta')}\circ D_{(0,\eta,\eta'),u_{\eta,\eta'}}\circ Q'_{\eta,\eta'}s-s
=\pi_{(0,\eta,\eta')}\left(
   D_{(0,\eta,\eta'),u_{\eta,\eta'}}\circ g_{\eta,\eta'}(w)-(s+t)\right)
\]
\begin{eqnarray*}
\parallel w \parallel_{W^{1,p}}&=&\parallel Q_{\eta,\eta'}\circ
\pi_{\eta,\eta'}\circ e_{\eta,\eta'}(s)\parallel_{W^{1,p}}\\
&\leq& C_2\parallel\pi_{\eta,\eta'}\circ e_{\eta,\eta'}(s)
\parallel_{L^p}\\
&\leq& C_3\parallel s\parallel_{L^p}
\end{eqnarray*}
So
\begin{eqnarray*}
&&\parallel \pi_{(0,\eta,\eta')}\circ D_{(0,\eta,\eta'),u_{\eta,\eta'}}\circ
Q'_{(0,\eta,\eta')}s-s \parallel
_{L^p\left(\kappa_{\eta,\eta'}^{-1}(N_{\ep_1}(\Si))\right)}\\
&\leq& C_4 \left(|\log(|\eta|+|\eta'|)|^{\frac{1}{p}-1}+
(|\eta|+|\eta'|)^\frac{1}{2p}\right)\parallel s \parallel_{L^p}\\
&\leq& C_5|\log(|\eta|+|\eta'|)|^{\frac{1}{p}-1} \parallel s \parallel_{L^p}.
\;\; \Box
\end{eqnarray*}

\bigskip

Let
\[
p_{(\xi,\eta,\eta')}:L^p(\Si_{(\xi,\eta,\eta')},\Lambda^{0,1}\Si_{(\xi,\eta,\eta')}\otimes 
u_{\eta,\eta'}^*TX )
\ra L^p(\Si_{(\xi,\eta,\eta')},\Lambda^{0,1}\Si_{(0,\eta,\eta')}\otimes 
  u_{\eta,\eta'}^*TX) 
\]
be the map determined by the bundle isomorphism 
$P\circ i:\Lambda^{0,1}\Si_{(\xi,\eta,\eta')}\ra\Lambda^{0,1}\Si_{(0,\eta,\eta')}$,
where $i: \Lambda^{0,1}\Si_{(\xi,\eta,\eta')}\ra
 \Lambda^1\Si_{(\xi,\eta,\eta')}\cong\Lambda^1\Si_{(0,\eta,\eta')}$
is the inclusion, and
$P:\Lambda^1\Si_{(0,\eta,\eta')}\ra\Lambda^{0,1}\Si_{(0,\eta,\eta')}$
is the projection. Let
\begin{eqnarray*}
&& Q'_{(\xi,\eta,\eta')}=g_{\eta,\eta'}\circ Q_{\eta,\eta'}
\circ \pi_{\eta,\eta'}\circ e_{\eta,\eta'}\circ(p_{(\xi,\eta,\eta')}|_{F_{(\xi,\eta,\eta')}}):
F_{(\xi,\eta,\eta')}\\
&\ra& W^{1,p}(\Si_{(\xi,\eta,\eta')},\bS_{(\xi,\eta,\eta')},
u_{\eta,\eta'}^* TX,(u_{\eta,\eta'}|_{\bS_{(\xi,\eta,\eta')}})^*TL)
\end{eqnarray*}
where $F_{(\xi,\eta,\eta')}=F_{(\lambda_{(\xi,\eta,\eta')},u_{\eta,\eta'})}$.
We have

\begin{pro} \label{close}
\[
 \parallel\left(\pi_{(\xi,\eta,\eta')}\circ D_{(\xi,\eta,\eta'),u_{\eta,\eta'}}
  \circ Q'_{(\xi,\eta,\eta')}\right) s-s
  \parallel_{L^p}
\leq C(|\xi|+|\log(|\eta|+|\eta'|)|^{\frac{1}{p}-1})
 \parallel s \parallel_{L^p}
\]
\end{pro}

\paragraph{Proof.}
We identify
\[
W^{1,p}(\Si_{(\xi,\eta,\eta')},\bS_{(\xi,\eta,\eta')},
u_{\eta,\eta'}^* TX,(u_{\eta,\eta'}|_{\bS_{(\xi,\eta,\eta')}})^*TL)
\]
with
\[
W^{1,p}(\Si_{(0,\eta,\eta')},\bS_{(0,\eta,\eta')},
u_{\eta,\eta'}^* TX,(u_{\eta,\eta'}|_{\bS_{(0,\eta,\eta')}})^*TL),
\]
and embed both
\[
L^p(\Si_{(\xi,\eta,\eta')},
\Lambda^{0,1}\Si_{(\xi,\eta,\eta')}\otimes u_{\eta,\eta'}^*TX)
\]
and 
\[
L^p(\Si_{(0,\eta,\eta')},
\Lambda^{0,1}\Si_{(0,\eta,\eta')}\otimes u_{\eta,\eta'}^*TX)
\]
into
\[
L^p(\Si_{(0,\eta,\eta')},
\Lambda^1\Si_{(0,\eta,\eta')}\otimes u_{\eta,\eta'}^*TX).
\]
We extend the domain of $Q'_{(0,\eta,\eta')}$ to
$L^p(\Si_{(0,\eta,\eta')},
\Lambda^{0,1}\Si_{(0,\eta,\eta')}\otimes u_{\eta,\eta'}^*TX)$,
and that of $Q'_{(\xi,\eta,\eta')}$ to
$L^p(\Si_{(0,\eta,\eta')},
\Lambda^{0,1}\Si_{(0,\eta,\eta')}\otimes u_{\eta,\eta'}^*TX)$.
In other words, we have
\[
Q'_{(0,\eta,\eta')}=g_{\eta,\eta'}\circ Q_{\eta,\eta'}
\circ \pi_{\eta,\eta'}\circ e_{\eta,\eta'},\;\;\;
Q'_{(\xi,\eta,\eta')}=Q'_{(0,\eta,\eta')}\circ p_{(\xi,\eta,\eta')}.
\]
With the above understanding, we have
\begin{eqnarray*}
&&\left(\pi_{(\xi,\eta,\eta')}\circ D_{(\xi,\eta,\eta'),u_{\eta,\eta'}}
  \circ Q'_{(\xi,\eta,\eta')}\right) s-s\\
&=&(\pi_{(\xi,\eta,\eta')}\circ D_{(\xi,\eta,\eta'),u_{\eta,\eta'}}
   -\pi_{(0,\eta,\eta')}  \circ D_{(0,\eta,\eta'),u_{\eta,\eta'}})
    \circ Q'_{(\xi,\eta,\eta')}s\\
  && +\left(\pi_{(0,\eta,\eta')}\circ D_{(0,\eta,\eta'),u_{\eta,\eta'}}
      \circ  Q'_{(0,\eta,\eta')}\right)p_{(\xi,\eta,\eta)}s
     -p_{(\xi,\eta,\eta)}s\\
  && +p_{(\xi,\eta,\eta)}s-s,
\end{eqnarray*}
where 
\[
\parallel \pi_{(\xi,\eta,\eta')}\circ D_{(\xi,\eta,\eta'),u_{\eta,\eta'}}
   -\pi_{(0,\eta,\eta')}  \circ D_{(0,\eta,\eta'),u_{\eta,\eta'}}\parallel
\leq C_1|\xi|,
\]
$\parallel p_{(\xi,\eta,\eta)}-Id \parallel \leq C_2 |\xi|$, and
$\parallel Q'_{(\xi,\eta,\eta')}\parallel \leq C_3|\xi|$
for all $(\xi,\eta,\eta')\in B_{\delta_2}\times D_d\times D'_{d'}$. 

\bigskip

We also have
\begin{eqnarray*}
&&\parallel \left(\pi_{(0,\eta,\eta')}\circ D_{(0,\eta,\eta'),u_{\eta,\eta'}}
      \circ  Q'_{(0,\eta,\eta')}\right)p_{(\xi,\eta,\eta)}s
     -p_{(\xi,\eta,\eta)}s\parallel\\
&\leq&C_4 \left(|\log(|\eta|+|\eta'|)|^{\frac{1}{p}-1}\right)
     \parallel p_{(\xi,\eta,\eta)}s\parallel_{L^p}\\
&\leq& C_4 \left(|\log(|\eta|+|\eta'|)|^{\frac{1}{p}-1}\right)
     \parallel s\parallel_{L^p}
\end{eqnarray*}
where $C_4$ is the constant $C$ in Proposition~\ref{close}.
Therefore,
\begin{eqnarray*}
&& \parallel
\left(\pi_{(\xi,\eta,\eta')}\circ D_{(\xi,\eta,\eta'),u_{\eta,\eta'}}
  \circ Q'_{(\xi,\eta,\eta')}\right) s- s\parallel_{L^p} \\
&\leq&(C_1 C_3+ C_4+ C_2)(|\xi|+|\log(|\eta|+|\eta'|)|^{\frac{1}{p}-1})
  \parallel s\parallel_{L^p}.\;\Box
\end{eqnarray*}

\begin{cor}\label{invertible}
There exist $\delta_2, d_1,\ldots,d_{l_0},d'_1,\ldots,d'_{l_1}>0$
such that for every $(\xi,\eta,\eta')\in B_{\delta_2}\times D_d\times
D'_{d'}$ there is a right inverse 
$Q_{(\xi,\eta,\eta')}$ of 
$\pi_{(\xi,\eta,\eta')}\circ D_{(\xi,\eta,\eta'),u_{\eta,\eta'}}$ such that
the operator norm $\parallel Q_{(\xi,\eta,\eta')} \parallel\leq C$
for some constant $C$.
\end{cor}

\paragraph{Proof.} By Proposition~\ref{closer}, there exist
 $\delta_2,d_1,\ldots,d_{l_0},d'_1,\ldots,d'_{l_1}>0$ sufficiently small
such that
\[
\parallel \left(\pi_{(\xi,\eta,\eta')}\circ D_{(\xi,\eta,\eta'),u_{\eta,\eta'}}\circ
 Q'_{(\xi,\eta,\eta')}\right)s-s \parallel_{L^p} \leq \half \parallel s \parallel_{L^p}. 
\]
for any $(\xi,\eta,\eta')\in B_{\delta_2}\times D_d\times D'_{d'}$

Let $A_{\xi,\eta,\eta'}=
\pi_{(\xi,\eta,\eta')}\circ D_{(\xi,\eta,\eta'),u_{\eta,\eta'}}\circ
Q'_{(\xi,\eta,\eta')}-I$, where $I$ is the identity map.
Then $\parallel A_{\xi,\eta,\eta'} \parallel\leq\half$,
so $I+A_{\xi,\eta,\eta'}$ is invertible and
$\parallel (I+A_{\xi,\eta,\eta'})^{-1} \parallel \leq 2$. Let
$Q_{(\xi,\eta,\eta')}=Q'_{(\xi,\eta,\eta')}\circ 
(I+A_{\xi,\eta,\eta'})^{-1}$,
then $Q_{(\xi,\eta,\eta')}$ has the desired properties. $\Box$

\bigskip

Let $\delta_2, d_1,\ldots,d_{l_0}, d'_1,\ldots, d'_{l_1}$ be chosen as in 
Corollary~\ref{invertible}. Then
$\pi_{(\xi,\eta,\eta')}\circ D_{(\xi,\eta,\eta'),u_{\eta,\eta'}}$
is surjective for 
$(\xi,\eta,\eta')\in B_{\delta_2}\times D_d\times D'_{d'}$.
We will construct a linear isomorphism
\[ 
i_{(\xi,\eta,\eta')}: \mathrm{Ker}(\pi\circ D_u)\longrightarrow
\mathrm{Ker}(\pi_{(\xi,\eta,\eta')}\circ D_{(\xi,\eta,\eta'),u_{\eta,\eta'}}).
\]

Given $w\in W^{1,p}(\Si,\bS,u^*TX,(u|_{\bS})^*TL)$, we will cut it off near
nodes to obtain 
\[
g'_{(0,\eta,\eta')}(w)\in W^{1,p}(\Si_{(0,\eta,\eta')},\bS_{(0,\eta,\eta')},
u_{\eta,\eta'}^*TX,(u_{\eta,\eta'}|_{\bS})^*TL).
\]

We first look at the local model of an interior node. Let
$u:B_{\ep_1,0}\ra X$ be a stable map. We have constructed
smooth maps $u_t:B_{\ep_1,t}\ra X$ for small $t\in\C$ such that
\[
u_t\left(z,\frac{t}{z}\right)=\left\{\begin{array}{ll}
 u\left(0,\frac{t}{z}\right))& \textup{ if }
\frac{r^2}{\ep_1}<|z|<\frac{r\sqrt{r}}{2},\\
 p &\textup{ if } r\sqrt{r}\leq |z| \leq \sqrt{r},\\
 u(z,0)& \textup{ if } 2\sqrt{r}<|z|<\ep_1.
\end{array}\right.
\]
where $p=u(0,0)$.

Let $r=\sqrt{|t|}$ as before, and set $s=(4r)^\frac{1}{3}$.
For $w\in W^{1,p}(B_{\ep_1,0},u^*TX,(u|_{\partial B_{\ep_1,0}})^*TL)$, define
$w_t\in W^{1,p}(B_{\ep_1,t},u^*TX,(u_t|_{\partial B_{\ep_1,t}})^*TL)$
by 
\[
w_t\left(z,\frac{t}{z}\right)=\left\{\begin{array}{ll}
  (1-\chi_s(\frac{t}{z}))w(0,\frac{t}{z})
  +\chi_s(\frac{t}{z}) P(z)w(0,0) & \textup{ if } 
  \frac{r^2}{\ep_1}<|z|\leq r\\
  (1-\chi_s(z))w(z,0)+\chi_s(z)P(z)w(0,0) & \textup{ if } r\leq|z|<\ep_1
\end{array}\right.
\]
where $\chi_s$ is the cutoff function in Lemma~\ref{cutoff}, and $P(z)$ is
the parallel transport along the unique length minimizing geodesic from
$p$ to $u_t(z,\frac{t}{z})$. We have
\[
w_t\left(z,\frac{t}{z}\right)=\left\{\begin{array}{ll}
  w(0,\frac{t}{z}) & \textup{ if } 
  \frac{r^2}{\ep_1}<|z|\leq 2^{-\frac{2}{3}} r^\frac{5}{3}\\
  P(z)w(0,0) & \textup{ if } \frac{r\sqrt{r}}{2}\leq|z|\leq 2\sqrt{r}\\
  w(z,0) &\textup{ if } (4r)^\frac{1}{3} \leq|z|<\ep_1 
\end{array}\right.
\]

We apply above construction to each interior node and similar
construction to each boundary node to obtain a linear map 
\begin{eqnarray*}
&&g'_{(0,\eta,\eta')}:
W^{1,p}(\Si,\bS,u^*TX,(u|_{\bS})^*TL)  \\
&&\ra W^{1,p}(\Si_{(0,\eta,\eta')},\bS_{(0,\eta,\eta')},
u^*_{\eta,\eta'}TX, (u_{\eta,\eta'}|_{\bS_{(0,\eta,\eta')}})^*TL),
\end{eqnarray*}
which can also be viewed as a map
\begin{eqnarray*}
&&g'_{(\xi,\eta,\eta')}: W^{1,p}(\Si,\bS,u^*TX,(u|_{\bS})^*TL)\\
&&\ra W^{1,p}(\Si_{(\xi,\eta,\eta')},\bS_{(\xi,\eta,\eta')},
u^*_{\eta,\eta'}TX, (u_{\eta,\eta'}|_{\bS_{(\xi,\eta,\eta')}})^*TL).
\end{eqnarray*}
The restriction of $g'_{(\xi,\eta,\eta')}$ to $\mathrm{Ker}(\pi\circ D_u)$
is injective by the unique continuity theorem (\cite{Ar}). Lemma~\ref{cutsection}
implies the following estimate for $w\in\mathrm{Ker}(\pi\circ D_u)$:
\[
 \parallel \pi_{(\xi,\eta,\eta')}\circ D_{u_{\eta,\eta'}}
  \circ g'_{(\xi,\eta,\eta')} w\parallel_{W^{1,p}}\leq 
  C(|\xi|+\log(|\eta|+|\eta'|)^{\frac{1}{p}-1})\parallel w \parallel_{W^{1,p}}.
\]

Let 
\[
i_{(\xi,\eta,\eta')}=(Id-Q_{(\xi,\eta,\eta')}\circ 
\pi_{(\xi,\eta,\eta')}\circ D_{(\xi,\eta,\eta'),u_{\eta,\eta'}})\circ
g'_{(\xi,\eta,\eta')}
\]
Then $i_{(\xi,\eta,\eta')}:\mathrm{Ker}(\pi\circ D_u)\ra
\mathrm{Ker}(\pi_{(\xi,\eta,\eta')}\circ D_{(\xi,\eta,\eta'),u_{\eta,\eta'}})$ 
is injective for $(\xi,\eta,\eta')$ sufficiently
small. It is actually a linear isomorphism since
\[
\dim \mathrm{Ker}(\pi_{(\xi,\eta,\eta')}\circ D_{(\xi,\eta,\eta'),u_{\eta,\eta'}}) 
=\mathrm{Ind} D_u +\dim E_\rho=\dim \mathrm{Ker}(\pi\circ D_u).
\]

We are now ready to find exact solutions near the approximate solution
$u_{\eta,\eta'}$. 

\begin{pro} \label{surjective}
There exist $\delta_2,d_1,\ldots, d_{l_0},d'_1,\ldots,d'_{l_1},\ep_1,\ep_2>0$
sufficiently small such that for all 
$(\xi,\eta,\eta')\in B_{\delta_2}\times D_d\times D_{d'}$,
if
\[ 
w\in\mathrm{Ker}
(\pi_{(\xi,\eta,\eta')}\circ D_{(\xi,\eta,\eta'),u_{\eta,\eta'}}),\;\;\; 
\parallel w \parallel_{W^{1,p}}\leq\ep_1,
\]
then there exists a unique
\[
h(\xi,\eta,\eta',w)\in L^p(\Si_{(\xi,\eta,\eta')},
\Lambda^{0,1}_{\Si_{(\xi,\eta,\eta')}}\otimes u_{\eta,\eta'}^*TX)
\] 
such that
\[
  \pi_{(\xi,\eta,\eta')}\circ\dbar_{J,\Si_{(\xi,\eta,\eta')}}
  \exp_{u_{\eta,\eta'}}
  (w+Q_{(\xi,\eta,\eta')} h(\xi,\eta,\eta',w))=0
\]
and
\[
  \parallel h(\xi,\eta,\eta',w)\parallel _{L^p}\leq \ep_2.
\]
\end{pro}

\paragraph{Proof.}
We assume that $(\xi,\eta,\eta')\in B_{\delta_2}\times D_d\times D'_{d'}$,
and $\parallel w \parallel_{W^{1,p}}\leq \ep_1$, where $\delta_2$,
$d=(d_1,\ldots,d_{l_0})$, $d'=(d'_1,\ldots,d'_{l_1})$, and $\ep_1$
will be determined later.

We will use Newton's method to find $h(\xi,\eta,\eta',w)$ as in
\cite[Theorem 3.3.4]{MS}. For convenience, write $v$ for 
$u_{\eta,\eta'}$, $Q$ for $Q_{(\xi,\eta,\eta')}$, 
$\dbar$ for $\dbar_{J,\Si_{(\xi,\eta,\eta')}}$. 
Set $h_0=0$, and 
\[
  h_{n+1}=h_n-P_n\circ\pi\circ\dbar\exp_v (w+Q h_n)
\]
where $P_n$ is the parallel transport along the geodesic
$t\in[0,1] \mapsto\exp_v((1-t)(w+Q h_n))$. Let $D_n$ denote the linearization
of $\dbar_{J,\Si_{(\xi,\eta,\eta')}}$ at $v_n=\exp_v(w+Qh_n)$, and 
write $\pi_n$ for $\pi_{(\lambda_{(\xi,\eta,\eta')},v_n)}$. We have
\begin{eqnarray*}
&& P_{n+1}\circ\pi\circ\dbar\exp_v(w+Q h_{n+1})\\
&=&P_{n+1}\circ \pi\circ\dbar\exp_v(w+Q h_n- Q\circ P_n\circ\pi\circ\dbar\exp_v(w+Q h_n))\\
&=&P_n\circ\pi\circ\dbar\exp_v(w+Qh_n)\\
 && -P_n\circ\pi_n\circ D_n\circ (d\exp_v)(w+Qh_n)
   (Q\circ P_n\circ\pi\circ\dbar\exp_v(w+Q h_n))\\
&& +R(\pi\circ\dbar\exp_v(w+Q h_n)),
\end{eqnarray*}
where 
\[
\parallel R(\pi\circ\dbar\exp_v((w+Q h_n)\parallel_{L^p}
\leq C_1\parallel \pi\circ\dbar\exp_v((w+Q h_n))\parallel_{L^p} ^2,    
\]
\[
   \parallel P_n\circ\pi_n\circ D_n\circ (d\exp_v)(w+Qh_n)
   -\pi\circ D_{(\xi,\eta,\eta'),v} \parallel\leq C_2 (|w|+|Q h_n|).
\]
Therefore,
\begin{eqnarray*}
  && \parallel \pi\circ\dbar\exp_v(w+Q h_{n+1})\parallel_{L^p}\\ 
  &\leq& C_3\left(\parallel w\parallel_{C^0}+\parallel Q h_n\parallel_{C^0}
       +\parallel \pi\circ\dbar\exp_v(w+Q h_n)\parallel_{L^p}\right)
        \parallel \pi\circ\dbar\exp_v(w+Q h_n)\parallel_{L^p}.
\end{eqnarray*}

We also have 
\[ 
  h_n=-\sum_{k=0}^{n-1}P_k\circ\pi\circ\dbar\exp_v(w+Q h_k),  
\]
so
\begin{eqnarray*}
 && \parallel\pi\circ\dbar\exp_v(w+Q h_{n+1})\parallel_{L^p}\\
  &\leq& C_3\left(\parallel w\parallel_{W^{1,p}}+
   \sum_{k=0}^n\parallel\pi\circ\dbar\exp_v(w+Q h_k)\parallel_{L^p}\right)
  \parallel \pi\dbar\exp_v(w+Q h_n)\parallel_{L^p}.
\end{eqnarray*}

 Let $a_n= \parallel\pi\circ\dbar\exp_v(w+Q h_n)\parallel_{L^p}\geq 0$ and 
$b=\parallel w \parallel_{W^{1,p}}\geq 0$, so that
\[ 
a_{n+1}\leq C_3\left(b+\sum_{k=0}^n a_k\right)a_n
\]
We will show that $a_0, b\leq \frac{1}{6C_3} \Rightarrow a_{n+1}\leq\half a_n$.
We prove by induction. For $n=0$, we have
\[ a_1\leq C_3(b+a_0)a_0\leq
C_3\left(\frac{1}{6C_3}+\frac{1}{6C_3}\right)a_0=\frac{1}{3}a_0\leq\half a_0.
\]
Now suppose that $a_{n+1}\leq a_n$ for $n=0,1,\ldots,m$. Then
\[ 
\sum_{k=0}^{m+1}a_k \leq 2 a_0,
\]
so 
\[
  a_{m+2}\leq C_3(b+2a_0)a_{m+1}\leq
  C_3\left(\frac{1}{6C_3}+\frac{1}{3C_3}\right)a_{m+1}
  =\half a_{m+1}.
\]

Let $D=D_{(\xi,\eta,\eta'),u_{\eta,\eta'}}$. Then
\[
\pi\circ\dbar\exp_v w=\pi\circ\dbar v+ \pi\circ D w +R(w),
\]
where $\parallel R(w)\parallel_{L^p}\leq 
C_4\parallel w\parallel_{C^0}\parallel w\parallel_{W^{1,p}}$, and
$\pi\circ D w=0$.
The proof of Lemma~\ref{almosthol} can be modified
to show that 
\[
\parallel \pi\circ\dbar v \parallel_{L^p}\leq 
C_5\left(|\xi|+(|\eta|+|\eta'|)^{\frac{1}{2p}}\right).
\]
So
\begin{eqnarray*}
a_0&=&\parallel \pi\circ\dbar\exp_v w\parallel_{L^p}\\
   &\leq& C_5\left(|\xi|+(|\eta|+|\eta'|)^{\frac{1}{2p}}\right)
     + C_6\parallel w\parallel_{W^{1,p}}^2\\
   &\leq& C_5\left(\delta_2+(|d|+|d'|)^{\frac{1}{2p}}\right)+C_6\ep_1^2,
\end{eqnarray*}
which is less than $\frac{1}{6C_3}$ for sufficiently small
$\delta_2,d,d',\ep_1$. $b=\parallel w\parallel_{W^{1,p}}\leq\ep_1$, which is also less than
$\frac{1}{6C_3}$ for small $\ep_1$. 
So if $\delta_2,d,d',\ep_1$ are sufficiently small, $h_n$ converges uniformly in
$L^p$ norm, and the limit  $h(\xi,\eta,\eta',w)$ satisfies
\[
\pi\circ\dbar\exp_v(w+Q h(\xi,\eta,\eta',w))=0,
\]
\[
\parallel h(\xi,\eta,\eta',w) \parallel_{L^{p}}\leq 2a_0\leq
2C_5(\delta_2+(|d|+|d'|)^{\frac{1}{2p}})+2C_6\ep_1^2.
\]
This proves the existence of $h(\xi,\eta,\eta',w)$.

We now turn to the uniqueness of $h(\xi,\eta,\eta',w)$. Write $s_1$ for 
$h(\xi,\eta,\eta',w)$, and suppose that $s_2$ satisfies 
\[
\pi\circ\dbar\exp_v(w+Q s_2)=0
\]
and
\[
\parallel s_2 \parallel_{L^{p}}\leq 2C_5(\delta_2+(|d|+|d'|)^{\frac{1}{2p}})+2C_6\ep_1^2.
\]
Let $P_1$, $P_2$ denote the parallel transports along the geodesics
$t\in[0,1]\mapsto\exp_v((1-t)(w+Qs_1)), \exp_v((1-t)(w+Qs_2))$, respectively.
Let $D'$ denote the linearization of $\dbar_{J,\Si_{(\xi,\eta,\eta')}}$
at $v'=\exp_v(w+Qs_1)$, and let $\pi'$ denote 
$\pi_{(\lambda_{(\xi,\eta,\eta')},v')}$. We have
\begin{eqnarray*}
0&=& P_2\circ\pi\circ \dbar\exp_v (w+Qs_2)\\
 &=& P_2\circ\pi\circ\dbar\exp_v(w+ Qs_1 +Q(s_2-s_1))\\
 &=& P_1\circ\pi\circ\dbar\exp_v(w+ Qs_1) 
     +P_1\circ\pi'\circ D'\circ (d\exp_v)(w+ Qs_1)\circ Q(s_2-s_1)\\
   && +R(s_2-s_1)\\
  &=&P_1\circ\pi'\circ D'\circ (d\exp_v)(w+ Qs_1)\circ Q(s_2-s_1)+R(s_2-s_1),
\end{eqnarray*}
\[
s_1-s_2=\left(P_1\circ\pi'\circ D'\circ (d\exp_v)(w+ Qs_1)-\pi\circ
      D\right)\circ Q(s_2-s_1)+R(s_2-s_1),
\]
where 
\[ 
\parallel P_1\circ\pi'\circ D'\circ(d\exp_v)-\pi\circ D \parallel
 \leq C_2 (|w|+|Qs_1|),
\]
and
\[
\parallel R(s_2-s_1) \parallel_{L^p} \leq 
\parallel s_2-s_1 \parallel_{L^p}^2.
\]
So
\begin{eqnarray*}
\parallel s_1-s_2\parallel_{L^p}&\leq& 
  C_7(\parallel w\parallel_{W^{1,p}} + 
 \parallel s_1\parallel_{L^p} +\parallel s_2\parallel_{L^p})
 \parallel s_1-s_2 \parallel_{L^p}\\
 &\leq& C_7 \left(\ep_1+4C_5\left(\delta_2+(|d|+|d'|)^\frac{1}{2p}\right)
        +4C_6\ep_1^2\right)\parallel s_1-s_2 \parallel_{L^p},
\end{eqnarray*} 
where  
\[
C_7 \left(\ep_1+4C_5\left(\delta_2+(|d|+|d'|)^\frac{1}{2p}\right)
  +4C_6\ep_1^2\right)\leq\half
\]
for sufficiently small $\delta_2,d,d',\ep_1$. We conclude that
$s_1=s_2$. 

The proposition is true if
$\delta_2,d_1,\ldots,d_{l_0},d'_1,\ldots,d'_{l_1},
\ep_1,\ep_2>0$ are chosen such that
\[
  C_5\left(\delta_2+(|d|+|d'|)^\frac{1}{2p}\right)+C_6\ep_1^2 \leq\frac{1}{6C_3},\;\;\;
  \ep_1\leq \frac{1}{6C_3},
\]
\[
 C_7\left(\ep_1+4C_5\left(\delta_2+(|d|+|d'|)^\frac{1}{2p}\right)+4C_6\ep_1^2\right)
  \leq\frac{1}{2},
\]
and $\ep_2=2C_5(\delta_2+(|d|+|d'|)^{\frac{1}{2p}})+2C_6\ep_1^2$. $\Box$

\subsubsection{Kuranishi neighborhood} \label{homeo}

Let $\delta_2,d_1,\ldots,d_{l_0}$, $d_1',\ldots,d_{l_1}'>0$ be chosen as in 
Proposition~\ref{surjective}. Let
$V_{\rho,\mathrm{map}}\subset H_{\rho,\mathrm{map}}=\mathrm{Ker}(\pi\circ D_u)$
be a neighborhood of the origin such that
\[
\parallel i_{(\xi,\eta,\eta')} w \parallel_{W^{1,p}} <\ep_1\;\;
\textup{ for all } w\in V_{\rho,\mathrm{map}}, 
(\xi,\eta,\eta')\in B_{\delta_2}\times D_d\times D'_{d'}.
\]

Write $\tilde{B}$ for $B_{\delta_2}\times D_d\times D'_{d'}$, as at the
beginning of Section~\ref{gluing}. Define a map
\begin{eqnarray*}
\psi: \tilde{B}\times V_{\rho,\mathrm{map}}
&\longrightarrow& W_B\\
(\xi,\eta,\eta',w)&\mapsto& [(\lambda_{(\xi,\eta,\eta')},u_{(\xi,\eta,\eta',w)})]
\end{eqnarray*}
where $W_B$ is defined as at the beginning of Section~\ref{gluing}, and
\[
 u_{(\xi,\eta,\eta',w)}=\exp_{u_{\eta,\eta'}}
 \left(i_{(\xi,\eta,\eta')}w+Q_{(\xi,\eta,\eta')}
 h(\xi,\eta,\eta',i_{(\xi,\eta,\eta')}w)\right).
\]
Then
\[
   \pi_{(\lambda_{(\xi,\eta,\eta')},u_{(\xi,\eta,\eta',w)})}\circ 
   \dbar_{J,\Si_{(\xi,\eta,\eta')}} u_{(\xi,\eta,\eta',w)}=0,
\]
so $(\xi,\eta,\eta',w)\mapsto \dbar_{J,\Si_{(\xi,\eta,\eta')}}
u_{(\xi,\eta,\eta',w)}$ defines a map 
$s: \tilde{B}\times V_{\rho,\mathrm{map}}\ra E_\rho$
such that $\psi(s^{-1}(0))\subset\MXL$. Actually,
$\psi(s^{-1}(0))$ contains a neighborhood of $\rho$ in $\MXL$.
To see this, note that any $\rho'\in\MXL$ which is 
sufficiently close to $\rho$ in $C^\infty$ topology can be written
in the form
\[
  \rho'=(\lambda_{(\xi,\eta,\eta')}, \exp_{u_{\eta,\eta'}}(w))
\]   
where $w\in W^{1,p}(\Si_{(\xi,\eta,\eta')},\bS_{(\xi,\eta,\eta')},
u_{\eta,\eta'}^* TX, (u_{\eta,\eta'}|_{\bS_{(\xi,\eta,\eta')}})^* TL)$
is small. There exist unique 
$w_0\in \mathrm{Ker}(\pi_{(\xi,\eta,\eta')}\circ D_{(\xi,\eta,\eta'),u_{\eta,\eta'}})$
and $h\in L^p(\Si_{(\xi,\eta,\eta')},
\Lambda^{0,1}\Si_{(\xi,\eta,\eta')}\otimes u_{\eta,\eta'}^* TX)$ such
that $w=w_0+Q_{(\xi,\eta,\eta')}h$. We may further assume that $w$ is so small such that
$w_0=i_{(\xi,\eta,\eta')}w_1$ for some $w_1\in V_{\rho,\mathrm{map}}$, and 
$\parallel h\parallel_{L^p}\leq \ep_2$. Since
\[ 
  \dbar_{J,\Si_{(\xi,\eta,\eta')}}\exp_{u_{\eta,\eta'}}
  (w_0+Q_{(\xi,\eta,\eta')}h)=0,
\]
we have $h=h(\xi,\eta,\eta',w_0)$ by the uniqueness part of Proposition~\ref{surjective}.
So $\rho'=\psi(\xi,\eta,\eta',w_1)$, and $s(\xi,\eta,\eta',w_1)=0$.

Now assume that $\lambda$ is stable, so that $\domaut$ is finite, and
$V'_{\rho,\mathrm{map}}=V_{\rho,\mathrm{map}}$. Then there is an isomorphism
$\phi:\lambda_{(\xi_1,\eta_1,\eta'_1)}\ra \lambda_{(\xi_2,\eta_2,\eta'_2)}$
if and only if $(\xi_2,\eta_2,\eta'_2)=
(\phi'\cdot\xi_1,\phi'\cdot\eta_1,\phi'\cdot\eta'_1)$ for some $\phi'\in\domaut$,
and the action of $\phi'$ on the universal family restricts to $\phi$ on 
$\lambda_{(\xi_1,\eta_1,\eta'_1)}$.
We may choose $\tilde{B}$, $V_{\rho,\mathrm{map}}$ small enough such that if 
$u_{(\xi,\eta,\eta',w_1)}=u_{(\phi\cdot\xi,\phi\cdot\eta,\phi\cdot\eta',w_2)}\circ\phi$
for $(\xi,\eta,\eta')\in\tilde{B}$, $w_1,w_2\in V_{\rho,\mathrm{map}}$, 
$\phi\in \domaut$, then $\phi\in \mapaut$. 

Given $\phi\in\mapaut$ and $w\in V_{\rho,\mathrm{map}}$, let 
$\phi_{(\xi,\eta,\eta')}w$ be the unique vector in $V_{\rho,\mathrm{map}}$
such that $u_{(\xi,\eta,\eta',w)}\circ\phi^{-1}=
u_{(\phi\cdot\xi,\phi\cdot\eta,\phi\cdot\eta',\phi_{(\xi,\eta,\eta')}w)}$.
Then $\mapaut$ acts on $\tilde{B}\times V_{\rho,\mathrm{map}}$ by
$\phi\cdot(\xi,\eta,\eta',w)=(\phi\cdot\xi,\phi\cdot\eta,\phi\cdot\eta',
\phi_{(\xi,\eta,\eta')}w)$. From the above discussion,
$\psi(\xi_1,\eta_1,\eta'_1,w_1)=\psi(\xi_2,\eta_2,\eta'_2,w_2)$ if and only
if $(\xi_2,\eta_2,\eta'_2,w_2)=\phi\cdot(\xi_1,\eta_1,\eta'_1,w_1)$ for 
some $\phi\in \mapaut$. 
Let $V_\rho$ be an $\mapaut$-invariant neighborhood
of the origin in $\tilde{B}\times V_{\rho,\mathrm{map}}\subset 
H_{\rho,\mathrm{domain}}\times H_{\rho,\mathrm{map}}=H_\rho$.
Let $s_\rho: V_\rho\ra E_\rho$ be the restriction of $s$. Then 
$s_\rho$ is $\mapaut$-equivariant. The restriction of $\psi$ gives a
continuous map $\psi_\rho:s_\rho^{-1}(0)\ra \MXL$ such that
$s_\rho^{-1}(0)/\mapaut\ra \MXL$ is injective. It is actually a homeomorphism
onto a neighborhood of $\rho$ in $\MXL$ since $C^\infty$ topology and
$W^{1,p}$ topology are equivalent on $\MXL$.
This completes the proof of Theorem~\ref{chart} when the domain $\lambda$
of $\rho$ is stable.

\begin{rem} \label{nonsmooth}
The section $s$ in the above construction is only continuous, not smooth. 
It is not hard to show that it is smooth within the stratum, but its dependence
on $\eta,\eta'$ is not even $C^1$ (see Lemma~\ref{almosthol}).
This is because the smooth structure of $V_\rho$  is canonical within
the stratum but dependent on our particular gluing construction
in the direction $(\eta,\eta')$ transversal to the stratum.
\end{rem}

If $\lambda=\domain$ is not stable, we add minimal number of
marked points to obtain a stable marked bordered Riemann surface
\[
\tilde{\lambda}=(\Si,\bB,\tilde{\bp};\tilde{\bq}^1,\ldots,
\tilde{\bq}^h)
\] 
where $\tilde{\bp}=(p_1,\ldots,p_{n+\hat{n}})$,
$\tilde{\bq}^i=(q^i_1,\ldots,q^i_{m^i+\hat{m}^i})$,
$p_{n+1},\ldots,p_{n+\hat{n}}$ are additional interior marked points,
and $q^i_{m^i+1},\ldots,q^i_{m^i+\hat{m}^i}$ are additional marked points
on $B^i$. Note that when counting the above minimal number, an interior 
marked point counts twice, while a boundary marked point counts once. We have 
$$
2\hat{n}+\hat{m}^1+\cdots \hat{m}^h=\dim_\R H_{\rho,\mathrm{aut}}.
$$

Let $\tilde{\rho}=(\tilde{\lambda},u)\in\SXL$. Then
$H_{\tilde{\rho},\mathrm{domain}}=H_{\rho,\mathrm{domain}}$ and 
$H_{\tilde{\rho},\mathrm{aut}}=0$.
If the additional marked points are chosen in $K_{3\sqrt{|d|+|d'|}}(\Si)$,
then the construction of $\lambda_{(\xi,\eta,\eta')}$ also yields
deformation $\tilde{\lambda}_{(\xi,\eta,\eta')}$ of 
$\tilde{\lambda}$ for $(\xi,\eta,\eta')\in 
B_{\delta_2}\times D_d\times D'_{d'}=\tilde{B}\subset 
H_{\tilde{\rho},\mathrm{domain}}=H_{\rho,\mathrm{domain}} $.
Both $\mathrm{Aut}\,\tilde{\lambda}$ and $\mapaut$ are subgroups
of $\domaut$, and $\mathrm{Aut}\,\tilde{\rho}=
\mathrm{Aut}\,\tilde{\lambda}\cap \mapaut$.
We may choose $E_{\tilde{\rho}}=E_\rho$, so that
$H_{\tilde{\rho},\mathrm{map}}= H_{\rho,\mathrm{map}}$, 
$H_{\tilde{\rho}}=H_\rho$.

There is a map
\[ 
F:\SXL\ra\MXL
\]
defined by forgetting $p_{n+1},\ldots,p_{n+\hat{n}},q_{m+1},\ldots,q_{m+\hat{m}}$ and
then contracting non-stable components which are mapped to points.
We have $F(\tilde{\rho})=\rho$.
We will construct a multi-valued map $A$ from a neighborhood $U_\rho$ of $\rho$ in $\MXL$ to $\SXL$
such that $F\circ A$ is the identity map, and $A(\rho)=\tilde{\rho}$.

The additional marked points are on
non-stable components, where $u$ is not a constant.
Let $2N$ be the dimension of $X$, as before.
For $j=n+1,\ldots,n+\hat{n}$, we may assume that there is a geodesic ball 
in $X$ centered at $u(p_j)$ such that its intersection with the image
of $u$ is an embedded holomorphic disc $D_j$, and
$u^{-1}(D_j)\ra D_j$ is a trivial cover. $u^{-1}(D_j)$ might consist
of more than one connected components if $u$ is not injective. 
Let $B_j\subset X$ be an embedded $(2N-2)$-dimensional ball which is 
the image of $N_{p_j}$ under $\exp_{u(p_j)}$, 
where $N_{p_j}$ is a small ball centered at the origin 
in the orthogonal complement of $u_* (T_{p_j}\Si)$ in $T_{u(p_j)}X$.
Then $D_j$ and $B_j$ intersect orthogonally at $u(p_j)$. 
We choose $N_{p_j}$ sufficiently small such that $B_j$ intersects
the image of $u$ at a single point $u(p_j)$.

Similarly, for $j^i=m^i+1,\ldots,m^i+\hat{m}^i$, 
we may assume that there is a geodesic ball 
in $X$ centered at $u(q_{j^i})$ such that its intersection with the image
of $u$ is an embedded holomorphic half disc $D^+_{j^i}$, and
$u^{-1}(D^+_{j^i})\ra D^+_{j^i}$ is a trivial cover. 
Let $B'_{j^i}\subset L$ be an embedded $(N-1)$-dimensional ball which is
the image of $N'_{q_{j^i}}$ under $\exp_{u(q_{j^i})}$, 
where $N'_{q_{j^i}}$ is a small ball centered at the origin
in the orthogonal complement of $u_* (T_{q_{j^i}}\bS)$ in $T_{u(q_{j^i})}L$. 
Then $I_{j^i}=D^+_{j^i}\cap L\cong[0,1]$ and $B'_{j^i}$ intersect orthogonally
in $L$ at $u(q_{j^i})$.
We choose $N'_{q_{j^i}}$ sufficiently small such that $B'_{j^i}$ intersects
the image of $u$ at a single point $u(q_{j^i})$.

In a small neighborhood of $\rho$ in $\MXL$, intersecting with
$B_j$ for $j=n+1,\ldots,n+\hat{n}$ and 
$B'_{j^i}$ for $j^i=m^i+1,\ldots,m^i+\hat{m}^i$ determines additional marked
points and gives the desired multi-valued map
$A:U_\rho\ra \SXL$ for some neighborhood of $\rho$ in $\MXL$.
$A$ is single-valued if $\mathrm{Aut}\;\tilde{\rho}=\mapaut$. Let
\begin{eqnarray*}
\tilde{U}'_\rho&=&\left\{\tilde{\rho}'=\anothermap\in
F^{-1}(U_\rho)\mid\right.\\
&&\left. u'(p_j)\in B_j\textup{ for }j=n+1,\ldots,n+\hat{n},\;
 u'(q_{j^i})\in B'_{j^i}\textup{ for }j^i=m^i+1,\ldots,m^i+\hat{m}^i\right\},
\end{eqnarray*}
and let $\tilde{U}_\rho$ be the connected component containing
$\tilde{\rho}$. Then $\tilde{U}_\rho=A(U_\rho)$, and 
the fiber of $\tilde{U}_\rho\ra  U_\rho$ is finite.
 
Let $\hat{W}$ be the (real) codimension 
$(2\hat{n}+\hat{m}^1+\cdots+\hat{m}^h)$ subspace of 
$W^{1,p}(\Si,\bS,u^*TX, (u|_{\bS})^*TL)$ defined by
\[ 
\hat{W}=\left\{w\in W^{1,p}(\Si,\bS, u^*TX, (u|_{\bS})^*TL)\left| 
\begin{array}{ll}
w(p_j)\in T_{u(p_j)} B_j&\textup{for }j=n+1,\ldots,n+\hat{n}\\
w(q_{j^i})\in T_{u(q_{j^i})} B'_{j^i}&\textup{for }j^i=m^i+1,\ldots,m^i+\hat{m}^i
\end{array}\right.\right\}
\]
Then $H_{\rho,\mathrm{aut}}\cap \hat{W}=\{0\}$, so 
\[
W^{1,p}(\Si,\bS, u^*TX, (u|_{\bS})^*TL)=H_{\rho,\mathrm{aut}}\oplus \hat{W}.
\]

Let 
\begin{eqnarray*}
\hat{W}_{(\xi,\eta,\eta')}&=&\left\{
w\in W^{1,p}(\Si_{(\xi,\eta,\eta')},\bS_{(\xi,\eta,\eta')}, 
u_{\eta,\eta'}^*TX,(u_{\eta,\eta'}|_{\bS_{(\xi,\eta,\eta')}})^*TL)\right | \\
&&w(p_j)\in T_{u(p_j)} B_j\textup{ for }j=n+1,\ldots,n+\hat{n}\\
&&\left. w(q_{j^i})\in T_{u(q_{j^i})} B'_{j'}\textup{ for }
j^i=m^i+1,\ldots,m^i+\hat{m}^i 
  \right\}.
\end{eqnarray*}
Then
\[
 W^{1,p}(\Si_{(\xi,\eta,\eta')},\bS_{(\xi,\eta,\eta')}, 
u_{\eta,\eta'}^*TX,(u_{\eta,\eta'}|_{\bS_{(\xi,\eta,\eta')}})^*TL)
=i_{(\xi,\eta,\eta')}H_{\rho,\mathrm{aut}}\oplus \hat{W}_{(\xi,\eta,\eta')}
\]
for $(\xi,\eta,\eta')$ sufficiently small.
Let 
\[
p_{(\xi,\eta,\eta')}:W^{1,p}(\Si_{(\xi,\eta,\eta')},\bS_{(\xi,\eta,\eta')}, 
u_{\eta,\eta'}^*TX,(u_{\eta,\eta'}|_{\bS_{(\xi,\eta,\eta')}})^*TL)\ra
\hat{W}_{(\xi,\eta,\eta')}
\]
be the projection. We have
$\pi_{(\xi,\eta,\eta')}\circ D_{(\xi,\eta,\eta'),u_{\eta,\eta'}}\circ
p_{(\xi,\eta,\eta')}=
\pi_{(\xi,\eta,\eta')}\circ D_{(\xi,\eta,\eta'),u_{\eta,\eta'}}$ since
$i_{(\xi,\eta,\eta')}H_{\rho,\mathrm{aut}}\subset 
\mathrm{Ker}(\pi_{(\xi,\eta,\eta')}\circ
D_{(\xi,\eta,\eta'),u_{\eta,\eta'}})$.
By replacing $Q_{(\xi,\eta,\eta')}$ with 
$p_{(\xi,\eta,\eta')}\circ Q_{(\xi,\eta,\eta')}$, we may assume that
$\mathrm{Im}Q_{(\xi,\eta,\eta')}\subset \hat{W}_{(\xi,\eta,\eta')}$.

Let 
$\hat{H}_{\tilde{\rho},\mathrm{map}}=
H_{\tilde{\rho},\mathrm{map}}\cap \hat{W}\cong H'_{\rho,\mathrm{map}}$.
We may modify the isomorphism
$ i_{(\xi,\eta,\eta')}:H_{\rho,\mathrm{map}}\ra
\mathrm{Ker}(\pi_{(\xi,\eta,\eta')}\circ D_{(\xi,\eta,\eta'),u_{\eta,\eta'}})$ such that
\[
i_{(\xi,\eta,\eta')}\hat{H}_{\tilde{\rho},\mathrm{map}}=
\mathrm{Ker}(\pi_{(\xi,\eta,\eta')}\circ
D_{(\xi,\eta,\eta'),u_{\eta,\eta'}})\cap \hat{W}_{(\xi,\eta,\eta')}.
\]

For $(\xi,\eta,\eta',w)\in \tilde{B}\times V_{\tilde{\rho},\mathrm{map}}
=\tilde{B}\times V_{\rho,\mathrm{map}} $, define
$\tilde{\psi}(\xi,\eta,\eta',w)=(\tilde{\lambda}_{(\xi,\eta,\eta')},u_{(\xi,\eta,\eta',w)})$.
For $w\in \hat{W}$, we have
$i_{(\xi,\eta,\eta')}w+Q_{(\xi,\eta,\eta')}
h(\xi,\eta,\eta',i_{(\xi,\eta,\eta')}w)\in \hat{W}_{(\xi,\eta,\eta')}$, so
\[
u_{(\xi,\eta,\eta',w)}=\exp_{u_{\eta,\eta'}}(i_{(\xi,\eta,\eta')}w+Q_{(\xi,\eta,\eta')}
h(\xi,\eta,\eta',i_{(\xi,\eta,\eta')})\in \hat{W}_{(\xi,\eta,\eta')}
\] 
satisfies $u_{(\xi,\eta,\eta',w)}(p_j)\in B_j$ for $j=n+1,\ldots,n+\hat{n}$,
and $u_{(\xi,\eta,\eta',w)}(q_{j'})\in B'_{j'}$ for $j=n+1,\ldots,n+\hat{n}$.
So $\tilde{\psi}(\xi,\eta,\eta',w)\in \tilde{U}_\rho$ if 
$w\in\hat{W}$, $s(\xi,\eta,\eta',w)=0$, and $\xi,\eta,\eta',w$ are
sufficiently small.

Conversely, if $(\tilde{\lambda}',u')$ is a stable map near $\tilde{\rho}$
in $\tilde{U}_\rho$, then $\tilde{\lambda}'=\lambda_{(\xi,\eta,\eta')}$ for 
some $(\xi,\eta,\eta')\in \tilde{B}$, and
$u'=\exp_{u_{\eta,\eta'}}(w)$ for some $w\in \hat{W}_{(\xi,\eta,\eta')}$.
There exist unique $w_0\in \mathrm{Ker}(\pi_{\xi,\eta,\eta'}\circ 
D_{(\xi,\eta,\eta'),u_{\eta,\eta'}})\cap \hat{W}_{(\xi,\eta,\eta')}$ and
$h\in L^p(\Si_{(\xi,\eta,\eta')},\Lambda^{0,1}\Si_{(\xi,0,0)}\otimes u_{\eta,\eta'}^*TX)$   
such that $w=w_0+Q_{(\xi,\eta,\eta')}h$. We have $h=h(\xi,\eta,\eta',w_0)$ since
\[
  \dbar_{J,\Si_{(\xi,\eta,\eta')}} \exp_{u_{\eta,\eta'}}(w_0+Q_{(\xi,\eta,\eta')}h)=0.
\] 
Let $w_1=i_{(\xi,\eta,\eta')}^{-1}w_0\in \hat{H}_{\tilde{\rho},\mathrm{map}}$. 
Then $u'=u(\xi,\eta,\eta',w_1)$, and $s(\xi,\eta,\eta',w_1)=0$.
Let $\hat{V}_{\rho,\mathrm{map}}=
V_{\rho,\mathrm{map}}\cap \hat{H}_{\tilde{\rho},\mathrm{map}}$. Then
$\mathrm{Aut}\,\tilde{\rho}$ acts on 
$\tilde{B}\times \hat{V}_{\tilde{\rho},\mathrm{map}}$. Let 
$\hat{V}_{\tilde{\rho}}$ be an $\mathrm{Aut}\,\tilde{\rho}$ invariant
neighborhood of the origin in $\tilde{B}\times
\hat{V}_{\tilde{\rho},\mathrm{map}}$.
It corresponds to a neighborhood $V'_\rho$ of the origin in 
$\tilde{B}\times H'_{\rho,\mathrm{map}}$ under the isomorphism
$\hat{H}_{\tilde{\rho},\mathrm{map}}\cong H'_{\rho,\mathrm{map}}$.
We have the following commutative diagram
\[
\begin{CD}
(s^{-1}(0)\cap \hat{V}_{\tilde{\rho}})/\mathrm{Aut}\,\tilde{\rho}
 @>\hat{\psi}>> U_{\tilde{\rho}}\\
@VVV @VV{F}V\\
(s^{-1}(0)\cap V'_\rho)/\mapaut @>\psi'>> U_\rho
\end{CD}
\]
where $\hat{\psi}$ and $\psi'$ are injective. $\psi'$ is a homeomorphism
onto its image.

\subsection{Transition functions} \label{transition}

For each $\rho\in \cM=\MXL$ and each choice of $E_\rho$,
we have constructed a Kuranishi neighborhood 
$(V'_\rho, E_\rho,\mapaut,\psi_\rho,s_\rho)$.

\begin{rem} \label{choices}
A different choice of $E_\rho$ yields a different, but
equivalent Kuranishi neighborhood in the sense of Definition~\ref{same}. 
Actually, if $E_{1,\rho}, E_{2,\rho}\subset L^p(\Si,\Lambda^{0,1}\otimes u^*
TX)$ are two different choices, set 
$E_\rho=E_{1,\rho}+E_{2,\rho}$. Then
$(V'_\rho,E_\rho,\mapaut, \psi_\rho,s_\rho)$ can serve as the
$(V_p,E_p,\Gamma_p,\psi_p,s_p)$ in Definition~\ref{same}.
\end{rem}

In this section, we will modify the Kuranishi neighborhoods
we constructed so that we can construct transition functions between them.
We will follow \cite[Section 15]{FO} closely.
We refer the reader to \cite[Section 15]{FO} for more details.

$\cM$ is compact in $C^\infty$ topology, so there exist
$\rho_1,\ldots,\rho_l\in \cM$ such that 
\[
\{ U'_i=\psi_{\rho_i}(s_{\rho_i}^{-1}(0)) \mid i=1,\ldots,l \}
\]
is an open cover of $\cM$, and there is $U_i\subset\subset U'_i$ 
such that $\{ U_i\mid i=1,\ldots,l \}$ is still an open
cover of $\cM$. $\cM$ is compact and Hausdorff, so the closure
$K_i$ of $U_i$ is compact for $i=1,\ldots,l$. 

Let $E_i\ra V_i=V_{\rho_i}$ be the
obstruction bundle constructed from 
$E_{\rho_i}$. We may choose $\rho_i$ and $E_{\rho_i}$ such
that if 
\[ 
\rho=\map \in K_{i_1}\cap\ldots\cap K_{i_k}
\]
and $\rho \notin K_i$ if $i\neq i_k$, 
then the subspace $\tilde{E}_\rho$ of $L^p(\Si,\Lambda^{0,1}\Si\otimes u^*TX)$  
spanned by $(E_{i_1})_\rho,\ldots, (E_{i_k})_\rho$ is actually a direct sum 
$\tilde{E}_\rho= (E_{i_1})_\rho\oplus\cdots\oplus (E_{i_k})_\rho$. 
We use $\tilde{E}_\rho$ to construct a Kuranishi neighborhood
$(\tilde{V}_\rho,\tilde{E}_\rho,\mapaut,\tilde{\psi}_\rho,\tilde{s}_\rho)$
as in Section~\ref{gluing}, \ref{homeo}.
We shrink $\tilde{V}_\rho$ such that $\tilde\psi_\rho(\tilde{s}_\rho^{-1}(0))\subset 
U_{i_1}\cap\cdots\cap U_{i_k}\cap K_{j_1}^c\cap\cdots\cap K_{j_{l-k}}^c$, where
$\{ i_1,\ldots,i_k,j_1,\ldots,j_{l-k} \}=\{1,\ldots,l\}$.

Suppose that $\rho'=\tilde{\psi}_\rho(\xi,\eta,\eta',w)$, where 
$(\xi,\eta,\eta',w)\in \tilde{V}_\rho$, $\tilde{s}_\rho(\xi,\eta,\eta',w)=0$.
Then $\mapaut'$ can be identified with the stabilizer $(\mapaut)_{(\xi,\eta.\eta',w)}$
of the action $\mapaut$ on $\tilde{V}_\rho$. This gives a monomorphism
$h_{\rho\rho'}:\mapaut'\ra\mapaut$.
Without loss of generality, we may assume that $\rho'\in K_1\cap\cdots K_k$, and
$\rho'\notin K_i$ for $i=k+1,\ldots,l$. We have $\rho\in K_1\cap\cdots\cap K_k$.
We may assume that $\rho\in K_1\cap\cdots\cap K_{k'}$, and $\rho\notin K_i$
for $i=k'+1,\ldots, l$, where $k'\geq k$. It follows from the construction
that there is an $\mapaut'$-invariant neighborhood $V_{\rho\rho'}$ 
of $0=\psi_{\rho'}^{-1}(\rho')$ in $\tilde{V}_{\rho'}$ such that
we have the following commutative diagram:
\[
\begin{CD}
\tilde{E}_{\rho'}|_{V_{\rho\rho'}} @>{\hat{\phi}_{\rho\rho'}}>> \tilde{E}_{\rho}\\
@VVV @VVV\\
V_{\rho\rho'} @>{\phi_{\rho\rho'}}>> \tilde{V}_{\rho}
\end{CD}
\] 
where $\hat{\phi}_{\rho\rho'}$ is induced by the inclusion
$E_1\oplus\cdots\oplus E_k\ra E_1\oplus\cdots\oplus E_{k'}$.
Both $\hat{\phi}_{\rho\rho'}$ and $\phi_{\rho\rho'}$ are 
$h_{\rho\rho'}$-equivariant embedding of codimension 
$\mathrm{rank}E_{k+1}+\cdots+\mathrm{rank} E_{k'}$. 

The Kuranishi neighborhoods
$(\tilde{V}_\rho,\tilde{E}_\rho,\mapaut,\tilde{\psi}_\rho,\tilde{s}_\rho)$
satisfy the properties listed in Definition~\ref{nbd}, 
and the transition functions
$(V_{\rho\rho'},\hat{\phi}_{\rho\rho'},\phi_{\rho\rho'},h_{\rho\rho'})$
satisfy the properties listed in Definition~\ref{Kuranishi}.
From now on, we will write $(V'_\rho,E_\rho,\mapaut,\psi_\rho,s_\rho)$
for the
$(\tilde{V_\rho},\tilde{E}_\rho,\mapaut,\tilde{\psi}_\rho,\tilde{s}_\rho)$ 
we just constructed. 

\begin{rem}
The $\hat{\phi}_{\rho\rho'},\phi_{\rho\rho'}$ in the above construction 
are smooth when restricted to a stratum. It is possible to refine the above
modification of Kuranishi neighborhoods such that 
$\hat{\phi}_{\rho\rho'},\phi_{\rho\rho'}$ 
are smooth (see \cite[Section 15]{FO}).
Such an refinement is artificial since our smooth structure in directions
transversal to a stratum is not natural, as discussed in
Remark~\ref{nonsmooth}. For simplicity of exposition, we will assume 
the smoothness of $\hat{\phi}_{\rho\rho'},\phi_{\rho\rho'}$ 
in Section~\ref{constructchain}, though such an assumption
is not absolutely necessary for our purpose.
\end{rem}

By Remark~\ref{choices}, the equivalence class of the Kuranishi structure
constructed above is independent of various choices in our construction. 

\subsection{Orientation} \label{orientation}

Recall that the orientation bundle of the Kuranishi structure is the real
line orbibundle obtained by gluing $\det(TV'_\rho)\times \det(E_\rho)^{-1}$.
In our case,
\[
(\det(TV'_\rho)\otimes \det(E_\rho)^{-1})_{(\xi,\eta,\eta',w)}
\cong \det(\mathrm{Ind}D_{(\xi,\eta,\eta'),u_{\eta,\eta'}})\otimes
      \det\left(\mathrm{Ind}(T_{\Si_{(\xi,\eta,\eta')}}, 
       T_{\bS_{(\xi,\eta,\eta')}})\right)^{-1} 
\]
where $\mathrm{Ind}(D)$ denotes the virtual real vector space
$\mathrm{Ker}(D)-\mathrm{Coker}(D)$, as in Section \ref{virtualdim}.

\begin{tm}\label{spin}
 The Kuranishi structure constructed above
 is orientable if $L$ is spin or if $h=1$ and $L$ is relatively spin
 (i.e., $L$ is orientable and $w_2(TL)=\alpha|_L$ for some
 $\alpha\in H^2(X,\Z_2)$). 
\end{tm}

To orient the index of the linearization $D\dbar_J$ of $\dbar_J$ at a stable
map, we need the following generalization of \cite[Proposition 21.3]{FO3}.
\begin{lm}\label{pinch}
 Let $(E,E_\R)\ra (\Si,\bS)$ be a Riemann-Hilbert bundle over a
 prestable bordered Riemann surface without boundary nodes
 $\Si$, and let $E_\R$ be a totally real subbundle of $E|_{\bS}$.
 Then an ordering of the connected components of 
 $\bS$ and a trivialization of $E_\R$ determine an orientation of
 $\mathrm{Ind}(E,E_\R)$, where
 $\mathrm{Ind}(E,E_\R)$ is defined as in \cite[Definition 3.4.1]{KL}.
\end{lm} 

\paragraph{Proof.} 

   Let $B^1,..., B^h$ be the ordered connected components of $\bS$. An
   isomorphism $E_\R\cong\bS\times\R^n$ is a collection of isomorphisms
   $E_\R|_{b^i}\cong B^i\times \R^n$. Let $\ep>0$ be such that 
   $A_i=B(B^i,\ep)$, the collar neighborhood of $B^i$ in $\Si$ of
   radius $\ep$ w.r.t. some admissible metric on $\Si$, are disjoint. 
   By tensoring with $\C$, we have 
   trivializations $E|_{B^i}\cong B^i\times \C^n$. By deforming the 
   Hermitian connection, we may assume that the connection is flat on 
   $A=\bigcup_{i=1}^h B(B^i,\frac{2}{3}\epsilon)$ and there are parallel sections
   $s_1, ..., s_h$ on $E|_A$ such that for $x\in\bS$, $s_i(x)$ corresponds 
   to $(x, e_i)$ under the isomorphism $E_\R\cong\bS\times\R^n$, where
   $e_1,\ldots e_n$ are the standard basis of $\R^n$.      

   The boundary of $A_i=B(B^i,\frac{1}{2})$ is the disjoint union of two
   circles, $B^i$ and $(B')^i$. We shrink $(B')^i$ to obtain a family
   of prestable bordered Riemann surface $\Si_t$, $t\in[0,1]$, such
   that $\Si_1=\Si$, $\Si_t$ are homeomorphic to $\Si$, and $\Si_0$ is
   obtained from $\Si$ by shrinking each $(B')^i$ to a point.
   $\Si_0=C\cup D^1\cup\ldots\cup D^h$ is a prestable bordered Riemann
   surfaces, where $C$ is a complex algebraic curve
   of genus $g$, and $D^i$ is a disc which intersects $C$ at an interior
   node on $\Si_0$, for $i=1,\ldots,h$.

   We extend $(E,E_\R)$ to a family of Riemann-Hilbert bundles
   $(E(t),E_\R(t))\ra
   (\Si_t,\bS_t)$ such that $s_1,\ldots,s_h$ extend to a neighborhood of 
   $\cup_{t\in[0,1]}\bS_t \textup{ in } \cup_{t\in[0,1]}\Si_t$ and
   give a holomorphic trivialization of $E(t)$ in a collar neighborhood
   of $\bS_t$ and a trivialization of $E_\R(t)$. In particular, they
   are defined on $D^i\subset \Si_0$ to give an identification 
   $(E(0),E_\R(0))|_{D_i}=(\C^n,\R^n)$. 
   
   We use the notation in \cite[Section 3.4]{KL}.
   $\mathrm{Ind}(E(t),E_\R(t))$ is a family of virtual real vector spaces
   over $[0,1]$. We have 
\[
   \mathrm{Ind}(E,E_\R)=\mathrm{Ind}(E(1),E_\R(1))\cong
   \mathrm{Ind}(E(0),E_\R(0)),
\]
  so it suffices to orient $\mathrm{Ind}(E(0),E_\R(0))$. We have a long exact
  sequence
\begin{eqnarray*} 
   0&\ra& H^0(\Si_0,\bS_0,E(0),E(0)_\R)\ra 
   H^0(C,F)\oplus\bigoplus_{i=1}^h H^0(D^i,\pa D^i,\C^n,\R^n)\stackrel{e}{\ra} \C^n\\
   &\ra& H^1(\Si_0,\bS_0,E(0),E(0)_\R)\ra 
   H^1(C,F)\oplus\bigoplus_{i=1}^h H^1(D^i,\pa D^i,\C^n,\R^n)\ra 0
   \end{eqnarray*}
where $F=E(0)|_{C}$ is a holomorphic vector bundle of degree
$\half\mu(E,E_\R)$, and $e$ is given by
\begin{eqnarray*}
H^0(C,F)\oplus\bigoplus_{i=1}^h 
H^0(D^i,\partial D^i,\C^n,\R^n)&\ra& \C^n\\
(\xi_0,\xi_1,...,\xi_h)&\mapsto&
(\xi_0(p_1)-\xi_1(0),...,\xi_0(p_h)-\xi_h(0))
\end{eqnarray*}
$p_i\in C$ and $0\in D_i$ are identified to form an interior node
of $\Si_0$. 
$E_\R\simeq \R^n$ gives the orientation on 
$H^0(D^i,\pa D^i,\C^n,\R^n)\simeq \R^n$ via the evaluation map.
$H^1(D^i,\pa D^i,\C^n,\R^n)=0$. 
$H^0(C,F), H^1(C,F)$ and $\C^n$ are complex vector spaces, thus
canonically oriented. Therefore,
\[
   \mathrm{Ind}(E(0),E_\R(0))= H^0(\Si_0,\bS_0,E(0),E(0)_\R)-
                               H^1(\Si_0,\bS_0,E(0),E(0)_\R)
\] 
is oriented. This orientation depends on the trivialization of $E_\R$ and
the ordering of connected components of $\bS$, since the trivialization
of $E_\R$ determines the orientation on
each $H^0(D^i,\pa D^i,\C^n,\R^n)\cong \R^n$ and the ordering of connected 
components of $\bS$ determines the ordering of these $h$ copies of $\R^n$.
$\Box$

\paragraph{\bf Proof of Theorem~\ref{spin}}
It suffices to show that the orientation bundle is trivial when
restricted to each loop $\gamma: S^1\ra \MXL$, $\gamma(t)=\rho_t$.
From the construction of Kuranishi structure we see that we may deform 
$\gamma$ to a family of smooth maps $\tilde{\rho}_t\in V_{\rho_t}$ such that 
the domain of $\tilde{\rho}_t$ are smooth bordered Riemann surfaces. 
We first assume that these domains are stable. 
The tangent bundle of $\Bar{M}_{(g,h),(n,m)}$ is orientable 
by Theorem~\ref{thm:Morientable}, so
it suffices to show that the index bundle $\mathrm{Ind}\,D\dbar$ is
orientable along the loop $\tilde{\gamma}(t)=\tilde{\rho}_t$.
Note that $\mapaut$ preserves the orientation of  $\mathrm{Ind}\,D\dbar$
since it does not permute boundary components of the domain.

Let $\Phi:(\Si,\bS)\times S^1\ra (X,L)$ be given by
$\Phi(z,t)=u_t(z)$, where $u_t$ is the map for $\rho_t$.
Since the connected components of $\bS$ are ordered, by
Lemma~\ref{pinch} it suffices to show that
$(\Phi|_{\bS\times S^1})^* TL$ is stably trivial. 

We first assume that $L$ is relatively spin, i.e,
$L$ is orientable and $w_2(TL)=\alpha|_L$ for some $\alpha\in H^2(X,\Z_2)$. 
Choose a cellular decomposition on $X$ such that $X^{(2)}\cap L=L^{(2)}$.
There exists a real orientable vector bundle of rank $2$
over $X^{(3)}$ such that $w_2(V)=\alpha|_{X^{(3)}}\in H^2(X^{(3)},\Z_2)$.
Then  
\[
w_2((TL\oplus V)|_{L^{(2)}})=0\in H^2(L^{(2)},\Z_2), 
\]
so $(TL\oplus V)|_{L^{(2)}}$ is spin, thus stably trivializable on $L^{(2)}$.

We write
\[ 
  (TL\oplus V)|_{L^{(2)}}\oplus \R^k=\R^{N+k+2}.
\]
Let $\tilde{\Phi}$ be homotopic to $\Phi$ such that
$\tilde{\Phi}(\bS\times S^1)\subset L^{(2)}$. It suffices to show that
$(\tilde{\Phi}|_{\bS\times S^1})^* TL$ is stably trivial.
We have
\[
 (\tilde{\Phi}|_{\bS\times S^1})^* TL \oplus
 (\tilde{\Phi}|_{\bS\times S^1})^* V \oplus \R^k=\R^{N+k+2}
\]
\[
\begin{CD}
 (\tilde{\Phi}|_{\bS\times S^1})^* V @>>> \tilde{\Phi}^*V\\
 @VVV @VVV\\
 \coprod_{i=1}^{h} R_i\times S^1=\bS\times S^1
 @>i>> \Si\times S^1
\end{CD}
\]

We may view $V$ as a complex line bundle. We need to show that
$(\tilde{\Phi}|_{\bS\times S^1})^* V$ is trivial, or equivalently,
$n_i=(\deg\tilde{\Phi}|_{R_i\times S^1})^* V=0$ for
$i=1,\ldots,h$. We have
\begin{eqnarray*}
\sum_{i=1}^h n_i&=&\deg (\tilde{\Phi}|_{\bS\times S^1})^* V \\
 &=& (i_*[\bS\times S^1])\cap c_1(\tilde{\Phi}^*V)\\
 &=&\pa (i_*[\Si\times S^1])\cap c_1(\tilde{\Phi}^*V)\\
 &=&0.  
\end{eqnarray*}
So $(\tilde{\Phi}|_{\bS\times S^1})^* V$ is trivial if $h=1$. 
For $h>1$, we assume
that $L$ is spin, so that we may take $V=\R^2$, the trivial bundle. 

We finally consider nonstable cases:
\begin{enumerate}
\item $(g,h)=(0,1)$, $n=1$, $\vm=(0)$.
\item $(g,h)=(0,1)$, $n=0$, $\vm=(1)$.  
\item $(g,h)=(0,1)$, $n=0$, $\vm=(2)$.
\item $(g,h)=(0,2)$, $n=0$, $\vm=(0,0)$.
\end{enumerate}

Cases (2) and (3) are treated in \cite{FO3}. For (1) the domain
is isomorphic to the unit disc with one interior marked point at the
origin, and the automorphism group $U(1)=\{e^{i\theta}\mid \theta\in\R\}$ of 
the domain can be oriented by $\frac{\pa}{\pa\theta}$. For Case 4,
the domain is isomorphic to an annulus
$\{z\in\C\mid 1\leq |z|\leq r\}$ for some $r\in (1,\infty)$, which
is oriented by $\frac{\pa}{\pa r}$, and automorphism group
$U(1)$ of the domain is oriented as above. 
$\Box$

\section{Virtual fundamental chain} \label{constructchain}

\subsection{Construction of virtual fundamental chain}

We consider the general setting in Section~\ref{defineku}.

\begin{df}
 Let $M$ be a Hausdorff space with a
 Kuranishi structure (with corners) 
\[
 \mathcal{K}=\left\{ (V_p,E_p,\Gamma_p,\psi_p,s_p):p\in M,
                    (V_{pq},\hat{\phi}_{pq},\phi_{pq},h_{pq}):q\in \psi_p(s_p^{-1}(0))
            \right\}. 
\]
 A Hausdorff topological space $W$ is an {\em ambient space} of $\mathcal{K}$ if 
 \begin{enumerate}
  \item $M$ is a subspace of $W$. 
  \item $\psi_p:V_p\ra W$, $\psi_p(x)\in M$ if and only if $s_p(x)=0$.
  \item $\psi_q=\psi_p\circ\phi_{pq}$.
  \item There is a subset $\pa W\subset W$ such that $\psi_p(x)\in\pa W$
        if and only if $x\in\pa V_p$, where $\pa V_p$ is the union of 
        corners of $V_p$. We take $\pa W=\emptyset$ if 
        $\mathcal{K}$ is a Kuranishi structure. We define $\pa M= M\cap\pa W$. 
  \item $V_p/\Gamma_p\ra W$ is injective.
   
\end{enumerate}
\end{df}

\begin{rem}
If $\mathcal{K}_1$ and $\mathcal{K}_2$ are equivalent Kuranishi structures
(with corners) in the sense of Definition~\ref{sameku}, and they have the
same ambient space $W$, we will implicitly assume that $W$ is also an ambient
space for the $\mathcal{K}$ in Definition~\ref{sameku}, and 
$\psi_{i,p}=\psi_p\circ\phi_i:V_{i,p}\ra W$ in Definition~\ref{same}. 
\end{rem}

\begin{ex}
$\cW=\onep$ is an ambient space of the Kuranishi structure with corners
on $\cM=\MXL$ constructed in Section~\ref{constructku}.
$\pa\cW\subset \cW$ is the subset corresponding to maps whose domain
has at least one boundary node. 
\end{ex}

Now assume that $M$ is a \emph{compact}, Hausdorff topological space
with an \emph{oriented} Kuranishi structure with corners $\cK$,
and that $\cK$ has an ambient space $W$.
Let $d$ be the virtual dimension of $\cK$. 
The {\em virtual fundamental chain}
we will construct is a singular $d$-chain 
$M_{\cK,\nu}\in \cS_d(W,\Q)$ such that 
$\pa M_{\cK,\nu}\in\mathcal{S}_{d-1}(\pa W,\Q)\subset
\cS_{d-1}(W,\Q)$, so it represents a relative
singular $d$-cycle $\bar{M}_{\mathcal{K},\nu}\in \mathcal{S}_d(W,\pa W,\Q)$.

$M$ is compact, so there exist 
finitely many $p_1,\ldots,p_l\in M$ such that
\begin{itemize}
\item $\{ U'_j=\psi_{p_j}(s_{p_j}^{-1}(0)) \mid j=1,\ldots, l \} $
      is an open cover of $M$.
\item There exists a $\Gamma_{p_j}$-invariant neighborhood $V_j$ of 
      $\psi_{p_j}^{-1}(p_j)$ in $V_{p_j}$ such that $V_j\subset\subset
      V_{p_j}$ and
$ \{ U_j=\psi_{p_j}(s_{p_j}^{-1}(0)\cap V_j) \mid j=1,\ldots, l \}$
is still an open cover of $M$.
\end{itemize}

Let $\beta_j: V_{p_j}\ra [0,1]$ be a smooth function with compact
support such that $\beta_j\equiv 1$ on the closure of $V_j$.
For any $\nu=(\nu_1,\ldots,\nu_l)$, where $\nu_i:V_{p_j}\ra E_{p_j}$
is a small continuous section, not necessarily $\Gamma_{p_j}$-equivariant,
we construct a section $\nu_p:\hat{V}_p\ra E_p|_{\hat{V}_p}$ 
for each $p\in X$, where $\hat{V}_p$ is an $\Gamma_p$-invariant 
neighborhood of $\psi_p^{-1}(p)$.

Given $p\in M$, if $p\in U_j$, then there is a $\Gamma_p$-invariant
neighborhood $V_{p_j p}$ of $\psi_p^{-1}(p)$ such that we have the
following commutative diagram:
\[
\begin{CD}
  E_p|_{V_{p_j p}} @>\hat{\phi}_{p_j p}>> E_{p_j}\\
  @VVV @VVV\\
  V_{p_j p} @>\phi_{p_j p}>> V_{p_j}
\end{CD}
\]
For any section $s:V_{p_j}\ra E_{p_j}$, let 
$\phi_{p_j p}^* s: V_{p_j p}\ra E_p|_{V_{p_j p}}$ be the 
unique section satisfying 
$\hat{\phi}_{p_j p}\circ \phi_{p_j p}^* s=s\circ\phi_{p_j p}$.
Let $\hat{V}_p=\cap_{p\in U'_j}V_{p_j p}$. We define 
$\nu_p=\sum_{p\in U'_j}\phi_{p_j p}^* (\beta_j \nu_j):\hat{V}_p\ra
E_p|_{\hat{V}_p}$. 

There exist $\hat{p}_1,\ldots,\hat{p}_{\hat{l}}\in M$ such that
$\{ \hat{U}_j=\psi_{\hat{p}_j}(s_{\hat{p}_j}^{-1}(0)\cap\hat{V}_{\hat{p}_j})
    \mid j=1,\ldots,\hat{l} \}$ is an open cover of $M$.   
We may choose $\nu=(\nu_1,\ldots,\nu_l)$ such that 
$s_{\hat{p}_j}+\nu_{\hat{p}_j}:
\hat{V}_{\hat{p}_j}\ra E_{\hat{p}_j}|_{\hat{V}_{\hat{p}_j}}$ is smooth and
intersects the zero section transversally for $j=1,\ldots,\hat{l}$. So 
$\hat{M}^\nu_j=(s_{\hat{p}_j}+\nu_{\hat{p}_j})^{-1}(0)$
is a $d$-dimensional submanifold of $\hat{V}_{\hat{p}_j}$, where $d$ is the
virtual dimension of the Kuranishi structure. The orientation
of the Kuranishi structure induces an orientation on $\hat{M}^\nu_j$.
If $\hat{V}_{\hat{p}_j}$ has corners, we may further require that
$\hat{M}^\nu_j$ intersect all the corners of $\hat{V}_{\hat{p}_j}$ 
transversally, i.e., $\hat{M}^\nu_j$ is a \emph{neat} submanifold of 
$\hat{V}_{\hat{p}_j}$ in the sense of \cite[Chapter 1, Section 4]{Hi}.
We call such $\nu=\{\nu_p:\hat{V}_p\ra E_p|_{\hat{V}_p}\mid p\in M\}$
a \emph{generic}  perturbation. Note that the difference between 
two generic perturbations is smooth. 

From our choice of $\nu_{\hat{p}_j}$, we have
\begin{equation} \label{compatible}
   \psi_{\hat{p}_j}(\hat{M}^\nu_j)\cap 
 \psi_{\hat{p}_{j'}}(\hat{V}_{\hat{p}_{j'}})
  =\psi_{\hat{p}_j}(\hat{V}_{p_j}) \cap \psi_{\hat{p}_{j'}}(\hat{M}^\nu_{j'}).
\end{equation}
for any $j,j'\in\{ 1,\ldots,\hat{l} \}$.

Choose a triangulation of $\hat{M}^\nu_j$ so that it 
becomes a simplicial complex, and all its corners
are subcomplexes. By (\ref{compatible}), we may assume that there is 
a compact subcomplex $K_j$ of $\hat{M}^\nu_j$ such that
\begin{itemize}
 \item $\cup_{j=1}^{\hat{l}}\psi_{\hat{p}_j}(K_j)=
        \cup_{j=1}^{\hat{l}}\psi_{\hat{p}_j}(\hat{M}^\nu_j)$.
 \item For any $j,j'\in\{1,\ldots,\hat{l}\}$,
       $\psi_{\hat{p}_j}^{-1}\left(\psi_{\hat{p}_j}(K_j)
        \cap \psi_{\hat{p}_{j'}}(K_{j'})\right)$ is a subcomplex of 
       the $(d-1)$-dimensional simplicial complex $\pa K_j$. 
\end{itemize}

Let $\{ \Delta^d_{j,\alpha} \mid \alpha=1,\ldots, N_j \}$ be the set
of $d$-simplices in the triangulation of $K_j$.
Each $\Delta^d_{j,\alpha}$ has an orientation induced by that of 
$\hat{M}^\nu_j$. The inclusion 
$i_{j,\alpha}:\Delta^d_{j,\alpha} \ra K_j\subset\hat{V}_{\hat{p}_j}$
can be viewed as a singular $d$-chain in $\hat{V}_{\hat{p}_j}$. We define
\begin{eqnarray*}
K^\nu_j&=&\sum_{\alpha=1}^{N_j} i_{j,\alpha}\in 
         \mathcal{S}_d(\hat{V}_{\hat{p}_j};\Z)\subset\mathcal{S}_d(\hat{V}_{\hat{p}_j};\Q)\\
M_{\mathcal{K},\nu}&=&\sum_{j=1}^{\hat{l}}\frac{1}{|\Gamma_{\hat{p}_j}|}
{\psi_{\hat{p}_j}}_* K^\nu_j\in\mathcal{S}_d(W;\Q)
\end{eqnarray*}
where $|\Gamma_{\hat{p}_j}|$ is the cardinality of $\Gamma_{\hat{p}_j}$.
It follows from our construction that 
$\pa M_{\mathcal{K},\nu}\in\mathcal{S}_{d-1}(\pa W)\subset \mathcal{S}_{d-1}(W)$,
so it represents a singular relative $d$-cycle
$\bar{M}_{\mathcal{K},\nu}\in \mathcal{S}_{d-1}(W,\pa W;\Q)$, which represents
a class $[M_{\mathcal{K},\nu}]^{\mathrm{rel}}\in H_d(W,\pa W;\Q)$. 

Note that up to subdivision, $M_{\mathcal{K},\nu}\in\mathcal{S}_d(W;\Q)$ depends
on $\nu=\{\nu_p:\hat{V}_p\ra E_p|_{\hat{V}_p} \mid p\in M\}$ 
but not on the choice of $\hat{p}_1,\ldots,\hat{p}_{\hat{l}}\in M$. 

\begin{pro} \label{perturb}
The class $[M_{\mathcal{K},\nu}]^{\mathrm{rel}}\in H_d(W,\pa W;\Q)$ 
is independent of the choice of 
$\nu=\{ \nu_p:\hat{V}_p\ra E_p|_{\hat{V}_p}\mid p\in M\}$.
So we may write $[M_{\mathcal{K}}]^\mathrm{rel}$ for this class. 
\end{pro}

\paragraph{Proof.}
We first observe that a Kuranishi structure with corners $\mathcal{K}$
on $M$ with ambient space $W$ gives rise to 
a Kuranishi structure with corners $\mathcal{K}\times [0,1]$ on 
$M\times[0,1]$ with ambient space $W\times[0,1]$. 
To see this, consider $(p,t)\in M\times [0,1]$. Let
$(V_p,E_p,\Gamma_p,\psi_p,s_p)$ be the Kuranishi neighborhood of $p$
assigned by the Kuranishi structure on $M$. Let
$V_{(p,t)}=V_p\times [0,1]$, and let $\pi_p:V_{(p,t)}\ra V_p$ be the
projection to the first factor. Let 
$E_{(p,t)}=\pi_p^* E_p \ra V_{(p,t)}$, and let 
$s_{(p,t)}=\pi_p^* s_p: V_{(p,t)}\ra E_{(p,t)}$.
We define $\psi_{(p,t)}=\psi_p\times id:
V_{(p,t)}=V_p\times[0,1]\ra W\times [0,1]$, where
$id$ is the identity map on $[0,1]$. Finally, let  $\Gamma_p$ act
on $[0,1]$ trivially. Then
$(V_{(p,t)},E_{(p,t)},\Gamma_p,\psi_{(p,t)}, s_{(p,t)})$ is a Kuranishi
neighborhood of $(p,t)$. The transition functions can be constructed
from those of the Kuranishi structure on $M$ in an obvious way.

Let $\nu=\{\nu_p:\hat{V}_p\ra E_p|_{\hat{V}_p} \mid p\in M\}$, 
$\nu'=\{\nu'_p:\hat{V}_p\ra E_p|_{\hat{V}_p} \mid p\in M\}$ 
be two choices of small generic perturbation of $\mathcal{K}$
in the above construction. There exists a generic perturbation 
\[
\mu=\{\mu_{(p,t)}:V_{(p,t)} \ra E_{(p,t)}\mid 
(p,t)\in M\times[0,1]\}
\]
of $\mathcal{K}\times[0,1]$ such that
\[
i_{p,0}^*\mu_{(p,\half)} =\nu_p :V_p\ra E_p,\; 
i_{p,1}^*\mu'_{(p,\half)}=\nu'_p:V_p\ra E_p,
\]
where $i_{p,t}:V_p\ra V_{(p,\half)}=V_p\times [0,1]$ 
is the inclusion $x\mapsto (x,t)$.

From the paragraph right before Proposition~\ref{perturb}, we may use
$\hat{p}_1,\ldots,\hat{p}_{\hat{l}}$ in the construction of
both $M^\nu$ and $M^{\nu'}$. Let
\[
Y_j=(s_{(\hat{p}_j,\half)}+\mu_{(\hat{p_j},\half)})^{-1}(0)\subset 
\hat{V}_{\hat{p}_j}\times[0,1].
\]
Then 
\[
Y_j\cap\pa(\hat{V}_{\hat{p}_j}\times[0,1])
=i_{\hat{p}_j,1}(\hat{M}^{\nu'}_j)\cup i_{\hat{p}_j,0}(\hat{M}^{\nu}_j)\cup 
 \left(Y_j\cap (\pa\hat{V}_{\hat{p}_j}\times[0,1])\right),
\]
where 
\[
\pa(\hat{V}_{\hat{p}_j}\times[0,1])
=i_{\hat{p}_j,1}(\hat{V}_{\hat{p}_j})\cup i_{\hat{p}_j,0}(\hat{V}_{\hat{p}_j})\cup
 (\pa\hat{V}_{\hat{p}_j}\times[0,1])
\]
is the union of corners of $\hat{V}_{\hat{p}_j}\times[0,1]$.
This gives rise to a singular $(d+1)$-chain $B_j$ in
$\hat{V}_{\hat{p}_j}\times[0,1]$ such that
\[ 
\pa B_j= (i_{\hat{p}_j,1})_* K^{\nu'}_j-(i_{\hat{p}_j,0})_* K^\nu_j +C_j+ D_j,
\] 
where $(i_{\hat{p}_j,1})_* K^{\nu'}_j$ comes from 
$i_{\hat{p}_j,1}(\hat{M}^{\nu'}_j)$,
$(i_{\hat{p}_j,0})_* K^\nu_j$ comes from 
$i_{\hat{p}_j,0}(\hat{M}^\nu_j)$, 
$D_j$ comes from $Y_j\cap(\pa\hat{V}_{\hat{p}_j}\times[0,1])$, and $C_j$
will get canceled:
\[ 
  \sum_{j=1}^{\hat{l}}(\psi_{\hat{p}_j}\times id)_* C_j =0\in
  \mathcal{S}_d(W\times[0,1]; \Q). 
\]
Let $\pi_j:\hat{V}_{\hat{p}_j}\times[0,1]\ra \hat{V}_{\hat{p}_j}$ be the
projection to the first factor. We have
\[
\pa({\pi_j}_* B_j)=K^{\nu'}_j-K^\nu_j+(\pi_j)_* C_j +(\pi_j)_* D_j.
\]

Let
\[
B= \sum_{j=1}^{\hat{l}} \frac{1}{|\Gamma_{\hat{p}_j}|}
     {\psi_{\hat{p}_j}}_* {\pi_j}_* B_j\in\mathcal{S}_{d+1}(W;\Q),
\]
and
\[
D= \sum_{j=1}^{\hat{l}} \frac{1}{|\Gamma_{\hat{p}_j}|}
   {\psi_{\hat{p}_j}}_* {\pi_j}_* D_j\in
   \mathcal{S}_{d}(\pa W;\Q)\subset \mathcal{S}_d(W;\Q).
\]
Then 
\[
\pa B= M_{\mathcal{K},\nu'}-M_{\mathcal{K},\nu} + D \in \mathcal{S}_d(W;\Q) 
\]
since 
\[ 
\sum_{j=1}^{\hat{l}} \frac{1}{|\Gamma_{\hat{p}_j}|}
{\psi_{\hat{p}_j}}_*{\pi_j}_* C_j
=\pi_* \left(\sum_{j=1}^{\hat{l}}\frac{1}{|\Gamma_{\hat{p}_j}|}
 (\psi_{\hat{p}_j}\times id)_* C_j\right) =0,
\]
where $\pi:W\times[0,1]\ra W$ is the projection to the first factor.
We have $\pa \bar{B}= \bar{M}_{\mathcal{K},\nu'}-
\bar{M}_{\mathcal{K},\nu}\in\mathcal{S}_d(W,\pa W;\Q)$, so
\[
  [M_{\mathcal{K},\nu'}]^{\mathrm{rel}}=[M_{\mathcal{K},\nu}]^{\mathrm{rel}}
  \in H_d(W,\pa W;\Q).\;\Box 
\]

\begin{pro}\label{invariance}
Let $\mathcal{K}_1$ and $\mathcal{K}_2$ be two equivalent Kuranishi
structures with corners on a compact Hausdorff topological space $M$. Let
$W$ be an ambient space of both $\mathcal{K}_1$ and $\mathcal{K}_2$.
Then $[M_{\mathcal{K}_1}]^{\mathrm{rel}}=[M_{\mathcal{K}_2}]^{\mathrm{rel}}$.
\end{pro}

\paragraph{Proof.} Let
\begin{eqnarray*}
 \mathcal{K}_1&=&\left\{ 
  (V_{1,p},E_{1,p},\Gamma_{1,p},\psi_{1,p},s_{1,p}):p\in M,
  (V_{1,pq}, \hat{\phi}_{1,pq},\phi_{1,pq},h_{1,pq}):
   q\in \psi_{1,p}(s_{1,p}^{-1}(0)) \right\}\\
 \mathcal{K}_2&=&\left\{ 
  (V_{2,p},E_{2,p},\Gamma_{2,p},\psi_{2,p},s_{2,p}):p\in M,
  (V_{2,pq}, \hat{\phi}_{2,pq},\phi_{2,pq},h_{2,pq}):
  q\in \psi_{2,p}(s_{2,p}^{-1}(0)) \right\}\\
 \mathcal{K}&=&\left\{
  (V_p,E_p,\Gamma_p,\psi_p,s_p): p\in M,
  (V_{pq},\hat{\phi}_{pq},\phi_{pq},h_{pq}): 
  q\in \psi_p(s_p^{-1}(0)) \right \}
\end{eqnarray*}
be as in Definition~\ref{sameku}.
 Let $\nu_i=\{\nu_{i,p}:V_{i,p}\ra E_{i,p}\mid p\in M\}$ be a generic perturbation
which can be used to define a virtual fundamental chain
$M_{\mathcal{K}_i,\nu_i}\in \mathcal{S}_d(W;\Q)$, for $i=1,2$. 
$\nu_i$ can be extended to a generic perturbation 
$\mu_i=\{\mu_{i,p}:V_p\ra E_p\mid p\in M\}$ such that 
$\hat{\phi}_i\circ \nu_{i,p}=\mu_{i,p}\circ\phi_i$. 
We have 
$M_{\mathcal{K}_i,\nu_i}=M_{\mathcal{K},\mu_i}\in \mathcal{S}_d(W;\Q)$. So
\[
 [M_{\mathcal{K}_1}]^\mathrm{rel}
=[M_{\mathcal{K},\mu_1}]^\mathrm{rel}
=[M_{\mathcal{K},\mu_2}]^\mathrm{rel}
=[M_{\mathcal{K}_2}]^\mathrm{rel}\in H_d(W,\pa W;\Q). \;\Box
\]

\begin{rem}\label{nocorner}
Let $M$ be a compact Hausdorff topological space
with an oriented Kuranishi structure $\mathcal{K}$,
and let $W$ be an ambient space of $\mathcal{K}$.
Then the above construction yields 
$M_{\mathcal{K},\nu}\in\mathcal{S}_d(W,\Q)$ such that 
$\pa M_{\mathcal{K},\nu}=0 \in \mathcal{S}_{d-1}(W,\Q)$, so it represents a class
$[M_{\mathcal{K},\nu}]\in H_d(W;\Q)$. The proof of Proposition~\ref{perturb}
shows that this class is independent of the choice of $\nu$, so we may write
$[M_{\mathcal{K}}]$ for this class. The proof of Proposition~\ref{invariance}
shows that if $\mathcal{K}'$ is another Kuranishi structure on $M$ such that
$\mathcal{K}'\sim\mathcal{K}$ and $W$ is also an ambient space of
$\mathcal{K}'$, then 
$[M_{\mathcal{K}}]=[M_{\mathcal{K}'}]\in H_d(W;\Q)$.
Let $i:M\ra W$ be the inclusion. Then 
$[M_{\mathcal{K}}]=i_*[M]\in H_d(W;\Q)$, where $i_*[M]$ is 
defined in \cite[Section 6]{FO}.
\end{rem}

\begin{ex}
Let $X$ be a Calabi-Yau $3$-fold, and let $\beta\in H_2(X;\Z)$.
Let $\overline{M}_{g,0}(X,\beta)$ denote the moduli space of stable maps
$f$ from a genus $g$ prestable curve $C$ to $X$ such that $f_*[C]=\beta$.
Then $\overline{M}_{g,0}(X,\beta)$ has an oriented  Kuranishi structure
\cite{FO}.  The virtual dimension of is $0$ for any 
$g\geq 0$ and $\beta\in H_2(X;\Z)$.
This Kuranishi structure has an ambient space $W^{1,p}_{g,0}(X,\beta)$, the
moduli space of stable $W^{1,p}$ maps from a genus $g$ prestable curve
$C$ to $X$ such that  $f_*[C]=\beta$. Then 
$[\overline{M}_{g,0}(X,\beta)_{\mathcal{K}}]\in H_0(W^{1,p}_{g,0}(X,\beta);\Q)$,
and $\deg [\overline{M}_{g,0}(X,\beta)_{\mathcal{K}}]\in\Q$ is some
Gromov-Witten invariant \cite[Section 17]{FO}.  
\end{ex}

\begin{ex} \label{trivial}
Let $X$ be a Calabi-Yau $3$-fold, $L$ be a special Lagrangian submanifold.
$\MXL$ is empty for $\mu\neq 0$.
The tangent bundle of $L$ is trivial, so $L$ is spin.
We have constructed an orientable Kuranishi structure with corners $\mathcal{K}$ on
$\cM=\Bar{M}_{(g,h),(0,\vec{0})}(X,L\mid \beta,\vga,0)$.
The equivalence class of $\cK$ is independent of choices in our construction.
The virtual dimension of $\cK$ is  $0$ for all
$g,h,\beta,\vga$.
The Kuranishi structure has an ambient space 
$\cW=W^{1,p}_{(g,h),(0,\vec{0})}(X,L\mid \beta,\vga,0)$.
$[\cM_{\cK}]^{\mathrm{rel}}\in H_0(\cW,\pa\cW;\Q)=0$,
so we cannot get a nontrivial number from $[\cM_{\cK}]^{\rel}$. 
\end{ex}

Most enumerative predictions about holomorphic curves with Lagrangian
boundary conditions concern the special case in Example~\ref{trivial}.
Our goal is to give a rigorous mathematical definition of these
highly nontrivial enumerative numbers. From Example~\ref{trivial}
we know that the relative cycle $[M_{\cK}]^{\rel}\in H_d(W,\pa W;\Q)$
is not the right object for our purpose. We want to remember the virtual 
fundamental \emph{chain} $M_{\cK,\nu} \in \cS_d(W,\Q)$ 
which contains more information. For example, when the virtual 
dimension $d=0$, the $0$-chain $M_{\cK,\nu}$ represents a class 
$[M_{\cK,\nu}]\in H_0(W;\Q)$ 
since any $0$-chain is a cycle. Let $\nu,\nu'$ be two perturbations.
From the proof of Proposition~\ref{invariance},
\[
\pa B=M_{\cK,\nu}-M_{\cK,\nu'}+ D,  
\]
where $B\in \cS_1(W,\Q)$, $D\in\cS_0(\pa W,\Q)$.
So $[M_{\cK,\nu}]$ \emph{depends on} $\nu$. 
We want to impose extra constraint on $\nu$ such that
if $\nu$ and $\nu'$ both satisfy the constraint 
then the above $D$ is zero. 
In particular, if $\nu$ and $\nu'$ satisfy
the same boundary conditions in the sense that 
$\nu_{\hat{p}_j}=\nu'_{\hat{p}_j}$ on $\pa \hat{V}_{\hat{p}_j}$ for 
all $j=1,\ldots,\hat{l}$, then $D=0$, so  
$[M_{\cK,\nu}]=[M_{\cK,\nu'}]\in H_0(W;\Q)$, and
$\deg[M_{\cK,\nu}]=\deg[M_{\cK,\nu'}]\in\Q$.

\subsection{$S^1$ action}

\begin{df} \label{circleku}
Let $M$ be a Hausdorff topological space, and let 
$\hat{\varrho}:S^1\times M\ra M$ be a continuous $S^1$ action.
A Kuranishi structure with corners 
\[
 \cK=\left\{ (V_p,E_p,\Gamma_p,\psi_p,s_p):p\in M,
  (V_{pq},\hat{\phi}_{pq},\phi_{pq},h_{pq}):q\in \psi_p(s_p^{-1}(0))\right\}. 
\]
on M is {\em $\hat{\varrho}$-equivariant on the boundary}
if the Kuranishi neighborhood $(V_p,E_p,\Gamma_p,\psi_p,s_p)$ of any 
$p\in\pa M$ is $\hat{\varrho}$-equivariant in the sense that
 \begin{enumerate}
\item There is a free continuous $S^1$ action on $V_p$ which commutes with
      the action of $\Gamma_p$ and leaves $\pa V_p$ invariant.
\item $E_p\ra V_p$ is an $S^1$-equivariant vector bundle.
\item $s_p: V_p\ra E_p$ is an
     $S^1$-equivariant section.
\item $\psi_p:s_p^{-1}(0)\ra M$ is $S^1$-equivariant, where
      $S^1$ acts on $M$ by $\hat{\varrho}$.
\end{enumerate}
\end{df}

\begin{rem}
Let $\mathcal{K}$ be as above, and let  $\pa \mathcal{K}$ be the Kuranishi
structure with corners on $\pa M$ defined in Remark~\ref{boundaryku}.
Then the quotients of the Kuranishi neighborhoods and transition functions 
of $\pa \mathcal{K}$ by the free $S^1$ action define
a Kuranishi structure $\pa \mathcal{K}/S^1$ of virtual dimension
$d-2$ on $\pa M/S^1$.
\end{rem}

\begin{rem}
Let $M$ be a Hausdorff topological space with a continuous $S^1$-action
$\hat{\varrho}:S^1\times M\ra M$. Then there exist an 
$\hat{\varrho}$-equivariant Kuranishi structure with corners only if
the action leaves $\pa M$ invariant and has no fixed point on $\pa M$.
\end{rem}

\begin{df} \label{samecircleku}
Let $M$ be a Hausdorff topological space with a continuous $S^1$-action
$\hat{\varrho}:S^1\times M\ra M$.
Let 
\[
  \mathcal{K}_1=\left\{ 
  (V_{1,p},E_{1,p},\Gamma_{1,p},\psi_{1,p},s_{1,p}):p\in M,
  (V_{1,pq}, \hat{\phi}_{1,pq},\phi_{1,pq},h_{1,pq}):
   q\in \psi_{1,p}(s_{1,p}^{-1}(0)) \right\}
\]
and 
\[
  \mathcal{K}_2=\left\{ 
  (V_{2,p},E_{2,p},\Gamma_{2,p},\psi_{2,p},s_{2,p}):p\in M,
  (V_{2,pq}, \hat{\phi}_{2,pq},\phi_{2,pq},h_{2,pq}):
  q\in \psi_{2,p}(s_{2,p}^{-1}(0)) \right\}
\]
be two Kuranishi structures with corners on $M$ which 
are $\hat{\varrho}$-equivariant on the boundary. $\mathcal{K}_1$ and
$\mathcal{K}_2$ are are {\em $\hat{\varrho}$-equivalent} if
\begin{itemize}
\item  There is another Kuranishi structure
\[
  \mathcal{K}=\left\{
  (V_p,E_p,\Gamma_p,\psi_p,s_p): p\in M,
  (V_{pq},\hat{\phi}_{pq},\phi_{pq},h_{pq}): 
  q\in \psi_p(s_p^{-1}(0)) \right \}
\]
on $M$ which is $\varrho$-equivariant on the boundary such that for all $p\in M$, 
$(V_{1,p},E_{1,p},\Gamma_{1,p},\psi_{1,p},s_{1,p})$,
$(V_{2,p},E_{2,p},\Gamma_{2,p},\psi_{2,p},s_{2,p})$, and
$(V_p,E_p,\Gamma_p,\psi_p,s_p)$ satisfy the relation described in
Definition~\ref{same}.
\item  The $\phi_i$ and $\hat{\phi}_i$ in Definition~\ref{same} are
$S^1$-equivariant.
\end{itemize}
In this case, we write $\mathcal{K}_1\stackrel{\hat{\varrho}}{\sim}\mathcal{K}_2$.
\end{df} 

Note that $\mathcal{K}_1\stackrel{\hat{\varrho}}{\sim}\mathcal{K}_2
\Rightarrow \mathcal{K}_1\sim\mathcal{K}_2$.

\begin{df}
Let $M$ be a Hausdorff topological space with a continuous 
$S^1$ action $\hat{\varrho}:S^1\times M\ra M$. Let
\[
\mathcal{K}=\left\{
  (V_p,E_p,\Gamma_p,\psi_p,s_p): p\in M,
  (V_{pq},\hat{\phi}_{pq},\phi_{pq},h_{pq}): 
  q\in \psi_p(s_p^{-1}(0)) \right \}
\]
be a Kuranishi structure with corners on $M$ which
is $\hat{\varrho}$-equivariant on the boundary. 
An ambient space $W$ of $\mathcal{K}$ is {\em $\hat{\varrho}$-equivariant}
if $\hat{\rho}$ extends to an action $S^1\times W\ra W$, and
$\psi_p:V_p\ra W$ is $S^1$-equivariant. 
\end{df}

We now assume that $M$ is a \emph{compact}, Hausdorff topological space
with a continuous $S^1$ action $\hat{\varrho}:S^1\times M\ra M$. Suppose that
$\mathcal{K}$ is an \emph{oriented} Kuranishi structure with corners of virtual
dimension \emph{zero} on $M$. We further assume that $\mathcal{K}$ is 
$\hat{\varrho}$-equivariant on the boundary, and $W$ is  
an $\hat{\varrho}$-equivariant ambient space of $\mathcal{K}$.

In this case, a generic perturbation 
$\nu=\{ \nu_p:V_p\ra E_p\mid p\in M \}$ of $\mathcal{K}$ 
is a perturbation such that $s_p+\nu_p$ 
intersects the zero section transversally at isolated points and 
is nowhere zero on $\pa V_p$. Let $\bar{p}$ denote the equivalent class of 
$p\in \pa M$ in $\pa M/S^1$.
The virtual dimension of 
\[
\pa\mathcal{K}/S^1=\left\{ (
(\bar{V}_{\bar{p}},\bar{E}_{\bar{p}},\Gamma_p,\bar{\psi}_{\bar{p}},\bar{s}_{\bar{p}}:
\bar{p}\in \pa M/S^1,\; (\bar{V}_{\bar{p}\bar{q}},\bar{\phi}_{\bar{p}\bar{q}},
\hat{\bar{\phi}}_{\bar{p}\bar{q}},\bar{h}_{\bar{p}\bar{q}}):\bar{q}\in
\bar{\psi}_{\bar{p}}(\bar{s}_{\bar{p}}^{-1}(0))\right\}
\]
is $-2$, so a generic perturbation $\bar{\nu}=\{ \bar{\nu}_{\bar{\rho}}:
\bar{V}_{\bar{p}}\ra\bar{E}_{\bar{p}}\mid p\in \pa M/S^1 \}$ of $\pa\mathcal{K}/S^1$
is a perturbation such that 
$\bar{s}_{\bar{p}}+\bar{\nu}_{\bar{p}}$ is nowhere zero for all 
$\bar{p}\in\pa M/S^1$. Let 
$Q_p:\pa V_p\ra \bar{V}_{\bar{p}}=\pa V_p/S^1$
be the natural projection. Then
$Q^*\bar{\nu}=\{ (Q^*\bar{\nu})_p=Q^*\bar{\nu}_{\bar{p}}:\pa V_p
\ra E|_{\pa V_p}\mid p\in\pa M\}$ is a generic perturbation of
$\pa \mathcal{K}$. We call such a perturbation 
an $\hat{\varrho}$-equivariant perturbation. A generic perturbation 
of $\pa \mathcal{K}$ can always be extended to a generic perturbation
of $\mathcal{K}$, so there exist a generic perturbation of $\mathcal{K}$
whose restriction to $\pa \mathcal{K}$ is $\hat{\varrho}$-equivariant. 

\begin{pro}\label{rhoperturb}
Let $M$ be a compact Hausdorff topological space with a continuous
$S^1$ action $\hat{\varrho}:S^1\times M\ra M$. Let $\mathcal{K}$ be
a Kuranishi structure with corners of virtual dimension $0$ on $M$ 
which is $\hat{\varrho}$-equivariant on the boundary.
Let $W$ be a $\hat{\varrho}$-equivariant ambient space of $\mathcal{K}$.
Let $\nu=\{\nu_p:V_p\ra E_p\mid p\in M\}$ be a generic perturbation
such that $\nu|_{\pa \mathcal{K}}=\{\nu_p|_{\pa V_p}:
\pa V_p\ra E_p|_{\pa V_p}\mid p\in \pa M\}$ is $\hat{\varrho}$-equivariant.
Then $[M_{\mathcal{K},\nu}]\in H_0(W;\Q)$ does not depend on the
perturbation $\nu$, so we may write $[M_{\mathcal{K}}^{\hat{\varrho}}]$ for 
this class. If $\mathcal{K}'\stackrel{\hat{\varrho}}{\sim}\mathcal{K}$, 
and $W$ is also a $\hat{\varrho}$-equivariant ambient space of $\mathcal{K}'$, then
$[M_{\mathcal{K}'}^{\hat{\varrho}}]=[M_{\mathcal{K}}^{\hat{\varrho}}]\in
 H_0(W;\Q)$.
\end{pro}
\paragraph{Proof.}Let
\[
\nu=\{\nu_p:V_p\ra E_p\mid p\in M\},\; \nu'=\{\nu'_p:V_p\ra E_p\mid p\in M\}
\]
be two generic perturbations of $\mathcal{K}$ such that
$\nu|_{\pa\mathcal{K}}=Q^*\bar{\nu}$,
$\nu'|_{\pa\mathcal{K}}=Q^*\bar{\nu}'$,
where 
\[
\bar{\nu}=\{\bar{\nu}_{\bar{p}}:\bar{V}_{\bar{p}}\ra \bar{E}_{\bar{p}}\mid
\bar{p}\in \pa M/S^1\},\;
\bar{\nu}'=\{\bar{\nu}'_{\bar{p}}:\bar{V}_{\bar{p}}\ra \bar{E}_{\bar{p}}\mid 
\bar{p}\in \pa M/S^1\}
\]
are perturbations for $\pa\mathcal{K}/S^1$ such that 
$\bar{s}_{\bar{p}}+\bar{\nu}_{\bar{p}}$,
$\bar{s}_{\bar{p}}+\bar{\nu}'_{\bar{p}}$
are nowhere zero for all $\bar{p}\in\pa M/S^1$. 
There exists a generic perturbation 
\[
\bar{\mu}=\{\bar{\mu}_{(\bar{p},t)}:V_{(\bar{p},t)}\ra E_{(\bar{p},t)}\mid
 (\bar{p},t)\in \pa M/S^1\times[0,1] \}
\]
such that 
\[
i_{\bar{p},0}^*\mu_{(\bar{p},\half)}=\bar{\nu}_p:
\bar{V}_{\bar{p}}\ra \bar{E}_{\bar{p}},\;
i_{\bar{p},1}^*\mu_{(\bar{p},\half)}=\bar{\nu}'_p:
\bar{V}_{\bar{p}}\ra \bar{E}_{\bar{p}},
\]
where $i_{\bar{p},t}:V_{\bar{p}}\ra
V_{(\bar{p},\half)}=V_{\bar{p}}\times[0,1]$ 
is the inclusion $x\mapsto (x,t)$.
The virtual dimension of $\pa\mathcal{K}/S^1\times[0,1]$ is $-1$, so 
$\bar{\mu}$ being generic simply means 
$\bar{s}_{(\bar{p},t)}+\bar{\mu}_{(\bar{p},t)}$ is nonzero for all
$(\bar{p},t)\in\pa M/S^1\times [0,1]$.
The pull back $Q^*\bar{\mu}$ of $\bar{\mu}$ is a generic perturbation of 
$\pa\mathcal{K}\times[0,1]$ such that
$s_{(p,t)}+(Q^*\bar{\mu})_{(p,t)}$ is nowhere zero for any
$(p,t)\in\pa M\times[0,1]$. 
We may extend $Q^*\bar{\mu}$ to a generic perturbation $\mu$ of
$\mathcal{K}\times [0,1]$ such that
\[
  i^*_{p,0}\mu_{(p,\half)}=\nu_p :V_p\ra E_p,\;
  i^*_{p,1}\mu_{(p,\half)}=\nu'_p:V_p\ra E_p
\]
as in the proof of Proposition~\ref{perturb}. We proceed as
the proof of Proposition~\ref{perturb}, and use the notation there.
$Y_j$ does not intersect $\pa\hat{V}_{\hat{p}_j}\times [0,1]$, so $D=0$, and
the $1$-chain $B\in\mathcal{S}_1(W;\Q)$ satisfies
\[
\pa B=M_{\mathcal{K},\nu'}-M_{\mathcal{K},\nu}\in \mathcal{S}_0(W;\Q).
\]
So we have $[M_{\mathcal{K},\nu'}]=[M_{\mathcal{K},\nu}]\in H_0(W;\Q)$. 

The last statement can be proved as Proposition~\ref{invariance}.
$\Box$

\subsection{Invariants for an $S^1$-equivariant pair} \label{circle}

Let $(X,\omega)$ be a compact symplectic manifold together with an
$\omega$-tame almost complex structure $J$, and $L$ be a Lagrangian submanifold.
We assume that $L$ is spin or that $h=1$ and $L$ is relative spin
so that $\cM=\MXL$ has an orientable Kuranishi structure.
We fix an orientation by choosing a stable trivialization of $TL$
or $TL\oplus V$ on the $2$-skeleton $L^{(2)}$ of $L$, where $V$ 
is chosen as in the proof of Theorem~\ref{spin}.
$\cW=\onep$, is an ambient space of the Kuranishi structure on $\cM$.

\begin{df}
An \emph{admissible} $S^1$ action on $(X,L)$ is an $S^1$ action
$\varrho:S^1\times X\ra X$ such that
\begin{itemize}
\item $\varrho$ preserves $J$ and $L$.
\item The restriction of $\varrho$ to $L$ is \emph{free}.
\end{itemize}
\end{df}

We now assume that there is an admissible $S^1$ action on $(X,L)$.
Given $t=e^{i\theta}\in S^1$, let $f_t: X\ra X$ be the $J$-holomorphic
diffeomorphism given by $x \mapsto t\cdot x$.
We have an $S^1$ action $\hat{\varrho}:S^1\times\cW\ra\cW$ given by
$(t,[(\lambda,u)])\mapsto t\cdot[(\lambda,u)]=[(\lambda,f_t\circ u)]$.
The action preserves $\cM\subset \cW$
since $f_t$ is $J$-holomorphic for all $t\in S^1$. 

If $\map$ represents a fixed point of the $S^1$ action, then 
$\Si=C\cup D_1\cup\cdots D_h$, where $C$ is a genus $g$ prestable curve,
and $D_k$ is a disc which intersects $C$ at an interior node for $k=1,\ldots,h$.
$u(\pa D_k)$ is an orbit of the $S^1$ action on $L$.
Let $\cW^{S^1}$ and $\cM^{S^1}$ denote the fixed loci of the $S^1$ action on
$\cW$ and $\cM$, respectively. Then $\cM^{S^1}=\cW^{S^1}\cap \cM$,
$\cW^{S^1}\cap \pa\cW=\emptyset$, and $\cM^{S^1}\cap \pa\cM=\emptyset$. 

\begin{tm}\label{existku}
Let $(X,L)$ be as above. Then there exists an oriented Kuranishi structure
on $\cM=\MXL$ which is $\varrho$-equivariant on the boundary. $\cW=\onep$
is a $\hat{\varrho}$-equivariant ambient space of $\cK$.
\end{tm}

\paragraph{Proof.} Let 
\[
\cK=\left\{ (V'_\rho,E_\rho,\mapaut, \psi_\rho, s_\rho):\rho\in\cM,
(V_{\rho\rho'},\hat{\phi}_{\rho\rho'},\phi_{\rho\rho'},h_{\rho\rho'}),
\rho'\in \psi_\rho(s^{-1}(0)) \right\} 
\]
be a Kuranishi structure constructed as in Section~\ref{constructku}.
For $\rho\in\pa\cM$, we will modify the Kuranishi neighborhood of $\rho$ such that
1.-4. in Definition~\ref{circleku} hold. 

Given $\rho=(\lambda,u)\in \pa\cM$, let $H_{\rho,\mathrm{domain}}$, 
$H_{\rho,\mathrm{map}}=\mathrm{Ker}(\pi\circ D_u)$ be defined as in  
Section~\ref{constructku}. Recall that $\phi\in \domaut\mapsto u\circ \phi^{-1}$
induces an inclusion $H_{\rho,\mathrm{aut}}\subset H_{\rho,\mathrm{map}}$.
Similarly, $t\in S^1\mapsto f_t\circ u$ induces an inclusion 
$H_{\rho,\mathrm{circle}}=T_1 S^1\subset H_{\rho,\mathrm{map}}$, where $T_1 S^1\cong\R$ is the 
tangent line of $S^1=\{e^{i\theta}\mid \theta\in\R\}$ at $1\in S^1$.
Note that $H_{\rho,\mathrm{circle}}\subset H_{\rho,\mathrm{aut}}$ if and only
if $\rho\in \cM^{S^1}$, which is excluded since we consider $\rho\in\pa \cM$. 
We may choose $H'_{\rho,\mathrm{map}}$ such that 
$H_{\rho,\mathrm{circle}}\subset H'_{\rho,\mathrm{map}}$ and
$H_{\rho,\mathrm{map}}=H_{\rho,\mathrm{aut}}\oplus H'_{\rho,\mathrm{map}}$.
Choose a subspace $H''_{\rho,\mathrm{map}}$ in $H'_{\rho,\mathrm{map}}$ such
that $H'_{\rho,\mathrm{map}}=H_{\rho,\mathrm{circle}}\oplus
H''_{\rho,\mathrm{map}}$. 

Let $V''_{\rho,\mathrm{map}}$ be a small neighborhood of $0$ in 
$H''_{\rho,\mathrm{map}}$, and $\ep>0$ be small.
Then for sufficiently small 
$(\xi,\eta,\eta')\in B_{\delta_2}\times D_d\times D'_{d'}$ and
$w'\in H'_{\rho,\mathrm{map}}$, there are unique 
$w''\in V''_{\rho,\mathrm{map}}$ and $\theta\in(-\ep,\ep)$ such that
\[
u_{(\xi,\eta,\eta',w')}=f_{e^{i\theta}}\circ u_{(\xi,\eta,\eta',w'')}
\equiv u_{(\xi,\eta,\eta',w'',e^{i\theta})},
\]
where the notation is the same as in Section~\ref{constructku}.
Both $\mapaut$ and $S^1$ act on 
$B_{\delta_2}\times D_d\times D_d'\times V''_{\rho,\mathrm{map}}\times S^1$ such that 
$u_{\phi\cdot(\xi,\eta,\eta',w'',t)}=u_{(\xi,\eta,\eta',w'',t)}\circ\phi^{-1}$ 
and $t'\cdot(\xi,\eta,\eta',w'',t)=(\xi,\eta,\eta',w'',t't)$ for
$\phi\in\mapaut$, $t'\in S^1$, $(\xi,\eta,\eta',w'',t)\in
B_{\delta_2}\times D_d\times D_d'\times V''_{\rho,\mathrm{map}}\times S^1$.
Note that the action  of $\mapaut$ commutes with that of $S^1$.
So we may choose a neighborhood $V''_\rho$ of $0$ in 
$B_{\delta_2}\times D_d\times D_d'\times V''_{\rho,\mathrm{map}}$ such that
$V''_\rho\times S^1$ is invariant under the action of $\mapaut$.
There is an $S^1$ equivariant map
\begin{eqnarray*}
\tilde{\psi}_\rho: V''_\rho\times S^1 &\longrightarrow &\cW\\
(\xi,\eta,\eta',w'',t)&\mapsto& 
[(\lambda_{(\xi,\eta,\eta')},u_{(\xi,\eta,\eta',w'',t)})]
\end{eqnarray*}
$(V''_\rho\times S^1)/\mapaut \ra\cW$ is injective if and only if
$S^1_\rho=\{1\}$, where $S^1_\rho$ is the stabilizer at $\rho$
of the $S^1$ action. In general, $S^1_\rho$ is a finite group
because $\rho\in\pa M$ is not a fixed point of the $S^1$ action.

  The automorphism of $\bar{\rho}\in\pa M/S^1$ is
\[
\mathrm{Aut}\,\bar{\rho}=
\{ \phi\in\domaut\mid u\circ \phi=f_t\circ u \textup{ for some }t\in S^1 \}, 
\]
and
\[
S^1_\rho=\{ t\in S^1\mid  u\circ \phi=f_t\circ u \textup{ for some }\phi\in\domaut\}.
\]
$\mapaut$ is a normal subgroup of $\mathrm{Aut}\,\bar{\rho}$, and
we have an exact sequence
\[
 1\ra \mapaut\ra\mathrm{Aut}\,\bar{\rho}\ra S^1_\rho\ra 1.
\]
The action of $\mapaut$ on $V''_\rho\times S^1$ extends to 
$\mathrm{Aut}\,\bar{\rho}$, and 
$(V''_\rho\times S^1)/\mathrm{Aut}\,\bar{\rho}\ra\cW$ is injective.

Let $q_\rho: V''_\rho\times S^1\ra V''_\rho$ be the projection to the first
factor, and let $E\ra V''_\rho$ be the restriction of the obstruction
bundle over $V'_\rho$. Let $\tilde{V}_\rho=V''_\rho\times S^1$, 
$\tilde{E}_\rho=q_\rho^* E\ra \tilde{V}_\rho$, 
$\tilde{s}_\rho=q_\rho^* (s_\rho|_{V''_\rho})$. Then
$(\tilde{V}_\rho,\tilde{E}_\rho,\mathrm{Aut}\,\bar{\rho},\tilde{\psi}_\rho,\tilde{s}_\rho)$
is a Kuranishi neighborhood of $\rho\in\pa\cM$ in $\cM$ which satisfies
1.-4. in Definition~\ref{circleku}. 

Finally, we proceed as in Section~\ref{transition} to construct new
transition functions. $\Box$

\begin{rm}\label{rhochoices}
The $\hat{\varrho}$-equivalence class of the Kuranishi structure with corners 
which is $\hat{\varrho}$-equivariant on the boundary constructed in the
above proof does not depend on various choices.
\end{rm}

Now consider the case $m=n=0$, and the virtual dimension
\[
d=\mu+(N-3)(2-2g-h)=0.
\]
We define the {\em Euler characteristic}
of $\cM=\Bar{M}_{(g,h),(0,\vec{0})}(X,L|\beta,\vga,\mu)$ to be
\[
\chi_{(g,h)}(X,L,\varrho|\beta,\vga,\mu)=
\deg[\cM_{\cK}^{\hat{\varrho}}]\in\Q,
\]
where $\mathcal{K}$ is a Kuranishi structure with corners which
is $\hat{\varrho}$-equivariant on the boundary, constructed in the proof of 
Theorem~\ref{existku}. This number is well-defined by Proposition~\ref{rhoperturb}.

$\chi_{(g,h)}(X,L,\varrho|\beta,\vga,\mu)$
is an invariant for the equivariant pair $(X,L,\varrho)$, but \emph{not}
an invariant for the pair $(X,L)$. In other words, it is possible that
\[ 
\chi_{(g,h)}(X,L,\varrho_1|\beta,\vga,\mu)\neq
\chi_{(g,h)}(X,L,\varrho_2|\beta,\vga,\mu) 
\]
for two different admissible $S^1$ actions $\varrho_1,\varrho_2$ on $(X,L)$.

\subsection{Multiple covers of the disc}

We consider the special case studied in \cite{OV,KL,LS}.
Let $(z, u, v)$ and $(\tilde{z}, \tilde{u}, \tilde{v})$ be the two charts
of $\OO$, related by $(\tilde{z},\tilde{u},\tilde{v})=(\frac{1}{z},zu,zv)$.
Let $X$ be the total space of $\OO$. There is an antiholomorphic
involution
\begin{eqnarray*} 
A:X&\longrightarrow&X\\
   (z,u,v)&\longmapsto&(\frac{1}{\bar{z}},\bar{z}\bar{v},\bar{z}\bar{u})
\end{eqnarray*} 
in terms of the first chart. The fixed locus $L=X^{A}$ is a special
Lagrangian submanifold of the noncompact Calabi-Yau
3-fold. For any integer $a$, let $\varrho_a:S^1\times X\ra X$ be the $S^1$
action on $X$ defined by $(e^{i\theta},(z,u,v))\mapsto
(e^{i\theta}z, e^{i(a-1)\theta} u, e^{-ia\theta}v)$
on the first chart. Then $\rho_a$ is an admissible $S^1$ action on $(X,L)$.
Let $D^2=\{(z,0,0)\mid |z|\leq 1\}$ be a disc in the first chart, oriented
by the complex structure. Then $\pa D^2\subset L$.
Let $\beta=[D^2]\in H_2(X,L;\Z)\cong H_2(\PP,S^1;\Z)$ and 
$\gamma=[\pa D^2]\subset H_1(L;\Z)\cong H_1(S^1;\Z)$.

By Schwartz reflection principle (see e.g \cite[Section 3.3.2]{KL}),
any nonconstant holomorphic map $f:(\Si,\bS)\ra(X,L)$ can be extended to a
nonconstant holomorphic map $f_\C:\Si_\C\ra X$, whose image must lie in $\PP$.
So we have
\[
\cM=\Bar{M}_{(g,h),(0,\vec{0})}(X,L\mid d \beta,(n_1\gamma,\ldots,n_h\gamma),0)=
\Bar{M}_{(g,h),(0,\vec{0})}(\PP, S^1\mid d\beta,(n_1\gamma,\ldots,n_h\gamma),2d)
\]
as topological spaces. Therefore, 
$\Bar{M}_{(g,h),(0,\vec{0})}(X,L\mid d \beta,(n_1\gamma,\ldots,n_h\gamma),0)$
is compact in $C^\infty$ topology. 
Note that the virtual dimension of the Kuranishi structure corners is $0$ for
$\Bar{M}_{(g,h),(n,\vec{0})}(X,L\mid d \beta, (n_1\gamma,\ldots,n_h\gamma),0)$
and $2(d+2g+h-2)$ for
$\Bar{M}_{(g,h),(n,\vec{0})}(\PP, S^1\mid d \beta,(n_1\gamma,\ldots,n_h\gamma),0)$.
In particular, these two Kuranishi structures on $\cM$  are
not equivalent. The following numbers are defined:
\[
C(g,h|d;n_1,\ldots,n_h|a)=
\chi_{(g,h)}(X,L,\varrho_a\mid d\beta,(n_1\gamma,\ldots,n_h\gamma),0)\in\Q,
\]
where $n_1,\ldots, n_h$ are positive integers, and $d=n_1+\cdots+n_h$.

Let $\pi:\cC\ra \Bar{M}_{g,h}$ be the universal family, $\omega_\pi$ be
the relative dualizing sheaf, and $s_i: \Bar{M}_{g,h}\ra\cC$ be the
section corresponding to the $i$-th marked point.
Let $\E=\pi_*\omega_\pi$ be the Hodge bundle on $\Bar{M}_{(g,h)}$,
and $\psi_i=c_1(s_i^*\omega_\pi)$. 

\begin{con}\label{localize}
Let $a$ be a positive integer. Then
\begin{eqnarray*}
&&(-1)^{d-h}C(0;h|d;n_1,\ldots,n_h|a)=C(0;h|d;n_1,\ldots,n_h|1-a)\\
&=&(a(1-a))^{h-1}
   \prod_{i=1}^h\left(\begin{array}{c}n_i a -1\\n_i-1 \end{array}\right)d^{h-3}.
\end{eqnarray*}

For $g>0$,
\begin{eqnarray*}
&&(-1)^{d-h}C(g;h|d;n_1,\ldots,n_h|a)=C(g;h|d;n_1,\ldots,n_h|1-a)\\
&=&(a(1-a))^{h-1}
 \prod_{i=1}^h\left(\begin{array}{c}n_i a -1\\n_i-1\end{array}\right)\cdot\\
&& \int_{\left(\bar{M}_{g,h}\right)_{U(1)}}\frac{c_g(\E^\vee(\lambda))
 c_g(\E^\vee((a-1)\lambda))c_g(\E^\vee(-a\lambda))\lambda^{2h-3}}
 {\prod_{i=1}^h (\lambda-n_i\psi_i)}.
\end{eqnarray*}
\end{con} 

The above formulae for $C(g;h|d;n_1,\ldots,n_h|a)$ are calculated in \cite{KL} 
by localization techniques using the $S^1$ action $\varrho_a$. Actually, 
the definition of $\chi_{(g,h)}(X,L,\rho \mid\beta,\vga,\mu)$
is inspired by R. Bott's interpretation of the computations in \cite{KL}.
The localization formula, and in particular the proof of
Conjecture~\ref{localize}, is left to future work.

Finally, the assumption of the existence of an admissible $S^1$ action is too
restrictive. The $S^1$ action disappears when we perturb the almost complex structure $J$
or the Lagrangian submanifold $L$, so the invariant is not even defined for
other almost complex structures, and it is not clear in which sense 
$\chi_{(g,h)}(X,L,\rho \mid\beta,\vga,\mu)$
is an ``invariant''. It is desirable to find a natural way to impose 
boundary conditions for the general case.


\begin{thebibliography}{AAAA}

\bibitem[Ab]{A}  W.\ Abikoff, 
 {\em The real analytic theory of Teichm\"{u}ller space},
 Lecture Notes in Mathematics, Vol. 820. Springer-Verlag, Berlin-New York,
 1980.

\bibitem[Ar]{Ar} N.\ Aronszajn, 
 {\em A unique continuation theorem for solutions of elliptic partial differential equations or
 inequalities of second order},
 J. Math. Pures Appl. (9) {\bf 36} (1957), 235--249.
\bibitem[AG]{AG} N.L.\ Alling and N.\ Greenleaf, 
 {\em Foundations of The theory of Klein Surfaces},
 Lecture Notes in Mathematics, Vol. 219, Springer-Verlag, Berlin-New York,
1971.

\bibitem[AKV]{AKV} M.\ Aganagic, A.\ Klemm, and C.\ Vafa,
 {\em Disc Instantons, Mirror Symmetry, and the Duality Web\/},
 Z. Naturforsch. A 57 (2002), no. 1-2, 1--28
  
\bibitem[AV]{AV} M.\ Aganagic and C.\ Vafa, 
 {\em Mirror Symmetry, D-Branes and Counting Holomorphic Discs\/},
 preprint, 2000, hep-th/0012041.
  
\bibitem[BEGG]{BEGG}
 E. Bujalance, J.J.\ Etayo, J.M.\ Gamboa, and  G. Gromadzki,
 {\em Automorphism Groups of Compact Bordered Klein Surfaces},
 Lecture Notes in Mathematics, Vol. 1439, Springer-Verlag, Berlin-New York, 1990.
  
\bibitem[BF]{BF} K.\ Behrend and B.\ Fantechi,
 {\em The intrinsic normal cone},
 Invent. Math. {\bf 128} (1997), no. 1, 45--88.

\bibitem[CdGP]{CdGP} P.\ Candelas, X.\ de la Ossa, P.\ Green and L.\ Parkes, 
 {\em A pair of Calabi-Yau manifolds as an exactly soluble superconformal 
 field theory\/}, 
 Nuclear Physics {\bf B359} (1991), 21--74 and in {\sl Essays on Mirror Manifolds\/}
   (S.-T.\ Yau, ed.), International Press, Hong Kong 1992, 31--95.
 
\bibitem[FO]{FO} K.\ Fukaya and K.\ Ono,
 {\em Arnold conjecture and Gromov-Witten invariant},
 Topology \textbf{38} (1999), no. 5, 933--1048.
 
\bibitem[FO$_3$]{FO3} 
 K.\ Fukaya, Y.-G.\ Oh, H.\ Ohta, and K.\ Ono, 
 {\em Lagrangian intersection Floer theory--anomaly and obstruction}, 
 preprint  2000.
 
\bibitem[Gi]{G} A.\ Givental, 
 {\em Equivariant Gromov-Witten Invariants\/}, 
 Internat.\ Math.\ Res.\ Notices (1996), 613--663.
 
\bibitem[GV]{gov} R.\ Gopakumar and C.\ Vafa,
 {\em M-Theory and Topological Strings--II\/},
 preprint 1998, hep-th/9812127.

\bibitem[GP]{GP} T.\ Graber and R.\ Pandharipande,
 {\em Localization of virtual classes}, 
 Invent. Math. \textbf{135} (1999), 487-518
 
\bibitem[GZ]{GZ} T.\ Graber, E.\ Zaslow
 {\em Open-String Gromov-Witten Invariants: Calculations and a
 Mirror "Theorem"}, hep-th/0109075.

\bibitem[Gr]{Gr} M.\ Gromov,
 {\em Pseudoholomorphic curves in symplectic manifolds},
 Invent. Math. \textbf{82} (1985), no. 2, 307--347.

\bibitem[Ha]{H} R.\ Hartshorne,
 {\em Algebraic geometry},
 Graduate Texts in Mathematics, 52, Springer-Verlag,
 New York-Heidelberg, 1977.
 
\bibitem[Hi]{Hi} M.W.\ Hirsch
 {\em Differential topology}, 
 Graduate Texts in Mathematics, 33, Springer-Verlag, 
 New York-Heidelberg, 1976.

\bibitem[HM]{HM} J.\ Harris and I.\ Morrison,
 {\em Moduli of curves}, 
 Graduate Texts in Mathematics, 187, Springer-Verlag, New York, 1998.


\bibitem[IS1]{IS1} S.\ Ivashkovich and V.\ Shevchishin,
{\em Gromov compactneww theorem for $J$-complex curves with boundary}
Internat. Math. Res. Notices 2000,  no. 22, 1167--1206.                   

\bibitem[IS2]{IS2}S.\ Ivashkovich and V.\ Shevchishin, 
{\em Holomorphic structure on the space of Riemann surfaces with marked 
boundary},
Tr. Mat. Inst. Steklova  \textbf{235} (2001),  
Anal. i Geom. Vopr. Kompleks. Analiza, 98--109;  
translation in  Proc. Steklov Inst. Math.  2001,  no. 4 (235), 91--102.

\bibitem[IS3]{IS3} S.\ Ivashkovich, V.\ Shevchishin, 
{\em Reflection principle and $J$-complex curves with boundary on totally 
real immersions},
Commun. Contemp. Math. \textbf{4} (2002), no. 1, 65--106.

\bibitem[Kl]{Klein} F.\ Klein,
 {\em \"{U}ber Realit\"{a}tsverh\"{a}ltnisse bei der einem beliebigen Geschlechte
    zugeh\"{o}rigen Normalkurve der $\phi$},
  Math. Ann., Vol. 42, 1892. 

\bibitem[Ko]{Ko} K.\ Kodaira,
 {\em Complex manifolds and deformation of complex structures},
 Grundlehrender Mathematischen Wissenschaften 
 [Fundamental Principles of Mathematical Sciences], 283.
 Springer-Verlag, New York, 1986. 
 
\bibitem[KL]{KL} S.\ Katz and C.C.\ Liu,
 {\em Enumerative Geometry of Stable Maps with Lagrangian Boundary
 Conditions  and Multiple Covers of the Disc}, 
 Adv. Theor. Math. Phys. {\bf 5} (2001), no. 1, 1--49.

\bibitem[LLY]{LLY} B.\ Lian, K.\ Liu and S.-T.~Yau,
 {\em Mirror principle I}, 
 Asian J.\ Math.\ {\bf 1} (1997), 729--763.
   
\bibitem[LS]{LS}  J.\ Li and Y.S.\ Song, 
{\em Open string instantons and relative stable morphisms},
 Adv. Theor. Math. Phys. 5 (2001), no. 1, 67--91.

\bibitem[LT1]{LT1} J.\ Li and G.\ Tian,
 {\em  Virtual moduli cycles and Gromov-Witten invariants of algebraic varieties},
 J. Amer. Math. Soc. 11 (1998), no. 1, 119--174.

\bibitem[LT2]{LT2} J.\ Li and G.\ Tian,
 {\em Virtual moduli cycles and Gromov-Witten invariants of general
 symplectic manifolds},
 Topics in symplectic $4$-manifolds (Irvine, CA, 1996), 47--83, 
 First Int. Press Lect. Ser., I, Internat. Press, Cambridge, MA, 1998.

\bibitem[LT]{LT} G.\ Liu and G.\ Tian,
 {\em Floer homology and Arnold conjecture}, 
 J. Differential Geom. {\bf 49} (1998), no. 1, 1--74.
 
\bibitem[LMV]{LMV} J.M.F.\ Labastida, M.\ Mari\~{n}o and C.\ Vafa
 {\em Knots, links and branes at large N\/},
 Jour.\ High En.\ Phys.\ 0011 (2000) 007.
  
\bibitem[MS]{MS} D.\ McDuff and D.\ Salamon,
{\em $J$-holomorphic curves and quantum cohomology}
University Lecture Series, 6. 
American Mathematical Society, Providence, RI, 1994. 
  
\bibitem[MV]{MV} M.\ Mari\~{n}o and C.\ Vafa
 {\em Framed knots at large N},
 hep-th/0108064.
  

\bibitem[OV]{OV} H.\ Ooguri and C.\ Vafa,
 {\em Knot invariants and topological strings},
 Nuclear Phys. \textbf{B577} (2000), 69-108.

\bibitem[PW]{PW} T.\ Parker and J.G.\ Wolfson, 
{\em Pseudo-holomorphic maps and bubble trees},
J. Geom. Anal. \textbf{3}  (1993),  no. 1, 63--98. 

\bibitem[Sa]{Sa} I.\ Satake, 
 {\em The Gauss-Bonnet theorem for $V$-manifolds},
  J. Math. Soc. Japan \textbf{9} (1957), 464-492.    

\bibitem[Se]{S} M.\ Sepp\"{a}l\"{a}, 
 {\em Moduli spaces of stable real algebraic curves},
  Ann. Sci. \'{E}cole norm.Sup. (4) \textbf{24} (1991).

\bibitem[Sie1]{Si} B.\ Siebert,
 {\em Gromov-Witten invariants of general symplectic manifolds},
 dg-ga/9608005 

\bibitem[Sie2]{Si2} B.\ Siebert, 
 {\em Symplectic Gromov-Witten invariants},
 New trends in algebraic geometry (Warwick, 1996), 375--424, 
 London Math. Soc. Lecture Note Ser., 264, Cambridge Univ. Press,
 Cambridge, 1999.

\bibitem[Sil]{Sil} R.\ Silhol, 
 {\em Real algebraic surfaces},
 Lecture Notes in Mathematics, Vol. 1392. Springer-Verlag, Berlin, 1989

\bibitem[St]{St} J.D. Stasheff,
{\em Homotopy associativity of $H$-spaces I, II},
Trans. Amer. Math. Soc. \textbf{108} (1963), 
275--292,  293--312. 

\bibitem[SS]{SS} M.\ Sepp\"{a}l\"{a} and R.\ Silhol,
 {\em Moduli spaces for real algebraic curves and real abelian varieties},
 Math. Z. \textbf{201} (1989), 151-165.

 
\bibitem[SU]{SU} J.\ Sacks and K.\ Uhlenbeck,
 {\em The existence of minimal immersions of $2$-spheres},
 Ann. of Math. (2) \textbf{113} (1981), no. 1, 1--24.

\bibitem[Th]{Th}  W.P. Thurston,
{\em The geometry and topology of three-manifolds},     
Princeton mimeographed lecture notes, 1979.

\bibitem[We]{W} G.\ Weichhold,
 {\em \"{U}ber symmetrische Riemannsche Fl\"{a}chen und die
 Periodizit\"{a}tsmodulen der zugerh\"{o}rigen Abelschen Normalintegrale
 erstes Gattung },
 Leipziger dissertation, 1883.

\bibitem[Wolf]{Wolf} J.G.\ Wolfson, 
{\em Gromov's compactness of pseudo-holomorphic curves and symplectic geometry},J. Differential Geom. \textbf{28}  (1988),  no. 3, 383--405. 

\bibitem[Wolp]{Wo} S.\ Wolpert, 
 {\em  On the Weil-Petersson geometry of the moduli space of curves},
  Amer. J. Math. \textbf{107} (1985), no. 4, 969--997.
  
\bibitem[Ye]{Ye} R.\ Ye,
 {\em Gromov's compactness theorem for pseudo holomorphic curves},
 Trans. Amer. Math. Soc. \textbf{342} (1994), no. 2, 671--694.

\end{thebibliography}
\end{document}